\documentclass[12pt,leqno]{amsbook}
\usepackage[all]{xy}

\swapnumbers
\theoremstyle{plain}
\newtheorem{thm}[subsection]{Theorem}
\newtheorem*{thmProlog}{Theorem}
\newtheorem{lem}[subsection]{Lemma}
\newtheorem{fact}[subsection]{Fact}
\newtheorem{prp}[subsection]{Proposition}
\newtheorem{obv}[subsection]{Observation}

\title[Koszul duality of operads]{Koszul duality of operads\\ and homology of partition posets}
\author{Benoit Fresse}
\date{November 13, 2002}
\address{Lab. J.A. Dieudon\-n\'e\\
Uni\-ver\-sit\'e de Nice-Sophia-Anti\-polis\\
Parc Val\-rose\\
F-06108 Nice Cedex 02\\
France}
\email{fresse@math.unice.fr}
\urladdr{http://www-math.unice.fr/\~{ }fresse}
\begin{document}

% Notation

\newcommand{\K}{\mathbb{K}}
\newcommand{\F}{\mathbb{F}}
\newcommand{\N}{\mathbb{N}}
\newcommand{\Z}{\mathbb{Z}}
\newcommand{\Q}{\mathbb{Q}}
\newcommand{\Func}{\mathcal{F}}
\newcommand{\Pol}{\mathcal{P}}
\newcommand{\Coef}{\mathcal{M}}
\newcommand{\dg}{\mathop{dg}}
\newcommand{\gr}{\mathop{gr}}
\newcommand{\s}{\mathop{s}}
\newcommand{\Hom}{\mathop{\mathrm{Hom}}\nolimits}
\newcommand{\Char}{\mathop{Char}}
\renewcommand{\lim}{\mathop{\mathrm{lim}}}
\newcommand{\colim}{\mathop{\mathrm{colim}}}
\newcommand{\coend}{\mathop{\mathrm{coend}}}
\newcommand{\im}{\mathop{\mathrm{im}}\nolimits}
\renewcommand{\ker}{\mathop{\mathrm{ker}}\nolimits}
\newcommand{\coker}{\mathop{\mathrm{coker}}\nolimits}
\newcommand{\one}{\mathbf{1}}
\newcommand{\sgn}{\mathop{\mathrm{sgn}}\nolimits}
\renewcommand{\L}{\mathcal{L}}
\newcommand{\A}{\mathcal{A}}
\newcommand{\C}{\mathcal{C}}
\renewcommand{\P}{\mathcal{P}}
\newcommand{\G}{\mathcal{G}}
\newcommand{\M}{\mathcal{M}}
\newcommand{\E}{\mathcal{E}}
\newcommand{\LMod}[1]{\mathop{{#1}\mathrm{Mod}}\nolimits}
\newcommand{\RMod}[1]{\mathop{\mathrm{Mod}{#1}}\nolimits}
\newcommand{\LAlg}[1]{\mathop{{#1}\mathrm{Alg}}\nolimits}
\newcommand{\In}{\mathop{In}}
\newcommand{\Out}{\mathop{Out}}
\newcommand{\T}{\mathcal{T}}
\newcommand{\Pn}{\mathop{\mathcal{P}n}}
\newcommand{\Sh}{\mathop{\mathcal{S}h}}
\newcommand{\Tor}{\mathop{\mathrm{Tor}}\nolimits}
\newcommand{\EZ}{\mathop{\mathrm{EZ}}}
\newcommand{\Ho}{\mathop{\mathrm{Ho}}}

\frontmatter

% Front page (Plain TeX code)
\thispagestyle{empty}
%\pdfbookmark[0]{Title}{Title}
\begin{center}
\vbox to 20mm{ }
\Large\textbf{\uppercase{Koszul duality of operads\\
\smallskip
and\\
\smallskip
homology of partition posets\\}}
\medskip
\textit{by}\\
\medskip
\textrm{Benoit Fresse}\\
\medskip
\textrm{November 13, 2002}
\end{center}
\vfill
\noindent\textbf{Address:} Benoit Fresse\\
\indent Lab. J.A. Dieudonn\'e\\
\indent Universit\'e de Nice-Sophia-Anti\-polis\\
\indent Parc Valrose\\
\indent F-06108 Nice Cedex 02\\
\indent France\\
\textbf{E-mail:} fresse@math.unice.fr\\
\textbf{URL:} {http://www-math.unice.fr/\~{ }fresse/}

\newpage
\cleardoublepage
\tableofcontents

\mainmatter

%\input{Prolog}
% 30/1/2003

\chapter*{Prolog}\label{Prolog}

We consider the set of partitions of $\{1,\ldots,r\}$
equipped with the partial order $\leq$
defined by the refinement of partitions
(for instance, we have $\{1,3\},\{2,4\}\leq\{1\},\{3\},\{2,4\}$).
The single set $\{1,\ldots,r\}$ forms the smallest partition of $\{1,\ldots,r\}$
and
the collection $\{1\},\ldots,\{r\}$
forms the largest partition.
We mention that the components of a partition are not assumed to be ordered
(for instance, the collections $\{1,3\},\{2,4\}$ and $\{2,4\},\{1,3\}$ define the same partition).
We observe also that the set of partitions is equipped with an action of the symmetric group $\Sigma_r$,
since a permutation $w: \{1,\ldots,r\}\,\rightarrow\,\{1,\ldots,r\}$ maps a partition of $\{1,\ldots,r\}$
to another partition.

We consider the simplicial set $\bar{K}(r)$
which consists of sequences of partitions
$\lambda_0\leq\lambda_1\leq\cdots\leq\lambda_n$
such that $\lambda_0 = \{1,\ldots,r\}$
and $\lambda_n = \{1\},\ldots,\{r\}$.
The face $d_i$ is given by the omission of $\lambda_i$.
The degeneracy $s_j$ is given by the repetition of $\lambda_j$.
We have the following result:

\begin{thmProlog}
We consider the reduced homology
of the simplicial set $\bar{K}(r)$
with coefficients in $\K = \Z$, $\Q$ or $\F_p$.
This homology is a graded $\Sigma_r$-module
since the set of partitions is equipped with an action of the symmetric group.
We have
\[ \tilde{H}_*(\bar{K}(r),\K)
= \begin{cases}
\L(r)^{\vee}\otimes\sgn_r, & \text{if $*=r-1$}, \\
0, & \text{otherwise}, \end{cases}\]
where $\L(r)$ denotes the $r$th component of the Lie operad $\L$.
\end{thmProlog}

The Lie operad $\L$ consists precisely of a sequence of representations $\L(r)$
which are given by the components of weight $(1,\ldots,1)$
in the free Lie algebras
$\L(x_1,\ldots,x_r)$.
In the theorem,
we consider tensor products of dual representations $\L(r)^{\vee}$
with the classical signature representations $\sgn_r$.

Here are bibliographical references about this theorem.
The character of the representation $\tilde{H}_{r-1}(\bar{K}(r),\K)$ is determined
by R. Stanley in \cite{Stanley} and by P. Hanlon in \cite{Hanlon}.
The relationship with the character of the Lie operad
is pointed out by A. Joyal in \cite{Joyal}.
These results about characters are refined by A. Robinson and S. Whitehouse in \cite{RobinsonWhitehouse}
(see also S. Whitehouse \cite{Whitehouse} for the integral version).

One can deduce from results of A. Bj\"orner that the (reduced) homology of partition posets
vanishes in degree $*\not=r-1$
(\emph{cf}. \cite{Bjorner}).
Then, in \cite{Barcelo},
H. Barcelo defines an isomorphism of representations $\tilde{H}_{r-1}(\bar{K}(r),\K)\simeq\L(r)^{\vee}\otimes\sgn_r$,
based on the Lyndon basis of the Lie operad
(\emph{cf}. M. Lothaire \cite{Lothaire}, C. Reutenauer \cite{Reutenauer}).
This result is improved by P. Hanlon and M. Wachs
in \cite{HanlonWachs}.
Namely,
these authors define a natural morphism
(which does not involve the choice of a basis of the Lie operad)
from the dual of the Lie operad
to the chain complex of the partition poset
\[\L(r)^{\vee}\otimes\sgn_r\,\rightarrow\,C_{r-1}(\bar{K}(r)).\]
This morphism fixes a representative of the homology class
associated to a given element of the Lie operad.
In addition,
P. Hanlon and M. Wachs generalize the theorem above
and give a relationship between partition posets
and structures of Lie algebras with $k$ary brackets.

On the other hand,
a topological proof of the theorem above, based on calculations of F. Cohen (\emph{cf}. \cite{Cohen}),
is given by G. Arone and M. Kankaanrinta
in \cite{AroneKankaanrinta}.
In connection with this result,
we should mention that an article of G. Arone and M. Mahowald (\emph{cf}. \cite{AroneMahowald})
sheds light on the importance of partition posets
in homotopy theory.
Namely,
these authors prove that the Goodwillie tower of the identity functor on topological spaces
is precisely determined by partition posets.

\chapter*{Introduction}

In this article,
we would like to point out that the theorem of the prolog
occurs a by-product of the theory of \emph{Koszul operads}.
For that purpose,
we improve results of V. Ginzburg and M. Kapranov (\emph{cf}. \cite{GinzburgKapranov})
in several directions.
More particularly,
we extend the Koszul duality of operads to operads defined over a field of positive characteristic
(or over a ring).
In addition,
we obtain more conceptual proofs of theorems
of V. Ginzburg and M. Kapranov.

We recall that an operad is an algebraic structure
which parame\-tri\-zes a set of operations.
Originally,
the notion of an operad is introduced by P. May in a topological context
in order to model the stucture of an iterated loop space
(\emph{cf}. \cite{MayLoop}).
To be precise,
the definition of an operad arises on one hand
from the notion of an \emph{algebraic theory}
(\emph{cf}. J. Boardman and R. Vogt \cite{BoardmanVogt},
F. Lawvere \cite{Lawvere},
S. Mac Lane \cite{MacLaneProp})
and on the other hand from the structure of a \emph{monad}
(also called \emph{triple}, \emph{cf}. J. Beck \cite{Beck}, S. Eilenberg and C. Moore \cite{EilenbergMoore}).
In particular,
we should mention that J. Boardman and R. Vogt consider certain \emph{prop} structures,
which are equivalent to operads,
in their work on iterated loop spaces (\emph{cf}. \cite{BoardmanVogt}).
We refer to Adams's book \cite{Adams}
for a nice survey of this subject.

In this article,
we are concerned with operads in a category of modules over a ground ring $\K$.
In this context,
an operad $P$ consists of a sequence of representations $P(r)$ of the symmetric groups $\Sigma_r$
together with composition products
\[P(r)\otimes P(n_1)\otimes\cdots\otimes P(n_r)\,\rightarrow\,P(n_1+\cdots+n_r).\]
In general,
an operad $P$ is associated to a category of algebras, called $P$-algebras,
for which the elements of $P(r)$
represent multilinear operations on $r$ variables.
The composition products of $P$
determine the composites of these operations.
For instance,
there are operads denoted by $\C$, respectively $\L$,
associated
to the category of commutative and associative algebras,
respectively to the category of Lie algebras.
(To be precise, in this article, we consider commutative algebras without units.)
The $\K$-module $\C(r)$ is generated by the single monomial $p(x_1,\ldots,x_r) = x_1\cdot\,\cdots\,\cdot x_r$
in $r$ commutative variable
which has degree $1$ in each variable.
Consequently, this $\K$-module $\C(r)$ is the trivial representation of the symmetric group $\Sigma_r$.
The $\K$-module $\L(r)$ is generated by the Lie monomials in $r$ variables
which have degree $1$ in each variable.

Classically,
a sequence of representations of the symmetric groups $M(r)$
is equivalent to a functor $S(M): \LMod{\K}\,\rightarrow\,\LMod{\K}$
whose expansion
\[S(M)(V) = \bigoplus_{r=0}^{\infty} (M(r)\otimes V^{\otimes r})_{\Sigma_r}\]
consists of coinvariant modules
$S_r(M)(V) = (M(r)\otimes V^{\otimes r})_{\Sigma_r}$.
For instance,
the functor associated to the commutative operad $S(\C): \LMod{\K}\,\rightarrow\,\LMod{\K}$
can be identified with the augmentation ideal of the symmetric algebra
$S(\C)(V) = \bigoplus_{r=1}^{\infty} (V^{\otimes r})_{\Sigma_r}$.
Similarly,
the functor associated to the Lie operad $S(\L): \LMod{\K}\,\rightarrow\,\LMod{\K}$
is given by the structure of a free Lie algebra.
In general,
the structure formed by a sequence of representations of the symmetric groups $M(r)$
is called a \emph{$\Sigma_*$-module}
(or a \emph{symmetric module} in plain English).

One observes that classical complexes associated to commutative algebras and Lie algebras
are functors determined by $\Sigma_*$-modules.
On one hand,
we consider the \emph{cotriple construction} of a commutative algebra $A$
which is a chain complex $C^{cotriple}_*(A)$
such that
\[C^{cotriple}_*(A) = \underbrace{S(\C)\circ\cdots\circ S(\C)}_{*-1}(A)\]
(\emph{cf}. M. Barr \cite{BarrHarrison}, J. Beck \cite{Beck}).
We prove precisely that the cotriple construction is the functor associated
to the chain complex of the partition posets $C_*(\bar{K}(r))$.
We have explicitly
$C^{cotriple}_*(A) = \bigoplus_{r=0}^\infty (C_*(\bar{K}(r))\otimes A^{\otimes r})_{\Sigma_r}$.

On the other hand,
we consider the \emph{Harrison complex} $C^{Harr}_*(A)$
(\emph{cf}. D. Harrison, \cite{Harrison})
which has the coinvariant module
\[C^{Harr}_{r-1}(A) = (\L(r)^{\vee}\otimes\sgn_r\otimes A^{\otimes r})_{\Sigma_r}\]
in degree $r-1$
(see paragraph~\ref{KoszulCommutativeHomology} for more precisions).
We define a comparison morphism
\[C^{Harr}_*(A)\,\rightarrow\,C^{cotriple}_*(A).\]
In positive characteristic,
the cotriple homology
differs from the Harrison homology.
In particular,
the cotriple homology of a symmetric algebra $H^{cotriple}_*(S(\C)(V))$ vanishes for general reasons
while the Harrison homology $H^{Harr}_*(S(\C)(V))$ does not
(\emph{cf}. M. Barr \cite{BarrHarrison}, D. Harrison \cite{Harrison}, S. Whitehouse \cite{Whitehouse}).
Consequently,
the comparison morphism above does not define a quasi-isomorphism of chain complexes.
Nevertheless,
we deduce from the Koszul duality of operads
that this comparison morphism
is induced by a quasi-isomorphism of $\Sigma_*$-modules
$\L(r)^{\vee}\otimes\sgn_r\,\xrightarrow{\sim}\,C_*(\bar{K}(r))$.
(In fact, we recover exactly the chain equivalence of M. Wachs and P. Hanlon.)
The theorem of the prolog follows from this assertion.

Similarly,
for Lie algebras,
we have, on one hand, the \emph{cotriple complex}
\[C^{cotriple}_*(\G) = \underbrace{S(\L)\circ\cdots\circ S(\L)}_{*-1}(\G)\]
and, on the other hand, the \emph{Chevalley-Eilenberg complex} $C^{CE}_*(\G)$,\linebreak
which can be defined by the formula
\[C^{CE}_{r-1}(\G) = (\C(r)^{\vee}\otimes\sgn_r\otimes\G^{\otimes r})_{\Sigma_{r}}\]
(\emph{cf}. M. Barr \cite{BarrCE}, J.-L. Koszul \cite{Koszul}).
One can deduce from classical arguments that the cotriple homology
agrees with the Chevalley-Eilenberg homology
(\emph{cf}. M. Barr \cite{BarrCE}).
We make this assertion more precise.
Namely, in the case of Lie algebras,
we have a comparison morphism
\[C^{CE}_*(\G)\,\rightarrow\,C^{cotriple}_*(\G)\]
which gives rise to an isomorphism from the Chevalley-Eilenberg homology
to the cotriple homology for all Lie algebras $\G$
which are projective over the ground ring $\K$.
In fact,
one purpose of this article is to point out that the comparison morphisms
\[C^{Harr}_*(A)\,\rightarrow\,C^{cotriple}_*(A)\qquad\text{and}\qquad C^{CE}_*(\G)\,\rightarrow\,C^{cotriple}_*(\G)\]
are determined by quasi-isomorphisms of $\Sigma_*$-modules in both cases,
although the first one does not give a quasi-isomorphism of chain complexes
for all commutative algebras $A$.

More generally,
we observe that $\Sigma_*$-modules have better homological properties than associated functors
and this makes results about coefficients easier.
Therefore,
we introduce an operation on $\Sigma_*$-modules $M,N\mapsto M\circ N$
which corresponds to the composition product of functors.
We have explicitly $S(M\circ N) = S(M)\circ S(N)$.
The cotriple constructions above are determined by composites $\Sigma_*$-modules,
since
$S(\C)\circ\cdots\circ S(\C)(V) = S(\C\circ\cdots\circ\C)(V)$
and
$S(\L)\circ\cdots\circ S(\L)(V) = S(\L\circ\cdots\circ\L)(V)$.
Similarly,
the comparison morphisms are induced by embeddings of $\Sigma_*$-modules
\[\L(r)^{\vee}\otimes\sgn_r\,\hookrightarrow\,\underbrace{\C\circ\cdots\circ\C}_{r-1}(r)
\quad\text{and}
\quad\C(r)^{\vee}\otimes\sgn_r\,\hookrightarrow\,\underbrace{\L\circ\cdots\circ\L}_{r-1}(r)\]
respectively.

Let us introduce the general constructions for operads
which give rise to these comparison morphisms.
The idea is to generalize classical structures from algebras to operad.
The tensor product of $\K$-modules is replaced by the composition product of $\Sigma_*$-modules.
The bar construction of algebras is replaced by the cotriple construction of operads
(which corresponds to the cotriple construction of monads
considered by J. Beck in \cite{Beck}
and by P. May in \cite{MayLoop}).
Explicitly,
for a given operad $P$,
we let $\bar{C}_*(P)$ denote the chain complex
such that
\[\bar{C}_*(P) = \underbrace{P\circ\cdots\circ P}_{*-1}.\]
We analyze the structure of this chain complex.
First,
we observe that the operadic bar construction $\bar{B}_*(P)$
introduced by E. Getzler and J. Jones in \cite{GetzlerJones}
and by V. Ginzburg and M. Kapranov in \cite{GinzburgKapranov}
is identified with a subcomplex of $\bar{C}_*(P)$.
Then,
we assume that $P$ is equipped with a weight grading $P = \bigoplus_{s=0}^{\infty} P_{(s)}$.
In this case,
the bar construction $\bar{B}_*(P)$ is equipped with an induced grading
$\bar{B}_*(P) = \bigoplus_{s=0}^{\infty} \bar{B}_*(P)_{(s)}$.
Furthermore,
we have $\bar{B}_*(P)_{(s)} = 0$ if $*>s$.
Therefore,
the homology modules $\bar{K}_s(P) = H_s(\bar{B}_*(P)_{(s)})$
form a subcomplex of $\bar{B}_*(P)$.
Finally,
we obtain embeddings of $\Sigma_*$-modules
\[\bar{K}_*(P)\,\hookrightarrow\,\bar{B}_*(P)\,\hookrightarrow\,\bar{C}_*(P).\]
The definition of $\bar{K}_*(P)$ generalizes a construction of V. Ginzburg and M. Kapranov.
To be precise,
these authors consider operads for which the weight grading coincides with the operadic grading.
More explicitly,
any operad $P$ can be equipped with a weight grading such that $P_{(s)}(r) = P(r)$ if $s = r+1$
and $P_{(s)}(r) = 0$ otherwise.
In this context,
V. Ginzburg and M. Kapranov
determine the $\Sigma_*$-modules $\bar{K}_*(P)$
associated to the commutative operad $P = \C$ and to the Lie operad $P = \L$.
We obtain precisely
$\bar{K}_*(\C)(r) = \L(r)^{\vee}\otimes\sgn_r$
and
$\bar{K}_*(\L)(r) = \C(r)^{\vee}\otimes\sgn_r$.
Hence,
the embeddings above provide
the required comparison morphisms.

Here is how we deduce the theorem of the prolog from Koszul duality arguments.
We introduce constructions with coefficients $C_*(L,P,R)$ and $K_*(L,P,R)$
such that
\[C_*(I,P,P) = \bar{C}_*(P)\circ P\qquad\text{and}\qquad K_*(I,P,P) = \bar{K}_*(P)\circ P.\]
Then, we generalize a classical comparison argument from algebras to operads.
Namely,
we prove that a morphism of $\Sigma_*$-modules $\phi: M\,\rightarrow\,N$,
where $M = \bar{M}\circ P$ and $N = \bar{N}\circ P$,
induces a quasi-isomorphism
$\bar{\phi}: \bar{M}\,\rightarrow\,\bar{N}$
if the $\Sigma_*$-modules $M$ and $N$ are both acyclic chain complexes.
(Let us mention that, in positive characteristic,
this property holds for $\Sigma_*$-modules
but not for functors.)
Consequently,
we conclude that the embedding $\bar{K}_*(P)\,\hookrightarrow\,\bar{C}_*(P)$
is a quasi-isomorphism
if and only if the complex $K_*(I,P,P) = \bar{K}_*(P)\circ P$
is acyclic.

Since $S(K_*(I,P,P)) = S(\bar{K}_*(P))\circ S(P)$,
the complex of $\Sigma_*$-mo\-du\-les $K_*(I,P,P)$
is associated to the Harrison complex of the symmetric algebra $C^{Harr}_*(S(\C)(V))$ in the case $P = \C$
and to the Chevalley-Eilenberg complex of the free Lie algebra $C^{CE}_*(S(\L)(V))$ in the case $P = \L$.
We mention that the latter complex is acyclic
while the former is not.
Nevertheless,
we prove that the associated complexes of $\Sigma_*$-modules $K_*(I,\C,\C)$ and $K_*(I,\L,\L)$
are both acyclic.
To be precise,
the complex $K_*(I,\L,\L)$
is acyclic because this property holds for the associated functor $S(K_*(I,\L,\L))(V) = C^{CE}_*(S(\L)(V))$.
Then,
we deduce that the complex $K_*(I,\C,\C)$
is acyclic by Koszul duality
between the commutative operad and the Lie operad.

As mentioned above,
our definition of the Koszul construction $P\mapsto\bar{K}(P)$
generalizes the definition of V. Ginzburg and M. Kapranov.
In fact,
there are Koszul duality results for generalizations
of the Lie operad
in which one considers Lie brackets
in more than $2$ variables
(\emph{cf}. E. Getzler \cite{Getzler}, A.V. Gnedbaye \cite{Gnedbaye},
A.V. Gnedbaye and M. Wambst \cite{GnedbayeWambst},
P. Hanlon and M. Wachs \cite{HanlonWachs},
Y. Manin \cite{Manin}).
The theory of V. Ginzburg and M. Kapranov
does not work for these examples,
but the generalized construction gives the right results.

Finally,
we would like to mention that similar Koszul duality results seem to occur for the derived functors
\begin{align*} & \Tor_*(\C)(V) = \bigoplus_{r=0}^\infty \Tor^{\Sigma_r}_*(\C(r),V^{\otimes r}) \\
\text{and}\qquad & \Tor_*(\L)(V) = \bigoplus_{r=0}^\infty \Tor^{\Sigma_r}_*(\L(r),V^{\otimes r})
\end{align*}
associated to the commutative and Lie operads.
Hints are given, on one hand, by the calculations of G. Arone and M. Mahowald
(\emph{cf}. \cite{AroneMahowald},
see also G. Arone and M. Kankaanrinta \cite{AroneKankaanrinta})
and, on the other hand, by theorems obtained by A.K. Bousfield and E.B. Curtis in the context of unstable algebras
(\emph{cf}. \cite{BousfieldCurtis}).

Here are a few indications about the contents of this article.
We survey operad structures in section~\ref{Operads}.
We generalize results of classical homology theory from modules over an algebra to modules over an operad
in section~\ref{ModuleComplexes}.
We recall the main properties of the Bar construction $\bar{B}_*(P)$ in section~\ref{ReducedBar}.
We introduce the cotriple construction $\bar{C}_*(P)$
and the comparison morphism $\bar{B}_*(P)\,\hookrightarrow\,\bar{C}_*(P)$
in section~\ref{BarCoefficients}.
We devote section~\ref{KoszulDuality} to the Koszul construction $\bar{K}_*(P)$.
We go back to partition posets and applications in the last section
of this article.
The main sections of the article (\ref{ModuleComplexes}-\ref{BarCoefficients})
have independent introductions.

\renewcommand{\thesubsection}{\thechapter.\arabic{subsection}}
%\input{Conventions}
% 30/1/2003

\chapter*{Conventions}

\subsection{Notation}\label{Conventions}
We fix a commutative ground ring $\K$.
We work within the category of $\K$-modules, denoted by $\LMod{\K}$.
We let $V^{\vee} = \Hom_\K(V,\K)$ denote the dual $\K$-module of any $V\in\LMod{\K}$.

The group of permutations of $\{1,\ldots,r\}$
is denoted by $\Sigma_r$.
If $M$ is a $\Sigma_r$-module,
then $M^{\vee}$ is the dual $\K$-module of $M$
and is equipped with an unsigned action of $\Sigma_r$.
Explicitly, if $f: M\,\rightarrow\,\K$,
then we set
$w\cdot f(x) = f(w^{-1}\cdot x)$,
for all $w\in\Sigma_r$.

\subsection{Symmetric monoidal categories}
We recall that a symmetric monoidal category
consists of a category $\M$
together with an associative bifunctor $\otimes: \M\times\M\,\rightarrow\,\M$
which possesses a symmetry isomorphism
$c(X,Y): X\otimes Y\,\rightarrow\,Y\otimes X$
and a unit $\one\in\M$.

In this article,
we are concerned with the following basic symmetric monoidal categories:
the category of $\K$-modules, denoted by $\LMod{\K}$;
the category of differential graded $\K$-modules, denoted by $\dg\LMod{\K}$;
the category of simplicial $\K$-modules, denoted by $\s\LMod{\K}$.
In the next paragraphs,
we recall the definition of these categories
and some related conventions.

\subsection{Differential graded modules}\label{DGMod}
A differential graded $\K$-mo\-du\-le (a \emph{dg-module}, for short)
denotes a lower $\N$-graded $\K$-module $V = \bigoplus_{*\in\N} V_*$
equipped with a differential $\delta: V_*\,\rightarrow\,V_{*-1}$
which decreases degrees by $1$.
A dg-module is equivalent to a chain complex
\[\cdots\,\xrightarrow{\delta}\,V_d\,\xrightarrow{\delta}\,\cdots\,\xrightarrow{\delta}\,V_1\,\xrightarrow{\delta}\,V_0.\]
In general,
the homology of this chain complex is denoted by $H_*(V)$.
But,
in the case of a non canonical differential $\delta: V_*\,\rightarrow\,V_{*-1}$,
we adopt the notation $H_*(V,\delta)$.
The notation $|v| = d$ indicates the degree of a homogeneous element $v\in V_d$.

A morphism $f: V\,\rightarrow\,V'$ is homogeneous of lower degree $|f| = d$
if we have $f(V_*)\subset V'_{*+d}$.
The dual of a dg-module $V$,
denoted by $V^{\vee}$,
is the dg-module
generated by homogeneous morphisms $f: V\,\rightarrow\,\K$.
More precisely,
we assume that the ground field $\K$
is equivalent to a dg-module concentrated in degree $0$.
Accordingly,
a homogeneous morphism of lower degree $|f| = d$
is equivalent to a map $f: V_{-d}\,\rightarrow\,\K$
and $V^{\vee}$ is the dg-module
such that $(V^{\vee})_{d} = (V_{-d})^{\vee}$.

A morphism of dg-modules $f: V\,\rightarrow\,V'$
is a homogeneous morphism of degree $0$
such that $\delta(f(v)) = f(\delta(v))$,
for all $v\in V$.
A quasi-isomorphism $f: V\,\xrightarrow{\sim}\,V'$
is a morphism of dg-modules
which induces a homology isomorphism $f_*: H_*(V)\,\xrightarrow{\simeq}\,H_*(V')$.

\subsection{The tensor product of dg-modules}\label{DGTensor}
The category of dg-modules, which we denote by $\dg\LMod{\K}$,
is equipped with the structure of a symmetric mo\-noi\-dal category.
We just consider the classical tensor product of dg-modules.
By definition,
we have $(V\otimes W)_d = \bigoplus_{i+j = d} V_i\otimes W_j$
and the differential of a homogeneous tensor $v\otimes w\in V\otimes W$
is given by the formula
$\delta(v\otimes w) = \delta(v)\otimes w + (-1)^{|v|} v\otimes\delta(w)$.
We recall that the symmetry isomorphism $c(V,W): V\otimes W\,\rightarrow\,W\otimes V$ involves a sign.
Precisely,
for $v\otimes w\in V\otimes W$,
we have the formula $c(V,W)(v\otimes w) = (-1)^{|v| |w|} w\otimes v$.
In general,
a commutation of homogeneous elements of respective degrees $d$ and $e$
is supposed to produce a sign $\pm = (-1)^{d e}$
as in the definition of the symmetry isomorphism.
In this article, we may denote this sign by $\pm$ without more specification.

For instance,
suppose given homogeneous morphisms
$f: V\,\rightarrow\,V'$ and $g: W\,\rightarrow\,W'$.
According to Koszul,
we have a homogeneous morphism
$f\otimes g: V\otimes W\,\rightarrow\,V'\otimes W'$
defined by the formula
$(f\otimes g)(v\otimes w) = (-1)^{|v| |g|} f(v)\otimes g(w)$,
for all $v\otimes w\in V\otimes W$.

\subsection{Graded modules}\label{GradedMod}
We have an inclusion of symmetric mo\-noi\-dal categories
$(\LMod{\K},\otimes)\,\subset\,(\dg\LMod{\K},\otimes)$.
Precisely,
a $\K$-module is equivalent to a dg-module
which has only one component in degree $0$.
We have a similar inclusion for the category of graded modules
$(\gr\LMod{\K},\otimes)\,\subset\,(\dg\LMod{\K},\otimes)$.
By definition,
a graded $\K$-module is a dg-module equipped with a trivial differential $\delta = 0$.
Because of this definition,
we assume the same commutation rule for graded modules
as for dg-modules.
In section~\ref{KoszulDuality},
we introduce the notion of a weight
in order to make distinct a grading
which does not contribute to the sign of a symmetry isomorphism.

\subsection{Simplicial modules}\label{SMod}
We assume the classical definition of a simplicial $\K$-module,
as a sequence of $\K$-modules $V_n$
equipped with faces $d_i: V_{n}\,\rightarrow\,V_{n-1}$
and degeneracies $s_j: V_{n}\,\rightarrow\,V_{n+1}$.
The tensor product of simplicial $\K$-modules $V\otimes W$ is defined by $(V\otimes W)_n = V_n\otimes W_n$
for all dimensions $n\in\N$.
The face $d_i: (V\otimes W)_{n}\,\rightarrow\,(W\otimes V)_{n-1}$
(respectively, the degeneracy $s_j: (V\otimes W)_{n}\,\rightarrow\,(W\otimes V)_{n+1}$)
satisfies the identity
$d_i(v\otimes w) = d_i(v)\otimes d_i(w)$
(respectively, $s_j(v\otimes w) = s_j(v)\otimes s_j(w)$),
for all $v\otimes w\in V_n\otimes W_n$.
The symmetry isomorphism $c(V,W): V\otimes W\,\rightarrow\,W\otimes V$
verifies $c(V,W)(v\otimes w) = w\otimes v$, for all $v\otimes w\in V_n\otimes W_n$.

We define the normalized chain complex of a simplicial $\K$-module
by the formula
\[N_n(V) = V_n/s_{0}(V_{n-1})+\cdots+s_{n-1}(V_{n-1}).\]
The differential of $N_*(V)$ is given by the alternate sum
\[\delta(v) = \sum_{i=0}^{n} (-1)^i d_i(v),\]
for all $v\in V_n$.
One observes that the normalization functor
\[N_*: \s\LMod{\K}\,\rightarrow\,\dg\LMod{\K}\]
is weak-monoidal,
since the classical Eilenberg-MacLane morphism defines a functorial and associative quasi-isomorphism
\[N_*(V)\otimes N_*(W)\,\xrightarrow{\sim}\,N_*(V\otimes W)\]
which commutes with symmetry isomorphisms
(\emph{cf}. S. Mac Lane \cite[Chapter 8]{MacLaneHomology}).

\subsection{Multiple dg-modules}\label{BiDGMod}
A bi-graded module $V$ is equipped with a splitting $V = \bigoplus_{s,t\in\N} V_{s t}$.
We refer to the module $V_{s t}$
as the homogeneous component of $V$
of horizontal degree $s$ and vertical degree $t$.
In the context of a bi-graded module,
the notation $|v|$ refers to the total degree $|v| = s+t$
of a homogeneous element $v\in V_{s t}$.
A bi-dg-module $V$ is a bi-graded module equipped with commuting differentials
such that
$\delta_h: V_{s t}\,\rightarrow\,V_{s-1 t}$
and
$\delta_v: V_{s t}\,\rightarrow\,V_{s t-1}$.
To be precise,
the commutation relation reads $\delta_h\delta_v + \delta_v\delta_h = 0$,
because $\delta_h$ and $\delta_v$ have total degree $-1$.
We generalize the commutation rule of tensors in the context of bi-dg-modules.
We assume precisely that signs are determined by total degrees.

\subsection{The spectral sequence of a bicomplex}\label{BicomplexSpectralSequence}
We consider the total complex $(V,\delta)$ of a bi-dg-module $(V,\delta_h,\delta_v)$.
The homogeneous component of degree $d$ of this dg-module is the sum $\bigoplus_{s+t=d} V_{s t}$.
The differential of the total complex is the sum $\delta = \delta_h+\delta_v$.
We have a spectral sequence
$I^r(V)\,\Rightarrow\,H_*(V,\delta)$
such that
\begin{align*} & (I^0_{s t},d^0) = (V_{s t},\delta_v),
\qquad (I^1_{s t},d^1) = (H_t(V_{s *},\delta_v),\delta_h) \\
& \text{and}\qquad I^2_{s t} = H_s(H_t(V_{* *},\delta_v),\delta_h). \end{align*}
We consider also the spectral sequence
$II^r(V)\,\Rightarrow\,H_*(V,\delta)$
such that
\begin{align*} & (II^0_{s t},d^0) = (V_{s t},\delta_h),
\qquad (II^1_{s t},d^1) = (H_s(V_{* t},\delta_h),\delta_v) \\
& \text{and}\qquad II^2_{s t} = H_t(H_s(V_{* *},\delta_h),\delta_v). \end{align*}

\subsection{Chain complexes of dg-modules}\label{DGModComplex}
In sections~\ref{ReducedBar} and~\ref{BarCoefficients},
we consider chain complexes
\[\cdots\,\xrightarrow{\partial}\,C_d(X)
\,\xrightarrow{\partial}\,C_{d-1}(X)
\,\xrightarrow{\partial}\,\cdots\,\xrightarrow{\partial}\,C_0(X),\]
where each component $C_d(X)$ is a dg-module $C_d(X) = (C_d(X)_*,\delta)$.
We refer to $\delta: C_d(X)_*\,\rightarrow\,C_d(X)_{*-1}$
as the internal differential of $C_d(X)$.
We assume that the external differential $\partial: C_d(X)\,\rightarrow\,C_{d-1}(X)$
is a homogeneous morphism of degree $-1$ and commutes with internal differentials.
We have explicitly $\partial(C_d(X)_*)\subset C_{d-1}(X)_{*-1}$ and $\delta\partial + \partial\delta = 0$.
In this situation,
the chain complex $(C_d(X),\partial)$ is equivalent to a bi-dg-module
which has the $\K$-module $C_s(X)_{s+t}$ in bidegree $(s,t)$.

In general,
such chain complexes are associated to a dg-object $X$ equipped with extra algebraic structures.
The map $X\mapsto C_d(X)$ defines a functor in $X$.
The internal differential $\delta: C_d(X)_*\,\rightarrow\,C_d(X)_{*-1}$
is induced by the differential of $X$.
The external differential $\partial: C_d(X)_*\,\rightarrow\,C_{d-1}(X)_{*-1}$
is determined by the algebraic structure of $X$.
The external subscript refers always to the internal degree of the dg-module $C_d(X)_*$
(and is equivalent to the total degree of the associated bi-dg-module).

\renewcommand{\thesection}{\thechapter.\arabic{section}}
\renewcommand{\thesubsection}{\thesection.\arabic{subsection}}
%\input{Operads}
% 30/1/2003

\chapter{Composition products and operad structures}\label{Operads}

\section{Introduction: monads and operads}

The purpose of this section is to recall the definition of an operad
and to survey some classical constructions
related to this structure.
In particular, we give in paragraph \ref{ClassicalOperadExamples}
the definition of the classical operads $\C$, $\A$ and $\L$,
which are associated to commutative and associative algebras,
associative algebras and Lie algebras respectively.

\subsection{The category of functors}
We let $\Func = \Func(\LMod{\K})$
denote the category of functors on $\K$-modules
$S: \LMod{\K}\,\rightarrow\,\LMod{\K}$.
We equip this category $\Func$
with the composition product of functors
$\circ: \Func\times\Func\,\rightarrow\,\Func$.
We observe that the composition product of functors is associative and has a unit $I\in\Func$
which is supplied by the identity functor $I(V) = V$.
One concludes that the composition product provides the category of functors $\Func$
with the structure of a (non-symmetric) monoidal category
(\emph{cf}. S. Mac Lane \cite{MacLaneCat}).

\subsection{Monads}\label{MonoidMonads}
According to S. Mac Lane \cite{MacLaneCat},
a \emph{monad} is a mo\-noid in the monoidal category of functors $(\Func,\circ)$.
Consequently,
a monad consists of a functor $S\in\Func$
equipped with a product $\mu: S\circ S\,\rightarrow\,S$ and a unit $\eta: I\,\rightarrow\,S$
that verify the classical associativity and unit relations.
Explicitly,
the following classical diagrams are supposed to commute:
\begin{align*} & \text{\textbf{a)} Associativity relation:} & \qquad & \text{\textbf{b)} Unit relations:} \\
& \begin{aligned}
\xymatrix{ S\circ S\circ S\ar[r]^{\mu\circ S}\ar[d]_{S\circ\mu} &
S\circ S\ar[d]^{\mu} \\ S\circ S\ar[r]^{\mu} & S \\ }
\end{aligned} & &
\begin{aligned}
\xymatrix{I\circ S\ar[r]^{\eta\circ S}\ar[dr]_{=} &
S\circ S\ar[d]^{\mu} &
S\circ I\ar[l]_{S\circ\eta}\ar[dl]^{=} \\
& S & \\ }
\end{aligned}\end{align*}
A monad morphism $\phi: S\,\rightarrow\,S'$ is a functor morphism
which preserves monad products.

\subsection{Algebras over a monad}\label{MonadAlgebras}
An \emph{algebra over a monad} $S$
is a $\K$-module $A$ equipped with a left monad action $\rho: S(A)\,\rightarrow\,A$
such that the following diagrams commute:
\begin{equation*}
\begin{aligned}
\xymatrix{ S(S(A))\ar[r]^{S(\rho)}\ar[d]_{\mu(A)} &
S(A)\ar[d]^{\rho} \\ S(A)\ar[r]^{\rho} & A }
\end{aligned}
\qquad\text{and}\qquad
\begin{aligned}
\xymatrix{ A\ar[r]^{\eta(A)}\ar[dr]_{=} & S(A)\ar[d]^{\rho} \\ & A \\ }
\end{aligned}
\end{equation*}
A morphism of $S$-algebras $\phi: A\,\rightarrow\,A'$ is a $\K$-module morphism
which commutes with monad actions.

The $\K$-module $S(V)$, where $V\in\LMod{\K}$,
represents the free $S$-algebra generated by $V$.
More precisely,
the monad product $\mu(V): S(S(V))\,\rightarrow\,S(V)$
defines a canonical monad action on $S(V)$.
In addition,
the monad unit induces a $\K$-module morphism
$\eta(V): V\,\rightarrow\,S(V)$
that satisfies the classical universality property.
Explicitly,
a $\K$-module morphism $\phi: V\,\rightarrow\,A'$, where $A'$ is an $S$-algebra,
has a unique extension
\[\xymatrix{ V\ar[rd]_{\eta(V)}\ar[rr]^\phi & & A' \\ & S(V)\ar@{-->}[ru]_{\tilde{\phi}} & \\ }\]
such that $\tilde{\phi}: S(V)\,\rightarrow\,A'$
is a morphism of $S$-algebras.
Equivalently,
we have an adjunction relation
\[\Hom_{\LAlg{S}}(S(V),A') = \Hom_{\LMod{\K}}(V,A').\]

One deduces from this relation that classical algebra categories,
which possess free objects,
are determined by monads
(\emph{cf}. S. Eilenberg and C. Moore \cite{EilenbergMoore}, S. Mac Lane \cite{MacLaneCat}).
For instance,
the category of associative algebras,
the category of associative and commutative algebras
and the category of Lie algebras
are associated to monad structures.
In paragraph~\ref{ClassicalOperadExamples},
we deduce the construction of such monads
from operad structures.

\subsection{The classical definition of an operad}\label{ClassicalOperadDefinition}
A \emph{$\Sigma_*$-module} $M$ is a sequence $M(r)$, $r\in\N$,
such that $M(r)$ is a representation of the symmetric group $\Sigma_r$.
A morphism of $\Sigma_*$-modules $f: M\,\rightarrow\,M'$
consists of a sequence of representation morphisms $f: M(r)\,\rightarrow\,M'(r)$.
A $\Sigma_*$-module $M$ has an associated functor $S(M): \LMod{\K}\,\rightarrow\,\LMod{\K}$
defined by the formula
\[S(M)(V) = \bigoplus_{r=0}^{\infty} (M(r)\otimes V^{\otimes r})_{\Sigma_r}.\]

According to P. May (\emph{cf}. \cite{MayLoop}),
an \emph{operad} is a $\Sigma_*$-module $P$
such that the functor $S(P): \LMod{\K}\,\rightarrow\,\LMod{\K}$
is equipped with the structure of a monad.
An operad morphism $\phi: P\,\rightarrow\,P'$ is a morphism of $\Sigma_*$-modules
such that the functor morphism $S(\phi): S(P)\,\rightarrow\,S(P')$
defines a morphism of monads.

In the next paragraph,
we interpret an operad element $p\in P(r)$ as a multilinear operation
in $r$ variables $p = p(x_1,\ldots,x_r)$.
The action of a permutation $w_*: P(r)\,\rightarrow\,P(r)$ is equivalent to a permutation of variables.
We have explicitly $w_*(p)(x_1,\ldots,x_r) = p(x_{w(1)},\ldots,x_{w(r)})$.

The monad product $S(P)\circ S(P)\,\rightarrow\,S(P)$ is determined by \emph{composition products}
\[P(r)\otimes P(n_1)\otimes\cdots\otimes P(n_r)\,\rightarrow\,P(n_1+\cdots+n_r),\]
defined for all $r\geq 1$ and all $n_1,\ldots,n_r\in\N$,
that satisfy natural $\Sigma_*$-invariance properties.
The operad composition products are equivalent to a composition process for multilinear operations.
Therefore,
the composite of an operation $p\in P(r)$ with $q_1\otimes\cdots\otimes q_r\in P(n_1)\otimes\cdots\otimes P(n_r)$
is also denoted by
\[p(q_1,\ldots,q_r)\in P(n_1+\cdots+n_r).\]
The associativity property of a monad product is equivalent to natural associativity relations
of operad composition products (\emph{cf}. P. May \cite{MayLoop, MayOp}).

The monad unit $I\,\rightarrow\,S(P)$
is determined by an element $1\in P(1)$ (the operad unit of $P$)
that represents an identity operation $1(x_1) = x_1$.
The unit properties of the monad product are also equivalent
to natural identities for operad composition products.
We have explicitly $1(p) = p$ and $p(1,\ldots,1) = p$, for all $p\in P(r)$.

\subsection{The classical definition of an algebra over an operad}\label{OperadAlgebras}
An \emph{algebra over an operad} $P$ is an algebra over the associated monad $S(P)$.
In the context of operads,
a monad action $S(P)(A)\,\rightarrow\,A$ is determined by $\Sigma_r$-invariant products
\[P(r)\otimes A^{\otimes r}\,\rightarrow\,A,\]
defined for all $r\in\N$.
Equivalently,
an element $p\in P(r)$ determines a multilinear operation $p: A^{\otimes r}\,\rightarrow\,A$.
Therefore,
the image of a tensor $a_1\otimes\cdots\otimes a_r\in A$ under the operad action is denoted by
\[p(a_1,\ldots,a_r)\in A.\]
The invariance relation reads
$w_*(p)(a_1,\ldots,a_r) = p(a_{w(1)},\ldots,a_{w(r)})$,
for all $w\in\Sigma_r$.

As hinted by the notation, the composite of $p: A^{\otimes r}\,\rightarrow\,A$
with $q_1: A^{\otimes n_1}\,\rightarrow\,A,\ldots,q_r: A^{\otimes n_r}\,\rightarrow\,A$
is supposed to agree with
$p(q_1,\ldots,q_r): A^{\otimes n_1+\cdots+n_r}\,\rightarrow\,A$,
the operation supplied by the composition product of $P$.
The operad unit $1\in P(1)$ is supposed to give the identity operation $1(a_1) = a_1$.

\subsection{Augmented and connected operads}\label{ConnectedOperads}
The identity functor $I\in\Func$ is associated to a $\Sigma_*$-module, also denoted by $I\in\M$.
We have explicitly $I(r) = \K$ for $r=1$ and $I(r) = 0$ for $r\not=1$.
We observe that this $\Sigma_*$-module $I$ gives a initial object in the category of operads,
since the unit of an operad $P$ is equivalent to a morphism $\eta: I\,\rightarrow\,P$.
Concretely, the operad $I$ contains no more elements than multiples of an operad unit $1\in I(1)$.
Accordingly, an $I$-algebra is nothing more than a $\K$-module.

An operad equipped with a morphism $\epsilon: P\,\rightarrow\,I$
is called an \emph{augmented operad}.
An augmented operad has an \emph{augmentation ideal}, denoted by $\tilde{P}$,
defined by $\tilde{P}(r) = \ker(\epsilon: P(r)\,\rightarrow\,I(r))$.
We mention that the augmentation $\epsilon: P\,\rightarrow\,I$
is a retraction of the unit morphism $\eta: I\,\rightarrow\,P$,
since an operad morphism is supposed to preserve units.
Hence, we have a splitting $P = I\oplus\tilde{P}$.

An operad $P$ is \emph{connected} if we have $P(0) = 0$ and $P(1) = \K\,1$.
We observe that a connected operad has a unique augmentation $\epsilon: P\,\rightarrow\,I$
and a canonical augmentation ideal, which satisfies
\[\tilde{P}(r) = \begin{cases} 0, & \text{for $r = 0,1$}, \\
P(r), & \text{for $r\not= 0,1$}. \end{cases}\]

\subsection{Partial composition products}\label{PartialComposites}
We call \emph{partial composition products} the operad operations
\[\circ_i: P(m)\otimes P(n)\,\rightarrow\,P(m+n-1)\]
defined by the relation $p\circ_i q = p(1,\ldots,q,\ldots,1)$
(the operation $q\in P(n)$ is put at the $i$th entry of $p\in P(m)$).
We have equivalently
\begin{multline*}(p\circ_i q)(x_1,\ldots,x_{m+n-1}) \\
= p(x_1,\ldots,x_{i-1},q(x_i,\ldots,x_{i+n-1}),x_{i+n},\ldots,x_{m+n-1}).\end{multline*}
Thus,
the index of an entry of $q(x_1,\ldots,x_n)$
is increased by $i-1$
in the composite operation $p\circ_i q$;
the index of an entry of $p(x_1,\ldots,x_m)$
such that $k>i$ is increased by $n-1$;
the index of an entry of $p(x_1,\ldots,x_m)$
such that $k<i$ is fixed.
In case of an augmented operad,
the partial composition products are determined by restrictions
\[\circ_i: \tilde{P}(m)\otimes\tilde{P}(n)\,\rightarrow\,\tilde{P}(m+n-1),\]
because we have the unit relations $p\circ_i 1 = p$ and $1\circ_1 p = p$,
for all $p\in P(m)$.

According to M. Markl (\emph{cf}. \cite{MarklModels}),
the partial composition products
satisfy relations
which are equivalent to the associativity property
of operad composition products.
Explicitly:
if we assume $p\in P(m)$, $i\in\{1,\ldots,m\}$, $q\in P(n)$, $j\in\{1,\ldots,n\}$ and $r\in P(s)$,
then we have $(p\circ_i q)\circ_{i-1+j} r = p\circ_i(q\circ_j r)$;
if we assume $p\in P(m)$, $i_1,i_2\in\{1,\ldots,m\}$ (together with $i_1<i_2$),
$q_1\in P(n_1)$ and $q_2\in P(n_2)$,
then we obtain $(p\circ_{i_1} q_1)\circ_{i_2+n_1-1} q_2 = (p\circ_{i_2} q_2)\circ_{i_1} q_1$.

Conversely,
the composition products
$P(r)\otimes P(n_1)\otimes\cdots\otimes P(n_r)\,\rightarrow\,P(n_1+\cdots+n_r)$
are determined by partial products.
For instance,
because of the associativity and unit relations,
we obtain the formula
$p(q_1,\ldots,q_r) = (\cdots((p\circ_r q_r)\circ_{r-1} q_{r-1})\circ_{r-2}\cdots q_2)\circ_1 q_1$.

\subsection{Indexing by finite sets}\label{OperadIndexing}
We may assume that the entries of an operation $p = p(x_1,\ldots,x_r)$
are indexed by any set $I = \{i_1,\ldots,i_r\}$
in bijection with $\{1,\ldots,r\}$.
For this purpose,
we introduce the module
\[P(I) = \bigoplus_{i_*: \{1,\ldots,r\}\,\rightarrow\,I} P(r)/\equiv.\]
The sum ranges over bijections $i_*: \{1,\ldots,r\}\,\rightarrow\,I$.
Hence,
an element of $P(I)$ consists of a pair $(p,i_*)$,
where $p\in P(r)$ and $i_*: \{1,\ldots,r\}\,\rightarrow\,I$.
This element is supposed to represent the operation $p = p(x_{i_1},\ldots,x_{i_r})$.
We make the element $(w_*(p),i_{*})$, where $w\in\Sigma_r$ is a permutation,
equivalent to $(p,i_{w(*)})$,
where $i_{w(*)}: \{1,\ldots,r\}\,\rightarrow\,I$
denotes the composite of the bijection $i_*: \{1,\ldots,r\}\,\rightarrow\,I$
with the permutation $w: \{1,\ldots,r\}\,\rightarrow\,\{1,\ldots,r\}$.
We set explicitly
$(w_*(p),i_{*})\equiv (p,i_{w(*)})$.
We have equivalently
\[w_*(p)(x_{i_1},\ldots,x_{i_r}) = p(x_{i_{w(1)}},\ldots,x_{i_{w(r)}}).\]
We observe that the $\K$-modules $P(I)$ define a functor on the category of finite sets and bijections,
since a natural reindexing morphism $u_*: P(I)\,\rightarrow\,P(I')$
is associated to any bijection $u: I\,\rightarrow\,I'$.

We have a canonical isomorphism $P(n) = P(\{1,\ldots,n\})$.
Moreover,
the action of the symmetric group on $P(n)$
is identified with the morphisms $w_*: P(\{1,\ldots,n\})\,\rightarrow\,P(\{1,\ldots,n\})$
associated to permutations $w: \{1,\ldots,n\}\,\rightarrow\,\{1,\ldots,n\}$.

In this framework, the operad structure is equivalent to partial composition products
\[\circ_i: P(I)\otimes P(J)\,\rightarrow\,P(I\setminus\{i\}\amalg J)\]
which are associated to any element $i\in I$.
In a composite $p\circ_i q$,
we subsitute the variable $x_i$ of $p = p(x_{i_1},\ldots,x_{i_m})\in P(I)$
by the operation $q = q(x_{j_1},\ldots,q_{j_n})\in P(J)$.
There is no need to reindex the entries of the resulting operation,
which are identified with the variables
$\{x_{i_1},\ldots,x_{i_m}\}\setminus\{x_{i}\}\amalg\{x_{j_1},\ldots,x_{j_n}\}$.
Let us rewrite the relations of partial composition products:
for any $p\in P(I)$, $i\in I$, $q\in P(J)$, $j\in J$ and $r\in P(K)$,
we have $(p\circ_i q)\circ_j r = p\circ_i(q\circ_j r)$ in $P(I\setminus\{i\}\amalg J\setminus\{j\}\amalg K)$;
for any $p\in P(I)$, $i_1,i_2\in I$, $q_1\in P(J_1)$ and $q_2\in P(J_2)$,
we have $(p\circ_{i_1} q_1)\circ_{i_2} q_2 = (p\circ_{i_2} q_2)\circ_{i_1} q_1$
in $P(I\setminus\{i_1,i_2\}\amalg J_1\amalg J_2)$.

\subsection{The free operad}\label{FreeOperad}
The free operad generated by a $\Sigma_*$-module $M$,
denoted by $F(M)$,
is characterized by the classical universal property:
we have a natural inclusion morphism $M\subset F(M)$
and a morphism of $\Sigma_*$-modules $\phi: M\,\rightarrow\,P$, where $P$ is an operad,
has a unique extension
\[\xymatrix{ M\ar[rd]\ar[rr]^\phi & & P \\ & F(M)\ar@{-->}[ru]_{\tilde{\phi}} & \\ }\]
such that $\tilde{\phi}: F(M)\,\rightarrow\,P$ is an operad morphism.

We have an explicit realization of the free operad,
due to V. Ginz\-burg and M. Kapranov (\emph{cf}. \cite{GinzburgKapranov}),
which arises from the work of J. Boardman and R. Vogt
on homotopy invariant structures (\emph{cf}. \cite{BoardmanVogt}).
We recall this construction in section~\ref{FreeOperadStructure}.
Concretely, the free operad $F(M)$ is spanned by formal expressions
$(\cdots((x_1\circ_{i_2} x_2)\circ_{i_3}\cdots x_{l-1})\circ_{i_l} x_l$
which represent composites of generators
$x_1\in M(n_1),\ldots,x_l\in M(n_l)$.
In this formalism,
the definition of the partial composition product
$\circ_i: F(M)\otimes F(M)\,\rightarrow\,F(M)$
can be deduced from the associativity relations of paragraph~\ref{PartialComposites}.

We have a natural splitting $F(M) = \bigoplus_{r=0}^\infty F_{(r)}(M)$,
such that the module $F_{(r)}(M)$ consists of composites
$(\cdots((x_1\circ_{i_2} x_2)\circ_{i_3}\cdots x_{l-1})\circ_{i_l} x_l$
of weight $l = r$.
We have $F_{(0)}(M) = I$, $F_{(1)}(M) = M$
and we observe that the weight grading is preserved by composition products.
Explicitly,
if $p\in F_{(r)}(M)(m)$ and $q\in F_{(s)}(M)(n)$,
then we obtain $p\circ_i q\in F_{(r+s)}(M)(m+n-1)$.
The unit of the free operad $I\,\rightarrow\,F(M)$
is given by the inclusion morphism $I = F_{(0)}(M)\,\subset\,F(M)$.
The projection $F(M)\,\rightarrow\,F_{(0)}(M)$ provides the free operad
with a canonical augmentation $F(M)\,\rightarrow\,F_{(0)}(M) = I$.
The augmentation ideal $\tilde{F}(M)$ consists of components $F_{(r)}(M)$ such that $r\geq 1$.
The universal morphism $M\,\rightarrow\,F(M)$
is given by the inclusion morphism $M = F_{(1)}(M)\,\subset\,F(M)$.
We observe that the projection $F(M)\,\rightarrow\,F_{(1)}(M) = M$ identifies the $\Sigma_*$-module $M$
with the indecomposable quotient of the free operad.
By definition, the indecomposable quotient of $F(M)$ is the quotient of the augmentation ideal
by the partial composites $p\circ_i q\in\tilde{F}(M)(m+n-1)$
such that $p\in\tilde{F}(M)(m)$ and $q\in\tilde{F}(M)(n)$.

\subsection{Operad ideals and quotients}\label{OperadIdeals}
An \emph{operad ideal} $R\subset P$ consists of a sequence $R(n)$,
where $R(n)$ is a subrepresentation of $P(n)$,
such that
$x\circ_i q\in R(m+n-1)$, as long as $x\in R(m)$ and $q\in P(n)$,
and
such that
$p\circ_i y\in R(m+n-1)$, as long as $p\in P(m)$ and $y\in R(n)$.
The operad ideal generated by a collection of elements $R = (r_i)_i$
is the smallest ideal of $P$
such that $r_i\in R(n_i)$, for all $i$.

A quotient by an operad ideal $P/R(n) = P(n)/R(n)$ is equipped with the structure of an operad.
Moreover, an operad morphism $\phi: P\,\rightarrow\,P'$
induces a morphism on the quotient operad
$\bar{\phi}: P/R\,\rightarrow\,P'$
provided $R\subset\ker(\phi)$.
In particular,
if $R = (r_i)_i$ is generated by elements $r_i\in P(n_i)$,
then a morphism $\bar{\phi}: P/R\,\rightarrow\,P'$
is equivalent to an operad morphism $\phi: P\,\rightarrow\,P'$
such that $\phi(r_i) = 0$, for all $i$.

\subsection{Classical examples}\label{ClassicalOperadExamples}
As mentioned in paragraph~\ref{MonadAlgebras},
the structure of an associative and commutative algebra,
the structure of an associative algebra
and the structure of a Lie algebra
are associated to operads, denoted by $\C$, $\A$ and $\L$ respectively.
We recall the definition of these operads.
As in \cite[V. Ginzburg and M. Kapranov]{GinzburgKapranov},
we observe that $\C$, $\A$ and $\L$
can be defined by generators and relations.
Formally,
we define operads such that $P = F(M)/R$,
where $M$ is a $\Sigma_*$-module of generators,
and $R$ is an ideal of the free operad $F(M)$
generated by relations
$r_i\in F(M)$, $i\in I$.

\textbf{a)} The \emph{commutative operad} $\C$ is generated by a single operation $m = m(x_1,x_2)$
such that $m(x_2,x_1) = m(x_1,x_2)$.
Hence,
the $\Sigma_*$-module $M(r)$ is the trivial representation $M(2) = \K\,m$ for $r = 2$
and vanishes otherwise.
The ideal $R$ is generated by the associativity relation $r = m\circ_1 m - m\circ_2 m\in F(M)(3)$.
We have equivalently $r(x_1,x_2,x_3) = m(m(x_1,x_2),x_3) - m(x_1,m(x_2,x_3))$.
Accordingly,
an algebra over the commutative operad $\C$
is a $\K$-module $A$
equipped with an operation $m: A\otimes A\,\rightarrow\,A$
such that $m(a_2,a_1) = m(a_1,a_2)$, for all $a_1,a_2\in A$,
and which makes the composite operation
$r: A\otimes A\otimes A\,\rightarrow\,A$
vanish.
As a conclusion,
an algebra over the commutative operad $\C$
is nothing but a commutative and associative algebra (without unit)
in the classical sense.
In addition,
we deduce from the classical theory
that the free $\C$-algebra is given by the augmentation ideal of the symmetric algebra
\[\tilde{S}(V) = \bigoplus_{r=1}^{\infty} (V^{\otimes r})_{\Sigma_r}.\]
Consequently, the $\Sigma_r$-module $\C(r)$ is the trivial representation of the symmetric group.

\textbf{b)} The \emph{associative operad} $\A$ is generated by a single operation $m = m(x_1,x_2)$
which supports a free action of the symmetric group $\Sigma_2$.
Hence,
the $\Sigma_*$-module $M(r)$ is the regular representation $M(2) = \K[\Sigma_2]\,m$ for $r = 2$
and vanishes otherwise.
The ideal $R$ is generated by the associativity relation $r = m\circ_1 m - m\circ_2 m\in F(M)(3)$.
We have equivalently $r(x_1,x_2,x_3) = m(m(x_1,x_2),x_3) - m(x_1,m(x_2,x_3))$.
Accordingly,
an algebra over the associative operad $\A$
is a $\K$-module $A$
equipped with an operation $m: A\otimes A\,\rightarrow\,A$
which makes the composite operation
$r: A\otimes A\otimes A\,\rightarrow\,A$
vanish.
As a conclusion,
an algebra over the associative operad $\A$
is nothing but an associative algebra (without unit)
in the classical sense.
In addition,
we deduce from the classical theory
that the free $\A$-algebra is given by the augmentation ideal of the tensor algebra
\[\tilde{T}(V) = \bigoplus_{r=1}^{\infty} V^{\otimes r}.\]
Consequently, the $\Sigma_r$-module $\A(r)$ is the regular representation of the symmetric group.

\textbf{c)} The \emph{Lie operad} $\L$ is generated by a single operation $l = l(x_1,x_2)$
such that $l(x_2,x_1) = - l(x_1,x_2)$.
Hence,
the $\Sigma_*$-module $M(r)$ is the signature representation $M(2) = \K\,l$ for $r = 2$
and vanishes otherwise.
The ideal $R$ is generated by the Jacobi relation
$r = l\circ_1 l + (2,3,1)_*(l\circ_1 l) +(3,1,2)_*(l\circ_1 l)\in F_3(M)$.
We have equivalently
$r(x_1,x_2,x_3) = l(x_1,l(x_2,x_3)) + l(x_3,l(x_1,x_2)) + l(x_2,l(x_3,x_1))$.
An algebra over the Lie operad $\C$ is a $\K$-module $\G$
equipped with an operation $l: \G\otimes\G\,\rightarrow\,\G$
such that $l(x_2,x_1) = - l(x_1,x_2)$, for all $x_1,x_2\in\G$,
and which makes the composite operation
$r: \G\otimes\G\otimes\G\,\rightarrow\,\G$
vanish.
We observe that this operation is not supposed to verify the identity $[a,a] = 0$.
In order to obtain this result,
we have to change the definition of an algebra associated to an operad
(see proposition~\ref{DividedLieAlgebras}).
We refer to Reutenauer's monograph \cite{Reutenauer}
for descriptions of the free Lie algebra
and results about the representations $\L(r)$
associated to the Lie operad.

\section{Composition products and operad structures}

According to the classical definition of paragraph \ref{ClassicalOperadDefinition},
an operad is equivalent to a monad in the category of $\K$-modules.
We make this relationship more precise in this section.
For that purpose,
we introduce the structure of a $\Sigma_*$-module
which is nothing but the coefficient sequence of an analytic functor on $\K$-modules.
The category of $\Sigma_*$-modules is equipped with an operation,
called the composition product,
which reflects the composition of functors.
We have then a conceptual definition
of the structure of an operad.
Namely,
an operad is a monoid
with respect to the composition product
of $\Sigma_*$-modules.
In this context,
we can deduce the relationship between operads and monads
from the relationship between the composition product of $\Sigma_*$-modules
and the composition product of functors.
In the next sections,
we observe that the composition product of $\Sigma_*$-modules
has better homological properties than the composition product of functors.
These observations motivate our definition
of the structure of an operad.

\subsection{The functor associated to a $\Sigma_*$-module}\label{SymmetricFunctors}
We recall that a \emph{$\Sigma_*$-module} (a \emph{symmetric module} in plain English)
consists of a sequence $M(r)$, $r\in\N$,
such that $M(r)$ is a representation of the symmetric group $\Sigma_r$.
A morphism of $\Sigma_*$-modules $f: M\,\rightarrow\,M'$
is a sequence of representation morphisms $f: M(r)\,\rightarrow\,M'(r)$, $r\in\N$.
In general,
we denote the category of $\Sigma_*$-modules by $\LMod{\Sigma_*}$,
but, in this section, we adopt the shorter notation $\Coef = \LMod{\Sigma_*}$.

As mentioned in paragraph~\ref{ClassicalOperadDefinition},
a $\Sigma_*$-module $M$ determines a functor
$S(M): \LMod{\K}\,\rightarrow\,\LMod{\K}$.
The image of a $\K$-module $V\in\LMod{\K}$
under $S(M)$ is denoted either by $S(M)(V)$
or by $S(M,V)$.
We have explicitly $S(M,V) = \bigoplus_{r=0}^{\infty} S_r(M,V)$
where $S_r(M,V)$ is the module of coinvariants
\[S_r(M,V) = (M(r)\otimes V^{\otimes r})_{\Sigma_r}.\]
We denote by $x(v_1,\ldots,v_r)\in S_r(M,V)$
the element of $S_r(M,V)$ associated to the tensor
$x\otimes(v_1\otimes\cdots\otimes v_r)\in M(r)\otimes V^{\otimes r}$,
so that the equivariance relation reads
$w_*(x)(v_1,\ldots,v_r) = x(v_{w(1)},\ldots,v_{w(r)})$
for all $w\in\Sigma_r$.

We say that a $\Sigma_*$-module $M$ connected if $M(0) = 0$.
Clearly, this property hold if and only if the associated functor verifies $S(M,0) = 0$.
Conversely, a $\Sigma_*$-module $M$ such that $M(r) = 0$, for $r\not=0$,
determines a constant functor $S(M,V) \equiv M(0)$.

The dual of a $\Sigma_*$-module $M^{\vee}$
consists of the dual representations $M(r)^{\vee}$
of the modules $M(r)$.
We would like to recall that the modules $M(r)^{\vee}$ are equipped
with an unsigned action of the symmetric groups $\Sigma_r$
(see paragraph \ref{Conventions}).

\subsection{Operations on functors}
The category of functors is endowed with classical operations
\[\oplus: \Func\times\Func\,\rightarrow\,\Func,
\quad\otimes: \Func\times\Func\,\rightarrow\,\Func
\quad\text{and}\quad\circ: \Func\times\Func\,\rightarrow\,\Func.\]
We have explicitly
$(F\oplus G)(V) = F(V)\oplus G(V)$, $(F\otimes G)(V) = F(V)\otimes G(V)$
and $(F\circ G)(V) = F(G(V))$, for all $V\in\LMod{\K}$.
We note that the functor associated to a $\Sigma_*$-module $S(M)$ is the analogue of a power series.
Therefore, we aim to prove that such functors are preserved by sums $S(M)\oplus S(N)$,
products $S(M)\otimes S(N)$ and composites $S(M\circ N)$,
as in the classical framework.

We introduce monoidal category structures in order to obtain precise results.
We can forget about addition of functors.
We have observed that the composition product $\circ: \Func\times\Func\,\rightarrow\,\Func$
provides the category of functors with the structure of a monoidal category.
The unit is the identity functor $I(V) = V$.
Similarly, the tensor product $\otimes: \Func\times\Func\,\rightarrow\,\Func$
provides the category of functors with the structure of a symmetric monoidal category.
The unit is the constant functor $\one(V)\equiv\K$.
The symmetry isomorphism $c(F,G): F\otimes G\,\rightarrow\,G\otimes F$
is given by the symmetry isomorphism
of the tensor product of $\K$-modules
$c(F(V),G(V)): F(V)\otimes G(V)\,\rightarrow\,G(V)\otimes F(V)$.
We introduce such structures in the category of $\Sigma_*$-modules:

\begin{prp}\label{MonoidalFunc}
There are operations
\[\oplus: \Coef\times\Coef\,\rightarrow\,\Coef,
\quad\otimes: \Coef\times\Coef\,\rightarrow\,\Coef
\quad\text{and}\quad\circ: \Coef\times\Coef\,\rightarrow\,\Coef\]
such that $S(M\oplus N) = S(M)\oplus S(N)$, $S(M\otimes N) = S(M)\otimes S(N)$
and $S(M\circ N) = S(M)\circ S(N)$.
To be more precise:

\textbf{a)} We have an associative and symmetric functor $\otimes: \Coef\times\Coef\,\rightarrow\,\Coef$
(the tensor product of $\Sigma_*$-modules)
which provides the category of $\Sigma_*$-modules
with the structure of a symmetric monoidal category.
The unit is the $\Sigma_*$-module $\one\in\Coef$ such that
\[\one(r) = \begin{cases} \K, & \text{if $r = 0$}, \\
0, & \text{otherwise}. \end{cases}\]
The map $M\mapsto S(M)$ defines a monoidal functor
$S: (\Coef,\otimes)\,\rightarrow\,(\Func,\otimes)$.

\textbf{b)} We have an associative functor $\circ: \Coef\times\Coef\,\rightarrow\,\Coef$
(the composition product of $\Sigma_*$-modules)
which provides the category of $\Sigma_*$-modules
with the structure of a (non-symmetric) monoidal category.
The unit is the $\Sigma_*$-module $I\in\Coef$ such that
\[I(r) = \begin{cases} \K, & \text{if $r = 1$}, \\
0, & \text{otherwise}. \end{cases}\]
The map $M\mapsto S(M)$ defines a monoidal functor
$S: (\Coef,\circ)\,\rightarrow\,(\Func,\circ)$.
\end{prp}

We make explicit the expansion of the tensor product of $\Sigma_*$-modules in section~\ref{OperationsExpansions}
(as well as the expansion of the composition product).

\subsection{The $\Sigma_*$-module associated to a functor}\label{CrossEffectModule}
We observe that a $\Sigma_*$-module $M$ is determined by the associated functor $S(M)$.
To be precise,
given a functor $F\in\Func$,
we consider the $\Sigma_*$-module $C(F)\in\Coef$,
such that
\[C(F)(r) = \Hom_{\Func}(T^{(r)},F),\]
where $T^{(r)}: \LMod{\K}\,\rightarrow\,\LMod{\K}$
denotes the power functor $T^{(r)}(V) = V^{\otimes r}$.
The module $C(F)(r)$ is a representation of the symmetric group $\Sigma_r$,
since a permutation $w\in\Sigma_r$ determines a morphism
$w^*: T^{(r)}\,\rightarrow\,T^{(r)}$.
Accordingly, the map $F\mapsto C(F)$ gives a functor $C: \Func\,\rightarrow\,\Coef$.
One proves the following result (\emph{cf}. I. Macdonald \cite{Macdonald}):

\begin{prp}\label{MonoidalCoef}
The functor $C: \Func\,\rightarrow\,\Coef$
is left adjoint to $S: \Coef\,\rightarrow\,\Func$.
Moreover, if the ground ring $\K$ is an infinite field,
then the adjunction unit $\eta(M): M\,\rightarrow C(S(M))$
is an isomorphism.
\end{prp}

\subsection{Remarks}\label{PolynomialCrossEffect}
The first assertion is straighforward.
The last assertion does not hold if the ground ring $\K$ is not an infinite field.
For instance,
the calculations of P. Kraso\'n and N. Kuhn \cite{KrasonKuhn} prove that the $\Sigma_r$-modules $C(S)(r)$
which are associated to the symmetric algebra
\[S(V) = \bigoplus_{r=0}^{\infty} (V^{\otimes r})_{\Sigma_r}\]
differ from the trivial representations.
In fact,
in the case $\K = \F_p$,
the $\Sigma_r$-modules $C(F)(r)$ are connected
to the \emph{polynomial filtration} of $F$
(\emph{cf}. H.-W. Henn, J. Lannes and L. Schwartz \cite{HennLannesSchwartz}, N. Kuhn \cite{Kuhn}).
Precisely,
the $\Sigma_r$-module $C(F)(r)$ is given by the $r$th \emph{cross effect}
of the $r$th polynomial filtration layer $p_r F\subset F$
(\emph{cf}. L. Piriou \cite{Piriou} and L. Piriou and L. Schwartz \cite{PiriouSchwartz},
see also V. Franjou and L. Schwartz \cite{FranjouSchwartz}).

Therefore,
in order to obtain a more general result,
one has to introduce the category of \emph{strict polynomial functors} $\Pol$
(\emph{cf}. E. Friedlander and A. Suslin \cite{FriedlanderSuslin})
and to consider the $\Sigma_*$-module
such that
\[C(F)(r) = \Hom_{\Pol}(T^{(r)},F).\]
This $\Sigma_*$-module is connected to the top \emph{homogeneous cross effects} of $F$
(\emph{cf}. A.K. Bousfield \cite{BousfieldFunc, BousfieldOp})
and
verifies the relation $M = C(S(M))$,
for all ground rings $\K$.

\subsection{Remarks}\label{PolAdjoint}
The adjunction morphism $\epsilon(F): S(C(F))\,\rightarrow\,F$
is supplied by the evaluation morphism
\[\Hom_\F(T^{(r)},F)\otimes V^{\otimes r}\,\rightarrow\,F(V).\]
In positive characteristic,
there are (polynomial) functors $F\in\F$
such that $\epsilon: S(C(F))\,\rightarrow\,F$
is not an isomorphism.
Here is a significant example.
We replace coinvariants by invariants in the definition
of the functor $S_r(M,V) = (M(r)\otimes V^{\otimes r})_{\Sigma_r}$.
We consider the functor
$\Gamma_r(M): \LMod{\K}\,\rightarrow\,\LMod{\K}$
such that
$\Gamma_r(M,V) = (M(r)\otimes V^{\otimes r})^{\Sigma_r}$.
We have a canonical isomorphism $M(r)\,\rightarrow\,C(\Gamma_r(M))(r)$.
One may prove that the adjunction augmentation
$\epsilon(\Gamma_r(M)): S_r(M)\,\rightarrow\,\Gamma_r(M)$
is identified with the norm map,
which is not an isomorphism in general.

In characteristic $0$,
all polynomials functors have the form $F = S(M)$
and
such phenomena do not occur (\emph{cf}. I. Macdonald \cite{Macdonald}).

\smallskip
We mention the following properties of the cross-effect functor $C: \Func\,\rightarrow\,\Coef$:

\begin{prp}
We have the relations $C(F\oplus G) = C(F)\oplus C(G)$ and $C(F\otimes G) = C(F)\otimes C(G)$.
In addition, if the ground ring $\K$ is an infinite field,
then we have the relation $C(F\circ G) = C(F)\circ C(G)$.
\end{prp}

We just refer to L. Piriou \cite{Piriou} for a proof of the relation $C(F\otimes G) = C(F)\otimes C(G)$.
The relation $C(F\circ G) = C(F)\circ C(G)$ holds in full generality
if we consider the homogeneous cross effect
$C(F)(r) = \Hom_{\Pol}(T^{(r)},F)$.

\subsection{Operads}\label{MonoidOperads}
We define an operad $P$ to be a monoid
in the category of $\Sigma_*$-modules $\Coef$
together with the composition product $\circ: \Coef\times\Coef\,\rightarrow\,\Coef$.
According to this definition,
an operad consists of a $\Sigma_*$-module $P\in\Coef$
equipped with a product $\mu: P\circ P\,\rightarrow\,P$ and a unit $\eta: I\,\rightarrow\,P$
that verify the classical associativity and unit relations,
as in the definition of a monad structure (see paragraph~\ref{MonoidMonads}).

We observe that this definition,
borrowed from \cite[E. Getzler and J. Jones]{GetzlerJones} and \cite[V. Smirnov]{Smirnov},
agrees with the classical one,
because of the following result
which is an immediate consequence of proposition~\ref{MonoidalFunc}:

\begin{prp}\label{SymmetricMonads}
If $P$ is an operad (according to the definition of paragraph~\ref{MonoidOperads}),
then the functor $S(P)$ is a monad.
\end{prp}

Explicitly,
the operad product $\mu: P\circ P\,\rightarrow\,P$ yields a morphism
\[S(P)\circ S(P) = S(P\circ P)\,\xrightarrow{\mu}\,S(P)\]
that gives the functor $S(P)$ the structure of a monad.
Similarly, we deduce the next result from proposition~\ref{MonoidalCoef}:

\begin{prp}
We assume that the ground ring $\K$ is an infinite field.
If $S$ is a monad, then the $\Sigma_*$-module $C(S)$ is an operad.
\end{prp}

\subsection{Digression: divided power functors}\label{DividedPowerFunctors}
We observe that an operad $P$ determines monads that differ from the classical one $S(P)$.
For this purpose,
we change the definition of the functor associated to a $\Sigma_*$-module.

Suppose given a $\Sigma_*$-module $M$.
As in paragraph~\ref{PolAdjoint},
we replace coinvariants by invariants in the definition
of the functor
\[S_r(M,V) = (M(r)\otimes V^{\otimes r})_{\Sigma_r}.\]
Explicitly,
we consider the functor $\Gamma(M,V) = \bigoplus_{r=0}^{\infty}\Gamma_r(M,V)$,
such that
\[\Gamma_r(M,V) = (M(r)\otimes V^{\otimes r})^{\Sigma_r}.\]
We consider also the functor $\Lambda(M,V) = \bigoplus_{r=0}^{\infty}\Lambda_r(M,V)$,
defined by the image of the norm maps
$(M(r)\otimes V^{\otimes r})_{\Sigma_r}\,\rightarrow\,(M(r)\otimes V^{\otimes r})^{\Sigma_r}$.
Consequently,
we have functor morphisms
\[S(M,V)\,\rightarrow\,\Lambda(M,V)\,\rightarrow\,\Gamma(M,V)\]
that are not isomorphisms in general.

Nevertheless, we prove in \cite{FresseSimpAlg}) that the following assertion holds:

\begin{lem}\label{DividedPowerComposite}
We assume that $N$ is a connected $\Sigma_*$-module.
In this situation, the composition product of $\Sigma_*$-modules,
which is characterized by the relation $S(M\circ N) = S(M)\circ S(N)$,
satisfies also the identities
$\Lambda(M\circ N) = \Lambda(M)\circ\Lambda(N)$
and
$\Gamma(M\circ N) = \Gamma(M)\circ\Gamma(N)$.
\end{lem}

The property above does not hold if the $\Sigma_*$-module $N$ is not connected.
For instance, if we assume that $N$ is a constant $\Sigma_*$-module
(explicitly, we assume $N(0) = C$ and $N(r) = 0$ for $r>0$),
then the associated functors are constant
$S(N,V) = \Lambda(N,V) = \Gamma(N,V)\equiv C$.
Consequently,
we obtain $S(M,S(N,V))\equiv S(M,C)$.
Clearly,
this functor his associated to a constant $\Sigma_*$-module $M\circ N$,
which has $M\circ N(0) = S(M,C)$ and $M\circ N(r) = 0$ for $r>0$.
We obtain then $\Lambda(M\circ N,V)\equiv S(M,C)$ and $\Gamma(M\circ N,V)\equiv S(M,C)$.
But,
on the other hand,
we obtain $\Lambda(M,\Lambda(N,V))\equiv\Lambda(M,C)$ and $\Gamma(M,\Gamma(N,V))\equiv\Gamma(M,C)$.
Our remark follows.

\smallskip
The lemma above implies the following assertions:

\begin{prp}\label{DividedPowerMonads}
If $P$ is a connected operad,
then the functor $\Lambda(P)$ (respectively, $\Gamma(P)$)
is equipped with the structure of a monad.
Furthermore,
the morphisms $S(P)\,\rightarrow\,\Lambda(P)\,\rightarrow\,\Gamma(P)$ are monad morphisms.
\end{prp}

As an example,
we make explicit the structure of a $\Lambda(P)$-algebra
(respectively, of a $\Gamma(P)$-algebra)
in the case of the classical operads $P = \C$ and $P = \L$.
The next result is classical (\emph{cf}. H. Cartan \cite{Cartan}, N. Roby \cite{Roby}):

\begin{prp}\label{DividedCommutativeAlgebras}
We assume that the ground ring $\K$ is a field of characteristic $p>0$.
Recall that an algebra over the monad $S(\C)$
is an associative and commutative algebra
in the classical sense.
An algebra over the monad $\Lambda(\C)$
is an associative and commutative algebra
such that $x^p = 0$, for all elements $x$.
An algebra over the monad $\Gamma(\C)$
is a divided power algebra.
\end{prp}

The next result is proved in \cite{FresseSimpAlg}:

\begin{prp}\label{DividedLieAlgebras}
We assume that the ground ring $\K$ is a field of characteristic $p>0$.
Recall that an algebra over the monad $S(\L)$
is a Lie algebra $\G$ together with a Lie bracket $[-,-]: \G\otimes\G\,\rightarrow\,\G$
such that $[x,y] = - [y,x]$, for all $x,y\in\G$.
An algebra over the monad $\Lambda(\L)$
is a Lie algebra $\G$ together with a Lie bracket $[-,-]: \G\otimes\G\,\rightarrow\,\G$
such that $[x,x] = 0$.
An algebra over the monad $\Gamma(\L)$
is a restricted Lie algebra $\G$.
\end{prp}

\subsection{Cooperads}\label{Cooperads}
The structure of a \emph{cooperad} is essentially dual to the structure of an operad.
We consider such structures in chapters~\ref{ReducedBar} and~\ref{KoszulDuality}.

To be precise,
we define a cooperad $D$ to be a $\Sigma_*$-module
equipped with a coproduct $\nu: D\,\rightarrow\,D\circ D$
and with a counit $\epsilon: D\,\rightarrow\,I$
that verify the classical coassociativity and counit relations.

We assume in general that $D$ is a connected $\Sigma_*$-module.
In this case,
we deduce from lemma~\ref{DividedPowerComposite}
that the functor
\[\Gamma(D,V) = \bigoplus_{r=0}^{\infty}\bigl(D(r)\otimes V^{\otimes r}\bigr)^{\Sigma_r}\]
is equipped with the structure of a comonad.
Explicitly, the cooperad coproduct induces a morphism
\[\Gamma(D)\,\xrightarrow{\nu}\,\Gamma(D\circ D)\simeq\Gamma(D)\circ\Gamma(D).\]
We define a \emph{$D$-coalgebra} to be a coalgebra over the monad $\Gamma(D)$.
Explicitly,
a $D$-coalgebra is a $\K$-module $B$ together with a coproduct $\rho: B\,\rightarrow\,\Gamma(D,B)$
that verify natural coassociativity and counit relations.

\smallskip
We consider the dual $\Sigma_*$-module of a cooperad $D^{\vee}$ (\emph{cf}. paragraph~\ref{SymmetricFunctors}).
We have the following results:

\begin{prp}\label{CoalgebraDuality}
If $D$ is a cooperad, then the dual sequence $D(r)^{\vee}$ is equipped with the structure of an operad.
Moreover, the dual $\K$-module of a $D$-coalgebra $B^{\vee}$
is equipped with the structure of a $D^{\vee}$-algebra.

Dually, we are given a connected operad $P$.
We assume that the sequence $P(r)$ consists of finitely generated projective $\K$-modules.
In this case, the dual sequence $P(r)^{\vee}$ is equipped with the structure of a cooperad.
Moreover, if $A$ is a finitely generated projective $\K$-module,
then the dual $\K$-module $A^{\vee}$ is equipped with the structure of a $P^{\vee}$-coalgebra.
\end{prp}

As an example, a coalgebra over the dual cooperad of the commutative operad $\C^{\vee}$
is nothing but a nilpotent cocommutative and coassociative coalgebra $B$.
(The iterated coproduct $\Delta^{r}(b)\in B^{\otimes r}$ of a given element $b\in B$
is supposed to vanish for large $r\in\N$.)

The results of proposition~\ref{CoalgebraDuality} are consequences of the following lemma
(proved in section~\ref{OperationsExpansions}):

\begin{lem}\label{DualCompositeModule}
We are given $\Sigma_*$-modules $M$ and $N$.
We have a natural morphism $M^{\vee}\circ N^{\vee}\,\rightarrow\,(M\circ N)^{\vee}$,
which is an isomorphism provided that $N$ is connected
and the sequences $M(r)$ and $N(r)$
consist of finitely generated projective $\K$-modules.
\end{lem}

We deduce these assertions from lemma~\ref{ExpansionComposite}
about the expansion of the composition product of $\Sigma_*$-modules.

\section{The composition product of symmetric modules}\label{OperationsExpansions}

We define explicitly the composition product of $\Sigma_*$-modules in this section.
We deduce the nice properties of this operation from the observation~\ref{PartitionsPermutation}
and from the expansion given in lemma~\ref{ExpansionComposite}.

\subsection{Indexing by finite sets}
Suppose given a $\Sigma_*$-module $M$.
As in paragraph~\ref{OperadIndexing},
we define a functor on the category of finite sets and bijections $I\mapsto M(I)$
such that $M(r) = M(\{1,\ldots,r\})$.

We consider the module $M(I) = \bigoplus_{i_*: \{1,\ldots,r\}\,\rightarrow\,I} M(r)/\equiv$
(as in the construction of paragraph~\ref{OperadIndexing}).
The sum ranges over bijections $i_*: \{1,\ldots,r\}\,\rightarrow\,I$.
Hence,
an element of $M(I)$ consists of a pair $(x,i_*)$,
where $x\in M(r)$ and $i_*: \{1,\ldots,r\}\,\rightarrow\,I$.
We set $(w_*(x),i_{*})\equiv (x,i_{w(*)})$ for any permutation $w\in\Sigma_r$.
We recall that $i_{w(*)}: \{1,\ldots,r\}\,\rightarrow\,I$
denotes the composite of the bijection $i_*: \{1,\ldots,r\}\,\rightarrow\,I$
with the permutation $w: \{1,\ldots,r\}\,\rightarrow\,\{1,\ldots,r\}$.

A bijection $u: I\,\rightarrow\,I'$ induces a reindexing morphism $u_*: M(I)\,\rightarrow\,M(I')$.
Accordingly, the collection $M(I)$ defines a functor on the category of finite sets and bijections.
We have clearly the relation $M(r) = M(\{1,\ldots,r\})$.
Furthermore, the morphisms
$w_*: M(\{1,\ldots,r\})\,\rightarrow\,M(\{1,\ldots,r\})$
associated to permutations
$w: \{1,\ldots,r\}\,\rightarrow\,\{1,\ldots,r\}$
give the action of the symmetric group on $M(r)$.

\subsection{The coend definition of the functor associated to a $\Sigma_*$-module}
Similarly,
for a given $\K$-module $V$,
the tensor powers $r\mapsto V^{\otimes r}$ extends to a contravariant functor $I\mapsto V^{\otimes I}$
on the category of finite sets and bijections.
To be explicit,
an element of $V^{\otimes I}$
is equivalent to a tensor $v_{i_1}\otimes\cdots\otimes v_{i_r}\in V\otimes\cdots\otimes V$.
The morphism
$u^*: V^{\otimes I'}\,\rightarrow\,V^{\otimes I}$
induced by a bijection
$u: I\,\rightarrow\,I'$
verifies the relation
$u^*(v_{i'_1}\otimes\cdots\otimes v_{i'_r}) = w_{i_1}\otimes\cdots\otimes w_{i_r}$,
where $w_{i_k} = v_{u(i_k)}$.

The functor $S(M): \LMod{\K}\,\rightarrow\,\LMod{\K}$
agrees with the coend
\[S(M,V) = \int_I M(I)\otimes V^{\otimes I}.\]
The elements of this coend are represented by tensors
$x\otimes(v_{i_1}\otimes\cdots\otimes v_{i_r})\in M(I)\otimes V^{\otimes I}$,
together with relations
\[u_*(x)\otimes(v_{i'_1}\otimes\cdots\otimes v_{i'_r})
\equiv x\otimes u^*(v_{i'_1}\otimes\cdots\otimes v_{i'_r}),\]
which generalize the invariance relations of the module of coinvariants.

\subsection{The tensor product of $\Sigma_*$-modules}\label{SymmetricTensorProduct}
We define the tensor product $M\otimes N$ of $\Sigma_*$-modules $M$ and $N$.
We consider the functors on finite sets associated to $M$ and $N$.
We set
\[(M\otimes N)(n) = \bigoplus_{(I,J)} M(I)\otimes N(J).\]
The sum ranges over ordered pairs of sets $(I,J)$ such that $I\amalg J = \{1,\ldots,n\}$.

We have a symmetry isomorphism
\[c(M,N): M\otimes N\,\rightarrow\,N\otimes M\]
defined componentwise by the symmetry isomorphism of $\K$-modules
$M(I)\otimes N(J)\,\rightarrow\,N(J)\otimes M(I)$.
Clearly, the $\Sigma_*$-module $\one\in\Coef$
such that
\[\one(r) = \begin{cases} \K, & \text{if $r = 0$}, \\
0, & \text{otherwise}, \end{cases}\]
satifies the unit relations $M\otimes\one = M = \one\otimes M$.

The next result is a straighforward consequence of the coend definition
of the functor associated to a $\Sigma_*$-module.

\begin{lem}
The $\Sigma_*$-module $M\otimes N$ defined in the paragraph above,
verifies the relation $S(M\otimes N) = S(M)\otimes S(N)$.
Furthermore,
the symmetry isomorphism of $\Sigma_*$-modules
$c(M,N): M\otimes N\,\rightarrow\,N\otimes M$
corresponds to the symmetry isomorphism of functors
$c(M,N): S(M)\otimes S(N)\,\rightarrow\,S(N)\otimes S(M)$.
Explicitly,
the diagram
\[\xymatrix{ S(M\otimes N)\ar[rr]^{c(M,N)}\ar[d]_{=} & & S(N\otimes M)\ar[d]^{=} \\
S(M)\otimes S(N)\ar[rr]^{c(S(M),S(N))} & & S(N)\otimes S(M) \\ }\]
is commutative.
\end{lem}

\subsection{The composition product of $\Sigma_*$-modules}\label{CompositionProductDefinition}
We define of the composition product $M\circ N$ of $\Sigma_*$-modules $M$ and $N$.
We give E. Getzler, J. Jones and V. Smirnov credit
for the introduction of the composition product in the context of operads
(\emph{cf}. E. Getzler and J. Jones \cite{GetzlerJones}, V. Smirnov \cite{Smirnov}),
but this operation is also related to the \emph{plethysm}
in combinatorial representation theory
(\emph{cf}. I. Macdonald \cite{Macdonald}, see also A. Joyal \cite{Joyal}).

We follow the construction of \cite{FresseLie} and \cite{FresseSimpAlg}.
We observe that the tensor power of a $\Sigma_*$-module $N^{\otimes r}$
is endowed with an action of the symmetric group $\Sigma_r$,
because the tensor product of $\Sigma_*$-modules has a symmetry isomorphism.
Accordingly,
the $\Sigma_n$-modules $N^{\otimes r}(n)$
have an external $\Sigma_r$-action.
We set
\[(M\circ N)(n) = \bigoplus_{r=0}^\infty\bigl(M(r)\otimes N^{\otimes r}(n)\bigr)_{\Sigma_r}.\]
The next assertion is a straighforward consequence of this definition:

\begin{lem}
The $\Sigma_*$-module $M\circ N$ defined in the paragraph above
verifies the relation $S(M\circ N) = S(M)\circ S(N)$.
\end{lem}

Explicitly, we have a sequence of immediate identities:
\begin{align*} S(M,S(N,V))
& = \bigoplus_{r=0}^{\infty}\Bigl(M(r)\otimes S(N,V)^{\otimes r}\Bigr)_{\Sigma_r} \\
& = \bigoplus_{r=0}^{\infty}\Bigl(M(r)\otimes S(N^{\otimes r},V)\Bigr)_{\Sigma_r} \\
& = S\Bigl(\bigoplus_{r=0}^{\infty}\bigl(M(r)\otimes N^{\otimes r}\bigr)_{\Sigma_r},V\Bigr) = S(M\circ N,V).
\end{align*}

\subsection{On the tensor power of a $\Sigma_*$-module}
We need more insights in the structure of a composite $\Sigma_*$-module $M\circ N$.
Therefore, for our purposes,
we make the expansion of a tensor power $N^{\otimes r}$ more explicit.
We have readily
\[N^{\otimes r}(\{1,\ldots,n\}) = \bigoplus_{(J_1,\ldots,J_r)} N(J_1)\otimes\cdots\otimes N(J_r).\]
This sum ranges over $r$-tuples of sets $(J_1,\ldots,J_r)$ such that $J_1\amalg\cdots\amalg J_r = \{1,\ldots,n\}$.
Suppose given a permutation $w\in\Sigma_r$.
We observe that the associated morphism
$w^*: N^{\otimes r}\,\rightarrow\,N^{\otimes r}$
permutes the summands
$N(J_1)\otimes\cdots\otimes N(J_r)$
in the expansion of $N^{\otimes r}$.
Precisely,
a tensor
$y_1\otimes\cdots\otimes y_r\in N(J_1)\otimes\cdots\otimes N(J_r)$
is mapped to
$y_{w(1)}\otimes\cdots\otimes y_{w(r)}\in N(J_{w(1)})\otimes\cdots\otimes N(J_{w(r)})$.

If $N$ is a connected $\Sigma_*$-module,
then this action has nice properties,
because of the following classical result:

\begin{obv}\label{PartitionsPermutation}
The symmetric group $\Sigma_r$ acts freely
on ordered partitions $J_1\amalg\ldots\amalg J_r = \{1,\ldots,n\}$
such that $J_1,\ldots,J_r\not=\emptyset$.
\end{obv}

The discussion above leads to the following conclusion:

\begin{lem}\label{ExpansionComposite}
Let $M,N$ be $\Sigma_*$-modules.
In general,
the components of a composite $\Sigma_*$-module $M\circ N$
are given by the coinvariant formula
\[\bigl(M(r)\otimes N^{\otimes r}(n)\bigr)_{\Sigma_r}
= \Bigl\{\bigoplus_{(J_1,\ldots,J_r)} M(r)\otimes N(J_1)\otimes\cdots\otimes N(J_r)\Bigr\}_{\Sigma_r},\]
where we consider all $r$-tuples of sets $(J_1,\ldots,J_r)$
such that $J_1\amalg\cdots\amalg J_r = \{1,\ldots,n\}$.
If $N$ is connected, then these modules have an expansion
\[\bigl(M(r)\otimes N^{\otimes r}(n)\bigr)_{\Sigma_r}
= \bigoplus_{(J_1,\ldots,J_r)'} M(r)\otimes N(J_1)\otimes\cdots\otimes N(J_r),\]
which does not involve coinvariants,
and
which can be obtained by restricting the summation range to $r$-tuples of non-empty sets
$J_k = \{j_{k 1}<\cdots<j_{k n_k}\}$, $k = 1,\ldots,r$,
such that $j_{1 1}<j_{2 1}<\cdots<j_{r 1}$.
\end{lem}

\subsection{Composite elements}\label{CompositeElements}
Suppose given $x\in M(r)$ and $y_1\in N(n_1),\ldots,y_r\in N(n_r)$.
We let $x(y_1,\ldots,y_r)\in(M\circ N)(n)$ denote the element
defined by the tensor
$x\otimes(y_1\otimes\cdots\otimes y_r)\in M(r)\otimes(N(n_1)\otimes\cdots\otimes N(n_r))$
in the component of $(M\circ N)(n)$
associated to the canonical partition $J_1\amalg\cdots\amalg J_r = \{1,\ldots,n\}$,
where $J_k = \{n_1+\cdots+n_{k-1}+1<\cdots<n_1+\cdots+n_{k-1}+n_k\}$,
for $k = 1,\ldots,r$.

We observe that a tensor $x(y_1,\ldots,y_r)(v_1,\ldots,v_n)\in S(M\circ N,V)$
associated to such an element
is identified with a composite
\begin{align*} & x(w_1,\ldots,w_r)\in S(M,S(N,V)) \\
\text{where}\qquad & w_k = y_k(v_{n_1+\cdots+n_{k-1}+1},\ldots,v_{n_1+\cdots+n_{k-1}+n_k})\in S(N,V).
\end{align*}
Therefore,
the image under an operad product $\mu: P\circ P\,\rightarrow\,P$
of such a tensor $p(q_1,\ldots,q_r)\in (P\circ P)(n)$
agrees with the element denoted by $p(q_1,\ldots,q_r)\in P(n)$
in paragraph~\ref{ClassicalOperadDefinition}.

\begin{proof}[Proof of lemma~\ref{DualCompositeModule}]
As we assume that the sequences $M(r)$ and $N(r)$ consist of finitely generated projective $\K$-modules,
we have a natural isomorphism
\begin{multline*}\Bigl\{\bigoplus_{(J_*)} M(r)^{\vee}\otimes
N(J_1)^{\vee}\otimes\cdots\otimes N(J_r)^{\vee}\Bigr\}^{\Sigma_r} \\
\xrightarrow{\simeq}\,\Bigl\{\bigoplus_{(J_*)} M(r)\otimes
N(J_1)\otimes\cdots\otimes N(J_r)\Bigr\}_{\Sigma_r}^{\vee}.\end{multline*}
Furthermore, as long as $N$ is a connected $\Sigma_*$-module,
the norm map provides a natural isomorphism
\begin{multline*}\Bigl\{\bigoplus_{(J_*)} M(r)^{\vee}\otimes
N(J_1)^{\vee}\otimes\cdots\otimes N(J_r)^{\vee}\Bigr\}_{\Sigma_r} \\
\xrightarrow{\simeq}\,\Bigl\{\bigoplus_{(J_*)} M(r)^{\vee}\otimes
N(J_1)^{\vee}\otimes\cdots\otimes N(J_r)^{\vee}\Bigr\}^{\Sigma_r},\end{multline*}
because the symmetric group permutes freely the partitions $(J_*)$ involved in the sum.
Hence,
we obtain an isomorphism
$M^{\vee}\circ N^{\vee}\,\xrightarrow{\simeq}\,(M\circ N)^{\vee}$
and the conclusion of lemma~\ref{DualCompositeModule} follows.
\end{proof}

%\input{ModuleComplexes}
% 30/1/2003

\chapter{Chain complexes of modules over an operad}\label{ModuleComplexes}

\section{Summary}\label{SummaryModuleComplexes}

In this section,
we introduce the notion of a module over an operad.
In fact,
our aim is to develop classical homological algebra techniques
in the framework of operads.
Therefore,
we first generalize the definitions of chapter~\ref{Operads}
to the differential graded context.

\subsection{Chain complexes of $\Sigma_*$-module}\label{DGSigmaModule}
We consider $\Sigma_*$-objects in the category of dg-modules.
We assume precisely that $M = M(r)$, $r\in\N$, is a sequence of dg-modules
and that the symmetric group $\Sigma_r$ operates on $M(r)$
by morphisms of dg-modules.
If the ground ring $\K$ is not a field,
then we may consider $\Sigma_*$-modules $M$
which are \emph{projective as $\K$-modules}.
(More explicitly, we may assume that the sequence $M(r)$ consists of projective $\K$-modules.)
If the ground ring $\K$ is not a field of characteristic $0$,
then we may consider $\Sigma_*$-modules $M$
which are \emph{projective in the category of $\Sigma_*$-modules}.
(Explicitly,
a $\Sigma_*$-module $M$ is projective
if the sequence $M(r)$ consists of projective $\Sigma_r$-modules.)
We recall that a $\Sigma_*$-module $M$ is connected if $M(0) = 0$.

We observe that the definition of the direct sum $M\oplus N$,
of the tensor product $M\otimes N$
and of the composition product $M\circ N$
makes sense in any symmetric monoidal category
and, in particular, in the category of dg-$\Sigma_*$-modules
(see section~\ref{DGMod}).
For our purposes,
we just record that the differential of a composite element $x(y_1,\ldots,y_r)\in M\circ N(n)$,
where $x\in M(r)$ and $y_1\in N(n_1),\ldots,y_r\in N(n_r)$
(\emph{cf}. paragraph~\ref{CompositeElements}),
is given by the formula
\[\delta(x(y_1,\ldots,y_r)) = \delta(x)(y_1,\ldots,y_r)
+ \sum_{i=1}^r \pm x(y_1,\ldots,\delta(y_i),\ldots,y_r).\]

\subsection{Classical homological problems arising from functors associated to dg-$\Sigma_*$-modules}
We have also a functor
\[S(M): \dg\LMod{\K}\,\rightarrow\,\dg\LMod{\K}\]
associated to any dg-$\Sigma_*$-module $M$,
since the formula of paragraph~\ref{SymmetricFunctors}
\[S_r(M)(V) = (M(r)\otimes V^{\otimes r})_{\Sigma_r}\]
makes sense in the category of dg-modules.
Similarly,
the relations
\begin{align*} & S(M\oplus N,V) = S(M,V)\oplus S(N,V), \\
& S(M\otimes N,V) = S(M,V)\otimes S(N,V) \\
& \text{and}\qquad S(M\circ N,V) = S(M,S(N,V)), \end{align*}
hold in any symmetric monoidal category
and, in particular, in the category of dg-modules $\dg\LMod{\K}$.

Classically,
one observes that the functor $S(M): \dg\LMod{\K}\,\rightarrow\,\dg\LMod{\K}$
does not preserve quasi-isomorphisms,
unless $M$ is a projective $\Sigma_*$-module.
Similarly,
the functor $M\mapsto S(M)$ is additive and right exact,
but not exact.

We observe that the composition product of dg-$\Sigma_*$-modules
has nice homological properties despite these problems,
because a composite $\Sigma_*$-module $M\circ N$
such that $N$ is connected
has a good expansion (\emph{cf}. lemma~\ref{ExpansionComposite}).
In particular,
if the ground ring $\K$ is a field,
then we obtain the following result:

\begin{lem}\label{KunnethComposition}
We assume that $M$ and $N$ are dg-$\Sigma_*$-modules over a field $\K$ of characteristic $p\geq 0$.
If the $\Sigma_*$-module $N$ is connected (explicitly, if $N(0) = 0$),
then we have a natural isomorphism $H_*(M)\circ H_*(N)\,\simeq\,H_*(M\circ N)$.
\end{lem}

We analyze the differential structure of a composite $\Sigma_*$-module $M\circ N$
in detail in section~\ref{DGComposite}.

\subsection{Differential graded operads}\label{DGOperads}
Finally,
the notion an operad makes sense in the category of dg-modules.
In this context,
the composition product of a dg-operad $P$
is supposed to be a morphism of dg-$\Sigma_*$-modules
$P\circ P\,\rightarrow\,P$.
Equivalently,
we assume that a composite operation $p(q_1,\ldots,q_r)\in P(n_1+\cdots+n_r)$
satisfies the derivation relation
\[\delta(p(q_1,\ldots,q_r)) = \delta(p)(q_1,\ldots,q_r)
+ \sum_{i=1}^r \pm p(q_1,\ldots,\delta(q_i),\ldots,q_r).\]
We recall that an operad is connected if $P(0) = 0$ and $P(1) = \K\,1$.

Dually,
we have a notion of a dg-cooperad,
in which the coproduct is a morphism of dg-modules $D\,\rightarrow\,D\circ D$.

We observe that the homology of a dg-operad $H_*(P)$ is equipped with the structure of a graded operad.
If the ground ring $\K$ is a field, then the same is true for dg-cooperads,
because lemma~\ref{KunnethComposition} above allows to define a morphism
\[H_*(D)\,\rightarrow\,H_*(D\circ D)\,\simeq\,H_*(D)\circ H_*(D).\]

\subsection{Right modules over an operad}\label{RightModules}
A \emph{right module} over an operad $P$ is a $\Sigma_*$-module $L\in\LMod{\Sigma_*}$
equipped with a morphism
$\lambda: L\circ P\,\rightarrow\,L$
that satisfies the classical associativity relation with respect to the operad product
$\mu: P\circ P\,\rightarrow\,P$
and the classical unit relation with respect to the monad unit
$\eta: I\,\rightarrow\,P$.
Explicitly,
we assume that the following classical diagrams commute:
\begin{equation*}\begin{aligned}
\xymatrix{ L\circ P\circ P\ar[r]^{\lambda\circ P}\ar[d]_{L\circ\mu} &
L\circ P\ar[d]^{\lambda} \\ L\circ P\ar[r]^{\lambda} & L \\ }
\end{aligned}\qquad\text{and}
\qquad\begin{aligned}
\xymatrix{ L\circ P\ar[d]_{\lambda} & L\circ I\ar[l]_{L\circ\eta}\ar[dl]^{=} \\ L & & \\ }
\end{aligned}\end{equation*}
The category of right $P$-modules is denoted by $\RMod{P}$.

Like an operad product,
an operad action $\lambda: L\circ P\,\rightarrow\,L$
is equivalent to composition products
\[L(r)\otimes P(n_1)\otimes\cdots\otimes P(n_r)\,\rightarrow\,L(n_1+\cdots+n_r),\]
defined for all $r\geq 1$ and all $n_1,\ldots,n_r\in\N$.
The composite element
associated to $x\in L(r)$ and $q_1\in P(n_1),\ldots,q_r\in P(n_r)$
is denoted by $x(q_1,\ldots,q_r)\in L(n_1+\cdots+n_r)$.
There is also a partial composition product such that $x\circ_i q = x(1,\ldots,q,\ldots,1)$.

In the dg-context,
we assume that the operad action is given by a morphism of dg-modules
$\lambda: L\circ P\,\rightarrow\,L$.
Accordingly,
a composite element $x(q_1,\ldots,q_r)\in L(n_1+\cdots+n_r)$
satisfies the same derivation relation as a composite operation
in a dg-operad.

\subsection{Left modules over an operad}\label{LeftModules}
Symmetrically,
a \emph{left module} over an operad $P$
is a $\Sigma_*$-module $R\in\LMod{\Sigma_*}$
equipped with a left operad action
$\rho: P\circ R\,\rightarrow\,R$.
We observe that a $P$-algebra $A$
is equivalent to a left $P$-module $R$
such that
\[R(n) = \begin{cases} A, & \text{if $n=0$}, \\
0, & \text{otherwise}. \end{cases}\]
Conversely,
we may consider \emph{connected left $P$-modules} $R$,
which have $R(0) = 0$,
because connected $\Sigma_*$-modules behave better than constant modules
in regard to homological constructions.
The category of left $P$-modules is denoted by $\LMod{P}$.
The category of connected left $P$-modules is denoted by $\LMod{P}_0$.

A left operad action is equivalent to composition products
\[P(r)\otimes R(n_1)\otimes\cdots\otimes R(n_r)\,\rightarrow\,R(n_1+\cdots+n_r)\]
defined for all $r\geq 1$ and all $n_1,\ldots,n_r\in\N$.
The composite element
associated to $p\in P(r)$ and $y_1\in R(n_1),\ldots,y_r\in R(n_r)$
is denoted by $p(y_1,\ldots,y_r)\in R(n_1+\cdots+n_r)$.
In the dg-context,
a composite element $p(y_1,\ldots,y_r)\in R(n_1+\cdots+n_r)$
satisfies the same derivation relation as a composite operation
in a dg-operad.

A \emph{bimodule} over an operad $P$ is a $\Sigma_*$-module $L\in\LMod{\Sigma_*}$
equipped with both a right and left $P$-module structure
such that the right and left operad actions commute.
This notion of a bimodule agrees with the notion
introduced by M. Kapranov and Y. Manin
in \cite{KapranovManin},
but differs from a notion of a bimodule introduced by M. Markl
in \cite{MarklModels}.
In fact,
the $P$-bimodules in Markl's sense are abelian group objects
in the relative category of operads over $P$
(\emph{cf}. \cite{FresseLie})
and form an abelian category unlike the $P$-bimodules
in our sense.
These notions of bimodules are also considered (in a more general framework)
by H. Baues, M. Jibladze and A. Tonks
in \cite{BauesJibladzeTonks}.

\subsection{Relative composition products}\label{RelativeCompositionProducts}
The \emph{relative composition product} $L\circ_P R$
is similar to the classical tensor product
over an associative algebra.
Explicitly,
we are given a right $P$-module $L$ and a left $P$-module $R$.
The relative composition product $L\circ_P R$
is defined by the cokernel
\[\xymatrix{ *+<2mm>{L\circ P\circ R}\ar@<1.5mm>[r]^{\lambda\circ R}\ar@<-1.5mm>[r]_{L\circ\rho} &
*+<2mm>{L\circ R}\ar[r] & *+<2mm>{L\circ_P R.} \\ }\]
We have $P\circ_P R = R$ and $L\circ_P P = L$.
More generally,
in case of a free right $P$-module $L = M\circ P$ (see paragraph~\ref{FreeModules} below),
we obtain $(M\circ P)\circ_P R = M\circ R$.
Similarly, if $R = P\circ M$, then we obtain $L\circ_P (P\circ M) = L\circ M$.

\subsection{Extension and restriction functors}\label{ExtensionRestrictionFunctors}
Suppose given an operad morphism $\phi: P\,\rightarrow\,P'$.
In this situation, we have a \emph{restriction functor} $\phi^!: \RMod{P'}\,\rightarrow\,\RMod{P}$,
since the composite
\[L'\circ P\,\xrightarrow{L'\circ\phi}\,L'\circ P'\,\xrightarrow{\lambda'}\,L'\]
provides any right $P'$-module $L'$ with a right $P$-action.
We have also an \emph{induction functor} $\phi_!: \RMod{P}\,\rightarrow\,\RMod{P'}$
together with an adjunction relation
\[\Hom_{\RMod{P'}}(\phi_! L,L') = \Hom_{\RMod{P}}(L,\phi^! L').\]
Explicitly,
the morphism $\phi: P\,\rightarrow\,P'$ gives the operad $P'$ the structure of a left $P$-module.
One verifies easily that the right $P'$-module $\phi_!(L)$
is represented by the relative composition product
$\phi_!(L) = L\circ_P P'$.
Furthermore,
the adjunction unit $\eta(L): L\,\rightarrow\,\phi^!(\phi_!(L))$
can be identified with the morphism
$L\circ_P\phi: L\circ_P P\,\rightarrow\,\,L\circ_P P'$.

Symmetrically, for left modules over an operad,
we have a restriction functor $\phi^!: \LMod{P'}\,\rightarrow\,\LMod{P}$
and an extension functor $\phi_!: \LMod{P}\,\rightarrow\,\LMod{P'}$
given by the relative composition product $\phi_!(R) = P'\circ_P R$.
These functors verify the adjunction relation
\[\Hom_{\LMod{P'}}(\phi_! R,R') = \Hom_{\LMod{P}}(R,\phi^! R').\]
The adjunction unit $\eta(R): L\,\rightarrow\,\phi^!(\phi_!(R))$
is identified with the morphism $\phi\circ_P R: P\circ_P R\,\rightarrow\,\,P'\circ_P R$.

\subsection{Indecomposable quotients}\label{ModuleIndecomposableQuotient}
Suppose $P$ is an augmented operad,
so that we are given an operad morphism $\epsilon: P\,\rightarrow\,I$.
In this case,
the restriction functor $\epsilon^!: \RMod{I}\,\rightarrow\,\RMod{P}$
provides any $\Sigma_*$-module
with the structure of a right $P$-module,
because any $\Sigma_*$-module is endowed with a right action of the identity operad.
The relative composition product $\bar{L} = L\circ_P I$
is called the \emph{indecomposable quotient} of $L$.
The adjunction relation reads
\[\Hom_{\LMod{\Sigma_*}}(\bar{L},M) = \Hom_{\RMod{P}}(L,M).\]
The adjunction unit is identified with the morphism
$L\circ_P\epsilon: L\circ_P P\,\rightarrow\,\,L\circ_P I$.
This morphism is clearly surjective.

Symmetrically,
the indecomposable quotient of a left $P$-module
is defined by the relative composition product $\bar{R} = I\circ_P R$.
We have an adjunction relation
\[\Hom_{\LMod{\Sigma_*}}(\bar{R},M) = \Hom_{\LMod{P}}(R,M).\]
The adjunction unit is identified with the morphism
$\epsilon\circ_P R: P\circ_P R\,\rightarrow\,\,I\circ_P R$
and is clearly surjective.

\smallskip
The relative composition product of dg-modules $L\circ_P R$ is a dg-$\Sigma_*$-module.
But,
as in the classical context of modules over an associative algebra,
we observe that the functor
$-\circ_P R: \RMod{P}\,\rightarrow\,\LMod{\Sigma_*}$
does not preserve all quasi-isomorphisms of right $P$-modules
$\phi: L\,\xrightarrow{\sim}\,L'$.
Similarly,
the functor $L\circ_P -: \LMod{P}\,\rightarrow\,\LMod{\Sigma_*}$,
associated to a right $P$-module $L$,
does not preserve all quasi-isomorphisms of left $P$-modules
$\psi: R\,\xrightarrow{\sim}\,R'$.
We introduce a notion of a \emph{quasi-free right $P$-module}
(see paragraph~\ref{QuasiFreeRightModule})
for which these properties hold
(see theorems~\ref{QuasiFreeCompositionProduct} and~\ref{ExactnessCompositionProduct}).
We define also the symmetric notion of a \emph{quasi-free left $P$-module}
(see paragraph~\ref{QuasiFreeLeftModule}).

\subsection{Free modules over an operad}\label{FreeModules}
We observe that the composite $\Sigma_*$-module $L = M\circ P$ defines a free object
in the category of right $P$-modules.
Precisely,
the operad product induces a morphism
$M\circ\mu: M\circ P\circ P\,\rightarrow\,M\circ P$
which gives $L = M\circ P$ the structure of a right $P$-module
and
the operad unit induces a morphism of $\Sigma_*$-modules
$M\circ\eta: M\,\rightarrow\,M\circ P$
that satisfies the classical universal property.
Explicitly,
a morphism $\phi: M\,\rightarrow\,L'$, where $L'$ is a right $P$-module,
has a unique extension
\[\xymatrix{ M\ar[rd]\ar[rr]^\phi & & L' \\ & M\circ P\ar@{-->}[ru]_{\tilde{\phi}} & \\ }\]
such that $\tilde{\phi}: M\circ P\,\rightarrow\,L'$
is a morphism of right $P$-modules.
Equivalently,
we have an adjunction relation
\[\Hom_{\RMod{P}}(M\circ P,L') = \Hom_{\LMod{\Sigma_*}}(M,L').\]
In regard to the expansion of a composite $\Sigma_*$-module
(see lemma~\ref{ExpansionComposite}),
the universal morphism $M\circ\eta: M\,\rightarrow\,M\circ P$
identifies $x\in M(r)$
with the composite element $x(1,\ldots,1)\in M\circ P$.

Suppose $P$ is an augmented operad.
In that case,
the augmentation morphism induces a morphism of $\Sigma_*$-modules
$M\circ\epsilon: M\circ P\,\rightarrow\,M$
and gives a canonical retraction of the universal morphism
$M\circ\eta: M\,\rightarrow\,M\circ P$.
On the other hand,
the adjunction relation implies immediately that the $\Sigma_*$-module $M$
is isomorphic to the indecomposable quotient
$\bar{L} = (M\circ P)\circ_P I$.
One observes that the morphism $M\circ\epsilon: M\circ P\,\rightarrow\,M$
is identified with the canonical surjection
$L\,\rightarrow\,\bar{L}$.

Symmetrically,
the composite $\Sigma_*$-module $R = P\circ N$
represents the free left $P$-module generated by $N$.
The operad action is induced by the operad product.
The operad unit induces a canonical morphism $\eta\circ N: N\,\rightarrow\,P\circ N$
which identifies an element $x\in N(n)$ with the composite element $1(x)\in P\circ N(n)$.
If $P$ is an augmented operad,
then we have a morphism $\epsilon\circ N: P\circ N\,\rightarrow\,N$
which identifies the module $N$ with the indecomposable quotient
$\bar{R} = I\circ_P(P\circ N)$.

\subsection{Quasi-free right $P$-modules}\label{QuasiFreeRightModule}
In the dg-context,
a free right $P$-module is also realized by a composite $L = M\circ P$,
where $M$ is a dg-$\Sigma_*$-module.
The differential of $L$ is the canonical differential of the composite dg-$\Sigma_*$-module $L = M\circ P$.
The canonical inclusion $M\,\rightarrow\,M\circ P$
which identifies an element $x\in M(r)$ with a composite $x(1,\ldots,1)\in M\circ P(r)$
is a morphism of dg-$\Sigma_*$-modules.
The canonical projection $M\circ P\,\rightarrow\,M$
is also a morphism of dg-$\Sigma_*$-modules.
Hence,
the indecomposable quotient of the free right $P$-module $L = M\circ P$
is isomorphic to $M$
as a dg-$\Sigma_*$-module.

A \emph{quasi-free right $P$-module} $L$ is a dg-$\Sigma_*$-module such that $L = M\circ P$,
but whose differential,
denoted by $\delta_\theta: M\circ P\,\rightarrow\,M\circ P$,
differs from the canonical differential of a free right $P$-module
(which is still denoted by $\delta: M\circ P\,\rightarrow\,M\circ P$).
We have explicitly $\delta_\theta = \delta + d_\theta$,
where
$d_\theta: M\circ P\,\rightarrow\,M\circ P$
is a certain homogeneous morphism of degree $-1$.
We assume
\[d_\theta(x(q_1,\ldots,q_r)) = d_\theta(x)(q_1,\ldots,q_r),\]
so that $\delta_\theta: M\circ P\,\rightarrow\,M\circ P$ satisfies the derivation relation
\[\delta_\theta(x(q_1,\ldots,q_r)) = \delta_\theta(x)(q_1,\ldots,q_r)
+ \sum_{i=1}^r \pm x(q_1,\ldots,\delta(q_i),\ldots,q_r).\]
In addition,
we assume the relation $\delta d_\theta + d_\theta\delta + d_\theta{}^2 = 0$
(because this equation is equivalent to the identity $\delta_\theta{}^2 = 0$).

We observe that the morphism $d_\theta: M\circ P\,\rightarrow\,M\circ P$
is determined by its restriction to $M\subset M\circ P$.
Consequently,
the inclusion morphism $M\,\rightarrow\,M\circ P$ is not a morphism of dg-modules
as long as $d_\theta$ is non zero.
But,
the canonical projection $M\circ P\,\rightarrow\,M$
may be a morphism of dg-modules
provided that the composite
\[M\circ P\,\xrightarrow{d_\theta}\,M\circ P\,\xrightarrow{}\,M\]
vanishes.
In general,
we assume that this property holds,
because, in this situation,
the dg-module $M$ can be identified with the indecomposable quotient of $L = M\circ P$.

\subsection{Quasi-free left $P$-modules}\label{QuasiFreeLeftModule}
Symmetrically,
a composite $R = P\circ N$, where $N$ is a dg-$\Sigma_*$-module,
is a free object in the category of dg-left $P$-modules.
A \emph{quasi-free left $P$-module} $R$ is a dg-$\Sigma_*$-module such that $R = P\circ N$,
but whose differential,
denoted by $\delta_\theta: P\circ N\,\rightarrow\,P\circ N$,
differs from the canonical differential
of a free left $P$-module
(still denoted by $\delta: P\circ N\,\rightarrow\,P\circ N$).
We still assume $\delta_\theta = \delta + d_\theta$,
where
$d_\theta: P\circ N\,\rightarrow\,P\circ N$
is a certain homogeneous morphism of degree $-1$.
In the case of left modules,
the derivation relation is equivalent to the equation
\[d_\theta(p(y_1,\ldots,y_r)) = \sum_{i=1}^r \pm p(y_1,\ldots,d_\theta(y_i),\ldots,y_r)\]
for all $p\in P(r)$ and $y_1\in P\circ N(n_1),\ldots,y_r\in P\circ N(n_r)$.
We observe that the morphism
$d_\theta: P\circ N\,\rightarrow\,P\circ N$
is determined by its restriction to $N\subset P\circ N$.

The inclusion morphism $N\,\rightarrow\,P\circ N$ is not a morphism of dg-modules
as long as $d_\theta$ is non zero.
But,
the canonical projection $P\circ N\,\rightarrow\,N$
may be a morphism of dg-modules
provided that the composite
\[P\circ N\,\xrightarrow{d_\theta}\,P\circ N\,\xrightarrow{}\,N\]
vanishes.
In general,
we assume that this property holds,
so that the dg-module $N$
can be identified with the indecomposable quotient of $R = P\circ N$.

\smallskip
We prove the next theorems in section~\ref{ModuleSpectralSequences}.

\begin{thm}\label{QuasiFreeCompositionProduct}
Let $P$ be a connected dg-operad.

We assume that $\phi: L\,\xrightarrow{\sim}\,L'$ is a quasi-isomorphism of quasi-free right $P$-modules.
We have explicitly $L = \bar{L}\circ P$ (respectively, $L' = \bar{L}'\circ P$).
If the ground ring $\K$ is not a field,
then we assume in addition that the operad $P$
and the $\Sigma_*$-module $\bar{L}$ (respectively, $\bar{L}'$)
are projective as $\K$-modules.
The induced morphism
$\phi\circ_P R: L\circ_P R\,\rightarrow\,L'\circ_P R$,
where $R$ is any connected left $P$-module (which is projective as a $\K$-module),
is also a quasi-isomorphism.

We assume that $\psi: R\,\xrightarrow{\sim}\,R'$ is a quasi-isomorphism of connected quasi-free left $P$-modules.
We have explicitly $R = P\circ\bar{R}$ (respectively, $R' = P\circ\bar{R}'$),
where $\bar{R}(0) = 0$ (respectively, $\bar{R}'(0) = 0$).
If the ground ring $\K$ is not a field,
then we assume in addition that the $\Sigma_*$-module $\bar{R}$ (respectively, $\bar{R}'$)
is projective as a $\K$-module.
The induced morphism
$L\circ_P\psi: L\circ_P R\,\rightarrow\,L\circ_P R'$,
where $L$ is any right $P$-module,
is also a quasi-isomorphism.
\end{thm}

\begin{thm}\label{ExactnessCompositionProduct}
Let $P$ be a connected dg-operad.

We assume that $L$ is a quasi-free right $P$-module.
We have explicitly $L = \bar{L}\circ P$.
The functor $L\circ_P -: \dg\LMod{P}\,\rightarrow\,\dg\LMod{\Sigma_*}$
preserves all quasi-isomorphisms of connected left $P$-modules $\psi: R\,\xrightarrow{\sim}\,R'$
(provided $R$ and $R'$ are projective as $\K$-modules).

Symmetrically,
we assume that $R$ is a connected quasi-free left $P$-module.
We have explicitly $R = P\circ\bar{R}$, where $\bar{R}(0) = 0$.
If the ground ring $\K$ is not a field,
then we assume in addition that the $\Sigma_*$-module $\bar{R}$ is projective as a $\K$-module.
The functor $-\circ_P R: \dg\RMod{P}\,\rightarrow\,\dg\LMod{\Sigma_*}$
preserves all quasi-isomorphisms
of right $P$-modules $\phi: L\,\xrightarrow{\sim}\,L'$.
\end{thm}

We deduce the main results of this article from the next theorems
(proved in section~\ref{ProofModuleComparisonTheorems})
which are generalizations of classical comparison theorems
(\emph{cf}. H. Cartan \cite[Expos\'es 2-3]{Cartan},
J. McCleary \cite[Chapter 3]{McCleary},
E.C. Zeeman \cite{Zeeman})
in the context of operads.

\begin{thm}\label{RightModuleComparison}
We are given a morphism of connected dg-ope\-rads $\psi: P\,\rightarrow\,P'$.
We consider a morphism of right $P$-modules $\phi: L\,\rightarrow\,L'$.
We assume $L$ (respectively, $L'$) is a quasi-free module over $P$ (respectively, $P'$).
We have precisely $L = \bar{L}\circ P$ (respectively, $L' = \bar{L}'\circ P'$).
We equip the dg-$\Sigma_*$-module $L'$ with the structure of a right $P$-module by restriction of structure.
If the ground ring $\K$ is not a field,
then we assume in addition that the operad $P$ (respectively, $P'$)
and the $\Sigma_*$-module $\bar{L}$ (respectively, $\bar{L}'$)
are projective as $\K$-modules.

We identify the $\Sigma_*$-module $\bar{L}$ (respectively, $\bar{L}'$)
with the indecomposable quotient of $L$ (respectively, $L'$).
Hence,
the $\Sigma_*$-modules $\bar{L}$ and $\bar{L}'$ are equipped with a natural differential
and we have an induced morphism of dg-$\Sigma_*$-modules
$\bar{\phi}: \bar{L}\,\rightarrow\,\bar{L}'$.

\textbf{a)} If the morphisms $\psi: P\,\rightarrow\,P'$ and $\bar{\phi}: \bar{L}\,\rightarrow\,\bar{L}'$
are quasi-isomorphisms, then so is $\phi: L\,\rightarrow\,L'$.

\textbf{b)} If the morphisms $\psi: P\,\rightarrow\,P'$ and $\phi: L\,\rightarrow\,L'$
are quasi-isomorphisms, then so is $\bar{\phi}: \bar{L}\,\rightarrow\,\bar{L}'$.

\textbf{c)} We assume in addition $\bar{L}(0) = \bar{L}'(0) = 0$ and $\bar{L}(1) = \bar{L}'(1) = \K$.
If the morphisms $\phi: L\,\rightarrow\,L'$ and $\bar{\phi}: \bar{L}\,\rightarrow\,\bar{L}'$
are quasi-isomorphisms, then so is $\psi: P\,\rightarrow\,P'$.
\end{thm}

\begin{thm}\label{LeftModuleComparison}
We are given a morphism of connected dg-ope\-rads $\phi: P\,\rightarrow\,P'$.
We consider a morphism of left $P$-modules $\psi: R\,\rightarrow\,R'$.
We assume $R = P\circ\bar{R}$ (respectively, $R' = P\circ\bar{R}'$)
is a quasi-free module over $P$ (respectively, $P'$)
such that $\bar{R}(0) = 0$ (respectively, $\bar{R}'(0) = 0$).
We equip the dg-$\Sigma_*$-module $R'$ with the structure of a left $P$-module by restriction of structure.
If the ground ring $\K$ is not a field,
then we assume in addition that the $\Sigma_*$-module $\bar{R}$ (respectively, $\bar{R}'$)
is projective as a $\K$-module.

We identify the $\Sigma_*$-module $\bar{R}$ (respectively, $\bar{R}'$)
with the indecomposable quotient of $R$ (respectively, $R'$).
Hence,
the $\Sigma_*$-modules $\bar{R}$ and $\bar{R}'$ are equipped with a natural differential
and we have an induced morphism of dg-$\Sigma_*$-modules
$\bar{\phi}: \bar{R}\,\rightarrow\,\bar{R}'$.

\textbf{a)} If the morphisms $\phi: P\,\rightarrow\,P'$ and $\bar{\psi}: \bar{R}\,\rightarrow\,\bar{R}'$
are quasi-isomorphisms, then so is $\phi: R\,\rightarrow\,R'$.

\textbf{b)} If the morphisms $\phi: P\,\rightarrow\,P'$ and $\psi: R\,\rightarrow\,R'$
are quasi-isomorphisms, then so is $\bar{\psi}: \bar{R}\,\rightarrow\,\bar{R}'$.

\textbf{c)} We assume in addition $\bar{R}(1) = \bar{R}'(1) = \K$.
If the morphisms $\psi: R\,\rightarrow\,R'$ and $\bar{\psi}: \bar{R}\,\rightarrow\,\bar{R}'$
are quasi-isomorphisms, then so is $\phi: P\,\rightarrow\,P'$.
\end{thm}

\section[Digression: model categories of modules]{Digression: model categories of modules over an operad}

We interpret the comparison theorems above within the language of model categories.
In fact,
we would like to point out that our results are stronger than the assertions
which one expects from classical model structures.
First,
we mention the next statements,
which are similar to classical results
about modules over an associative algebra
(\emph{cf}. W. Dwyer and J. Spalinski \cite{DwyerSpalinski},
I. Gelfand and Y. Manin \cite{GelfandManin},
M. Hovey \cite{Hovey}).

\begin{thm}\label{ModelModuleCategories}
Let $P$ be a connected dg-operad.

The category of dg-right $P$-modules $\dg\RMod{P}$
is endowed with the structure of a model category
in which a weak-equivalence is a quasi-isomorphism of right $P$-modules.
A fibration is a morphism of right $P$-modules which is surjective in degree $*>0$.
A cofibration is a morphism of right $P$-modules which has the left lifting property with respect to fibrations.

The category of connected dg-left $P$-modules $\dg\LMod{P}_0$
is endowed with the structure of a model category
in which a weak-equivalence is a quasi-isomorphism of left $P$-modules.
A fibration is a morphism of left $P$-modules which is surjective in degree $*>0$.
A cofibration is a morphism of left $P$-modules which has the left lifting property with respect to fibrations.
\end{thm}

In case of a connected operad $P$,
the $\K$-module $0$ is the initial object in the category of right $P$-modules
and in the category of left $P$-modules.
We recall that a right $P$-module $L$ (respectively, a connected left $P$-module $R$)
is a cofibrant object
if the initial morphism $0\,\rightarrow\,L$ (respectively, $0\,\rightarrow\,R$)
is a cofibration.

\begin{lem}
A cofibrant object in the category of dg-right $P$-modules
is a retract of a quasi-free module $F = M\circ P$
such that $M$ is a complex of projective $\Sigma_*$-modules.

A cofibrant object in the category of dg-left $P$-modules
is a retract of a quasi-free module $F = P\circ N$
such that $N$ is a complex of projective $\Sigma_*$-modules.
\end{lem}

\subsection{Model structures for algebras over an operad}
In general, the definition of theorem~\ref{ModelModuleCategories}
does not provide the category of all dg-left $P$-modules $\dg\LMod{P}$
with a model structure.
We observe precisely that the category of dg-left $P$-modules $\dg\LMod{P}$ possesses a model structure
if and only if so does the category of $P$-algebras $\dg\LAlg{P}$.
Here are a few indications.

We consider the dg-$\K$-module $C = \K\,e\oplus\K\,b$
such that $|b| = 0$, $|e| = 1$ and $\delta(e) = b$.
On one hand, the free algebra $F = S(P,C)$ has the left lifting property
with respect to all surjective morphisms of $P$-algebras.
On the other hand,
the functor $S(P,-)$ does not preserves quasi-isomorphisms in general.
Accordingly,
the free algebra $S(P,C)$ is not acyclic.
Therefore,
we conclude that the category of all dg-left $P$-modules $\dg\LMod{P}$
is not equipped with the structure of a model category
in general.
We refer to V. Hinich \cite{Hinich} and M. Spitzweck \cite{Spitzweck}
for a discussion about such assertions.

We recall that a $P$-algebra $A$ is equivalent to a left $P$-module $R$
which has $R(0) = A$ and $R(n) = 0$ for $n>0$.
We identify the dg-$\K$-module $C$ with the dg-$\Sigma_*$-module $N$
which has $N(0) = C$ and $N(n) = 0$ for $n>0$.
We have $(P\circ N)(0) = S(P,C)$ and $(P\circ N)(n) = 0$ for $n>0$.
Therefore, the free module $F = P\circ N$ is acyclic if and only if the free algebra $S(P,C)$ is so.
On the other hand,
the free module $F = P\circ N$ has the left lifting property
with respect to all surjective morphisms of left $P$-modules.
We conclude that the category of all dg-left $P$-modules $\dg\LMod{P}$
is not equipped with the structure of a model category in general.

\subsection{Derived functors and quasi-free resolutions}\label{ModuleResolutions}
In positive characteristic,
the representations of the symmetric groups
are not all projective objects.
Therefore, the quasi-free modules over an operad
are not all cofibrant objects.

A \emph{cofibrant resolution} of a right $P$-module $L$
consists of a cofibrant module $F$
together with a surjective quasi-isomorphism $F\,\xrightarrow{\sim}\,L$.
A \emph{quasi-free resolution} of $L$
consists of a quasi-free module $L' = M'\circ P$
together with a surjective quasi-isomorphism $L'\,\xrightarrow{\sim}\,L$.
If the ground ring $\K$ is not a field,
then we consider \emph{$\K$-projective quasi-free resolutions}
in which the $\Sigma_*$-module $M'$ is supposed to be projective as a $\K$-module.
Any right $P$-module $L$
has a cofibrant resolution $F\,\xrightarrow{\sim}\,L$
such that $F$ is a quasi-free right $P$-module.
We have explicitly $F = M\circ P$, where $M$ is a complex of projective $\Sigma_*$-modules.
A cofibrant resolution is connected to any quasi-free resolution by a quasi-isomorphism $F\,\xrightarrow{\sim}\,L'$,
because the left lifting property provides the following diagram
with a fill in morphism
\[\xymatrix{ & L'\ar@{->>}[d]^{\sim} \\ F\ar[r]^{\sim}\ar@{-->}[ur] & L \\ }\]
We have similar observations for left $P$-modules.

The functor $-\circ_P R: \dg\RMod{P}\,\rightarrow\,\dg\LMod{\Sigma_*}$
does not preserve quasi-isomorphism of right $P$-modules.
Therefore,
we introduce the left derived functor of the relative composition product $-\circ^{\mathbb{L}}_P R$.
We replace a right $P$-module $L$ by a cofibrant resolution $F$
and we consider the composition product $F\circ_P R$.
We have by definition $L\circ^{\mathbb{L}}_P R = F\circ_P R$.
In general,
the dg-module $L\circ^{\mathbb{L}}_P R$ does not depend on the choice of a cofibrant resolution
(up to quasi-isomorphism).
In our context,
we obtain the following stronger result:

\begin{prp}\label{DerivedCompositionProduct}
Suppose given a connected dg-operad $P$, a right $P$-module $L$
and a connected left $P$-module $R$.
The relative composition product $L'\circ_P R$,
where $L'$ is any $\K$-projective quasi-free resolution of $L$,
is a representative of the derived composition product $L\circ^{\mathbb{L}}_P R$.
The relative composition product $L\circ_P R'$,
where $R'$ is any $\K$-projective quasi-free resolution of $R$,
is also equivalent to the derived composition product $L\circ^{\mathbb{L}}_P R$.
\end{prp}

\begin{proof}
We observe in paragraph~\ref{ModuleResolutions}
that the quasi-free module $L'$ is connected to a cofibrant resolution of $L$
by a quasi-isomorphism of right $P$-modules $F\,\xrightarrow{\sim}\,L'$.
We deduce from theorem~\ref{QuasiFreeCompositionProduct}
that the induced morphism
$F\circ_P R\,\xrightarrow{}\,L'\circ_P R$
is also a quasi-isomorphism.
Therefore,
we conclude that the composition product $L'\circ_P R$ is equivalent to $L\circ^{\mathbb{L}}_P R = F\circ_P R$.

In addition,
if $R'$ is a $\K$-projective quasi-free resolution of $R$,
then theorem~\ref{ExactnessCompositionProduct} provides a zig-zag of quasi-isomorphisms
\[L\circ_P R'\,\xleftarrow{\sim}\,L'\circ_P R'\,\xrightarrow{\sim}\,L'\circ_P R.\]
As a consequence, the composition product $L\circ_P R'$
is also equivalent to $L\circ^{\mathbb{L}}_P R$.
\end{proof}

\subsection{Remark}\label{NonClassicalModelStructure}
We would like to point out that the theorems of section~\ref{SummaryModuleComplexes}
could be deduced from non-classical model structures.
For instance,
we consider the category of right modules over an operad $P$.
We still assume that a morphism of right $P$-modules $\phi: L\,\rightarrow\,L'$
is a weak-equivalence (respectively, a fibration)
if $\phi: L\,\rightarrow\,L'$ is a weak-equivalence (respectively, a fibration)
in the category of dg-$\Sigma_*$-modules,
but we modify the definition of a fibration
in the category of dg-$\Sigma_*$-modules.

To be explicit,
we equip the category of dg-$\Sigma_*$-modules with the model structure
in which weak-equivalences are quasi-isomorphisms
and cofibrations are injective morphisms.
In this context,
a fibration of dg-$\Sigma_*$-modules is a morphism of dg-$\Sigma_*$-modules
which has the right-lifting property
with respect to acyclic cofibrations.

This definition makes all quasi-free modules cofibrant objects.

\section[On composite symmetric modules]{On composite symmetric modules in the differential graded framework}\label{DGComposite}

In this section,
we analyze the homology of composite $\Sigma_*$-modules.
We deduce our results from the nice expansion of the composition product
given in lemma~\ref{ExpansionComposite}.

\subsection{The differential structure of a composite $\Sigma_*$-module}\label{CompositionBicomplex}
We observe that a composite dg-$\Sigma_*$-module $M\circ N$ forms a bicomplex.
Precisely,
the bidegree of an element
\[x(y_1,\ldots,y_r)\in(M\circ N)(n),\]
such that $x\otimes y_1\otimes\cdots\otimes y_r\in M_d(r)\otimes N_{e_1}(n_1)\otimes\cdots\otimes N_{e_1}(n_r)$,
is defined by $(s,t) = (d,e_1+\cdots+e_r)$.
The horizontal differential $\delta_h: M\circ N\,\rightarrow\,M\circ N$
is induced by the differential of $M$.
The vertical differential $\delta_v: M\circ N\,\rightarrow\,M\circ N$
is induced by the differential of $N$.
We set explicitly
\begin{align*} \delta_h(x(y_1,\ldots,y_r)) &
= \delta(x)(y_1,\ldots,y_r) \\
\text{and}\qquad\delta_v(x(y_1,\ldots,y_r)) &
= \sum_{i=1}^r \pm x(y_1,\ldots,\delta(y_i),\ldots,y_r). \end{align*}
Hence, we have clearly
\[\delta(x(y_1,\ldots,y_r)) = \delta_h(x(y_1,\ldots,y_r)) + \delta_v(x(y_1,\ldots,y_r)).\]
Consequently,
as recalled in paragraph~\ref{BicomplexSpectralSequence}
we have a spectral sequence
$I^r(M\circ N)\,\Rightarrow\,H_*(M\circ N,\delta)$
such that
$I^1(M\circ N) = H_*(M\circ N,\delta_v)$
and
$I^2(M\circ N) = H_*(H_*(M\circ N,\delta_v),\delta_h)$.
Similarly,
we have a spectral sequence
$II^r(M\circ N)\,\Rightarrow\,H_*(M\circ N,\delta)$
such that
$II^1(M\circ N) = H_*(M\circ N,\delta_h)$
and
$II^2(M\circ N) = H_*(H_*(M\circ N,\delta_h),\delta_v)$.

We observe also that the splitting
$M\circ N = \bigoplus_{r=0}^{\infty} (M(r)\otimes N^{\otimes r})_{\Sigma_r}$
is compatible with the bicomplex structure.
Hence,
the spectral sequences above
have natural splittings
\begin{align*} I^r(M\circ N) & = \bigoplus_{r=0}^{\infty} I^r(M(r)\otimes N^{\otimes r})_{\Sigma_r} \\
\text{and}\qquad II^r(M\circ N) & = \bigoplus_{r=0}^{\infty} II^r(M(r)\otimes N^{\otimes r})_{\Sigma_r}. \end{align*}

We go back to lemma~\ref{ExpansionComposite}.
According to this statement,
the modules
$I^1(M(r)\otimes N^{\otimes r})_{\Sigma_r}$
and
$II^1(M(r)\otimes N^{\otimes r})_{\Sigma_r}$
have the following expansions:

\begin{fact}\label{CompositeHorizontalHomology}
We assume that $M$ is a $\Sigma_*$-module which is projective as a $\K$-module.
We assume that $N$ is a connected $\Sigma_*$-module.
We have
$I^1(M(r)\otimes N^{\otimes r})_{\Sigma_r}
= \bigoplus_{(J_1,\ldots,J_r)'} M(r)\otimes H_*(N(J_1)\otimes\cdots\otimes N(J_r))$.
\end{fact}

\begin{fact}\label{CompositeVerticalHomology}
We let $M$ be any $\Sigma_*$-module.
We assume that $N$ is a connected $\Sigma_*$-module which is projective as a $\K$-module.
We have
$II^1(M(r)\otimes N^{\otimes r})_{\Sigma_r}
= \bigoplus_{(J_1,\ldots,J_r)'} H_*(M(r))\otimes N(J_1)\otimes\cdots\otimes N(J_r)$.
\end{fact}

\subsection{Morphisms of composite dg-$\Sigma_*$-modules and spectral sequences}
Suppose given morphisms of dg-$\Sigma_*$-modules $\phi: M\,\rightarrow\,M'$ and $\psi: N\,\rightarrow\,N'$.
The composite $\phi\circ\psi: M\circ N\,\rightarrow\,M'\circ N'$
is a morphism of bicomplexes.
Consequently,
we obtain morphisms of spectral sequences
$I^r(\phi\circ\psi): I^r(M\circ N)\,\rightarrow\,I^r(M'\circ N')$
and
$II^r(\phi\circ\psi): II^r(M\circ N)\,\rightarrow\,II^r(M'\circ N')$.

\smallskip
We have the following result:

\begin{lem}\label{CompositionMorphismHomology}
We consider the morphism of spectral sequences
\[I^r(\phi\circ\psi): I^r(M\circ N)\,\rightarrow\,I^r(M'\circ N').\]
The $\Sigma_*$-modules $N$ and $N'$ are supposed to be connected.
If the ground ring $\K$ is not a field,
then we assume in addition that all $\Sigma_*$-modules $M$, $M'$, $N$ and $N'$
are projective as $\K$-modules.

Suppose
$\phi_*: H_*(M(r))\,\rightarrow\,H_*(M'(r))$
is an isomorphism for all $r\in\N$ such that $r<a$
and
$\psi_*: H_*(N(n))\,\rightarrow\,H_*(N'(n))$
is an isomorphism for all $n\in\N$ such that $n<b$.
Then,
the induced morphism
\[I^2(\phi\circ\psi): I^2(M(r)\otimes N^{\otimes r}(n))_{\Sigma_r}
\,\rightarrow\,I^2(M'(r)\otimes N^{\otimes r}(n))_{\Sigma_r}\]
is an isomorphism for all pairs $(r,n)$ such that $r<a$ and $n<b$.
In the case $r>1$, we have also an isomorphism for $n=b$.
\end{lem}

\begin{proof}
We deduce expansions of $I^2(M\circ N)$ and $I^2(M'\circ N')$
from the identity of paragraph~\ref{CompositeHorizontalHomology}.
In the case $r>1$ and $n\leq b$,
the tensor product
\[(M(r)\otimes N{}^{\otimes r}(n))_{\Sigma_r}
= \bigoplus_{(J_1,\ldots,J_r)'} M(r)\otimes N(J_1)\otimes\cdots\otimes N(J_r)\]
involves modules $N(J_i)\simeq N(n_i)$ such that $0<n_i<b$.
In fact,
since $N(0) = 0$,
we may assume $n_i>0$
and
the relation $n_1+\cdots+n_r = n\leq b$ implies $n_i<b$.
In the case $r=1$,
we have
\[(M(1)\otimes N^{\otimes 1}(n))_{\Sigma_1} = M(1)\otimes N(n).\]
Similar observations hold for the tensor product $(M'(r)\otimes N'{}^{\otimes r}(n))_{\Sigma_r}$.

If the morphisms $\psi: N(n)\,\rightarrow\,N'(n)$, where $n<b$,
are quasi-isomorphisms of projective $\K$-modules,
then the induced homology morphism
\[H_*(N(J_1)\otimes\cdots\otimes N(J_r))\,\rightarrow\,H_*(N'(J_1)\otimes\cdots\otimes N'(J_r))\]
is an isomorphism (provided $n_1,\ldots,n_r<b$).
If a morphism $\phi: M(r)\,\rightarrow\,M'(r)$
is a quasi-isomorphism of projective $\K$-modules,
then the tensor product of a graded $\K$-module
with $\phi: M(r)\,\rightarrow\,M'(r)$
gives also a homology isomorphism.

The lemma follows from these observations.
\end{proof}

We have a similar result for the second spectral sequence:

\begin{lem}
We consider the morphism of spectral sequences
\[II^r(\phi\circ\psi): II^r(M\circ N)\,\rightarrow\,II^r(M'\circ N').\]
The $\Sigma_*$-modules $N$ and $N'$ are supposed to be connected
and $\K$-pro\-jec\-tive (if the ground ring is not a field).

Suppose
$\phi_*: H_*(M(r))\,\rightarrow\,H_*(M'(r))$
is an isomorphism for all $r\in\N$ such that $r<a$
and
$\psi_*: H_*(N(n))\,\rightarrow\,H_*(N'(n))$
is an isomorphism for all $n\in\N$ such that $n<b$.
Then,
the induced morphism
\[II^2(\phi\circ\psi): I^2(M(r)\otimes N^{\otimes r}(n))_{\Sigma_r}
\,\rightarrow\,II^2(M'(r)\otimes N^{\otimes r}(n))_{\Sigma_r}\]
is an isomorphism for all pairs $(r,n)$ such that $r<a$ and $n<b$.
In the case $r>1$, we have also an isomorphism for $n=b$.
\end{lem}

\section{The spectral sequence of a quasi-free module}\label{ModuleSpectralSequences}

In this section,
we prove that a quasi-free right $P$-module is equip\-ped with a natural filtration
giving rise to a spectral sequence.
Precisely,
we obtain the following result:

\begin{prp}\label{RightModuleSpectralSequence}
We let $P$ be a connected dg-operad.
We suppose given a quasi-free right $P$-module $L = M\circ P$.
We assume that $M$ can be identified with the indecomposable quotient of $L$.
We consider a relative composition product $L\circ_P R$
where $R$ is any connected left $P$-module.
We have then $L\circ_P R = M\circ R$
(but the differential of $L\circ_P R$ differs from the canonical differential of the composite module $M\circ R$).

We have a spectral sequence of $\Sigma_*$-modules $E^r(L\circ_P R)\,\Rightarrow\,H_*(L\circ_P R)$.
To be precise,
the component $E^0_{s\,*-s}(L\circ_P R)$ of the $\Sigma_*$-modules $E^0(L\circ_P R)$
is generated by tensors $x(y_1,\ldots,y_r)\in(M\circ R)(n)$
such that $x\in M(r)_d$
where $d + r = s$.
Consequently,
for fixed $n\in\N$,
the $\Sigma_n$-module $E^r_{s t}(L\circ_P R)(n)$
lies in the quadrant $(s\geq 0,t\geq -n)$.

The differential $d^0: E^0(L\circ_P R)\,\rightarrow\,E^0(L\circ_P R)$
is induced by the differential of $R$.
The differential $d^1: E^1(L\circ_P R)\,\rightarrow\,E^1(L\circ_P R)$
is induced by the differential of $M$.
Hence,
the $E^2$-term of the spectral sequence is connected to the module $I^2(M\circ R)$
introduced in paragraph~\ref{CompositionBicomplex}.
We have precisely
\[E^2_{s t}(L\circ_P R)(n) = \bigoplus_{r=0}^{\infty} I^2_{s-r\,t+r}(M(r)\otimes R^{\otimes r}(n))_{\Sigma_r}.\]
In particular,
if the ground ring $\K$ is a field (of characteristic $p\geq 0$),
then we obtain $E^2(L\circ_P R) = H_*(M)\circ H_*(R)$.
\end{prp}

In the case $R = P$, we have $L\circ_P P = L$ and we obtain a spectral sequence $E^r(L)\,\rightarrow\,H_*(L)$.

\smallskip
We have a symmetrical result for quasi-free left $P$-modules.
We obtain precisely:

\begin{prp}\label{LeftModuleSpectralSequence}
We let $P$ be a connected dg-operad.
We suppose given a connected quasi-free left $P$-module $R = P\circ N$.
(We have then $N(0) = 0$.)
We assume that $N$ can be identified with the indecomposable quotient of $R$.
We consider a relative composition product $L\circ_P R$
where $L$ is any left $P$-module.
We have then $L\circ_P R = L\circ N$ (but the differential of $L\circ_P R$
differs from the canonical differential of the composite module $L\circ N$).

We have a spectral sequence of $\Sigma_*$-modules $E^r(L\circ_P R)\,\Rightarrow\,H_*(L\circ_P R)$.
To be precise,
the component $E^0_{s\,*-s}(L\circ_P R)$ of the $\Sigma_*$-module $E^0(L\circ_P R)$
is generated by elements $x(y_1,\ldots,y_r)\in(L\circ N)(n)$
such that $y_1\in N(n_1)_{e_1},\ldots,y_r\in N(n_r)_{e_r}$
where $(n_1 + e_1 - 1) + \cdots + (n_r + e_r - 1) = s$.
Consequently,
for fixed $n\in\N$,
the $\Sigma_n$-module $E^r_{s t}(L\circ_P R)(n)$
lies in the quadrant $(s\geq 0,t\geq 1 - n)$.

The differential $d^0: E^0(L\circ_P R)\,\rightarrow\,E^0(L\circ_P R)$
is induced by the differential of $L$.
The differential $d^1: E^1(L\circ_P R)\,\rightarrow\,E^1(L\circ_P R)$
is induced by the differential of $N$.
Hence,
the $E^2$-term of the spectral sequence is connected to the module $II^2(L\circ N)$
introduced in paragraph~\ref{CompositionBicomplex}.
We have precisely
\[E^2_{s t}(L\circ_P R)(n) = \bigoplus_{r=0}^{\infty} II^2_{t+n-r\,s+r-n}(L(r)\otimes N^{\otimes r}(n))_{\Sigma_r}.\]
In particular,
if the ground ring $\K$ is a field (of characteristic $p\geq 0$),
then we obtain $E^2(L\circ_P R) = H_*(L)\circ H_*(N)$.
\end{prp}

The purpose of the next paragraphs is to define the spectral sequence of a quasi-free right $P$-module.
The symmetrical constructions for quasi-free left $P$-modules are omitted.

\subsection{The differential of a quasi-free right $P$-module}\label{RightModuleDifferential}
We are given a right module $L$ over a connected dg-operad $P$
as in proposition~\ref{RightModuleSpectralSequence}.
We assume $L$ is quasi-free $L = M\circ P$
and is equipped with a differential $\delta_\theta: M\circ P\,\rightarrow\,M\circ P$
as in the discussion of paragraph~\ref{QuasiFreeRightModule}.

We still assume that the composite
\[M\circ P\,\xrightarrow{d_\theta}\,M\circ P\,\xrightarrow{}\,M\]
vanishes,
so that $M$ is identified with the indecomposable quotient $\bar{L} = L\circ_P I$.
We make this assumption more explicit.
We have
\[(M\circ P)(n) = \bigoplus_{r=0}^{\infty} (M(r)\otimes P^{\otimes r}(n))_{\Sigma_r}.\]
Furthermore,
since $P$ is supposed to be connected,
we obtain
$(M(r)\otimes P^{\otimes r}(n))_{\Sigma_r} = 0$ for $r>n$
and
$(M(n)\otimes P^{\otimes n}(n))_{\Sigma_n} \simeq M(n)$.
The canonical projection $M\circ P\,\rightarrow\,M$ is identified with the projection onto the summands
$(M(n)\otimes P^{\otimes n}(n))_{\Sigma_n}\subset M\circ P(n)$.

Finally,
the property above is equivalent to the relation
\[d_\theta(M(n))\subset\bigoplus_{r<n} (M(r)\otimes P^{\otimes r}(n))_{\Sigma_r}.\]
To be precise,
if $d_\theta(x)\in\bigoplus_{r<m} (M(r)\otimes P^{\otimes r}(m))_{\Sigma_r}$ for all $x\in M(m)$,
then the same property holds for composite elements
$x(q_1,\ldots,q_m)\in (M(m)\otimes P^{\otimes m}(n))_{\Sigma_m}$,
since we have the relation
\[d_\theta(x(q_1,\ldots,q_m)) = d_\theta(x)(q_1,\ldots,q_m)
\in\bigoplus_{r<m} (M(r)\otimes P^{\otimes r}(n))_{\Sigma_r}\]
(see paragraph~\ref{QuasiFreeRightModule}).

\subsection{The filtration of a quasi-free right $P$-module}\label{RightModuleFiltration}
We equip a quasi-free module with a filtration
$0 = F_{-1}(L)\subset F_0(L)\subset\cdots\subset F_s(L)\subset\cdots\subset L$.
Explicitly,
the module $F_s(L)\subset L$
consists of tensors $x(q_1,\ldots,q_r)\in M\circ P(n)$
such that $x\in M(r)_d$
together with $d + r\leq s$.

A homogeneous component $E^0_{s t}(L) = F_s(L)_{s+t}/F_{s-1}(L)_{s+t}$
can be identified to a submodule of $M\circ P$.
To be more precise,
we have a natural splitting $M\circ P = \bigoplus_{s,t} E^0_{s t}(L)$
determined by the bigrading of a composite dg-$\Sigma_*$-module
introduced in paragraph~\ref{CompositionBicomplex}.
Explicitly,
by definition of the filtration,
an element
$x(q_1,\ldots,q_r)\in(M\circ P)(n)$,
such that
$x\otimes q_1\otimes\cdots\otimes q_r\in M(r)_d\otimes P(n_1)_{e_1}\otimes\cdots\otimes P(n_r)_{e_r}$,
belongs to the module $E^0_{s t}(L(n))$,
where $s = d+r$ and $t = e_1+\cdots+e_r-r$.
Hence,
if we go back to the notation of paragraph~\ref{CompositionBicomplex},
then we obtain the formula
\[E^0_{s t}(L)(n) = \bigoplus_{r=0}^{\infty} I^0_{s-r\,t+r}(M(r)\otimes P^{\otimes r}(n))_{\Sigma_r}.\]

We observe that the differential
$\delta_\theta: M\circ P\,\rightarrow\,M\circ P$
preserves the filtration above.
To be more precise,
we split the canonical differential of $M\circ P$ into $\delta = \delta_M+\delta_P$
where
$\delta_M: M\circ P\,\rightarrow\,M\circ P$
is induced by the differential of $M$
and
$\delta_P: M\circ P\,\rightarrow\,M\circ P$
is induced by the differential of $P$.
The differential $\delta_M$ (respectively, $\delta_P$)
is identified with the horizontal (respectively, vertical) differential of paragraph~\ref{CompositionBicomplex}.
Consequently,
we obtain
$\delta_M(E^0_{s t}(L)(n))\subset E^0_{s-1 t}(L)(n)$
and
$\delta_P(E^0_{s t}(L)(n))\subset E^0_{s t-1}(L)(n)$.

We observe that the extra component of the differential $d_\theta: M\circ P\,\rightarrow\,M\circ P$
satisfies $d_\theta(F_s(L))\subset F_{s-2}(L)$.
More precisely,
this property is verified if we assume
$d_\theta(M(n))\subset\bigoplus_{r<n} (M(r)\otimes P^{\otimes r})_{\Sigma_r}$.
In fact,
in this situation,
the expansion of $d_\theta(x)\in(M\circ P)(n)_{d-1}$, where $x\in M(n)_d$,
involves elements $x'(q'_1,\ldots,q'_r)\in(M\circ P)(n)_{d-1}$
such that $x'\in M(r)_{d'}$
satisfies $r\leq n-1$ and $d'\leq |x'(q'_1,\ldots,q'_r)| = d-1$.

We deduce from all these observations that the filtration $F_s(L)\subset L$
gives rise to a spectral sequence $E^r(L)\,\Rightarrow\,H_*(L)$.

\subsection{The filtration of a relative composition product}
We observe that the module filtration layers of a quasi-free right $P$-module
are quasi-free right $P$-modules.
We have explicitly $F_s(L) = F_s(M)\circ P$,
where the module $F_s(M)\subset M$ consists of elements $x\in M(r)_d$
such that $d + r\leq s$.

Suppose given a left $P$-module $R$.
We equip the relative composition product $L\circ_P R$
with the induced filtration
\[0 = F_{-1}(L)\circ_P R\subset F_0(L)\circ_P R\subset\cdots\subset F_s(L)\circ_P R\subset\cdots\subset L\circ_P R.\]
Since the module $F_s(L)$ is quasi-free, we obtain $F_s(L)\circ_P R = F_s(M)\circ R$.

As in paragraph~\ref{RightModuleFiltration},
we split the canonical differential of $M\circ R$
into $\delta = \delta_M + \delta_R$
where $\delta_M: M\circ R\,\rightarrow\,M\circ R$
is induced by the differential of $M$
and $\delta_R: M\circ R\,\rightarrow\,M\circ R$
is induced by the differential of $R$.
The differential of a composite element $x(y_1,\ldots,y_r)\in(M\circ R)(n)_d$
in $L\circ_P R$
is the sum of the terms
\begin{align*} & \delta_M(x(y_1,\ldots,y_r)) = \delta(x)(y_1,\ldots,y_r), \\
& \delta_R(x(y_1,\ldots,y_r)) = \sum\nolimits_i \pm x(y_1,\ldots,\delta(y_i),\ldots,y_r) \\
& \text{and}\quad d_\theta(x)(y_1,\ldots,y_r)\in(L\circ_P R)(n). \end{align*}
We have clearly
$\delta_M(F_s(M)\circ R)_{d}\subset(F_{s-1}(M)\circ R)_{d-1}$,
$\delta_R(F_s(M)\circ R)_{d}\subset(F_s(M)\circ R)_{d-1}$.
Furthermore, if $x(y_1,\ldots,y_r)\in(F_s(M)\circ R)(n)_d$,
then we obtain $d_\theta(x)(y_1,\ldots,y_r)\in(F_{s-2}(L)\circ_P R)(n)_{d-1}$.

We conclude from these observations that the filtration $F_s(L\circ_P R) = F_s(L)\circ_P R$
gives rise to a spectral sequence $E^r(L\circ_P R)\,\Rightarrow\,H_*(L\circ_P R)$
of the form stated in proposition~\ref{RightModuleSpectralSequence}.

\begin{lem}\label{RightModuleMorphism}
We are given a morphism of connected dg-operads $\psi: P\,\rightarrow\,P'$.
We assume $L$ (respectively, $L'$) is a quasi-free module over $P$ (respectively, $P'$)
as in proposition~\ref{RightModuleSpectralSequence}.
We have precisely $L = M\circ P$ (respectively, $L' = M'\circ P'$).
We assume that the dg-$\Sigma_*$-module $M$ (respectively, $M'$)
can be identified with the indecomposable quotient of $L$ (respectively, $L'$)

We consider a morphism of right $P$-modules $\phi: L\,\rightarrow\,L'$.
(The dg-$\Sigma_*$-module $L'$ is equipped with the structure of a right $P$-module by restriction of structure.)
Such a morphism $\phi: L\,\rightarrow\,L'$
preserves the filtrations defined in paragraph~\ref{RightModuleFiltration}
and gives rise to a morphism of spectral sequences
$E^r(\phi): E^r(L)\,\rightarrow\,E^r(L')$.

We consider the induced morphism of dg-$\Sigma_*$-modules
$\bar{\phi}: \bar{L}\,\rightarrow\,\bar{L}'$,
where $\bar{L} \simeq M$
(respectively, $\bar{L}' \simeq M'$)
denotes the indecomposable quotient of $L$
(respectively, $L'$).
The morphism
$E^0(\phi): E^0(L)\,\rightarrow\,E^0(L')$
is identified with the composite
$\bar{\phi}\circ\psi: \bar{L}\circ P\,\rightarrow\,\bar{L}'\circ P'$.
\end{lem}

\begin{proof}
By definition of a $P$-module morphism,
we have the formula
\[\phi(x(q_1,\ldots,q_r)) = \phi(x)(\psi(q_1),\ldots,\psi(q_r)),\]
for all $x(q_1,\ldots,q_r)\in (M\circ P)(n)$.
Therefore,
we prove the relation $\phi(x)\in F_{r+d}(M'\circ P')(r)$, for $x\in M(r)_d$,
and the general case follows.

We have $(M'\circ P')(r) = \bigoplus_{s=0}^r (M'(s)\otimes P'{}^{\otimes s}(r))_{\Sigma_s}$.
Consequently,
the expansion of $\phi(x)\in(M'\circ P')(r)_d$, where $x\in M(r)_d$,
involves elements $x'(q'_1,\ldots,q'_{s})\in(M'\circ P')(r)_d$
such that $x'\in M'(s)_{d'}$
satisfies $s\leq r$ and $d'\leq |x'(q'_1,\ldots,q'_s)| = d$.
Hence,
we conclude $\phi(x)\in F_{r+d}(M'\circ P')(r)$.

We mention in paragraph~\ref{RightModuleDifferential}
that the canonical projection
$(M'\circ P')(n)\,\rightarrow\,M'(n)$
is given by the projection onto the summand $(M'(n)\otimes P'{}^{\otimes n}(n))_{\Sigma_n}$ of $(M'\circ P')(n)$.
In the expansion of $\phi(x)\in F_{n+d}(M'\circ P')(n)_d$
the elements
$x'(q'_1,\ldots,q'_{n'})\in (M'(n')\otimes P'{}^{\otimes n'}(n))_{\Sigma_{n'}}$
such that $n'<n$
lies in $F_{n-1+d}(M'\circ P')(n)$.
Therefore,
the projection of $\phi(x)\in F_{n+d}(M'\circ P')(n)_d$ onto $E^0_{n+d\,d}(L)(n)$
is represented by the image of $x\in M(n)_d$
under the composite
\[M(n)\,\subset\,(M\circ P)(n)\,\xrightarrow{\phi}\,(M'\circ P')(n)\,\xrightarrow{}\,M'(n),\]
which we identify with the morphism
$\bar{\phi}: \bar{L}\,\rightarrow\,\bar{L}'$.
\end{proof}

We deduce the theorems~\ref{QuasiFreeCompositionProduct} and~\ref{ExactnessCompositionProduct}
from the following lemmas.

\begin{lem}\label{RightModuleMorphismComposition}
We assume $L$ (respectively, $L'$) is a quasi-free right $P$-module
as in proposition~\ref{RightModuleSpectralSequence}.
We have precisely $L = M\circ P$ (respectively, $L' = M'\circ P$).
We assume that the dg-$\Sigma_*$-module $M$ (respectively, $M'$)
is identified with the indecomposable quotient of $L$ (respectively, $L'$).

We consider a relative composition product $-\circ_P R$,
where $R$ is any connected left $P$-module.
A morphism of right $P$-modules
$\phi: L\,\rightarrow\,L'$
induces a morphism of spectral sequences
$E^r(\phi\circ_P R): E^r(L\circ_P R)\,\rightarrow\,E^r(L'\circ_P R)$.

We consider the induced morphism of dg-$\Sigma_*$-modules
$\bar{\phi}: \bar{L}\,\rightarrow\,\bar{L}'$,
where $\bar{L}\simeq M$
(respectively, $\bar{L}'\simeq M'$)
denotes the indecomposable quotient of $L$
(respectively, $L'$).
The morphism
$E^0(\phi\circ_P R): E^0(L\circ_P R)\,\rightarrow\,E^0(L'\circ_P R)$
is identified with the composite
$\bar{\phi}\circ R: \bar{L}\circ R\,\rightarrow\,\bar{L}'\circ R$.
\end{lem}

\begin{lem}\label{LeftModuleMorphismComposition}
We assume $L$ is a quasi-free right $P$-module
as in proposition~\ref{RightModuleSpectralSequence}.
We have precisely $L = M\circ P$.
We assume that the dg-$\Sigma_*$-module $M$
is identified with the indecomposable quotient of $L$.

A morphism of left $P$-modules
$\psi: R\,\rightarrow\,R'$
induces a morphism of spectral sequences
$E^r(L\circ_P\psi): E^r(L\circ_P R)\,\rightarrow\,E^r(L\circ_P R')$.
The morphism
$E^0(L\circ_P\psi): E^0(L\circ_P R)\,\rightarrow\,E^0(L\circ_P R')$
is identified with the composite
$M\circ\psi: M\circ R\,\rightarrow\,M\circ R'$.
\end{lem}

\begin{proof}
The proofs of lemmas~\ref{RightModuleMorphismComposition}
and~\ref{LeftModuleMorphismComposition} involve similar arguments.
We give some indications in the case of lemma~\ref{RightModuleMorphismComposition}.

By lemma~\ref{RightModuleMorphism}, we have $\phi(F_s(L))\subset F_s(L')$.
Consequently,
the morphism $\phi\circ_P R: L\circ_P R\,\rightarrow\,L'\circ_P R$ preserves
the induced filtrations of the relative composition products.
Therefore,
we have a morphism of spectral sequences $E^r(\phi\circ_P R): E^r(L\circ_P R)\,\rightarrow\,E^r(L'\circ_P R)$.
We omit the proof of the last assertion
about the morphism $E^0(\phi\circ_P R): E^0(L\circ_P R)\,\rightarrow\,E^0(L'\circ_P R)$.
We have a similar result in lemma~\ref{RightModuleMorphism}
which follows from the same observations.
\end{proof}

We omit to state the analogues of lemmas~\ref{RightModuleMorphism},
\ref{RightModuleMorphismComposition} and~\ref{LeftModuleMorphismComposition}
for left modules, because these assertions are exactly the same
for left and right modules over an operad.

\begin{proof}[Proof of theorem~\ref{QuasiFreeCompositionProduct}]
Suppose given a quasi-isomorphism of quasi-free right $P$-modules $\phi: L\,\xrightarrow{\sim}\,L'$
as in theorem~\ref{QuasiFreeCompositionProduct}.
We consider the associated morphism of spectral sequences
$E^r(\phi\circ_P R): E^r(L\circ_P R)\,\rightarrow\,E^r(L'\circ_P R)$
supplied by lemma~\ref{RightModuleMorphismComposition}.
By the comparison theorem~\ref{RightModuleComparison},
we have also a quasi-isomorphism $\bar{\phi}: \bar{L}\,\xrightarrow{\sim}\,\bar{L}'$
on the indecomposable quotients.
Therefore,
by lemma~\ref{CompositionMorphismHomology},
the morphism $E^0(\phi\circ_P R): E^0(L\circ_P R)\,\rightarrow\,E^0(L'\circ_P R)$,
which is identified with the composition product
$\bar{\phi}\circ R: \bar{L}\circ R\,\xrightarrow{\sim}\,\bar{L}'\circ R$,
gives rise to a homology isomorphism.
More precisely,
we deduce from lemma~\ref{CompositionMorphismHomology}
that the morphism $E^2(\phi\circ_P R): E^2(L\circ_P R)\,\rightarrow\,E^2(L'\circ_P R)$
is an isomorphism.
We conclude immediately that the composition product
$\phi\circ_P R: L\circ_P R\,\rightarrow\,L'\circ_P R$
is a quasi-isomorphism.

The arguments are same in the symmetrical case of a composition product
$L\circ_P\psi: L\circ_P R\,\rightarrow\,L\circ_P R'$,
such that $\psi: R\,\xrightarrow{\sim}\,R'$ is a quasi-isomorphism
of quasi-free left $P$-modules.
\end{proof}

\begin{proof}[Proof of theorem~\ref{ExactnessCompositionProduct}]
Suppose given a quasi-free right $P$-module $L$
and a quasi-isomorphism of connected left $P$-modules $\psi: R\,\xrightarrow{\sim}\,R'$
as in theorem~\ref{ExactnessCompositionProduct}.
We consider the associated morphism of spectral sequences
$E^r(L\circ_P\psi): E^r(L\circ_P R)\,\rightarrow\,E^r(L\circ_P R')$
supplied by lemma~\ref{LeftModuleMorphismComposition}.
By lemma~\ref{CompositionMorphismHomology},
the morphism $E^0(L\circ_P\psi): E^0(L\circ_P R)\,\rightarrow\,E^0(L\circ_P R')$,
which is identified with the composition product
$\bar{L}\circ\psi: \bar{L}\circ R\,\xrightarrow{\sim}\,\bar{L}\circ R'$,
gives rise to a homology isomorphism.
The conclusion follows.

The arguments are similar in the symmetrical case of a composition product
$\phi\circ_P R: L\circ_P R\,\rightarrow\,L'\circ_P R$,
where $R$ is a quasi-free left $P$-module.
\end{proof}

\section{Proof of comparison theorems}\label{ProofModuleComparisonTheorems}

The purpose of this section is to give the proof of the comparison theorem~\ref{RightModuleComparison}.
The proof of theorem~\ref{LeftModuleComparison} is similar and is omitted.

\smallskip
We consider the morphism of spectral sequences
$E^r(\phi): E^r(L)\,\rightarrow\,E^r(L')$
given by lemma~\ref{RightModuleMorphism}.
As mentioned in the summary, we deduce theorem~\ref{RightModuleComparison}
from a generalization of classical comparison results
(\emph{cf. loc. cit.}).

\subsection{Proof of assertion \textbf{(a)}}
\emph{We assume $\psi: P\,\rightarrow\,P'$ and $\bar{\phi}: \bar{L}\,\rightarrow\,\bar{L}'$
are quasi-isomorphisms.}

\smallskip
In this case,
we deduce from lemma~\ref{CompositionMorphismHomology}
that $E^2(\phi): E^2(L)\,\rightarrow\,E^2(L')$
is an isomorphism.
The conclusion follows.

\subsection{Proof of assertion \textbf{(b)}}
\emph{We assume $\psi: P\,\rightarrow\,P'$ and $\phi: L\,\rightarrow\,L'$ are quasi-isomorphisms.}

\smallskip
By proposition~\ref{RightModuleSpectralSequence}
(and because $I^2(\bar{L}\circ P)$ lies in the first quadrant),
we obtain
$E^2_{s t}(L(n)) = \bigoplus_{-t\leq r\leq s} I^2_{s-r\,t+r}(\bar{L}(r)\otimes P^{\otimes r}(n))_{\Sigma_r}$.
Since $P$ is connected,
we have
$(\bar{L}(n)\otimes P^{\otimes n}(n))_{\Sigma_r} = \bar{L}(n)$
and
$(\bar{L}(r)\otimes P^{\otimes r}(n))_{\Sigma_r} = 0$ for $r>n$.
Consequently,
in the expansion of $E^2_{s t}(L(n))$,
we assume $r\leq n$.
Furthermore,
we have
\[I^2_{s-n\,t+n}(\bar{L}(n)\otimes P^{\otimes n}(n))_{\Sigma_n}
= \begin{cases} H_{s-n}(\bar{L}(n)) & \text{if $t=-n$}, \\
0, & \text{otherwise}. \end{cases}\]
We have similar results for the module $E^2_{s t}(L'(n))$.

The following assertions are immediate consequences of these observations.

\begin{fact}
We have $E^2_{s t}(L(n)) = E^2_{s t}(L'(n)) = 0$ for all bidegrees $(s,t)$ such that $t<-n$.
Moreover, for $t = -n$,
we obtain $E^2_{s\,-n}(L(n)) = H_{s-n}(\bar{L}(n))$ and $E^2_{s\,-n}(L'(n)) = H_{s-n}(\bar{L}'(n))$.
\end{fact}

We prove by induction that the morphism $\bar{\phi}: \bar{L}\,\rightarrow\,\bar{L}'$
induces an isomorphism $\bar{\phi}_*: H_*(\bar{L}(r))\,\rightarrow\,H_*(\bar{L}'(r))$
for all $r\in\N$.
Precisely,
we assume that $\bar{\phi}_*: H_*(\bar{L}(r))\,\rightarrow\,H_*(\bar{L}'(r))$
is an isomorphism in all degrees for $r<n$
and in degree $*<d$ for $r=n$.
We prove that $\bar{\phi}_*: H_*(\bar{L}(n))\,\rightarrow\,H_*(\bar{L}'(n))$
is also an isomorphism in degree $*=d$.

\begin{fact}
Under the induction assumption,
the morphism
\[E^2(\phi): E^2_{s t}(L(n))\,\rightarrow\,E^2_{s t}(L'(n))\]
is an isomorphism for all bidegrees $(s,t)$ such that $s<d+n$.
\end{fact}

\begin{proof}
If $r<n$,
then, by lemma~\ref{CompositionMorphismHomology},
the induction assumption implies that the morphism
\[I^2(\bar{\phi}\circ\psi): I^2_{s-r\,t+r}(\bar{L}(r)\otimes P{}^{\otimes r}(n))_{\Sigma_r}
\,\rightarrow\,I^2_{s-r\,t+r}(\bar{L}'(r)\otimes P'{}^{\otimes r}(n))_{\Sigma_r}\]
is an isomorphism for all $(s-r,t+r)$.
If $r=n$,
then we obtain an isomorphism for all bidegrees $(s-n,t+n)$ such that $s<d+n$,
because the components $I^2_{s-n\,t+n}$
are identified with homology modules of $\bar{L}(n)$
and $\bar{L}'(n)$.
We conclude that $E^2(\phi): E^2_{s t}(L(n))\,\rightarrow\,E^2_{s t}(L'(n))$
is an isomorphism for all $(s,t)$ such that $s<d+n$.
\end{proof}

\begin{fact}
The morphism $E^2(\phi): E^2_{s t}(L(n))\,\rightarrow\,E^2_{s t}(L'(n))$
is also an isomorphism in bidegree $(s,t)=(d+n,-n)$.
\end{fact}

\begin{proof}
We borrow a classical argument.
We consider the mapping cylinder of $\phi: L(n)\,\rightarrow\,L'(n)$
and the associated spectral sequence
$E^r(C_\phi)\,\Rightarrow\,H_*(C_\phi)$.
We have a long exact sequence
\[\cdots\,\rightarrow\,E^2_{s t}(L(n))\,\rightarrow\,E^2_{s t}(L'(n))
\,\rightarrow\,E^2_{s t}(C_\phi)\,\rightarrow\,E^2_{s-1 t}(L(n))\,\rightarrow\,\cdots\]
because $(E^1(C_\phi),d^1)$ is identified with the mapping cylinder of the morphism
$E^1(\phi): E^1(L(n))\,\rightarrow\,E^1(L'(n))$.

Since $E^2_{s t}(L(n)) = E^2_{s t}(L'(n)) = 0$ for $t<-n$,
we have $E^2_{s t}(C_\phi) = 0$ for $t<-n$.
In addition,
the previous assertion implies that the modules $E^2_{s t}(C_\phi)$ are zero for $s<d+n$.
Consequently,
we have $E^2_{d+n\,-n}(C_\phi) = E^\infty_{d+n\,-n}(C_\phi)$.
This module is zero since $\phi: L(n)\,\rightarrow\,L'(n)$
is supposed to be a quasi-isomorphism.

The conclusion follows.
\end{proof}

We conclude from the assertion above that the morphism $\psi_*: H_*(\bar{L}(n))\,\rightarrow\,H_*(\bar{L}'(n))$
is an isomorphism in degree $*=d$.
Hence,
this result achieves the induction process and the proof of assertion \textbf{(b)} of the theorem.

\subsection{Proof of assertion  \textbf{(c)}}
\emph{We assume $\phi: L\,\rightarrow\,L'$ and $\bar{\phi}: \bar{L}\,\rightarrow\,\bar{L}'$
are quasi-isomorphisms.}

\smallskip
We have
$E^2_{s t}(L(n)) = \bigoplus_{-t\leq r\leq s} I^2_{s-r\,t+r}(\bar{L}(r)\otimes P^{\otimes r}(n))_{\Sigma_r}$.
We assume $\bar{L}(0) = 0$ and $\bar{L}(1) = \K$.
Consequently,
in the expansion above,
we may assume $r\geq 1$.
Moreover,
for $r = 1$,
we obtain
\[I^2_{s-1\,t+1} = \begin{cases} H_t(P(n)), & \text{if $s=1$}, \\
0, & \text{otherwise}. \end{cases}\]
We have similar results for the module $E^2_{s t}(L'(n))$.

The following assertions are immediate consequences of these observations.

\begin{fact}
We have $E^2_{s t}(L(n)) = E^2_{s t}(L(n)) = 0$ for all bidegrees $(s,t)$ such that $s<1$.
Moreover, for $s = 1$, we obtain $E^2_{1 t}(L(n)) = H_t(P(n))$ and $E^2_{1 t}(L'(n)) = H_t(P'(n))$.
\end{fact}

We prove by induction that $\psi: P\,\rightarrow\,P'$
induces a homology isomorphism
$\psi_*: H_*(P(n))\,\rightarrow\,H_*(P'(n))$
for all $n\in\N$.
Precisely,
we assume that $\psi_*: H_*(P(m))\,\rightarrow\,H_*(P'(m))$
is an isomorphism in all degrees for $m<n$
and in degree $*<d$ for $m=n$.
We prove that $\psi_*: H_*(P(n))\,\rightarrow\,H_*(P'(n))$
is also an isomorphism in degree $*=d$.

\begin{fact}
Under the induction assumption,
the morphism
\[E^2(\phi): E^2_{s t}(L(n))\,\rightarrow\,E^2_{s t}(L'(n))\]
is an isomorphism for all bidegrees $(s,t)$ such that $t<d$.
\end{fact}

\begin{proof}
If $r>1$,
then, by lemma~\ref{CompositionMorphismHomology},
the induction assumption implies that the morphism
\[I^2(\bar{\phi}\circ\psi): I^2_{s-r\,t+r}(\bar{L}(r)\otimes P{}^{\otimes r}(n))_{\Sigma_r}
\,\rightarrow\,I^2_{s-r\,t+r}(\bar{L}'(r)\otimes P'{}^{\otimes r}(n))_{\Sigma_r}\]
is an isomorphism for all $(s-r,t+r)$.
If $r=1$,
then we obtain an isomorphism for all bidegrees $(s-1,t+1)$ such that $t<d$,
because the components $I^2_{s-1\,t+1}$
are identified with homology modules of $\bar{L}(n)$
and $\bar{L}'(n)$.
We conclude that $E^2(\phi): E^2_{s t}(L(n))\,\rightarrow\,E^2_{s t}(L'(n))$
is an isomorphism for $t<d$.
\end{proof}

As in the proof of assertion \textbf{(b)},
the next property follows from a straighforward inspection of the spectral sequence.

\begin{fact}
The morphism $E^2(\phi): E^2_{s t}(L(n))\,\rightarrow\,E^2_{s t}(L'(n))$
is also an isomorphism in bidegree $(s,t)=(1,d)$.
\end{fact}

We conclude from this assertion that $\psi_*: H_*(P(n))\,\rightarrow\,H_*(P'(n))$ is an isomorphism in degree $*=d$.
Hence,
the result above achieves the induction process
and the proof of assertion \textbf{(c)} of the theorem.

%\input{ReducedBar}
% 30/1/2003

\chapter{The reduced bar construction}\label{ReducedBar}

We recall the definition of the \emph{reduced bar construction},
which is introduced
by E. Getzler and J. Jones in \cite{GetzlerJones}
and
by V. Ginzburg and M. Kapranov in \cite{GinzburgKapranov}.

The general bar construction (defined in chapter~\ref{BarCoefficients})
is a complex of $\Sigma_*$-modules $B(L,P,R)$
associated to a dg-operad $P$
with coefficients in a right $P$-module $L$ and a left $P$-module $R$.
The reduced bar construction is the bar construction with trivial coefficients $\bar{B}(P) = B(I,P,I)$.

We adopt the conventions of E. Getzler and J. Jones \cite{GetzlerJones}.
Consequently,
the bar construction $\bar{B}(P)$ is the dual module of the quasi-free operad denoted by $C(P)$
in \cite[V. Ginzburg and M. Kapranov]{GinzburgKapranov}.

\section[Summary]{Summary: quasi-free operads and reduced bar constructions}

\subsection{The free operad in the dg-context}\label{dgFreeOperad}
We recall the structure of the free operad $F(M)$ in the context of a dg-$\Sigma_*$-module $M$.
We make explicit a realization of $F(M)$ in section~\ref{FreeOperadStructure}.
In the dg-context,
the free operad $F(M)$ forms a dg-operad.
Equivalently,
we have a canonical differential $\delta: F(M)\,\rightarrow\,F(M)$ induced by the internal differential of $M$.
The natural splitting
$F(M) = \bigoplus_{r=0}^\infty F_{(r)}(M)$
introduced in paragraph~\ref{FreeOperad}
is preserved by the differential.
We have explicitly $\delta(F_{(r)}(M))\subset F_{(r)}(M)$,
for all $r\in\N$.

Let us recall that $F_{(0)}(M) = I$.
The inclusion
$I = F_{(0)}(M)\,\subset\,F(M)$
gives the unit of the free operad.
The projection
$F(M)\,\rightarrow\,F_{(0)}(M) = I$
defines an augmentation morphism $F(M)\,\rightarrow\,I$.
Accordingly, the free operad has an augmentation ideal $\tilde{F}(M)\subset F(M)$
defined by $\tilde{F}(M) = \bigoplus_{r=1}^{\infty} F_{(r)}(M)$.

We have in addition $F_{(1)}(M) = M$.
The inclusion
$M = F_{(1)}(M)\,\rightarrow\,F(M)$
defines the universal morphism $M\,\rightarrow\,F(M)$.
The projection
$F(M)\,\rightarrow\,F_{(1)}(M) = M$
allows to identify the dg-module $M$
with the indecomposable quotient of $F(M)$.
Let us observe that the indecomposable quotient of a dg-operad
forms a dg-$\Sigma_*$-module.
The projection
$F(M)\,\rightarrow\,F_{(1)}(M) = M$
is a morphism of chain complexes.
Therefore,
the differential of indecomposable elements induced by the quotient morphism
agrees with the internal differential of $M$.

We observe that the free operad functor has nice homological properties
like the composition product of $\Sigma_*$-modules.
For instance,
if the ground ring $\K$ is a field,
then we obtain the following result:

\begin{prp}\label{FreeOperadHomology}
We assume that $M$ is a connected dg-$\Sigma_*$-mo\-du\-le
defined over a field $\K$ of characteristic $p\geq 0$.
We have an isomorphism of graded operads $F(H_*(M))\simeq H_*(F(M))$,
which extends the natural morphism of graded $\Sigma_*$-modules $H_*(M)\,\rightarrow\,H_*(F(M))$
induced by the inclusion morphism $M\,\rightarrow\,F(M)$.
\end{prp}

We deduce this statement from the expansion of the free operad.
We refer to paragraph~\ref{ReducedFreeOperad} for details about this result.

\subsection{The cofree cooperad}\label{CofreeCooperad}
There is also a cofree cooperad, denoted by $F^c(M)$,
which satisfies a classical universal property
dual to that of a free operad.
In fact,
because of our definition of a cooperad (\emph{cf}. paragraph~\ref{Cooperads}),
the functors $F(M)$ and $F^c(M)$ have the same expansion
and are isomorphic as $\Sigma_*$-modules.
In particular,
we have a natural splitting
$F^c(M) = \bigoplus_{r=0}^{\infty} F^c_{(r)}(M)$
such that
$F^c_{(0)}(M) = I$ and $F^c_{(1)}(M) = M$.
The projection
$F^c(M)\,\rightarrow\,F^c_{(0)}(M) = I$
gives the counit of $F^c(M)$.
The inclusion
$I = F^c_{(0)}(M)\,\rightarrow\,F^c(M)$
gives a coaugmentation morphism.
Accordingly,
the cofree operad has a coaugmentation coideal $\tilde{F}^c(M)$
defined by $\tilde{F}^c(M) = \bigoplus_{r=1}^{\infty} F^c_{(r)}(M)$.

We recall that the dual $\Sigma_*$-module $D^{\vee}$ of a cooperad $D$
is equipped with the structure of an operad
(\emph{cf}. proposition \ref{CoalgebraDuality}).
In the case of a cofree cooperad $D = F^c(M)$,
we obtain the following straighforward result:

\begin{prp}\label{FreeOperadDuality}
Suppose given a $\Sigma_*$-module $M$
such that the sequence $M(r)$ consists of finitely generated projective $\K$-modules.
We assume in addition $M(0) = M(1) = 0$.
In this situation,
the components of the cofree cooperad $F^c(M)(r)$ are finitely projective $\K$-modules
(and so are the components of the free operad $F(M)(r)$).
Furthermore, we have an operad isomorphism $F(M^{\vee})\simeq F^c(M)^{\vee}$
which extends the natural morphism of $\Sigma_*$-modules $M^{\vee}\,\rightarrow\,F^c(M)^{\vee}$
dual to the canonical projection $F^c(M)\,\rightarrow\,M$.
\end{prp}

\subsection{Operad derivations}\label{OperadDerivationDefinition}
An operad derivation is a homogeneous morphism $d: P\,\rightarrow\,P$
that satisfies the identity
\[d(p(q_1,\ldots,q_r)) = d(p)(q_1,\ldots,q_r)
+ \sum_{i=1}^r \pm p(q_1,\ldots,d(q_i),\ldots,q_r).\]
We have a dual notion for cooperads.

As an example,
the differential of a dg-operad $\delta: P\,\rightarrow\,P$ is a homogeneous derivation of degree $-1$
such that $\delta^2 = 0$.
We observe that a morphism $\delta+d: P\,\rightarrow\,P$,
where $d: P\,\rightarrow\,P$ is a homogeneous derivation of degree $-1$,
is still a derivation of $P$.
Accordingly,
the morphism $\delta+d: P\,\rightarrow\,P$ defines a differential
(not equivalent to the original one)
provided we have $(\delta+d)^2 = 0$
(or, equivalently, $\delta d + d\delta + d^2 = 0$).

We consider such differentials $\delta+d: P\,\rightarrow\,P$
in the case of a free operad $P = F(M)$.
In this situation,
we have the following result:

\begin{lem}\label{FreeOperadDerivation}
Let $P = F(M)$ be a free dg-operad.
For a given homogeneous morphism $\theta: M\,\rightarrow\,F(M)$,
there is one and only one homogeneous derivation $d_\theta: F(M)\,\rightarrow\,F(M)$
whose restriction to $M\subset F(M)$
agrees with $\theta: M\,\rightarrow\,F(M)$.
Furthermore,
if $\theta(M)\subset F_{(r)}(M)$,
then we have $d_\theta(F_{(s)}(M))\subset F_{(r+s-1)}(M)$,
for all $s>0$.

Dually,
let $D = F^c(M)$ be a cofree dg-cooperad.
For a given homogeneous morphism $\theta: F^c(M)\,\rightarrow\,M$,
there is one and only one homogeneous coderivation
$d_\theta: F^c(M)\,\rightarrow\,F^c(M)$
whose projection onto $M\subset F(M)$
agrees with $\theta: F^c(M)\,\rightarrow\,M$.
Furthermore,
if the morphism $\theta: F^c(M)\,\rightarrow\,M$
vanishes over components $F^c_{(s)}(M)\subset F^c(M)$, such that $s\not=r$,
then we have $d_\theta(F^c_{(s+r-1)}(M))\subset F^c_{(s)}(M)$,
for all $s>0$.
\end{lem}

We make explicit the derivation $d_\theta: F(M)\,\rightarrow\,F(M)$
determined by a homogeneous morphism $\theta: M\,\rightarrow\,F(M)$
in paragraph~\ref{TreeDerivation}
and the coderivation $d_\theta: F^c(M)\,\rightarrow\,F^c(M)$
in paragraph~\ref{TreeCoderivation}.
The proof of lemma~\ref{FreeOperadDerivation}
reduces then to straightfoward verifications (wich are omitted).

\subsection{Quasi-free operads}\label{QuasiFreeOperads}
A \emph{quasi-free operad} is a dg-operad $P$
such that $P = F(M)$,
but whose differential,
denoted by $\delta_\theta: F(M)\,\rightarrow\,F(M)$,
differs from the canonical differential of a free operad.
We have then $\delta_\theta = \delta + d_\theta$,
where $d_\theta: F(M)\,\rightarrow\,F(M)$ is a derivation.

By lemma~\ref{FreeOperadDerivation},
this derivation is determined by a homogeneous morphism $\theta: M\,\rightarrow\,F(M)$.
In particular,
we observe that the inclusion morphism $M\,\rightarrow\,F(M)$
is not a morphism of dg-modules
as long as $\theta: M\,\rightarrow\,F(M)$
is non zero.
But,
the projection morphism $F(M)\,\rightarrow\,M$ is a morphism of dg-modules
provided we have the relation
\[d_\theta(F(M))\subset\bigoplus_{r\geq 2} F_{(r)}(M).\]
In this situation,
the indecomposable quotient of the quasi-free operad $P = F(M)$
is isomorphic to the dg-$\Sigma_*$-module $M$.
By lemma~\ref{FreeOperadDerivation},
the property above holds if and only if we have $\theta(M)\subset\bigoplus_{r\geq 2} F_{(r)}(M)$.

\subsection{Quasi-cofree cooperads}
Dually,
a \emph{quasi-cofree cooperad} is a dg-cooperad $D$
such that $D = F^c(M)$,
but whose differential,
denoted by $\delta_\theta: F^c(M)\,\rightarrow\,F^c(M)$,
differs from the canonical differential of a cofree cooperad.
We have then $\delta_\theta = \delta + d_\theta$,
where $d_\theta: F^c(M)\,\rightarrow\,F^c(M)$
is a coderivation.
By lemma~\ref{FreeOperadDerivation},
this coderivation is determined by a homogeneous morphism $\theta: F^c(M)\,\rightarrow\,M$.
Consequently,
the projection morphism $F^c(M)\,\rightarrow\,M$
is not a morphism of dg-modules
as long as $\theta: F^c(M)\,\rightarrow\,M$
is non zero.
But,
the inclusion morphism $M\,\rightarrow\,F^c(M)$ is a morphism of dg-modules
provided $\theta: F^c(M)\,\rightarrow\,M$
vanishes over $F^c_{(1)}(M)\subset F^c(M)$.

\subsection{The reduced bar construction}\label{ReducedBarPpties}
We assume that $P$ is a connected dg-operad.
We consider the augmentation ideal of $P$,
defined by $\tilde{P}(r) = 0$ for $r = 0,1$ and $\tilde{P}(r) = P(r)$ for $r\geq 2$
(see paragraph~\ref{ConnectedOperads}).
The \emph{reduced bar construction} of $P$ is a quasi-cofree cooperad $\bar{B}(P) = F^c(M)$,
such that $M = \Sigma\tilde{P}$.
We recall that the suspension of $\tilde{P}(r)$ is the dg-module
such that $\Sigma\tilde{P}(r)_* = \tilde{P}(r)_{*-1}$
(see paragraph~\ref{Suspension}).
In particular, if $P$ is a non-graded operad,
then $\Sigma\tilde{P}(r)$ has the module $\tilde{P}(r)$ in degree $*=1$
and vanishes in degree $*\not=1$.

The reduced bar construction is equipped with the canonical differential of the cofree cooperad
$\delta: F^c(\Sigma\tilde{P})\,\rightarrow\,F^c(\Sigma\tilde{P})$,
which is induced by the internal differential of $P$,
and with the bar differential
$\beta: F^c(\Sigma\tilde{P})\,\rightarrow\,F^c(\Sigma\tilde{P})$,
which is determined by the operad composition product of $P$.
To be precise,
we observe that the partial composition products of an operad are equivalent to a homogeneous morphism
$\theta: F^c_{(2)}(\Sigma\tilde{P})\,\rightarrow\,\Sigma\tilde{P}$
of degree $-1$.
The bar differential is the coderivation $\beta = d_\theta$ determined by this morphism.
One proves that the bar differential verifies the identities $\beta^2 = 0$
and $\delta\beta + \beta\delta = 0$.
We give the details of this construction in section~\ref{ConstructionBarDifferential}.

We observe more precisely that the reduced bar construction defines a complex of dg-modules
\[\cdots\,\xrightarrow{\beta}\,\bar{B}_{s}(P)\,\xrightarrow{\beta}\,\bar{B}_{s-1}(P)\,\xrightarrow{\beta}
\,\cdots\,\xrightarrow{\beta}\,\bar{B}_{0}(P).\]
Explicitly,
we consider the grading $\bar{B}(P) = \bigoplus_{s=0}^\infty \bar{B}_{s}(P)$
such that $\bar{B}_{s}(P) = F^c_{(s)}(\Sigma\tilde{P})$.
We recall that the internal differential verifies the relation
\[\delta(F^c_{(s)}(\Sigma\tilde{P}))\subset F^c_{(s)}(\Sigma\tilde{P}).\]
By lemma~\ref{FreeOperadDerivation},
we have in addition
\[\beta(F^c_{(s)}(\Sigma\tilde{P}))\subset F^c_{(s-1)}(\Sigma\tilde{P}).\]
The conclusion follows.

\subsection{The reduced cobar construction}\label{ReducedCobarPpties}
Dually, we assume that $D$ is a connected dg-cooperad.
The \emph{reduced cobar construction} of $D$ is a quasi-free operad $\bar{B}^c(D) = F(M)$,
where $M = \Sigma^{-1}\tilde{D}$ is the desuspension of the coaugmentation coideal of $D$.

The reduced cobar construction is equipped with the canonical differential of the cofree cooperad
$\delta: F(\Sigma^{-1}\tilde{D})\,\rightarrow\,F(\Sigma^{-1}\tilde{D})$,
which is induced by the internal differential of $D$,
and with the cobar differential
$\beta: F(\Sigma^{-1}\tilde{D})\,\rightarrow\,F(\Sigma^{-1}\tilde{D})$,
which is determined by the cooperad coproduct of $D$.
We have $\beta^2 = 0$ and $\delta\beta + \beta\delta = 0$.

The reduced cobar construction defines a (negatively graded) complex of dg-modules
\[\bar{B}^c_{0}(D)\,\xrightarrow{\beta}\,\cdots
\,\xrightarrow{\beta}\,\bar{B}^c_{-s}(D)\,\xrightarrow{\beta}\,\bar{B}^c_{-s-1}(D)\,\xrightarrow{\beta}\,\cdots,\]
where $\bar{B}^c_{-s}(D) = F_{(s)}(\Sigma^{-1}\tilde{C})$,
since we have the relations
\begin{align*} & \delta(F_{(s)}(\Sigma^{-1}\tilde{D}))\subset F_{(s)}(\Sigma^{-1}\tilde{D}) \\
\text{and}\qquad & \beta(F_{(s)}(\Sigma^{-1}\tilde{D}))\subset F_{(s+1)}(\Sigma^{-1}\tilde{D}). \end{align*}

\smallskip
The next proposition is a consequence of proposition~\ref{FreeOperadDuality}:

\begin{prp}\label{LinearDualBarCooperad}
We assume that $P$ is a connected operad such that the sequence $P(r)$
consists of finitely generated projective $\K$-mo\-du\-les.
We claim that the cobar complex $(\bar{B}^c(P^{\vee}),\beta)$
is dual to the bar complex $(\bar{B}(P),\beta)$.
Accordingly,
we have a natural isomorphism of differential graded operads
$\bar{B}^c(P^{\vee})\simeq\bar{B}(P)^{\vee}$.
\end{prp}

To be precise,
one observes that the natural isomorphism of proposition~\ref{FreeOperadDuality}
\[F^c(\Sigma\tilde{P})^{\vee}\simeq F((\Sigma\tilde{P})^{\vee})\simeq F(\Sigma^{-1}\tilde{P}^{\vee})\]
makes the differential of the bar construction $B(P)$
dual to the differential of the cobar construction $B^c(P^{\vee})$.

\smallskip
In the topological context,
a construction of J. Boardman and R. Vogt,
called the \emph{$W$-construction} (\emph{cf}. \cite{BoardmanVogt}),
gives quasi-free resolutions for all topological operads.
In the differential graded framework,
one observes that the composite of the bar and cobar constructions
plays the role of the $W$-construction:

\begin{prp}[\emph{cf}. V. Ginzburg and M. Kapranov, \cite{GinzburgKapranov}]\label{CobarBarResolution}
Suppose given a connected dg-operad $P$.
If the ground ring $\K$ is not a field,
then we assume that $P$ is projective as a $\K$-module.
We consider the operad morphism
\[F(\Sigma^{-1}\tilde{F}^c(\Sigma\tilde{P}))\,\xrightarrow{}\,P\]
induced by the canonical projection
\[\Sigma^{-1} F^c(\Sigma\tilde{P})\,\xrightarrow{}\,\Sigma^{-1} F^c_{(1)}(\Sigma \tilde{P})\]
(recall that $\Sigma^{-1} F^c_{(1)}(\Sigma \tilde{P})\,\simeq\,\tilde{P}\,\subset\,P$).
This morphism defines a quasi-isomorphism of dg-operads
\[\bar{B}^c(\bar{B}(P))\,\xrightarrow{\sim}\,P.\]
\end{prp}

In section~\ref{ProofCobarBarResolution},
we observe that this proposition,
proved by V. Ginz\-burg and M. Kapranov in \emph{loco citato},
is an easy consequence
of the comparison theorems~\ref{RightModuleComparison} and~\ref{LeftModuleComparison}.

\section{Digression: the model category of operads}

We prove the following comparison theorem in section~\ref{ProofsComparisonQuasiFreeOperads}.

\begin{thm}\label{ComparisonQuasiFreeOperads}
We are given a morphism of dg-operads $\phi: P\,\rightarrow\,P'$.
We assume that $P$ (respectively, $P'$)
is a quasi-free operad $P = F(M)$ (respectively, $P' = F(M')$)
such that $M(0) = M(1) = 0$ (respectively, $M'(0) = M'(1) = 0$).
If the ground ring $\K$ is not a field,
then we assume in addition that $M$ and $M'$ are projective as $\K$-modules.
We assume that the morphism $\phi: P\,\rightarrow\,P'$ preserves augmentation ideals of free operads
(explicitly, $\phi(\tilde{F}(M))\subset\tilde{F}(M')$).
Furthermore,
we identify the $\Sigma_*$-module $M$ (respectively, $M'$)
with the indecomposable quotient of $P$ (respectively, $P'$).
Hence, in this situation,
the $\Sigma_*$-modules $M$ and $M'$ are equipped with a natural differential
and we have an induced morphism of dg-$\Sigma_*$-modules
$\bar{\phi}: M\,\rightarrow\,M'$.

The morphism of dg-operad $\phi: P\,\rightarrow\,P'$ is a quasi-isomorphism
if and only if the associated morphism of dg-$\Sigma_*$-modules $\bar{\phi}: M\,\rightarrow\,M'$
is so.
\end{thm}

The aim of this section is to interpret theorem~\ref{ComparisonQuasiFreeOperads}
within the language of model categories.
As for modules over an operad,
we would like to point out that this result is stronger than the assertion
which one expects from classical model structures.
We quote the following theorem,
which is similar to a classical result about the category of differential graded commutative algebras
in characteristic $0$
(\emph{cf}. D. Quillen \cite{QuillenBook, QuillenProc, QuillenRationalHomotopy}).

\begin{thm}[\emph{cf}. C. Berger and I. Moerdijk \cite{BergerMoerdijk}, V. Hinich, \cite{Hinich}]\label{ModelOperads}
The dg-operads $P$ which have $P(0) = 0$ forms a model category.
A weak-equivalence is a quasi-isomorphism of dg-operads.
A fibration is a morphism of dg-operads which is surjective in degree $*>0$.
A cofibration is a morphism of dg-operads which has the left lifting property with respect to fibrations.
\end{thm}

\subsection{Remark}\label{LiftingProperty}
In positive characteristic,
we have to restrict ourself to dg-operads $P$ such that $P(0) = 0$.
In order to justify this remark,
we consider a morphism $P\,\rightarrow\,P\vee F(M)$,
where $P\vee F(M)$ is the coproduct of an operad $P$ with a free operad $F(M)$.
Such a morphism has the left-lifting property with respect to fibrations
provided $M$ is an acyclic dg-$\Sigma_*$-module.
We prove that the morphism $P\,\rightarrow\,P\vee F(M)$ is not a weak-equivalence in general,
as long as $M(0)\not=0$.

To be precise,
we assume $M(0) = V$ and $M(r) = 0$ for $r>0$.
In this case,
we have $P\vee F(M) = P[V]$,
where $P[V]$ is the operad such that
\[P[V](n) = \bigoplus_{r=0}^\infty (P(r+n)\otimes V^{\otimes r})_{\Sigma_r}\]
(\emph{cf}. \cite{FresseLie}).
As an example,
in the case of the commutative operad $P = \C$,
the module $\C[V](n)$ is identified with the symmetric algebra $\C[V](n) = S(V)$
and has a non-trivial homology.
Our conclusion follows.

\smallskip
The identity operad $I$ is the initial object in the category of dg-operads.
We recall that a dg-operad $P$ is a cofibrant object if the initial morphism $I\,\rightarrow\,P$ is a cofibration.
We have the following result:

\begin{lem}[\emph{cf}. V. Hinich \cite{Hinich}]
A cofibrant object in the category of dg-operads
is a retract of a quasi-free operad $F = F(M)$
such that the sequence $M(r)$ consists of projective $\Sigma_r$-modules.
\end{lem}

In positive characteristic,
the representations of the symmetric\linebreak groups are not all projective objects.
Therefore,
the quasi-free operads are not all cofibrant objects
with respect to the classical model structure.
On the other hand,
there is a non-classical model structure
which makes all quasi-free operads cofibrant objects.
For that purpose,
one can adapt the definition of paragraph~\ref{NonClassicalModelStructure}.

\subsection{Cofibrant and quasi-free resolutions of dg-operads}\label{CofibrantOperadResolutions}
A \emph{cofibrant resolution} of a dg-operad $P$
consists of a cofibrant operad $Q$
together with a surjective quasi-isomorphism $Q\,\xrightarrow{\sim}\,P$.
A \emph{quasi-free resolution} of $P$
consists of a quasi-free operad $Q' = F(M')$
together with a surjective quasi-isomorphism $F(M')\,\xrightarrow{\sim}\,P$.
Any operad $P$
has a cofibrant resolution $Q\,\xrightarrow{\sim}\,L$
such that $Q$ is a quasi-free operad.
We have explicitly $Q = F(M)$,
where the sequence $M(r)$ consists of projective $\Sigma_r$-modules.
We would like to mention the following result:

\begin{prp}\label{OperadCofibrantResolutions}
Suppose given a connected dg-operad $P$.
If the ground ring $\K$ is not a field,
then we assume that $P$ is projective as a $\K$-module.
We let $\tilde{B}(P)$ denote the coaugmentation coideal of the bar construction of $P$.

If $Q = F(M)$ is a cofibrant quasi-free resolution of $P$
such that $M(0) = M(1) = 0$,
then we have a quasi-isomorphism of dg-$\Sigma_*$-modules
$M\,\xrightarrow{\sim}\,\Sigma^{-1}\tilde{B}(P)$.
Conversely,
if we are given a resolution of $\Sigma^{-1}\tilde{B}(P)$ by projective $\Sigma_*$-modules,
let $M\,\xrightarrow{\sim}\,\Sigma^{-1}\tilde{B}(P)$,
then we have a cofibrant quasi-free resolution of $P$ such that $Q = F(M)$.
\end{prp}

\begin{proof}
A cofibrant resolution $Q$ is connected to any quasi-free resolution
by a quasi-iso\-mor\-phism $Q\,\xrightarrow{\sim}\,Q'$,
because the left lifting property provides the following diagram
with a fill in morphism
\[\xymatrix{ & Q'\ar@{->>}[d]^{\sim} \\ Q\ar[r]^{\sim}\ar@{-->}[ur] & P \\ }\]
In this situation,
theorem~\ref{ComparisonQuasiFreeOperads} asserts that the induced morphism $M\,\rightarrow\,M'$
is a quasi-isomorphism.
In the case $Q' = \bar{B}^c(\bar{B}(P))$, we have $Q' = F(M')$,
where $M' = \Sigma^{-1}\tilde{B}(P)$.
Therefore,
we obtain a quasi-isomorphism $M\,\xrightarrow{\sim}\,\Sigma^{-1}\tilde{B}(P)$.

The second assertion of the lemma generalizes a known result
which holds in characteristic $0$
(\emph{cf}. M. Kontsevich and Y. Soibelman \cite{KontsevichSoibelman})
and is related to homotopy invariance principles
(\emph{cf}. J. Boardman and R. Vogt \cite{BoardmanVogt}).
We omit the proof of this assertion
which is mentioned as a remark.
\end{proof}

\section{The language of trees}\label{Trees}

We borrow the construction of the free operad from \cite[V. Ginz\-burg and M. Kapranov]{GinzburgKapranov}.
We need to make this construction very precise.
Therefore,
we introduce some conventions about tree structures

\subsection{On the structure of a tree}\label{TreeStructure}
An \emph{$n$-tree} is an oriented abstract tree
together with one outgoing edge (the \emph{root} of the tree)
and
$n$ ingoing edges (the \emph{entries} of the tree)
indexed by the set $\{1,\ldots,n\}$.

Formally,
an $n$-tree $\tau$ is specified by a set of vertices $V(\tau)$
and by a set of edges $E(\tau)$.
An edge $e\in E(\tau)$ is oriented
from a source $s(e)\in V(\tau)\amalg\{1,\ldots,n\}$
to a target $t(e)\in V(\tau)\amalg\{0\}$.
For any vertex $v\in V(\tau)$,
we assume that there is one and only one edge
$e\in E(\tau)$
such that $s(e) = v$.
Similarly, for $i\in\{1,\ldots,n\}$,
there is one and only one edge such that $s(e) = i$.
This edge is the $i$th entry of the tree.
There is also one and only one edge
such that $t(e) = 0$.
This edge is the root of the tree.
The entries and the root are the \emph{external edges} of $\tau$.
The \emph{internal edges} verify $s(e)\in V(\tau)$ and $t(e)\in V(\tau)$.

As long as we omit to fix the set of egdes (which is not important for our purposes),
the tree structure is equivalent to a partition
$V(\tau)\amalg\{1,\ldots,n\} = \coprod_{v\in V(\tau)\amalg\{0\}} I_v$.
Explicitly,
for $v\in V(\tau)\amalg\{0\}$,
we define
$I_v = \{\,s(e),\ \text{where}\ e\in E(\tau)\ \text{verifies}\ t(e) = v\,\}$,
so that there is an edge going from $v'$ to $v$
if and only if $v'\in I_v$.

\subsection{Example}
\begin{figure}
\[\xymatrix@M=0pt@R=10mm@C=10mm{ & & *+<3mm>{1}\ar[dr]|-*+<1mm>{e_6} & &
*+<3mm>{5}\ar[dl]|-*+<1mm>{e_7} & & \\
& *+<3mm>{3}\ar[drr]|-*+<1mm>{e_4} & &
*+<3mm>[o][F]{v_3}\ar[d]|-*+<1mm>{e_5} &
*+<3mm>{4}\ar[dr]|-*+<1mm>{e_9} & &
*+<3mm>{6}\ar[dl]|-*+<1mm>{e_{10}} \\
*+<3mm>{2}\ar[drrr]|-*+<1mm>{e_2} & & &
*+<3mm>[o][F]{v_2}\ar[d]|-*+<1mm>{e_3} & &
*+<3mm>[o][F]{v_4}\ar[dll]|-*+<1mm>{e_8} & \\
& & & *+<3mm>[o][F]{v_1}\ar[d]|-*+<1mm>{e_1} & & & \\
& & & *+<3mm>{0} & & & \\ }\]
\caption{}\label{TreeExample}\end{figure}
The tree of figure~\ref{TreeExample}
has the set of vertices
$V(\tau) = \{v_1,v_2,v_3,v_4\}$
and the set of edges
$E(\tau) = \{e_1,\ldots,e_{10}\}$.
For each vertex $v_i\in V(\tau)$, we obtain precisely:
\begin{gather*} I_{0} = \{v_1\},\quad I_{v_1} = \{2,v_2,v_4\}, \\
I_{v_2} = \{3,v_3\},
\quad I_{v_3} = \{1,5\}
\quad\text{and}\quad I_{v_4} = \{4,6\}.\end{gather*}
As an example,
the edge $e_4$ starts from $s(e_4) = 3$ and targets to $t(e_4) = v_2$.
Hence, we have $3\in I_{v_2}$.

In the next figures, we may omit the notation of edges and vertices.

\subsection{Tree isomorphisms}\label{TreeIsomorphisms}
An isomorphism of $n$-trees $f: \tau\,\rightarrow\,\tau'$
consists of bijections
\[f_V: V(\tau)\,\rightarrow\,V(\tau')\quad\text{and}\quad f_E: E(\tau)\,\rightarrow\,E(\tau')\]
which preserve the tree structure
(the source and the target of each edge).
In that situation,
we extend
$f_V: V(\tau)\,\rightarrow\,V(\tau')$
to a bijection
$f_V: V(\tau)\amalg\{0,1,\ldots,n\}\,\rightarrow\,V(\tau')\amalg\{0,1,\ldots,n\}$
such that $f_V(i) = i$ for $i\in\{0,1,\ldots,n\}$.
We have the relation $I_{f_V(v)} = f_V(I_{v})$,
for all $v\in V(\tau)\amalg\{0\}$.

\subsection{Indexing}\label{EntryIndexing}
We can assume
that the entries of an $n$-tree are indexed by a fixed set $I = \{i_1,\ldots,i_n\}$
with $n$ elements.
We obtain the category of $I$-trees.
In this context,
the tree structure is equivalent to a partition
\[V(\tau)\amalg I = \coprod_{v\in V(\tau)\amalg\{0\}} I_v.\]

We observe that a bijection $u: I\,\rightarrow\, I'$
induces a functor from the category of $I$-trees
to the category of $I'$-trees.
The tree $\tau' = u_*(\tau)$
associated to $\tau$
has the same set of vertices
and the same set of edges as $\tau$.
We just reindex the entries.
Precisely,
if $e$ is an entry of $\tau$,
then the source of $e$ in $\tau'$
is the element
$s'(e) = u(s(e))\in I'$.

In the case $I = \{1,\ldots,n\}$,
we obtain that a permutation $w\in\Sigma_n$
induces a functor from the category of $n$-trees to the category of $n$-trees.
Hence,
we conclude that the category of $n$-trees is equipped with an action of the symmetric group $\Sigma_n$.

\subsection{Subtrees}\label{SubTrees}
A subtree of $\tau$ is a tree $\sigma\subset\tau$
such that $V(\sigma)\subset V(\tau)$
and
$E(\sigma)\subset E(\tau)$.
To be precise,
we assume that an edge $e\in E(\tau)$ belongs to $E(\sigma)\subset E(\tau)$
if and only if $s(e)\in V(\sigma)$ or $t(e)\in V(\sigma)$.
The source (respectively, the target) of an edge in $\sigma$
is given
by the source (respectively, the target) of this edge in the tree $\tau$.
(Consequently, the structure of the subtree is determined by the vertex set $V(\sigma)\subset V(\tau)$.)

By definition,
the root of $\sigma$ is an edge $e\in E(\tau)$
that verifies $s(e)\in V(\sigma)$
but $t(e)\not\in V(\sigma)$.
We assume that there is one and only one edge $e\in E(\tau)$
which has this property.
Similarly,
an entry of $\sigma$ is an edge $e\in E(\tau)$
that verifies $t(e)\in V(\sigma)$
but $s(e)\not\in V(\sigma)$.
The internal edges of $\sigma$ verify $s(e)\in V(\sigma)$
and $t(e)\in V(\sigma)$.

The entries of $\sigma$ are naturally indexed by a subset of $V(\tau)\setminus V(\sigma)\amalg\{1,\ldots,n\}$
which we denote by $I_\sigma$.
To be precise,
the index of an entry, which is associated to an edge $e\in E(\tau)$,
is given by the element $s(e)\in V(\tau)\setminus V(\sigma)\amalg\{1,\ldots,n\}$.

As an example, the vertices $V(\sigma) = \{v_1,v_2\}$
of the tree $\tau$ of figure~\ref{TreeExample}
generate a subtree of $\tau$.
We have $E(\sigma) = \{e_1,e_2,e_4,e_5,e_8\}$
and the entries of $\sigma$ are indexed
by the set $\{2,3,v_3,v_4\}$.

\subsection{Quotient trees}\label{QuotientTrees}
We obtain a quotient tree $\tau/\sigma$
by collapsing a subtree $\sigma\subset\tau$ to a vertex.
We define explicitly
$V(\tau/\sigma) = V(\tau)\setminus V(\sigma)\amalg\{\sigma\}$.
The elements of $E(\tau/\sigma)$
are the edges $e\in E(\tau)$
such that $s(e)\not\in V(\sigma)$ or $t(e)\not\in V(\sigma)$.
(Thus, we remove the internal edges of the subtree $\sigma\subset\tau$.)
In the quotient tree $\tau/\sigma$,
the entries of $\sigma$ targets to the collapsed vertex $\sigma\in V(\tau/\sigma)$;
the root of $\sigma$ starts from the collapsed vertex $\sigma\in V(\tau/\sigma)$;
sources and targets of other edges $e\in E(\tau)\setminus E(\sigma)$ are unchanged.
In the example above,
we obtain the tree of figure~\ref{QuotientTreeExample}.
\begin{figure}
\[\xymatrix@M=0pt@R=10mm@C=10mm{ & & &
*+<3mm>{1}\ar[dr]|-*+<1mm>{e_6} & &
*+<3mm>{5}\ar[dl]|-*+<1mm>{e_7} &
*+<3mm>{4}\ar[dr]|-*+<1mm>{e_9} & &
*+<3mm>{6}\ar[dl]|-*+<1mm>{e_{10}} \\
*+<3mm>{2}\ar[drrrr]|-*+<1mm>{e_2} & &
*+<3mm>{3}\ar[drr]|-*+<1mm>{e_4} & &
*+<3mm>[o][F]{v_3}\ar[d]|-*+<1mm>{e_5} & & &
*+<3mm>[o][F]{v_4}\ar[dlll]|-*+<1mm>{e_8} & \\
& & & & *+<4mm>[o][F]{\sigma}\ar[d]|-*+<1mm>{e_1} & & & & \\
& & & & *+<3mm>{0} & & & & \\ }\]
\caption{}\label{QuotientTreeExample}\end{figure}

Observe that the entry set of the collapsed vertex $I_\sigma\subset V(\tau/\sigma)\amalg\{1,\ldots,n\}$
is the same as the index set
$I_\sigma\subset V(\tau)\setminus V(\sigma)\amalg\{1,\ldots,n\}$
defined in paragraph~\ref{SubTrees}.

\subsection{Blowing up vertices}\label{BlowingUpVertices}
We have an operation inverse of the quotient process.
Namely,
by blowing up a vertex in $\tau$,
we obtain a tree $\tau'$
such that $\tau = \tau'/\sigma'$.

Precisely,
suppose given a vertex $v\in V(\tau)$
and a tree $\sigma'$
whose entries are in bijection with the entries of $v$.
(Hence, we assume that $\sigma'$ is an $I_v$-tree.)
In order to obtain $\tau'$,
we replacing the vertex $v$ in $\tau$
by the tree $\sigma'$.
More formally, the set of vertices of $\tau'$ is the sum
$V(\tau') = V(\tau)\setminus\{v\}\amalg V(\sigma')$.
The set of edges is the sum
$E(\tau') = E(\tau)\amalg E(\sigma')/\sim$
in which:
the root of $\sigma'$ is identified with the egde $e\in E(\tau)$ such that $s(e) = v$;
the entry of $\sigma'$ indexed by $x\in I_v$ is identified with the edge $e\in E(\tau)$
such that $t(e) = v$ and $s(e) = x$.
We have clearly $\tau'/\sigma'\simeq\tau$.

\section{Trees and free operads}\label{FreeOperadStructure}

In this section,
we recall the construction of the free operad,
made precise within the formalism of trees.
This representation, which we borrow from \cite[V. Ginzburg and M. Kapranov]{GinzburgKapranov},
goes back to J. Boardman and R. Vogt \cite{BoardmanVogt}.

\subsection{The $\Sigma_*$-category of trees}
We let $T(n)$ denote the category of $n$-trees with isomorphisms as morphisms
(see paragraphs~\ref{TreeStructure} and~\ref{TreeIsomorphisms}).
We have a weight splitting
\[T(n) = \coprod_{r=0}^\infty T_{(r)}(n),\]
where the category $T_{(r)}(n)$ is generated by trees
which have $r$ vertices.

We consider also the category $T(I) = \coprod_{r=0}^\infty T_{(r)}(I)$
which consists of trees
whose entries are indexed by a given set with $n$ elements $I$
(see paragraph~\ref{EntryIndexing}).
We have then $T(n) = T(\{1,\ldots,n\})$.
A bijection $u: I\,\rightarrow\,I'$ induces a functor $u_*: T(I)\,\rightarrow\,T(I')$
(see paragraph~\ref{EntryIndexing}).
Furthermore,
we have clearly $u_*(T_{(r)}(I))\subset T_{(r)}(I')$.

We conclude that the sequence $T(n)$, $n\in\N$, defines a (graded) $\Sigma_*$-object
in the category of small categories.

\subsection{The operad of trees}\label{TreeOperad}
We equip the sequence of categories $T(n)$, $n\geq 0$, with an operad structure.
Explicitly,
we define a partial composition product
\[\circ_i: T_{(r)}(m)\times T_{(s)}(n)\,\rightarrow\,T_{(r+s)}(m+n-1).\]

The composite tree $\sigma\circ_i\tau\in T(m+n-1)$, where $\sigma\in T_{(r)}(m)$ and $\tau\in T_{(s)}(n)$,
is obtained by gluing the root of $\tau$ to the $i$th entry of $\sigma$.
More formally,
the set of vertices of $\sigma\circ_i\tau$ is the sum $V(\sigma\circ_i\tau) = V(\sigma)\amalg V(\tau)$.
The set of edges of $\sigma\circ_i\tau$ is the sum $E(\sigma\circ_i\tau) = E(\sigma)\amalg E(\tau)/\sim$
in which we identify the $i$th entry of $\sigma$, let $e\in E(\sigma)$,
with the root of $\tau$, let $f\in E(\tau)$.
In $\sigma\circ_i\tau$,
this edge starts from $s(f)\in V(\tau)$ and targets to $t(e)\in V(\sigma)$.
We modify the indices of the entries of $\sigma\circ_i\tau$
according to the rule of paragraph~\ref{PartialComposites}.
Precisely,
the edges $e\in E(\sigma)$ such that $s(e) = k<i$
are indexed by the same element $s(e) = k$ in the composite tree;
the edges $e\in E(\sigma)$ such that $s(e) = k>i$
are indexed by $s(e) = k+n-1$;
the edges $e\in E(\tau)$ such that $s(e) = k$
are indexed by $s(e) = k+i-1$.
The source of an edge $e\in E(\sigma)$ (respectively, $e\in E(\tau)$)
is unchanged in $\sigma\circ_i\tau$
as long as $e$ rises from a vertex $v\in V(\sigma)$
(respectively, $v\in V(\tau)$).
Similarly, the target of an edge is unchanged as long as this edge targets to a vertex.

\begin{figure}
\begin{equation*}\begin{aligned}
\xymatrix@M=0pt@R=6mm@C=6mm{ & *+<3mm>{2}\ar[dr] & &
*+<3mm>{3}\ar[dl] \\
*+<3mm>{1}\ar[dr] & & *+<3mm>[o][F]{u_2}\ar[dl] & \\
& *+<3mm>[o][F]{u_1}\ar[d] & & \\ & *+<3mm>{0} & & \\ }
\end{aligned}\circ_2\begin{aligned}
\xymatrix@M=0pt@R=6mm@C=6mm{ *+<3mm>{1}\ar[dr] & &
*+<3mm>{2}\ar[dl] \\
& *+<3mm>[o][F]{v_1}\ar[d] & \\ & *+<3mm>{0} & \\ }
\end{aligned}=\begin{aligned}
\xymatrix@M=0pt@R=6mm@C=6mm{ *+<3mm>{2}\ar[dr] & & *+<3mm>{3}\ar[dl] & \\
& *+<3mm>[o][F]{v_1}\ar[dr] & & *+<3mm>{4}\ar[dl] \\
*+<3mm>{1}\ar[dr] & & *+<3mm>[o][F]{u_2}\ar[dl] & \\
& *+<3mm>[o][F]{u_1}\ar[d] & & \\
& *+<3mm>{0} & & \\ }\end{aligned}\end{equation*}
\caption{}\label{GraftingTrees}\end{figure}
As an example,
in figure~\ref{GraftingTrees},
we display a partial composite $\sigma\circ_2\tau\in T_{(3)}(4)$,
where $\sigma\in T_{(2)}(3)$ and $\tau\in T_{(1)}(2)$.

\subsection{The treewise tensor module}\label{TreewiseTensors}
We are given a $\Sigma_*$-module $M$.
A $\K$-module $\tau(M)$ is associated to each $I$-tree $\tau\in T(I)$.
An isomorphism of $I$-trees $f: \tau\,\rightarrow\,\tau'$
induces an isomorphism of $\K$-modules
$f_*: \tau(M)\,\rightarrow\,\tau'(M)$.
A bijection $u: I\,\rightarrow\,I'$
induces an isomorphism of $\K$-modules
$u_*: \tau(M)\,\rightarrow\,\tau'(M)$,
where $\tau' = u_*(\tau)$.

We set precisely:
\[\tau(M) = \bigotimes_{v\in V(\tau)} M(I_v).\]
Equivalently,
an element of $\tau(M)$
is a tensor $\bigotimes_{v\in V(\tau)} x_v$,
where $x_v\in M(I_v)$.
Concretely,
we mark each vertex $v\in V(\tau)$
with an element $x_v\in M(I_v)$.
Furthermore,
by the very definition of $M(I_v)$,
we suppose given a bijection from the entries of $x_v\in M(I_v)$
to the entries of $v\in V(\tau)$.
Consequently,
an element of $\tau(M)$
is represented by a labeled tree as in figure~\ref{TreewiseTensorExample}.
\begin{figure}
\[\xymatrix@M=0pt@R=10mm@C=10mm{ & & *+<3mm>{1}\ar[dr] & &
*+<3mm>{5}\ar[dl] & & \\
& *+<3mm>{3}\ar[drr] & &
*+<3mm>{x_3}\ar[d] &
*+<3mm>{4}\ar[dr] & &
*+<3mm>{6}\ar[dl] \\
*+<3mm>{2}\ar[drrr] & & &
*+<3mm>{x_2}\ar[d] & &
*+<3mm>{x_4}\ar[dll] & \\
& & & *+<3mm>{x_1}\ar[d] & & & \\
& & & *+<3mm>{0} & & & \\ }\]
\caption{}\label{TreewiseTensorExample}\end{figure}
In this example,
we consider the tree $\tau$ introduced in figure~\ref{TreeExample}.
We have then
\[\tau(M) = M(\{2,v_2,v_4\})\otimes M(\{3,v_3\})\otimes M(\{1,5\})\otimes M(\{4,6\})\]
and the display is the representation of a tensor $x_1\otimes x_2\otimes x_3\otimes x_4\in\tau(M)$,
\begin{gather*}\text{where}
\quad x_1\in M(\{2,v_2,v_4\}),
\quad x_2\in M(\{3,v_3\}), \\
\quad x_3\in M(\{1,5\})
\quad\text{and}\quad x_4\in M(\{4,6\}).\end{gather*}
According to the next definitions,
this treewise tensor represents the composite element $((x_1\circ_{v_4} x_4)\circ_{v_2} x_2)\circ_{v_3} x_3$
in the free operad.

\subsection{The expansion of the free operad}
We consider the colimit
\[F(M)(n) = \colim_{\tau\in T(n)}\tau(M),\]
where elements $x'\in\tau'(M)$ and $x\in\tau(M)$
which correspond under an isomorphism of $n$-trees $f: \tau\,\rightarrow\,\tau'$
are identified.
Explicitly,
we have $x'\equiv x$ in $F(M)$
if $x'=f_*(x)$.

To be more general,
a module $F(M)(I) = \colim_{\tau\in T(I)}\tau(M)$ is defined
for each given set with $n$ elements $I$.
A bijection $u: I\,\rightarrow\,I'$ determines an isomorphism $u_*: F(M)(I)\,\rightarrow\,F(M)(I')$,
since we have a natural morphism $u_*: \tau(M)\,\rightarrow\,u_*(\tau)(M)$,
for each tree $\tau\in T(I)$.
In particular,
a permutation $w\in\Sigma_n$ induces a morphism $w_*: F(M)(n)\,\rightarrow\,F(M)(n)$.

We conclude that the sequence $F(M)(n)$, $n\in\N$, forms a $\Sigma_*$-module.

\subsection{The weight components of the free operad}\label{WeightTreewiseTensors}
We have clearly a splitting $F(M)(n) = \bigoplus_{r=0}^\infty F_{(r)}(M)$,
where $F_{(r)}(M)$ is the submodule of $F(M)$
which consists of the summands $\tau(M)\subset F(M)$
associated to a tree with $r$ vertices $\tau\in T_{(r)}(n)$.
(Hence, the module $F_{(r)}(M)$ consists of tensors of weight $r$ in $F(M)$.)

We have only one tree with no vertex $\tau = \downarrow$
and this tree has one entry.
Consequently,
we obtain $F_{(0)}(M)(n) = \K$, for $n = 1$, and $F_{(0)}(M)(n) = 0$, for $n\not=1$.
In addition, the module $F_{(0)}(M)(1)$ has a canonical generator,
which we denote by $1\in F_{(0)}(M)(1)$,
since it represents the unit of the free operad.

Similarly,
for a given set with $n$ elements $I = \{i_1,\ldots,i_n\}$,
there is exactly one $I$-tree $\sigma_I\in T(I)$
with one vertex $v\in V(\sigma_I)$.
This tree is associated to the canonical partition $V(\sigma_I)\amalg I = I_{0}\amalg I_{v}$,
where $I_{0} = \{v\}$ and $I_{v} = I$.
We have by definition $\sigma_I(M)(I) = M(I)$.
Consequently,
we have an isomorphism $F_{(1)}(M)(I)\simeq M(I)$,
which is functorial in regard to bijections $u: I\,\rightarrow\,I'$.
We conclude that the $\Sigma_*$-module $F_{(1)}(M)$ is isomorphic to $M$.

\begin{figure}
\[\xymatrix@M=0pt@R=10mm@C=10mm{ *+<3mm>{i_1}\ar[drr] &
*+<3mm>{i_2}\ar[dr] & \cdots & & *+<3mm>{i_n}\ar[dll] \\
& & *+<3mm>{x}\ar[d] & & \\
& & *+<3mm>{0} & & \\ }\]
\caption{}\label{OneVertexTree}\end{figure}
Graphically,
we identify an element $x\in M(I)$ with the labeled tree of figure~\ref{OneVertexTree}.

\subsection{The treewise composition process in the free operad}
The $\Sigma_*$-module $F(M)$ is a realization
of the free operad generated by $M$.
The universal morphism $M\,\rightarrow\,F(M)$
is given by the isomorphism $M = F_{(1)}(M)$
as mentioned in paragraphs~\ref{FreeOperad} and~\ref{dgFreeOperad}.
The operad product is deduced
from the operad structure of the category of trees
and preserves the weight grading of the free operad.

We define the partial composition product
\[\circ_i: F_{(r)}(M)(m)\otimes F_{(s)}(M)(n)\,\rightarrow\,F_{(r+s)}(M)(m+n-1).\]
In fact,
we have a canonical isomorphism
\[\sigma\circ_i\tau(M) = \sigma(M)\otimes\tau(M).\]
The partial composition product maps simply
an element of $\sigma(M)\otimes\tau(M)\subset F_{(r)}(M)(m)\otimes F_{(s)}(M)(n)$
to the associated element
in $\sigma\circ_i\tau(M)\subset F_{(r+s)}(M)(m+n-1)$.

As mentioned in paragraph~\ref{WeightTreewiseTensors},
the unit of the free operad is associated to the unique tree with no vertex $\tau = \downarrow$
and generates the weight component $F_{(0)}(M)\subset F(M)$.

\subsection{Derivations}\label{TreeDerivation}
We make explicit the derivation of the free operad $d_\theta: F(M)\,\rightarrow\,F(M)$
associated to a homogeneous morphism $\theta: M\,\rightarrow\,F(M)$.

The morphism $\theta: M\,\rightarrow\,F(M)$
is determined by maps
$\theta: M(I)\,\rightarrow\,\sigma'(M)$,
where $\sigma'\in T(I)$.
We define a map $d_\theta^{\tau,\tau'}: \tau(M)(n)\,\rightarrow\,\tau'(M)(n)$
for each $n$-tree $\tau'\in T(n)$
obtained by blowing up a vertex $v\in V(\tau)$ in the tree $\tau\in T(n)$.
The derivation $d_\theta: F(M)(n)\,\rightarrow\,F(M)(n)$
is the sum of these maps
$d_\theta^{\tau,\tau'}: \tau(M)(n)\,\rightarrow\,\tau'(M)(n)$.

We have by definition $\tau = \tau'/\sigma'$,
where $\sigma'\subset\tau'$ is a tree whose entries are indexed by $I_v$.
We set
\[d_\theta^{\tau,\tau'}\Bigl(\bigotimes_{v'} x_{v'}\Bigr)
= \theta(x_v)\otimes\Bigl\{\bigotimes_{v'\not=v} x_{v'}\Bigr\}.\]
More precisely, in this formula, we consider the element $\theta(x_v)\in\sigma'(M)$
associated to $x_v\in M(I_v)$
by the map $\theta: M(I_v)\,\rightarrow\,\sigma'(M)$.
Since $V(\tau') = V(\sigma')\amalg V(\tau)\setminus\{v\}$,
we have
\[\tau'(M) = \sigma'(M)\otimes\bigl\{\bigotimes_{v'\in V(\tau)\setminus\{v\}} M(I_{v'})\bigr\}.\]
Therefore,
the formula above returns an element of $\tau'(M)$.

\subsection{Coderivations}\label{TreeCoderivation}
If $M$ is connected,
then the treewise $\Sigma_*$-module $F(M)$
is also a realization of the cofree cooperad cogenerated by $M$,
denoted by $F^c(M)$
(see paragraph~\ref{CofreeCooperad}).
We make explicit the process
which gives the coderivation
$d_\theta: F^c(M)\,\rightarrow\,F^c(M)$
associated to a given homogeneous morphism
$\theta: F^c(M)\,\rightarrow\,M$.

A morphism $\theta: F^c(M)\,\rightarrow\,M$ is determined by maps
$\theta: \sigma(M)\,\rightarrow\,M(I)$, where $\sigma\in T(I)$.
We define a map $d_\theta^{\tau,\tau'}: \tau(M)(n)\,\rightarrow\,\tau'(M)(n)$
for each $n$-tree $\tau'\in T(n)$
obtained by collapsing a subtree $\sigma\subset\tau$.
The coderivation $d_\theta: F^c(M)(n)\,\rightarrow\,F^c(M)(n)$
is the sum of these maps
$d_\theta^{\tau,\tau'}: \tau(M)(n)\,\rightarrow\,\tau'(M)(n)$.
We set precisely
\[d_\theta^{\tau,\tau'}\Bigl(\bigotimes_{v} x_{v}\Bigr)
= \theta\Bigl(\bigotimes_{v\in V(\sigma)} x_{v}\Bigr)
\otimes\Bigl\{\bigotimes_{v\not\in V(\sigma)} x_{v}\Bigr\}.\]
To be more explicit,
we consider the element of $M(I_\sigma)$
associated to the tensor
\[\bigotimes_{v\in V(\sigma)} x_{v}\in\bigotimes_{v\in V(\sigma)} M(I_v) = \sigma(M)\]
by the map $\theta: \sigma(M)\,\rightarrow\,M(I_\sigma)$.
This element gives the label of the vertex $\sigma\in V(\tau/\sigma)$
which results from the collapsing process.
The labels of vertices $v\not\in V(\sigma)$ are unchanged in the quotient tree.

\section{Trees and reduced bar constructions}\label{ConstructionBarDifferential}

In this section,
we make explicit the differential $\beta: \bar{B}(P)\,\rightarrow\,\bar{B}(P)$
of the reduced bar construction of a connected dg-operad $P$.
We recall that the reduced bar construction is a quasi-cofree cooperad
such that $\bar{B}(P) = F^c(\Sigma\tilde{P})$
(see paragraph~\ref{ReducedBarPpties}).
To begin with,
we recall the definition of the suspension of dg-modules.

\subsection{Suspension of dg-modules}\label{Suspension}
The \emph{suspension} of a dg-mo\-du\-le is the dg-module $\Sigma V$
defined by the tensor product
$\Sigma V = \K e\otimes V$, where $\deg(e) = 1$.
Thus,
we have a canonical isomorphism $(\Sigma V)_d = V_{d-1}$.
The element associated to $v\in V_{d-1}$
is represented by the tensor product $e\otimes v\in(\F e)_1\otimes V_{d-1}$
and is also denoted by $\Sigma v\in(\Sigma V)_d$.
According to the definition of paragraph~\ref{DGTensor},
the differential of $\Sigma V$ is given by the formula $\delta(\Sigma v) = - \Sigma(\delta v)$, for all $v\in V$.
Finally,
our definition makes the suspension symbol $\Sigma$ equivalent to a homogeneous element of degree $1$.

\subsection{Suspensions and tensor products}\label{SuspensionTensorProduct}
We have a canonical isomorphism
$\sigma: \Sigma U\otimes\Sigma V\,\simeq\,\Sigma^2(U\otimes V)$
defined by a permutation of tensors.
We have explicitly
\[\sigma(\Sigma u\otimes\Sigma v) = \pm\Sigma^2(u\otimes v).\]
In this formula,
the sign $\pm$ is yielded by the permutation of the element $u\in U$ with the permutation symbol $\Sigma$
as in the definition of paragraph~\ref{DGTensor}.
On the other hand,
we have isomorphisms
\begin{align*} & c(\Sigma U,\Sigma V): \Sigma U\otimes\Sigma V\,\xrightarrow{}\,\Sigma V\otimes\Sigma U \\
\text{and}\quad & \Sigma^2 c(U,V): \Sigma^2(U\otimes V)\,\xrightarrow{}\,\Sigma^2(V\otimes U) \end{align*}
induced by the symmetry isomorphism of the tensor product.
Hence,
we obtain a diagram
\[\xymatrix{ \Sigma U\otimes\Sigma V\ar[r]^\sigma\ar[d]_{c(\Sigma U,\Sigma V)} &
\Sigma^2(U\otimes V)\ar[d]^{\Sigma^2 c(U,V)} \\
\Sigma V\otimes\Sigma U\ar[r]^\sigma & \Sigma^2(V\otimes U) }\]
Both composites $c(U,V)\cdot\sigma$ and $\sigma\cdot c(\Sigma U,\Sigma V)$
permute factors of the tensor product $\Sigma U\otimes\Sigma V$.
We observe that the morphism $\sigma\cdot c(\Sigma U,\Sigma V)$ permutes suspension symbols,
while the morphism $\Sigma^2 c(U,V)\cdot\sigma$ does not.
Therefore:

\begin{fact}
The composite isomorphisms above verify the relation
$\Sigma^2 c(U,V)\cdot\sigma = - \sigma\cdot c(\Sigma U,\Sigma V)$.
The difference of signs corresponds to a permutation of suspensions.
\end{fact}

\subsection{Partial composition products}
We recall that an operad structure is equivalent to partial composition products
\[\circ_i: \tilde{P}(I)\otimes\tilde{P}(J)\,\rightarrow\,\tilde{P}(I\setminus\{i\}\amalg J)\]
that satisfy the following identities.
We assume $p\in P(I)$, $i\in I$, $q\in P(J)$, $j\in J$ and $r\in P(K)$.
We have then
\begin{equation}\label{AssPartialComp}
(p\circ_i q)\circ_j r = p\circ_i(q\circ_j r).
\end{equation}
We assume $p\in P(I)$, $i_1,i_2\in I$, $q_1\in P(J_1)$ and $q_2\in P(J_2)$.
We have then
\begin{equation}\label{ComPartialComp}
(p\circ_{i_1} q_1)\circ_{i_2} q_2 = (p\circ_{i_2} q_2)\circ_{i_1} q_1.
\end{equation}
We refer to paragraphs~\ref{PartialComposites} and~\ref{OperadIndexing}.
In the dg-context,
the permutation of the elements $q_1\in\tilde{P}(J_1)$ and $q_2\in\tilde{P}(J_2)$
yields an additional sign $\pm$.

\begin{lem}\label{BarCochain}
A collection of partial composition products
\[\circ_i: \tilde{P}(I)\otimes\tilde{P}(J)\,\rightarrow\,\tilde{P}(I\setminus\{i\}\amalg J),\]
where $I$ and $J$ range over the category of finite sets,
is equivalent to a homogeneous morphism of degree $-1$
\[\theta: F^c_{(2)}(\Sigma\tilde{P})\,\rightarrow\,\Sigma\tilde{P}.\]

The associated coderivation $d_\theta: F^c(\Sigma\tilde{P})\,\rightarrow\,F^c(\Sigma\tilde{P})$
verifies the equation $\delta d_\theta + d_\theta \delta = 0$
if and only if the partial composition products are morphism of dg-modules
and the equation $d_\theta{}^2 = 0$
if and only the relations~\eqref{AssPartialComp} and~\eqref{ComPartialComp}
are satisfied.
\end{lem}

In the next paragraph,
we analyze the structure of the $\Sigma_*$-module $F^c_{(2)}(\Sigma\tilde{P})$
and
we make explicit the morphism $\theta: F^c_{(2)}(\Sigma\tilde{P})\,\rightarrow\,\Sigma\tilde{P}$
equivalent to partial composition products.
In paragraph~\ref{ReducedBarDifferential},
we make explicit the coderivation $d_\theta: F^c(\Sigma\tilde{P})\,\rightarrow\,F^c(\Sigma\tilde{P})$
associated to $\theta: F^c_{(2)}(\Sigma\tilde{P})\,\rightarrow\,\Sigma\tilde{P}$
according to the process of paragraph~\ref{TreeCoderivation}.
The relations $\delta d_\theta + d_\theta \delta = 0$ and $d_\theta{}^2 = 0$
follow from straighforward verifications
which are omitted.

\subsection{On trees with $2$ vertices}
We observe that the structure of an $n$-tree with $2$ vertices is equivalent
to a partition $\{1,\ldots,n\} = \{i_1,\ldots,i_{k-1}\}\amalg\{j_1,\ldots,j_l\}$.
To be explicit,
we let $V(\tau) = \{u,v\}$.
We assume that $u$ is the source of the root of $\tau$.
We set then $I_{u} = \{i_1,\ldots,i_{k-1}\}\amalg\{v\}$ and $I_{v} = \{j_1,\ldots,j_l\}$.
We obtain the tree represented in figure~\ref{TwoVertexTree}.
\begin{figure}
\[\xymatrix@M=0pt@R=6mm@C=6mm{ & & & *+<3mm>{j_{1}}\ar[drr] &
*+<3mm>{j_{2}}\ar[dr] & \cdots &
*+<3mm>{j_{l-1}}\ar[dl] &
*+<3mm>{j_{l}}\ar[dll] \\
*+<3mm>{i_{1}}\ar[drrr] & &
*+<3mm>{i_{2}}\ar[dr] &
\cdots &
*+<3mm>{i_{k-1}}\ar[dl] &
*+<3mm>[o][F]{v}\ar[dll] & & \\
& & & *+<3mm>[o][F]{u}\ar[d] & & & & \\
& & & *+<3mm>{0} & & & & \\ }\]
\caption{}\label{TwoVertexTree}\end{figure}

We deduce from this observation that, for a given $\Sigma_*$-module $M$,
we have
\[F^c_{(2)}(M)(n) = \bigoplus_{\{i_*\},\{j_*\}}
M(\{i_1,\ldots,i_{k-1}\}\amalg \{v\})\otimes M(\{j_1,\ldots,j_l\}.\]
The sum ranges over partitions
$\{1,\ldots,n\} = \{i_1,\ldots,i_{k-1}\}\amalg\{j_1,\ldots,j_l\}$.
Therefore,
in the case $M = \Sigma\tilde{P}$,
we conclude that a homogeneous morphism of degree $-1$
\[\theta: F^c_{(2)}(\Sigma\tilde{P})(n)\,\rightarrow\,\Sigma\tilde{P}(n)\]
is equivalent to a collection of partial composition products homogeneous of degree $0$
\begin{multline*}\circ_{v}: \tilde{P}(\{i_1,\ldots,i_{k-1}\}\amalg\{v\})\otimes\tilde{P}(\{j_1,\ldots,j_l\}) \\
\,\rightarrow\,P(\{i_1,\ldots,i_{k-1}\}\amalg\{j_1,\ldots,j_l\}).\end{multline*}

\begin{figure}
\begin{equation*}\begin{aligned}
\xymatrix@M=0pt@R=4mm@C=4mm{ & & *+<2mm>{j_1}\ar[dr] & \cdots & *+<2mm>{j_l}\ar[dl] \\
*+<2mm>{i_1}\ar[drr] & \cdots & *+<2mm>{i_{k-1}}\ar[d] & *+<3mm>{\Sigma q}\ar[dl] & \\
& & *+<3mm>{\Sigma p}\ar[d] & & \\
& & *+<2mm>{0} & & \\ }
\end{aligned}
\mapsto
\begin{aligned}
\xymatrix@M=0pt@R=4mm@C=4mm{ *+<2mm>{i_1}\ar[drr] & \cdots & *+<2mm>{i_{k-1}}\ar[d] &
*+<2mm>{j_1}\ar[dl] & \cdots & *+<2mm>{j_l}\ar[dlll] \\
& & *+<3mm>{\Sigma(p\circ_i q)}\ar[d] & & & \\
& & *+<2mm>{0} & & & \\ }
\end{aligned}\end{equation*}
\caption{}\label{DisplayEdgeContraction}\end{figure}
Graphically,
the morphism $\theta: F^c_{(2)}(\Sigma\tilde{P})\,\rightarrow\,\Sigma\tilde{P}$
is represented by the picture of figure~\ref{DisplayEdgeContraction}.

\subsection{The differential of the bar construction}\label{ReducedBarDifferential}
The differential of the bar construction $\beta: \bar{B}(P)\,\rightarrow\,\bar{B}(P)$
is the coderivation $d_\theta: F^c(\Sigma\tilde{P})\,\rightarrow\,F^c(\Sigma\tilde{P})$
supplied by lemma~\ref{BarCochain}.
We make explicit the components
$\beta^{\tau,\tau'}: \tau(\Sigma\tilde{P})\,\rightarrow\,\tau'(\Sigma\tilde{P})$
of this coderivation.
We just follow the process of paragraph~\ref{TreeCoderivation}.

We observe that a subtree $\sigma\subset\tau$ with $2$ vertices $V(\sigma) = \{u,v\}$
is determined by an edge $e\in E(\tau)$.
We have explicitly $u = t(e)$ and $v = s(e)$.
We denote this subtree, associated to $e\in E(\tau)$, by the notation $\sigma = \sigma_e$.
The quotient tree $\tau/\sigma_e$ is obtained by contracting the edge
$\xymatrix@M=0pt@R=10mm@C=10mm{*+<3mm>[o][F]{v}\ar[r]|*+<1mm>{e} & *+<3mm>[o][F]{u} \\ }$
to a vertex.
Therefore, the process $\tau\mapsto\tau/\sigma_e$ is also known as an \emph{edge contraction}.

The bar differential has a component
$\beta^{\tau,\tau'}: \tau(\Sigma\tilde{P})\,\rightarrow\,\tau'(\Sigma\tilde{P})$
for each edge contraction $\tau' = \tau/\sigma_e$.
The vertex $\sigma_e\in V(\tau')$ (which results from the contraction process)
is labeled by the partial composition product
$- \Sigma(p_{u}\circ_{v} p_{v})\in\Sigma\tilde{P}(I_{u}\setminus\{v\}\amalg I_{v})$,
where $\Sigma(p_{u})\in\Sigma\tilde{P}(I_{u})$
(respectively, $\Sigma(p_{v})\in\Sigma\tilde{P}(I_{v})$)
labels the vertex $u = t(e)\in V(\tau)$
(respectively, $v = s(e)\in V(\tau)$).
The elements $\Sigma p_{x}\in\Sigma\tilde{P}(I_{x})$
associated to the other vertices of the quotient tree
are unchanged.
The motivation of the additional sign comes from a chain morphism introduced in section~\ref{BarCoefficients}
(the levelization morphism).
To be precise,
the sign makes the bar differential agree with faces of a simplicial bar construction
(see lemma~\ref{LevelizationDifferential}).

As an example,
the bar differential of a tensor
$x_1\otimes x_2\otimes x_3\otimes x_4 = \Sigma p_1\otimes\Sigma p_2\otimes\Sigma p_3\otimes\Sigma p_4
\in\tau(\Sigma\tilde{P})$,
where $\tau$ is the tree of figure~\ref{TreeExample},
consists of the terms displayed in figure~\ref{BarDifferentialExample}.
\begin{figure}
\[\xymatrix@M=0pt@R=6mm@C=6mm{ & & &
*+<2mm>{1}\ar[dr] & &
*+<2mm>{5}\ar[dl] &
*+<2mm>{4}\ar[dr] & &
*+<2mm>{6}\ar[dl] \\
*+<2mm>{2}\ar[drrrr] & &
*+<2mm>{3}\ar[drr] & &
*+<2mm>{\Sigma p_3}\ar[d] & & &
*+<2mm>{\Sigma p_4}\ar[dlll] & \\
& & & & *+<2mm>{\Sigma(p_1\circ_{v_2} p_2)}\ar[d] & & & & \\
& & & & *+<2mm>{0} & & & & \\
& & & & & & & & \\
& & & *+<2mm>{3}\ar[dr] & *+<2mm>{1}\ar[d] & *+<2mm>{5}\ar[dl] &
*+<2mm>{4}\ar[dr] & & *+<2mm>{6}\ar[dl] \\
*+<2mm>{2}\ar[drrrr] & & & & *+<2mm>{\Sigma(p_2\circ_{v_3} p_3)}\ar[d] & & &
*+<2mm>{\Sigma p_4}\ar[dlll] & \\
& & & & *+<2mm>{\Sigma p_1}\ar[d] & & & & \\
& & & & *+<2mm>{0} & & & & \\
& & & & & & & & \\
& & & *+<2mm>{1}\ar[dr] & & *+<2mm>{5}\ar[dl] & & & \\
& & *+<2mm>{3}\ar[drr] & & *+<2mm>{\Sigma p_3}\ar[d] & & \\
*+<2mm>{2}\ar[drrrr] & & & &
*+<2mm>{\Sigma p_2}\ar[d] & & *+<2mm>{4}\ar[dll] & &
*+<2mm>{6}\ar[dllll] \\
& & & & *+<2mm>{\Sigma(p_1\circ_{v_4} p_4)}\ar[d] & & & & \\
& & & & *+<3mm>{0} & & & & \\ }\]
\caption{}\label{BarDifferentialExample}\end{figure}

\section[Comparison of quasi-free operads]{Comparison of quasi-free operads: proof of theorem~\ref{ComparisonQuasiFreeOperads}}
\label{ProofsComparisonQuasiFreeOperads}

In this section,
we aim to prove theorem~\ref{ComparisonQuasiFreeOperads}.
Therefore,
we are concerned with quasi-free operads $F = F(M)$
such that $M(0) = M(1) = 0$.
To begin with,
we analyze the expansion of the free operad $F(M)$
in this situation.

\subsection{Reduced trees and free operads}\label{ReducedFreeOperad}
An $n$-tree $\tau\in T(n)$ is reduced if all vertices $v\in V(\tau')$ have at least $2$ entries.
We observe that the categories of reduced $n$-trees, which are denoted by $\tilde{T}(n)$,
define a suboperad of the operad of trees.
In the context of reduced $n$-trees,
we have $\tilde{T}_{(r)}(n) = \emptyset$, for $r\geq n$.

Suppose given a $\Sigma_*$-module $M$ such that $M(0) = M(1) = 0$.
In that case,
the $\K$-module $\tau(M)$ associated to a non-reduced tree $\tau\in T(n)$ vanishes.
Consequently,
the free operad generated by $M$ has an expansion
\[F(M)(n) = \colim_{\tau\in\tilde{T}(n)}\tau(M),\]
such that the colimit ranges over the category of reduced trees $\tilde{T}(n)$.
Since
$\tilde{T}_{(r)}(n) = \emptyset$, for $r\geq n$,
we obtain in this case
$F_{(r)}(M)(n) = 0$, for $r\geq n$.

In addition,
we observe that a given reduced $n$-tree $\tau\in\tilde{T}(n)$ has no non-trivial automorphism.
Consequently,
the associated module $\tau(M)$ is a direct summand of $F(M)(n)$.
Furthermore,
if we fix a representative $\tau_i$ of each isomorphism class,
then we obtain a quotient-free expansion of the free operad
$F(M)(n) = \bigoplus_{\tau_i} \tau(M)$.

This last assertion holds in the more general case of a connected $\Sigma_*$-modules $M$
(and implies the result of lemma~\ref{FreeOperadHomology}).
We consider dg-$\Sigma_*$-modules $M$ which have $M(0) = M(1) = 0$
in order to guarantee the convergence of the spectral sequence
of lemma~\ref{QuasiFreeOperadSpectralSequence}.

\begin{lem}\label{QuasiFreeOperadSpectralSequence}
We consider a quasi-free operad $P = F(M)$, where $M$ is a connected dg-$\Sigma_*$-module.
We assume that the differential of $P$
is determined by a derivation $d_\theta: F(M)\,\rightarrow\,F(M)$
such that $d_\theta(F(M))\subset \bigoplus_{r\geq 2} F_{(r)}(M)$
(so that the dg-$\Sigma_*$-module $M$ can be identified with the indecomposable quotient of $P$).

We have a second quadrant spectral sequence of operads $E^r(P)\,\Rightarrow\,H_*(P)$
such that $E^0_{s t}(P) = F_{(-s)}(M)_{-s+t}$,
and where the differential $d^0: E^0_{s t}(P)\,\rightarrow\,E^0_{s t-1}(P)$
is induced by the differential of $M$.
In particular, if the ground ring $\K$ is a field,
then we obtain $E^1_{s *} = F_{(-s)}(H_*(M))$.

The spectral sequence converges strongly if we have $M(0) = M(1) = 0$.
\end{lem}

\begin{proof}
We consider the filtration $0\subset\cdots\subset F_s P\subset\cdots\subset F_1 P\subset F_0 P = P$
such that $F_s P = \bigoplus_{r\geq s} F_{(r)} M$,
and
so that
\[E^0_{s t}(P) = F_{(-s)}(M)_{-s+t}.\]
For fixed $n\in\N$,
we have $E^0_{s t}(P(n)) = F_{(-s)}(M)(n)_{-s+t} = 0$, for $s<-n$.
Hence,
in the case $M(0) = M(1) = 0$,
the spectral sequence is bounded in horizontal degree.
The convergence assertion follows from this fact.

We prove that the differential of the quasi-free operad
$\delta_\theta: F(M)\,\rightarrow\,F(M)$
preserves the filtration above.
On one hand,
the canonical differential of the free operad
$\delta: F(M)\,\rightarrow\,F(M)$
verifies
$\delta(F_{(r)} M)\subset F_{(r)} M$.
On the other hand,
by lemma~\ref{FreeOperadDerivation},
the derivation $d_\theta: F(M)\,\rightarrow\,F(M)$
induces a morphism such that
$d_\theta(F_s P)\subset F_{s-1} P$,
because we have by assumption
$\theta(M)\subset\bigoplus_{r\geq 2} F_{(r)}(M)$.

The result follows immediately from these observations.
\end{proof}

\begin{lem}
Suppose given a morphism of quasi-free operads $\phi: P\,\rightarrow\,P'$
where $P = F(M)$ and $P' = F(M')$ are as in lemma~\ref{QuasiFreeOperadSpectralSequence} above.
We assume that the morphism $\phi: P\,\rightarrow\,P'$ preserves augmentation ideals of free operads
(explicitly, $\phi(\tilde{F}(M))\subset\tilde{F}(M')$).
In this situation,
we have an induced morphism of spectral sequences
$E^r(\phi): E^r(P)\,\rightarrow\,E^r(P')$.

Since $M$ (respectively, $M'$) is identified with the indecomposable quotient of $P = F(M)$
(respectively, $P' = F(M')$),
we have a morphism of dg-$\Sigma_*$-modules $\bar{\phi}: M\,\rightarrow\,M'$
induced by $\phi: F(M)\,\rightarrow\,F(M')$.
We observe in addition that the morphism $E^0(\phi): E^0(P)\,\rightarrow\,E^0(P')$
is identified with the morphism of free operads $F(\bar{\phi}): F(M)\,\rightarrow\,F(M')$
associated to $\bar{\phi}: M\,\rightarrow\,M'$.
\end{lem}

\begin{proof}
We prove that the morphism $\phi: F(M)\,\rightarrow\,F(M')$
preserves the filtrations introduced in the proof of lemma~\ref{QuasiFreeOperadSpectralSequence}.
We mention in the definition of a free operad
that the weight component $F_{(r)}(M)$ is generated by composites
$\bigl(\cdots\bigl(\bigl(x_1\circ_{i_2} x_2\bigr)\circ_{i_3}\cdots\bigr)\circ_{i_r} x_r$
where $x_1\in M(n_1),\ldots,x_r\in M(n_r)$ (see paragraphs~\ref{FreeOperad} and~\ref{TreewiseTensors}).
An operad morphism $\phi: F(M)\,\rightarrow\,F(M')$
maps such a composite to the element
$\bigl(\cdots\bigl(\bigl(\phi(x_1)\circ_{i_2}\phi(x_2)\bigr)\circ_{i_3}\cdots\bigr)\circ_{i_r}\phi(x_r)\in F(M')$.
By assumption,
we have $\phi(x_i)\in\bigoplus_{s=1}^{\infty} F_{(s)}(M')$.
Consequently,
the composite element above verifies
\[\bigl(\cdots\bigl(\bigl(\phi(x_1)\circ_{i_2}\phi(x_2)\bigr)\circ_{i_3}\cdots\bigr)\circ_{i_r}\phi(x_r)
\in\bigoplus_{s=r}^{\infty} F_{(s)}(M').\]
Therefore, we obtain $\phi(F_r P)\subset F_r P'$.

The morphism $\bar{\phi}: M\,\rightarrow\,M'$ is the composite of $\phi: F(M)\,\rightarrow\,F(M')$
with the canonical injection $M\,\rightarrow\,F(M)$ and with the canonical projection $F(M')\,\rightarrow\,M'$.
Equivalently, we have
\[\phi(x_i)\equiv\bar{\phi}(x_i)\quad\bigl[\mathop{\rm mod}\,\bigoplus_{s=2}^{\infty} F_{(s)}(M')\bigr],\]
for all $x_i\in M(n_i)$.
As a consequence,
we obtain
\begin{multline*} \bigl(\cdots\bigl(\bigl(\phi(x_1)\circ_{i_2}\phi(x_2)\bigr)\circ_{i_3}\cdots
\bigr)\circ_{i_r}\phi(x_r) \\
\equiv\bigl(\cdots\bigl(\bigl(\bar{\phi}(x_1)\circ_{i_2}\bar{\phi}(x_2)\bigr)\circ_{i_3}\cdots
\bigr)\circ_{i_r}\bar{\phi}(x_r)
\quad\bigl[\mathop{\mathrm{mod}}\,\bigoplus_{s=r+1}^{\infty} F_{(s)}(M')\bigr]. \end{multline*}
Therefore,
we conclude that the morphism $E^0(\phi): E^0(P)\,\rightarrow\,E^0(P')$
is identified with the morphism of free operads $F(\bar{\phi}): F(M)\,\rightarrow\,F(M')$
associated to $\bar{\phi}: M\,\rightarrow\,M'$.
\end{proof}

\begin{proof}[Proof of theorem~\ref{ComparisonQuasiFreeOperads}]
We deduce theorem~\ref{ComparisonQuasiFreeOperads} from the result above.
If the morphism $\bar{\phi}: M\,\rightarrow\,M'$ is a quasi-isomorphism,
then the morphism $\phi: P\,\rightarrow\,P'$
induces an isomorphism $E^1(\phi): E^1(P)\,\rightarrow\,E^1(P')$
at the $E^1$-stage of the spectral sequence.
Therefore,
we conclude that $\phi: P\,\rightarrow\,P'$ is a quasi-isomorphism.

Conversely,
suppose given a quasi-isomorphism $\phi: P\,\rightarrow\,P'$.
We prove by induction that the morphism of dg-modules
$\bar{\phi}: M(n)\,\rightarrow\,M'(n)$
is a quasi-isomorphism
for all $n\in\N$.
Precisely,
we assume that $\bar{\phi}: M(r)\,\rightarrow\,M'(r)$ is a quasi-isomorphism for all $r<n$.

We observe that, for $s\geq 2$,
the component
\[E^0_{s t}(P(n)) = F_{(-s)}(M)(n)_{-s+t}\]
involves only modules $M(I)$,
such that $I$ has less than $n$ elements.
We have the same observation for $E^0_{s t}(P'(n))$.
Consequently,
we obtain an isomorphism $E^1(\phi): E^1_{s t}(P(n))\,\rightarrow\,E^1_{s t}(P'(n))$
in all bidegrees $(s,t)$ such that $s\geq 2$.
We have $E^0_{0 *}(P(n)) = E^0_{0 *}(P'(n)) = 0$ as long as $n\geq 2$.
Therefore,
the previous assertion implies that $E^1(\phi): E^1_{s t}(P(n))\,\rightarrow\,E^1_{s t}(P'(n))$
is also in isomorphism in degree $s=1$.
We conclude that $\bar{\phi}: M(n)\,\rightarrow\,M'(n)$ is a quasi-isomorphism,
since we have $E^0_{1 *}(P(n)) = F_{(1)}(M)(n) = M(n)$
and $E^0_{1 *}(P'(n)) = F_{(1)}(M')(n) = M'(n)$.
\end{proof}

%\input{CoefficientBar}
% 30/1/2003

\chapter{Bar constructions with coefficients}\label{BarCoefficients}

\section{Summary}\label{BarCoefficientsPpties}

We generalize the definition of the reduced bar construction $\bar{B}(P)$
in order to have a complex of dg-$\Sigma_*$-modules $B(L,P,R)$,
with coefficients in a right $P$-module $L$ and a left $P$-module $R$.
This complex is equipped with a natural augentation morphism
$\epsilon(L,P,R): B(L,P,R)\,\rightarrow\,L\circ_P R$
and satisfies nice homological properties
(see statements~\ref{QuasiFreeBarComplex},~\ref{BarModule} and~\ref{BarResolution} below).
As mentioned in chapter~\ref{ReducedBar},
the reduced bar construction $\bar{B}(P)$ is the bar construction $B(L,P,R)$
with trivial coefficients $L = R = I$.

The next facts~\ref{QuasiFreeBarComplex} and~\ref{BarModule} are immediate consequences of the definition,
given in section~\ref{DifferentialGradedBarConstruction}.

\begin{fact}\label{QuasiFreeBarComplex}
We have $B(L,P,R) = L\circ\bar{B}(P)\circ R$.
Consequently,
the chain complex $B(L,P,P)$ (respectively, $B(P,P,R)$) is a quasi-free right $P$-module
(respectively, a quasi-free left $P$-module).
\end{fact}

\begin{fact}\label{BarModule}
The dg-$\Sigma_*$-module $B(P,P,P)$ is a complex of $P$-bi\-mo\-dules.
Furthermore, we have $B(L,P,R) = L\circ_P B(P,P,P)\circ_P R$
and the augmentation morphism is identified with the relative composition product
$\epsilon(L,P,R) = L\circ_P\epsilon(P,P,P)\circ_P R$.
\end{fact}

The next assertion is proved in section~\ref{ProofsBarCoefficientsPpties}.

\begin{lem}\label{BarResolution}
We are given a connected dg-operad $P$ and a right $P$-module $L$
(respectively, a connected left $P$-module $R$).
If the ground ring $\K$ is not a field,
then we assume in addition that $P$ (respectively, $L$, $R$)
is projective as a $\K$-module.

The augmentation morphism $\epsilon: B(L,P,P)\,\rightarrow\,L$
(respectively, $\epsilon: B(P,P,R)\,\rightarrow\,R$)
is a quasi-isomorphism of right $P$-modules
(respectively, left $P$-modules).
In particular, the dg-module $B(I,P,P)$ (respectively, $B(P,P,I)$) is acyclic.
\end{lem}

We mention the next result,
which, according to proposition~\ref{DerivedCompositionProduct},
occurs as a corollary of the statements above.

\begin{prp}\label{BarHomologyInterpretation}
The chain complex $B(L,P,R)$ is a representative
of the derived composition product $L\circ^{\mathbb{L}}_P R$
provided the assumptions of lemma~\ref{BarResolution} are satisfied.
\end{prp}

The bar constructions $B(I,P,P)$ and $B(P,P,I)$
are defined by E. Getzler and J. Jones in \cite{GetzlerJones}.
These authors prove also that $B(I,P,P)$ and $B(P,P,I)$ form acyclic chain complexes.
We give different arguments for these results.
More precisely,
we deduce a chain homotopy for $B(I,P,P)$ from a simplicial construction
(introduced in the next paragraph).
In fact,
in this case,
we recover the chain homotopy of \cite{GetzlerJones}.
But,
in the case of the second complex $B(P,P,I)$,
our argument differs completely from \cite{GetzlerJones},
because we observe that the result about $B(I,P,P)$
and theorem~\ref{QuasiFreeCompositionProduct} (about derived composition products)
imply immediately that the complex $B(P,P,I)$ is acyclic
(\emph{cf}. section~\ref{ProofsBarCoefficientsPpties} for details).

\subsection{The simplicial bar construction}\label{SimplicialBarConstructionPpties}
We recall a classical construction
which supplies a simplicial dg-$\Sigma_*$-module $C(L,P,R)$
together with an augmentation $\epsilon(L,P,R): C(L,P,R)\,\rightarrow\,L\circ_P R$
like the differential graded bar construction $B(L,P,R)$.
We call this simplicial dg-$\Sigma_*$-module $C(L,P,R)$ the \emph{simplicial bar construction}.
We consider also the reduced simplicial bar construction $\bar{C}(P)$,
defined as the simplicial bar construction
with trivial coefficients $\bar{C}(P) = C(I,P,I)$.

We have by definition
\[C_d(L,P,R) = L\circ\underbrace{P\circ\cdots\circ P}_{d}\circ R.\]
The face
$d_i: C_d(L,P,R)\,\rightarrow\,C_{d-1}(L,P,R)$
is defined:
for $i = 0$,
by the right operad action $\lambda: L\circ P\,\rightarrow\,L$;
for $0<i<d$,
by the operad product $\mu: P\circ P\,\rightarrow\,P$ of the $i$th and $i+1$th copies of $P$;
for $i = d$, by the left operad action $\rho: P\circ R\,\rightarrow\,R$.
The degeneracy
$s_j: C_d(L,P,R)\,\rightarrow\,C_{d+1}(L,P,R)$
is given by the insertion of an operad unit $\eta: I\,\rightarrow\,P$
between the $j$th and the $j+1$th copies of $P$ in $C_{d}(L,P,R)$.
The augmentation
$\epsilon: C(L,P,R)\,\rightarrow\,L\circ_P R$
is given by the canonical projection
\[C_0(L,P,R) = L\circ R\,\rightarrow\,L\circ_P R.\]

A similar simplicial bar construction is defined in the context of monads by J. Beck
(\emph{cf}. J Beck \cite{Beck}, S. Mac Lane \cite{MacLaneCat})
and introduced by P. May in the study of iterated loop spaces
(\emph{cf}. \cite{MayLoop}).
In fact,
since the composition product is preserved by the correspondence $M\mapsto S(M)$,
the functor associated to the operadic simplicial bar construction
verifies
\[S(C_*(L,P,R)) = S(L)\circ S(P)\circ\cdots\circ S(P)\circ S(R)\]
and agrees with Beck's construction $C_*(S(L),S(P),S(R))$.
The simplicial bar construction $C_*(S(L),S(P),S(R))$ is also considered by P. May
(in the topological framework).
As mentioned previously,
in the differential graded framework,
the composition product of $\Sigma_*$-modules
behaves better in regard to homology
than the composition product of functors.
Therefore,
we are motivated to introduce simplicial constructions $C(L,P,R)$
at the level of $\Sigma_*$-modules.

We state the analogues of properties~\ref{QuasiFreeBarComplex},~\ref{BarModule} and~\ref{BarResolution} above
in the context of the simplicial bar construction.

\begin{fact}
We have $C_d(L,P,R) = L\circ\bar{C}_d(P)\circ R$.
Moreover,
the simplicial module $C(L,P,P)$ (respectively, $C(P,P,R)$)
is a quasi-free simplicial right $P$-module
(respectively, a quasi-free simplicial left $P$-module).
\end{fact}

We refer to Quillen's monograph \cite{QuillenBook}
for the definition of a quasi-free object in the simplicial context
(called \emph{almost free objects} in \emph{loc. cit.}).
We just make this definition explicit in the case of the simplicial right $P$-module $C(L,P,P)$.

We observe precisely that each module $C_d(L,P,P)$ is a free right $P$-module.
We have explicitly $C_d(L,P,P) = C_d(L,P,I)\circ P$
and we identify $C_d(L,P,I)$ with the module
$C_d(L,P,I)\circ I\subset C_d(L,P,I)\circ P = C_d(L,P,P)$.
The degeneracies
$s_j: C_d(L,P,P)\,\rightarrow\,C_{d+1}(L,P,P)$
verify
$s_j(C_d(L,P,I))\subset C_{d+1}(L,P,I)$.
We have also
$d_i(C_d(L,P,I))\subset C_{d-1}(L,P,I)$
for all faces $d_i: C_d(L,P,P)\,\rightarrow\,C_{d-1}(L,P,P)$
such that $i<d$.
We have symmetric observations in the case of the simplicial left $P$-module $C(P,P,R)$.

\begin{fact}
The augmentation $\epsilon: C(L,P,P)\,\rightarrow\,L$
(respectively, $\epsilon: C(P,P,R)\,\rightarrow\,R$)
is a quasi-isomorphism of simplicial right $P$-modules
(respectively, left $P$-modules).
In particular, the simplicial module $C(I,P,P)$ (respectively, $C(P,P,I)$) is acyclic.
\end{fact}

In fact,
the simplicial module $C(L,P,P)$
has a well known extra degeneracy
$s_{d+1}: C_{d}(L,P,P)\,\rightarrow\,C_{d+1}(L,P,P)$
which gives a contracting homotopy on the normalized chain complex of $C(L,P,P)$.
To be explicit,
this extra degeneracy $s_{d+1}: C_{d}(L,P,P)\,\rightarrow\,C_{d+1}(L,P,P)$
is given by the insertion of an operad unit $\eta: I\,\rightarrow\,P$
at the last position
of the composition product
$C_{d+1}(L,P,P) = L\circ P\circ\cdots\circ P\circ P$.
Similarly,
the simplicial module $C(P,P,R)$
has an extra degeneracy
$s_{-1}: C_{d}(P,P,R)\,\rightarrow\,C_{d+1}(P,P,R)$,
given by the insertion of an operad unit $\eta: I\,\rightarrow\,P$
at the first position
of the composition product
$C_{d+1}(P,P,R) = P\circ P\circ\cdots\circ P\circ R$.

\smallskip
We aim to compare the differential graded bar construction to the simplicial one.
We let $N(L,P,R) = N_*(C(L,P,R))$ denote the normalized chain complex of $C(L,P,R)$.
We obtain the following result:

\begin{thm}\label{Levelization}
We are given a dg-operad $P$, a right $P$-module $L$ and a left $P$-module $R$.
We have a natural morphism of dg-$\Sigma_*$-modules
\[\phi(L,P,R): B(L,P,R)\,\rightarrow\,N(L,P,R),\]
which we call the levelization morphism.

We assume that the dg-operad $P$ (respectively, the dg-$\Sigma_*$-module $L$, $R$) is connected.
If the ground ring $\K$ is not a field,
then we assume in addition that $P$ (respectively, $L$, $R$)
is projective as a $\K$-module.

In this context,
the levelization morphism $\phi(L,P,R): B(L,P,R)\,\rightarrow\,N(L,P,R)$
is injective and is a quasi-isomorphism.
In particular,
for trivial coefficients $L = R = I$,
the levelization morphism
\[\bar{\phi}(P): \bar{B}(P)\,\rightarrow\,\bar{N}(P)\]
provides a quasi-isomorphism from the differential graded reduced bar construction
to the simplicial one.
\end{thm}

We make the levelization morphism explicit in section~\ref{LevelizationProcess}.
We prove the theorem above in section~\ref{ProofsBarCoefficientsPpties}.
We would like to mention that S. Shnider and D. Van Osdol
prove in \cite{ShniderOsdol}
that the chain complexes $\bar{B}(P)$ and $\bar{N}(P)$
have the same homology
(see also M. Markl, S. Shnider and J. Stasheff \cite{MarklShniderStasheff}).
Let us mention that the proof of Shnider and Van Osdol depends on a K\"unneth formula
which does not hold if the ground ring $\K$ is not a field.

If we assume lemma~\ref{BarResolution},
then the result of S. Shnider and D. Van Osdol
occurs as a direct consequence
of the comparison theorems of chapter~\ref{ModuleComplexes}.
Namely,
the augmentation morphisms
$N(L,P,P)\,\rightarrow\,L$ and $B(P,P,R)\,\rightarrow\,R$
yield quasi-isomorphisms
\[L\circ_P B(P,P,R)\,\xleftarrow{\sim}\,N(L,P,P)\circ_P B(P,P,R)\,\xrightarrow{\sim}\,N(L,P,P)\circ_P R\]
because $N(L,P,P)$ and $B(P,P,R)$
are quasi-free modules
(see theorem \ref{ExactnessCompositionProduct}).
In \cite{Balavoine},
D. Balavoine uses close arguments for comparing operadic homology with cotriple homology
(but his comparison result holds only in characteristic $0$).

Finally,
the levelization process makes simply the relationship between $\bar{B}(P)$ and $\bar{N}(P)$
more precise.

\smallskip
In paragraph~\ref{AssociativeLevelization},
we give a few indications about the levelization morphism
in the example of the associative operad $P = \A$.
We recover a construction of A. Tonks
which relates Stasheff's \emph{associahedron} to Milgram's \emph{permutohedron}
(\emph{cf}. J.-L. Loday and M. Ronco \cite{LodayRonco},
S. Saneblidze and R. Umble \cite{SaneblidzeUmble},
A. Tonks \cite{Tonks}).
The case of the commutative operad $P = \C$,
which provides the main application of the article,
is handled in section~\ref{CommutativeDuality}.

To conclude this introduction,
we state properties of the normalized bar complex $N(L,P,R)$
that arise out of the structure of the simplicial bar construction $C(L,P,R)$.
In fact,
we deduce theorem~\ref{Levelization}
from the next assertion
and from the comparison theorems of chapter~\ref{ModuleComplexes}.

\begin{fact}\label{NormalizedRightModule}
The chain complex $N(L,P,P) = N_*(C(L,P,P))$
(respectively, the chain complex $N(P,P,R) = N_*(C(P,P,R))$)
is equipped with the structure of a right $P$-module
(respectively, left $P$-module).

We have in addition $N(L,P,R) = N(L,P,I)\circ R$.
Consequently,
the chain complex $N(L,P,P)$ forms a quasi-free right $P$-module
and we have $N(L,P,R) = N(L,P,P)\circ_P R$.
\end{fact}

\subsection{The normalized chain complex of a composite simplicial $\Sigma_*$-module}
In general,
for given simplicial $\Sigma_*$-modules $M$ and $N$,
we have a natural morphism
\[N_*(M)\circ N_*(N)\,\xrightarrow{}\,N_*(M\circ N),\]
which is induced by the classical Eilenberg-MacLane equivalence,
and that satisfies unit and associativity properties.
Consequently,
if $P$ is a simplicial operad,
then the associated complex $N_*(P)$
is equipped with the structure of a dg-operad.
Similarly, if $L$ (respectively, $R$)
is a simplicial right (respectively, left) module
over a simplicial operad $P$,
then the complex $N_*(L)$ (respectively, $N_*(R)$)
is equipped with the structure
of a differential graded right (respectively, left) module
over the dg-operad $N_*(P)$.

To be precise,
the Eilenberg-MacLane equivalence
induces a morphism of dg-modules
\begin{multline*}N_*(M(r))\otimes N_*(N)(J_1)\otimes\cdots\otimes N_*(N)(J_r) \\
\,\rightarrow\,N_*\bigl(M(r)\otimes N(J_1)\otimes\cdots\otimes N(J_r)\bigr),\end{multline*}
for all partitions $\{1,\ldots,n\} = J_1\amalg\cdots\amalg J_r$.
This morphism commutes with permutations by symmetry of the Eilenberg-MacLane equivalence.
Therefore, if we go back to the expansion of the composition product of $\Sigma_*$-modules,
then we obtain a natural morphism $N_*(M)\circ N_*(N)\,\xrightarrow{}\,N_*(M\circ N)$.
We would like to mention that,
in the case of connected $\Sigma_*$-modules,
we have a more precise result,
which we deduce from the expansion of the composition product given in lemma~\ref{ExpansionComposite}.
Namely:

\begin{lem}
We are given a simplicial $\Sigma_*$-module $M$
(respectively, a connected simplicial $\Sigma_*$-module $N$).
If the ground ring $\K$ is not a field,
then we assume in addition that $M$ (respectively, $N$)
is projective as a $\K$-module.

We have a unique natural morphism
\[N_*(M)\circ N_*(N)\,\xrightarrow{}\,N_*(M\circ N)\]
which extends the identity morphism
$N_0(M)\circ N_0(N)\,\xrightarrow{=}\,N_0(M\circ N)$
in degree $0$.
This morphism is a quasi-isomorphism.
\end{lem}

The first assertions of paragraph~\ref{NormalizedRightModule} are direct consequences
of the construction above.
We deduce the properties of the chain complex $N(L,P,P)$
from the following immediate result:

\begin{fact}
If $N$ is a (discrete) $\Sigma_*$-module,
then we have a canonical isomorphism
$N_*(M)\circ N\,\xrightarrow{\simeq}\,N_*(M\circ N)$.
\end{fact}

To be explicit,
in the situation of paragraph~\ref{NormalizedRightModule},
the natural morphism
\[N_*(C(L,P,P))\circ P\,\rightarrow\,N_*(C(L,P,P)\circ P)\,\rightarrow\,N_*(C(L,P,P))\]
restricts to isomorphisms
\[N_d(C(L,P,I))\circ P\,\xrightarrow{\simeq}\,N_d(C(L,P,I)\circ P)\,\xrightarrow{\simeq}\,N_d(C(L,P,P)).\]
We conclude immediately that $N(L,P,P)$ is a quasi-free right $P$-mo\-du\-le.
We obtain similarly
$N(L,P,R) = N(L,P,I)\circ R = N(L,P,P)\circ_P R$.

Let us observe that such properties do not hold for the chain complex $N(P,P,R)$
associated to a left $P$-module.

\section{Composite symmetric modules and trees with levels}\label{TreewiseComposite}

In sections~\ref{SimplicialBarConstruction} and~\ref{DifferentialGradedBarConstruction},
we introduce trees equipped with level structures
in order to represent the bar constructions
with coefficients $C(L,P,R)$ and $B(L,P,R)$.
In this section,
we warm up by relating the composition product of $\Sigma_*$-modules $M\circ N$
to tree structures.

\subsection{Trees with $2$ levels}\label{TwoLevelTrees}
An \emph{$n$-tree with $2$ levels} is an $n$-tree $\tau$
equipped with a level map
$l: V(\tau)\,\rightarrow\,\{0,1\}$
which satisfies the following properties:
if $e$ is an entry of the tree,
then we have $l(t(e)) = 1$;
if $e$ is an internal edge of the tree,
then we have $l(s(e)) = 1$ and $l(t(e)) = 0$;
if $e$ is the root of the tree,
then we have $l(s(e)) = 0$.
We represent a tree with $2$ levels in figure~\ref{DisplayTwoLevelTree}.
\begin{figure}
\[\xymatrix@M=0pt@R=5mm@C=5mm{ *+<3mm>{i_{1 1}}\ar[dr] & *{\cdots} &
*+<3mm>{i_{1 n_1}}\ar[dl] &
*+<3mm>{i_{2 1}}\ar[dr] & *{\cdots} &
*+<3mm>{i_{2 n_2}}\ar[dl] &
*{\cdots} &
*+<3mm>{i_{r 1}}\ar[dr] & *{\cdots} &
*+<3mm>{i_{r n_r}}\ar[dl] \\
*+<3mm>{1}\ar@{.}[r] &
*+<3mm>[o][F]{v_1}\ar@{.}[r]\ar[drrrr] &
\ar@{.}[r] & \ar@{.}[r] &
*+<3mm>[o][F]{v_2}\ar@{.}[r]\ar[dr] &
& *{\cdots} & \ar@{.}[r] &
*+<3mm>[o][F]{v_r}\ar@{.}[r]\ar[dlll] & \\
*+<3mm>{0}\ar@{.}[r] &
\ar@{.}[r] & \ar@{.}[r] & \ar@{.}[r] & \ar@{.}[r] &
*+<4mm>[o][F]{v_0}\ar@{.}[r]\ar[d] &
\ar@{.}[r] & \ar@{.}[r] & \ar@{.}[r] & \\
& & & & & *+<3mm>{0} & & & & \\ }\]
\caption{}\label{DisplayTwoLevelTree}\end{figure}

The level structure is equivalent to a partition $V(\tau) = V_0(\tau)\amalg V_1(\tau)$,
where $V_i(\tau)\subset V(\tau)$ consists of vertices $v\in V(\tau)$
such that $l(v) = i$.
Clearly,
the set $V_0(\tau)$ is reduced to the vertex $v_0 = s(e)$,
where $e$ is the root of $\tau$.
Moreover,
we have $V_1(\tau) = I_{v_0}$ and $\{1,\ldots,n\} = \coprod_{v\in V_1(\tau)} I_v$.

An isomorphism of trees with $2$ levels is an isomorphism of trees $f: \tau\,\rightarrow\,\tau'$
such that $f_V: V(\tau)\,\rightarrow\,V(\tau')$
preserves level partitions
$V(\tau) = V_0(\tau)\amalg V_1(\tau)$.
Accordingly,
as long as we omit to fix edge sets,
an isomorphism $f: \tau\,\rightarrow\,\tau'$
is equivalent to a bijection $f_1: V_1(\tau)\,\rightarrow\,V_1(\tau')$
such that
$I_{f_1(v)} = I_v$, for all $v\in V_1(\tau)$.

We let $T_{[0,1]}(n)$ denote the category formed by $n$-trees with $2$ levels
together with isomorphisms as morphisms.

\subsection{Trees with $2$ levels and composite $\Sigma_*$-modules}\label{TreewiseComposition}
We are given $\Sigma_*$-modules $M$ and $N$.
As in paragraph~\ref{TreewiseTensors} (about the construction of the free operad),
we associate to each $n$-tree with $2$ levels $\tau\in T_{[0,1]}(n)$
a $\K$-module $\tau(M,N)$.
An isomorphism of $n$-trees with $2$ levels $f: \tau\,\rightarrow\,\tau'$
induces an isomorphism of $\K$-modules
$f_*: \tau(M,N)\,\rightarrow\,\tau'(M,N)$.

To be explicit,
if $V_0(\tau) = \{v_0\}$ and $V_1(\tau) = \{v_1,\ldots,v_r\}$ (as in figure~\ref{DisplayTwoLevelTree}),
then we set
\[\tau(M,N) = M(I_{v_0})\otimes N(I_{v_1})\otimes\cdots\otimes N(I_{v_r}).\]
Equivalently:
we mark the vertex of $V_0(\tau)$ with an element $x\in M(I_{v_0})$;
we mark the vertices of $V_1(\tau)$ with elements $y_1\in N(I_{v_1}),\ldots,y_r\in N(I_{v_r})$.
Accordingly,
if $I_{v_1} = \{1,\ldots,n_1\}$, $I_{v_2} = \{n_1+1,\ldots,n_1+n_2\}$
and $I_{v_r} = \{n_1+\cdots+n_{r-1}+1,\ldots,n_1+\cdots+n_{r-1}+n_r\}$,
then the resulting labeled tree represents the composite element
$x(y_1,\ldots,y_r)\in(M\circ N)(n)$.

As in the construction of the free operad,
we consider the colimit
\[C(M,N)(n) = \colim_{\tau\in T_{[0,1]}(n)}\tau(M,N),\]
where elements $x'\in\tau'(M,N)$ and $x\in\tau(M,N)$
which correspond under an isomorphism of $n$-trees with levels $f: \tau\,\rightarrow\,\tau'$
are identified.

\begin{lem}\label{CoendTreewiseComposite}
The construction above supplies an expansion
of the composition product of $\Sigma_*$-mo\-du\-les.
We have precisely $(M\circ N)(n) = \colim_{\tau\in T_{[0,1]}(n)}\tau(M,N)$.
\end{lem}

\begin{proof}
We borrow the notation of paragraph~\ref{TreewiseComposition}.
We observe in paragraph~\ref{TwoLevelTrees},
that an isomorphism of trees with $2$ levels $f: \tau\,\rightarrow\,\tau'$,
where $V_0(\tau) = \{v_0\}$ and $V_1(\tau) = I_{v_0} = \{v_1,\ldots,v_r\}$
(respectively, $V_0(\tau') = \{v'_0\}$ and $V_1(\tau') =  I_{v'_0} = \{v'_1,\ldots,v'_r\}$),
is equivalent to a bijection
$f_1: \{v_1,\ldots,v_r\}\,\rightarrow\,\{v'_1,\ldots,v'_r\}$.
By making explicit the induced morphisms
$f_*: M(I_{v_0})\otimes\bigl(N(I_{v_1})\otimes\cdots\otimes N(I_{v_r})\bigr)
\,\rightarrow\,M(I_{v'_0})\otimes\bigl(N(I_{v'_1})\otimes\cdots\otimes N(I_{v'_r})\bigr)$,
we verify readily that the module $C(M,N)(n)$ is given by a coend
\[C(M,N)(I) = \int_{\{v_1,\ldots,v_r\}}
M(\{v_1,\ldots,v_r\})\otimes\bigl(N(I_{v_1})\otimes\cdots\otimes N(I_{v_r})\bigr),\]
in which we consider the category of finite sets
$V_1(\tau) = \{v_1,\ldots,v_r\}$
and bijections
$f_1: \{v_1,\ldots,v_r\}\,\rightarrow\,\{v'_1,\ldots,v'_r\}$.
According to lemma~\ref{ExpansionComposite},
the composition product of $\Sigma_*$-module is given by an equivalent coend
in which we restrict ourself to permutations
of the single set $\{1,\ldots,r\}$.
This observation achieves the proof of lemma~\ref{CoendTreewiseComposite}.
\end{proof}

\section{The simplicial bar construction}\label{SimplicialBarConstruction}

Recall that the $d$-dimensional component of the simplicial bar construction
is given by the composition product
$C_d(L,P,R) = L\circ P\circ\cdots\circ P\circ R$.
By generalizing the construction of section~\ref{TreewiseComposite},
we represent an element of this composite $\Sigma_*$-module
by a labeled tree equipped with $d+2$ levels.
The purpose of this section is to make this representation precise.

\subsection{Trees with levels}\label{LevelTrees}
We consider $n$-trees $\tau$ equipped with $d+2$ levels indexed by the integers $0,1,\ldots,d,d+1$.
To be precise,
the level structure is defined by a map $l: V(\tau)\,\rightarrow\,\{0,\ldots,d+1\}$
endowed with the following properties:
if $e$ is an entry of the tree,
then we have $l(t(e)) = d+1$;
if $e$ is an internal edge of the tree,
then we have $l(t(e)) = l(s(e))-1$;
if $e$ is the root of the tree,
then we have $l(s(e)) = 0$.
The level structure is also equivalent to a partition $V(\tau) = \coprod_{i=0}^{d+1} V_i(\tau)$,
where $V_i(\tau)\subset V(\tau)$ consists of vertices $v\in V(\tau)$
such that $l(v) = i$.
The set $V_0(\tau)$ is reduced to the source of the root of the tree.
Moreover, we have $V_i(\tau) = \coprod_{v\in V_{i-1}(\tau)} I_v$, for $i = 1,\ldots,d+1$,
and $\{1,\ldots,n\} = \coprod_{v\in V_{d+1}(\tau)} I_v$.
We observe that the structure of the tree $\tau$ is determined by these partitions.
We give an example of a tree with $3$ levels in figure~\ref{DisplayLevelTree}
\begin{figure}
\[\xymatrix@M=0pt@R=6mm@C=6mm{ & *+<3mm>{2}\ar[d] &
*+<3mm>{4}\ar[dr] & *+<3mm>{1}\ar[d] & *+<3mm>{6}\ar[dl] & *+<3mm>{3}\ar[dr] & & *+<3mm>{5}\ar[dl] \\
*+<3mm>{2}\ar@{.}[r] & *+<3mm>[o][F]{w_1}\ar@{.}[rr]\ar[dr] & & *+<3mm>[o][F]{w_2}\ar@{.}[rrr]\ar[dl] & & &
*+<3mm>[o][F]{w_3}\ar@{.}[r]\ar[d] & \\
*+<3mm>{1}\ar@{.}[rr] & & *+<3mm>[o][F]{v_1}\ar@{.}[rrrr]\ar[drr] & & & & *+<3mm>[o][F]{v_2}\ar@{.}[r]\ar[dll] & \\
*+<3mm>{0}\ar@{.}[rrrr] & & & & *+<3mm>[o][F]{u_1}\ar@{.}[rrr]\ar[d] & & & \\
& & & & *+<3mm>{0} & & & & \\ }\]
\caption{}\label{DisplayLevelTree}\end{figure}

An isomorphism of $n$-trees with levels
is an isomorphism of trees $f: \tau\,\rightarrow\,\tau'$
such that $f_V: V(\tau)\,\rightarrow\,V(\tau')$
preserves level partitions
$V(\tau) = \coprod_{i=0}^{d+1} V_i(\tau)$.
Consequently,
an isomorphism of $n$-trees with levels
is equivalent to bijections $f_i: V_i(\tau)\,\rightarrow\,V_i(\tau')$
such that
$I_{f_i(v)} = f_{i+1}(I_v)$, for all $v\in V_i(\tau)$,
and
$I_{f_{d+1}(v)} = I_v$, for all $v\in V_{d+1}(\tau)$.
We let $T_{[0,\ldots,d+1]}(n)$ denote the category formed by $n$-trees with $d+2$ levels
together with isomorphisms as morphisms.

\subsection{Trees with levels and bar simplices}\label{OperadBarSimplices}
We generalize the construction of paragraph~\ref{TreewiseComposition}.
We are given a dg-operad $P$, a right $P$-module $L$
and a left $P$-module $R$.
A $\K$-module $\tau(L,P,R)$ is associated to each tree with levels $\tau\in T_{[0,\ldots,d+1]}(n)$.
An isomorphism of trees with levels $f: \tau\,\rightarrow\,\tau'$
yields an isomorphism of $\K$-modules
$f_*: \tau(L,P,R)\,\rightarrow\,\tau'(L,P,R)$.
We set explicitly
\[\tau(L,P,R) = \Bigl\{\bigotimes_{v\in V_{0}(\tau)} L(I_v)\Bigr\}
\otimes\bigotimes_{i=1}^d\Bigl\{\bigotimes_{v\in V_i(\tau)} P(I_v)\Bigr\}
\otimes\Bigl\{\bigotimes_{v\in V_{d+1}(\tau)} R(I_v)\Bigr\}.\]
Equivalently, a vertex $v\in V(\tau)$ of level $i = l(v)$ is labeled:
by an element $x_v\in L(I_v)$, for $i = 0$;
by an operation $p_v\in P(I_v)$, for $i\in\{1,\ldots,d\}$;
and by an element $y_v\in R(I_v)$, for $i = d+1$.
We assume that the tensor product $\tau(L,P,R)$ is ordered according
to the level of vertices.

As an example,
if $\tau$ is the tree with levels displayed in figure~\ref{DisplayLevelTree},
then we consider tensors
$x_1\otimes p_1\otimes p_2\otimes y_1\otimes y_2\otimes y_3\in\tau(L,P,R)$.
The vertex $u_1$ is labeled by an element of $L$,
namely $x_1\in L(\{v_1,v_2\}$,
the vertices $v_1$ and $v_2$ by operations of $P$,
respectively $p_1\in P(\{w_1,w_2\})$ and $p_2\in P(\{w_3\})$,
and the vertices $w_1$, $w_2$ and $w_3$ by elements of $R$,
respectively $y_1\in R(\{2\})$, $y_2\in R(\{4,1,6\})$ and $y_3\in R(\{3,5\})$.

The composite $\Sigma_*$-module $C_d(L,P,R) = L\circ P\circ\cdots\circ P\circ R$
is equivalent to the colimit
\[C_d(L,P,R)(n) = \colim_{\tau\in T_{[0,\ldots,d+1]}(n)}\tau(L,P,R),\]
where elements $x'\in\tau'(L,P,R)$ and $x\in\tau(L,P,R)$
which correspond under an isomorphism of $n$-trees with levels
$f: \tau\,\rightarrow\,\tau'$
are identified.

\subsection{Unit vertices and the case of trivial coefficients}
A \emph{unit vertex} in a tree is a vertex
which has exactly one entry
and
which can be labeled only by a unit element
as in the following figure
$\xymatrix@M=0pt@R=6mm@C=6mm{ \ar[r] & *+<3mm>[o][F]{1}\ar[r] & \\ }$.

We consider the bar construction with coefficients in trivial modules $L = R = I$.
We observe that the $\K$-module $\tau(I,P,I)$ vanishes unless the vertices at levels $0$ and $d+1$
are all unit vertices.
We can remove these unit vertices,
so that the tree with $d+2$ levels $\tau$ is equivalent to a tree $\bar{\tau}$
with $d$ levels indexed by the integers $1,\ldots,d$.
Accordingly,
in the case of trivial coefficients,
we consider the category $T_{[1,\ldots,d]}(n)$
which consists of $n$-trees equipped with $d$ levels indexed by the integers $1,\ldots,d$.
For all $\bar{\tau}\in T_{[1,\ldots,d]}(n)$,
we form the tensor products
$\bar{\tau}(P) = \bigotimes_{i=1}^d\bigl\{\bigotimes_{v\in V_i(\tau)} P(I_v)\bigr\}$.
The module $C_d(I,P,I)(n)$ is represented by the colimit
\[C_d(I,P,I)(n) = \colim_{\bar{\tau}\in T_{[1,\ldots,d]}(n)}\bar{\tau}(P).\]
We obtain similarly
\begin{align*}
& C_d(L,P,I)(n) = \colim_{\tau\in T_{[0,\ldots,d]}(n)}\tau(L,P) \\
\text{and}\qquad & C_d(I,P,R)(n) = \colim_{\tau\in T_{[1,\ldots,d+1]}(n)}\tau(P,R).
\end{align*}
We remove level $d+1$ from tree structures in the first case
and level $0$ in the second case.

\subsection{Degeneracies}
Recall that the degeneracy
\[s_j: C_d(L,P,R)\,\rightarrow\,C_{d+1}(L,P,R)\]
is obtained
by inserting a unit morphism $\eta: I\,\rightarrow\,P$ between the $j$th and $j+1$th operad factors
in $C_d(L,P,R) = L\circ P\circ\cdots\circ P\circ R$.

Accordingly,
the image of a tensor $\bigotimes_v x_v\in\tau(L,P,R)$ under a degeneracy
$s_j: C_d(L,P,R)\,\rightarrow\,C_{d+1}(L,P,R)$
is obtained by inserting unit vertices
between the $j$th and $j+1$th levels of $\tau$.
More formally,
we consider a tree with levels $s_j(\tau)$
such that
\[V_k(s_j(\tau)) = \begin{cases} V_k(\tau), & \text{for $k = 0,\ldots,j,j+1$}, \\
V_{k-1}(\tau), & \text{for $k = j+2,\ldots,d+2$}. \end{cases}\]
We have
$s_j(\bigotimes_v x_v)\in s_j(\tau)(L,P,R)$.
To be more precise,
a vertex $v\in V_k(s_j(\tau))$ of level $k\not=j+1$
corresponds to a vertex of $\tau$
and has unchanged entries $I_v\subset V_{k+1}(s_j(\tau))$
and an unchanged label $x_v\in P(I_v)$ in $s_j(\tau)$.
Since we have a bijection $V_{j+1}(s_j(\tau))\simeq V_{j+2}(s_j(\tau))$,
a vertex $v\in V_k(s_j(\tau))$ of level $k=j+1$
is connected to a vertex of level $j+2$
and is labeled by an operad unit $1\in P(1)$.

As an example,
the $s_1$-degeneracy
of the tensor $x_1\otimes p_1\otimes p_2\otimes y_1\otimes y_2\otimes y_3\in\tau(L,P,R)$,
introduced in paragraph~\ref{OperadBarSimplices},
is represented in figure~\ref{DegenerateLevelTree}.
\begin{figure}
\[\xymatrix@M=0pt@R=6mm@C=6mm{ & *+<3mm>{2}\ar[d] &
*+<3mm>{4}\ar[dr] & *+<3mm>{1}\ar[d] & *+<3mm>{6}\ar[dl] & *+<3mm>{3}\ar[dr] & & *+<3mm>{5}\ar[dl] \\
*+<3mm>{3}\ar@{.}[r] & *+<3mm>{y_1}\ar@{.}[rr]\ar[d] & & *+<3mm>{y_2}\ar@{.}[rrr]\ar[d] & & &
*+<3mm>{y_3}\ar@{.}[r]\ar[d] & \\
*+<3mm>{2}\ar@{.}[r] & *+<3mm>{1}\ar@{.}[rr]\ar[dr] & & *+<3mm>{1}\ar@{.}[rrr]\ar[dl] & & &
*+<3mm>{1}\ar@{.}[r]\ar[d] & \\
*+<3mm>{1}\ar@{.}[rr] & & *+<3mm>{p_1}\ar@{.}[rrrr]\ar[drr] & & & & *+<3mm>{p_2}\ar@{.}[r]\ar[dll] & \\
*+<3mm>{0}\ar@{.}[rrrr] & & & & *+<3mm>{x_1}\ar@{.}[rrr]\ar[d] & & & \\
& & & & *+<3mm>{0} & & & & \\ }\]
\caption{}\label{DegenerateLevelTree}\end{figure}

\subsection{Faces and level contractions}
We make explicit the process (called level contraction)
which gives the $d_i$-face of an element of $C_d(L,P,R)$
represented by a labeled tree.
We recall that a face $d_i: C_d(L,P,R)\,\rightarrow\,C_{d-1}(L,P,R)$
is induced by the product $\mu: P\circ P\,\rightarrow\,P$
of the $i$th and $i+1$th copies of the operad $P$
in the composite $C_d(L,P,R) = L\circ P\circ\cdots\circ P\circ R$.
In the case $i = 0$ (respectively, $i = d$)
we just replace the operad product $\mu: P\circ P\,\rightarrow\,P$
by the operad action $\lambda: L\circ P\,\rightarrow\,L$
(respectively, $\rho: P\circ R\,\rightarrow\,R$).

Accordingly,
the image of a tensor $\bigotimes_v x_v\in\tau(L,P,R)$
under a face $d_i: C_d(L,P,R)\,\rightarrow\,C_{d-1}(L,P,R)$
is obtained by composing the elements of levels $i$ and $i+1$ in $\tau$
according to the tree structure.
More formally,
we consider a tree with levels $d_i(\tau)$
such that
\[V_k(d_i(\tau)) = \begin{cases} V_k(\tau), & \text{for $k = 0,\ldots,i$}, \\
V_{k+1}(\tau), & \text{for $k = i+1,\ldots,d$}. \end{cases}\]
We have
$d_i(\bigotimes_v x_v)\in d_i(\tau)(L,P,R)$.
To be more precise,
a vertex $v\in V_k(d_i(\tau))$ of level $k\not=i$
has unchanged entries $I_v\subset V_{k+1}(d_i(\tau))$
and an unchanged label $x_v\in P(I_v)$ (respectively, $x_v\in L(I_v)$, $x_v\in R(I_v)$)
in $d_i(\tau)$.
We consider a vertex $v_0\in V_k(d_i(\tau))$ of level $k = i$.
In the original tree $\tau$,
we have $I_{v_0} = \{v_1,\ldots,v_r\}$,
where $v_1,\ldots,v_r\in V_{i+1}(\tau)$.
In the tree $d_i(\tau)$,
we connect the vertex $v_0$ to the vertices
of $I_{v_1}\amalg\cdots\amalg I_{v_r}\subset V_{i+2}(\tau) = V_{i+1}(d_i(\tau))$.
Accordingly, in $d_i(\tau)$,
we have $I_{v_0} = I_{v_1}\amalg\cdots\amalg I_{v_r}$.
The operad product $P(\{v_1,\ldots,v_r\})\otimes(P(I_{v_1})\otimes\cdots\otimes P(I_{v_1}))
\,\rightarrow\,P(I_{v_1}\amalg\cdots\amalg I_{v_r})$
associates a composite
$x_{v_0}(x_{v_1},\ldots,x_{v_r})\in P(I_{v_1}\amalg\cdots\amalg I_{v_r})$
to the elements which label the vertices $v_0,v_1,\ldots,v_r$
in $\tau(L,P,R)$.
We mark the vertex $v$ with this composite $x_{v_0}(x_{v_1},\ldots,x_{v_r})\in P(I_{v_1}\amalg\cdots\amalg I_{v_r})$
in order to obtain the element $d_i(\bigotimes_v x_v)\in d_i(\tau)(L,P,R)$.
In cases $i = 0$ and $i = d$,
we just replace the operad product above by the operad actions
$L(\{v_1,\ldots,v_r\})\otimes(P(I_{v_1})\otimes\cdots\otimes P(I_{v_1}))
\,\rightarrow\,L(I_{v_1}\amalg\cdots\amalg I_{v_r})$
and
$P(\{v_1,\ldots,v_r\})\otimes(R(I_{v_1})\otimes\cdots\otimes R(I_{v_1}))
\,\rightarrow\,R(I_{v_1}\amalg\cdots\amalg I_{v_r})$.

As an example,
the faces of the tensor $x_1\otimes p_1\otimes p_2\otimes y_1\otimes y_2\otimes y_3\in\tau(L,P,R)$
introduced in paragraph~\ref{OperadBarSimplices}
are represented in figure~\ref{FacesLevelTree}.
\begin{figure}
\[\xymatrix@M=0pt@R=5mm@C=5mm{ & *+<3mm>{2}\ar[dr] &
*+<3mm>{4}\ar[d] & *+<3mm>{1}\ar[dl] & *+<3mm>{6}\ar[dll] & *+<3mm>{3}\ar[dr] & & *+<3mm>{5}\ar[dl] \\
*+<3mm>{1}\ar@{.}[rr] & & *+<3mm>{p_1(y_1,y_2)}\ar@{.}[rrrr]\ar[drr] & & & & *+<3mm>{p_2(y_3)}\ar@{.}[r]\ar[dll] & \\
*+<3mm>{0}\ar@{.}[rrrr] & & & & *+<3mm>{x_1}\ar@{.}[rrr]\ar[d] & & & \\
& & & & *+<3mm>{0} & & & & \\
& & & & & & & & \\
& *+<3mm>{2}\ar[d] & *+<3mm>{4}\ar[dr] & *+<3mm>{1}\ar[d] & *+<3mm>{6}\ar[dl] & *+<3mm>{3}\ar[dr] & & *+<3mm>{5}\ar[dl] \\
*+<3mm>{1}\ar@{.}[r] & *+<3mm>{y_1}\ar@{.}[rr]\ar[drrr] & & *+<3mm>{y_2}\ar@{.}[rrr]\ar[dr] & & &
*+<3mm>{y_3}\ar@{.}[r]\ar[dll] & \\
*+<3mm>{0}\ar@{.}[rrrr] & & & & *+<3mm>{x_1(p_1,p_2)}\ar@{.}[rrr]\ar[d] & & & \\
& & & & *+<3mm>{0} & & & & \\ }\]
\caption{}\label{FacesLevelTree}\end{figure}

\section{The differential graded bar construction}\label{DifferentialGradedBarConstruction}

In this section,
we give the details of the definition of the differential graded bar construction
with coefficients $B(L,P,R)$.
Since $B(L,P,R) = L\circ\bar{B(P)}\circ R$,
an element of $B(L,P,R)$ is represented by a tree with $3$ levels
whose vertices are labeled by elements of $L$, $\bar{B}(P)$, $R$.
Recall that an element of $\bar{B}(P)(I)$
is represented by a labeled tree $\tau$
whose entries are in bijection with the elements of $I$.
Consequently,
if this element labels a vertex $v$ in a tree with $3$ levels,
then the entries of $\tau$ are in bijection with the entries of $v$.
Therefore,
we can replace the vertex $v$ by the labeled tree $\tau$
in the construction of $B(L,P,R)$.
We make this construction precise in the next paragraphs.
We define the differential of $B(L,P,R)$ in paragraph~\ref{CoefficientBarDifferential}.

\subsection{Composite trees}\label{CompositeTrees}
We consider trees, called \emph{composite\linebreak trees},
which are equipped with a lower level, a main level
and an upper level.
Precisely,
a composite $n$-tree $\tau$ is equipped with a map $l: V(\tau)\,\rightarrow\,\{l,p,r\}$
which verifies the following properties:
we have $l(v) = r$ if and only if the vertex $v$ is the target of an entry of the tree;
we have $l(v) = l$ if and only if the vertex $v$ is the source of the root of the tree;
otherwise, we have $l(v) = p$.
Accordingly,
let $e$ be an internal edge of $\tau$.
If $l(s(e)) = r$, then we have either $l(t(e)) = p$ or $l(t(e)) = l$.
If $l(t(e)) = l$, then we have either $l(s(e)) = p$ or $l(s(e)) = r$.
We may have also $l(s(e)) = l(t(e)) = p$.

Consequently,
the vertex set $V(\tau)$ is equipped with a partition
$V(\tau) = V_l(\tau)\amalg V_p(\tau)\amalg V_r(\tau)$
where $V_x(\tau)\subset V(\tau)$
consists of vertices $v\in V(\tau)$
of level $l(v) = x$.
By definition,
a vertex $v\in V(\tau)$ is connected to an entry of the tree if and only if $v\in V_r(\tau)$.
In this case,
we have $I_v\subset\{1,\ldots,n\}$ (the entries of $v\in V_r(\tau)$ are all entries of the tree).

An example of a composite tree is represented in figure~\ref{DisplayCompositeTree}
\begin{figure}
\[\xymatrix@M=0pt@R=6mm@C=6mm{ *+<3mm>{2}\ar[dr] & & *+<3mm>{4}\ar[dl] &
*+<3mm>{6}\ar[d] &
*+<3mm>{1}\ar[d] &
*+<3mm>{3}\ar[d] &
*+<3mm>{5}\ar[d] \\
*+<3mm>{R}\ar@{.}[r] & *+<3mm>[o][F]{w_1}\ar[dr]\ar@{.}[rr] & &
*+<3mm>[o][F]{w_2}\ar[dl]\ar@{.}[r] &
*+<3mm>[o][F]{w_3}\ar[dr]\ar@{.}[r] &
*+<3mm>[o][F]{w_4}\ar[d]\ar@{.}[r] &
*+<3mm>[o][F]{w_5}\ar[dl] \\
P\save\left\{\vbox to 8mm{ }\right.\restore & &
*+<3mm>[o][F]{v_1}\ar[dr] & & & *+<3mm>[o][F]{v_2}\ar[dll] & \\
*+<3mm>{L}\ar@{.}[rrr] & & & *+<3mm>[o][F]{u_1}\ar@{.}[rrr]\ar[d] & & & \\
& & & *+<3mm>{0} & & & \\ }\]
\caption{}\label{DisplayCompositeTree}\end{figure}

An isomorphism of composite trees is an isomorphism of $n$-trees
which preserves level structures.
We denote the category of composite $n$-trees by $T_{\{l,p,r\}}(n)$.
We consider also reduced composite trees $\tau\in\tilde{T}_{\{l,p,r\}}(n)$
in which vertices $v\in V_p(\tau)$
are supposed to have at least $2$ entries.

\subsection{The expansion of the bar construction with coefficients}\label{ExpansionCoefficientBar}
We are given a connected dg-operad $P$, a right $P$-module $L$ and a left $P$-module $R$.
We associate a dg-module $\tau(L,P,R)$ to each composite tree $\tau\in T_{\{l,p,r\}}(n)$.
An isomorphism of composite trees $f: \tau\,\rightarrow\,\tau'$
induces an isomorphism of dg-modules $f_*: \tau(L,P,R)\,\rightarrow\,\tau'(L,P,R)$.
We set precisely
\[\tau(L,P,R) = \Bigl\{\bigotimes_{v\in V_l(\tau)} L(I_v)\Bigr\}
\otimes\Bigl\{\bigotimes_{v\in V_p(\tau)} (\Sigma\bar{P})(I_v)\Bigr\}
\otimes\Bigl\{\bigotimes_{v\in V_R(\tau)} R(I_v)\Bigr\}.\]
Equivalently,
the unique vertex $v\in V_l(\tau)$ is labeled by an element $x\in L(I_v)$,
the vertices $v\in V_p(\tau)$ are labeled by elements $\Sigma p\in\Sigma\tilde{P}(I_v)$,
and the vertices $v\in V_r(\tau)$ by elements $y\in R(I_v)$.

According to the introduction of this section,
the composite module $B(L,P,R)(n) = L\circ\bar{B}(P)\circ R$
is equivalent to the colimit
\[B(L,P,R)(n) = \colim_{\tau\in T_{\{l,p,r\}}(n)}\tau(L,P,R),\]
where elements $x'\in\tau'(L,P,R)$ and $x\in\tau(L,P,R)$
which correspond under an isomorphism of composite $n$-trees
$f: \tau\,\rightarrow\,\tau'$
are identified.

As in the context of the simplicial bar construction,
if $R = I$ is a trivial module,
then we can remove the upper level from tree structures,
so that
$B(L,P,I)(n) = \colim_{\tau\in T_{\{l,p\}}(n)}\tau(L,P)$.
Similarly,
if $L = I$,
then we have an expansion
$B(I,P,R)(n) = \colim_{\tau\in T_{\{p,r\}}(n)}\tau(P,R)$,
which involves composite trees without lower level.

\subsection{The differential of the bar construction with coefficients}
\label{CoefficientBarDifferential}
The bar differential $\beta: B(L,P,R)\,\rightarrow\,B(L,P,R)$ is defined componentwise.
We suppose given an element
\[\Bigl\{\bigotimes_{v\in V_l(\tau)} x_v\Bigr\}
\otimes\Bigl\{\bigotimes_{v\in V_l(\tau)} \Sigma p_v\Bigr\}
\otimes\Bigl\{\bigotimes_{v\in V_l(\tau)} y_v\Bigr\}\in\tau(L,P,R).\]

{\bf a)} We consider a subtree $\sigma_e\subset\tau$ associated to an internal edge $e\in E(\tau)$
as the definition of the differential of the reduced bar construction
(see paragraph~\ref{ReducedBarDifferential}).
We have then $V(\sigma_e) = \{u,v\}$, where $u = t(e)$ and $v = s(e)$.
We assume $l(v) = l(u) = p$.
The bar differential has a component
$\beta^{\tau,\tau'}: \tau(L,P,R)\,\rightarrow\,\tau'(L,P,R)$,
such that $\tau' = \tau/\sigma_e$.
As in the definition of the differential of the reduced bar construction,
we consider the homogeneous morphism
\[\Sigma\tilde{P}(I_{u})\otimes\Sigma\tilde{P}(I_{v})
\,\rightarrow\,\Sigma\tilde{P}(I_{u}\setminus\{v\}\amalg I_{v}),\]
induced by the partial composition product of $P$.
The collapsed vertex of $\tau/\sigma_e$
is labeled by the partial composite
$- \Sigma(p_{u}\circ_{v}p_{v})\in\tilde{P}(I_{u}\setminus\{v\}\amalg I_{v})$
(see paragraph~\ref{ReducedBarDifferential}).
The vertices $x\not=u,v$ of $\tau/\sigma_e$
have unchanged label.

{\bf b)} We consider a subtree $\sigma_e\subset\tau$ as above,
but such that $l(v) = p$ and $l(u) = l$.
The bar differential has a component
$\beta^{\tau,\tau'}: \tau(L,P,R)\,\rightarrow\,\tau'(L,P,R)$,
such that $\tau' = \tau/\sigma_e$.
In this case,
we just consider the homogeneous morphism of degree $-1$
\[L(I_{u})\otimes\Sigma\tilde{P}(I_{v})\,\rightarrow\,L(I_{u}\setminus\{v\}\amalg I_{v})\]
determined by the right module structure.
Accordingly,
the collapsed vertex of $\tau/\sigma_e$
is labeled by the partial composite
$x_{u}\circ_{v}p_{v}\in L(I_{u}\setminus\{v\}\amalg I_{v})$.

{\bf c)} The bar differential has also a component
\[\beta^{\tau,\tau/\sigma}: \tau(L,P,R)\,\rightarrow\,\tau/\sigma(L,P,R),\]
for each subtree $\sigma\subset\tau$
such that $V(\sigma) = \{u\}\amalg\{v_1,\ldots,v_r\}$
where $u\in V_p(\tau)$, $I_{u} = \{v_1,\ldots,v_r\}$
and $v_1,\ldots,v_r\in V_r(\tau)$.
Explicitly,
the operad action provides a homogeneous morphism of degree $-1$
\[\Sigma\tilde{P}(I_{u})\otimes R(I_{v_1})\otimes\cdots\otimes R(I_{v_1})
\,\rightarrow\,R(I_{v_1}\amalg\cdots\amalg I_{v_r})\]
which associates a composite
$p_{u}(y_{v_1},\ldots,y_{v_r})\in R(I_{v_1}\amalg\cdots\amalg I_{v_r})$
to the elements
$\Sigma p_{u}\in\Sigma\tilde{P}(I_{v_0})$ and $y_{v_1}\in R(I_{v_1}),\ldots,y_{v_r}\in R(I_{v_r})$
which label the vertices $u,v_1,\ldots,v_r$
in $\tau(L,P,R)$.
The collapsed vertex of the quotient tree $\tau/\sigma$
is labeled by the composite element
$- p_{u}(y_{v_1},\ldots,y_{v_r})\in R(I_{v_1}\amalg\cdots\amalg I_{v_r})$.
A vertex $x\not=u,v_1,\ldots,v_r$ of $\tau$
has an unchanged label in $\tau/\sigma$.

\subsection{On the components of the differential of the bar construction}\label{TwistingBarDifferential}
We recall that the chain complex $B(P,P,P)$ is both a quasi-free right $P$-module
and a quasi-free left $P$-module,
since we have by definition $B(P,P,P) = P\circ\bar{B}(P)\circ P$.
Therefore,
a construction of the differential of $B(P,P,P)$
can be deduced from observations of paragraphs~\ref{QuasiFreeRightModule} and~\ref{QuasiFreeLeftModule}
(about quasi-free modules).
More explicitly,
the differential of $B(P,P,P)$ is determined by the internal differential of the reduced bar construction
$\beta: \bar{B}(P)\,\rightarrow\,\bar{B}(P)$
and by homogeneous morphisms
\[\theta_R: \bar{B}(P)\,\rightarrow\,\bar{B}(P)\circ P
\qquad\text{and}
\qquad\theta_L: \bar{B}(P)\,\rightarrow\,P\circ\bar{B}(P).\]

We recall that elements of the reduced bar construction $\bar{B}(P) = B(I,P,I)$
are identified with elements of the bar construction with coefficients $B(P,P,P)$
which have unit vertices in lower and upper levels.
We have similar identifications for the elements of $\bar{B}(P)\circ P = B(I,P,P)$
and $P\circ\bar{B}(P) = B(P,P,I)$.
The homogeneous morphism
$\theta_R: \bar{B}(P)\,\rightarrow\,\bar{B}(P)\circ P$
is the restriction to $\bar{B}(P)\subset B(P,P,P)$
of a bimodule derivation
$d_{\theta_R}: P\circ\bar{B}(P)\circ P\,\rightarrow\,P\circ\bar{B}(P)\circ P$
and consists of components of the differential $\beta: B(P,P,P)\,\rightarrow\,B(P,P,P)$
defined in paragraph~\ref{CoefficientBarDifferential}{\bf .c)}.
The homogeneous morphism
$\theta_L: \bar{B}(P)\,\rightarrow\,P\circ\bar{B}(P)$
is the restriction of a bimodule derivation
$d_{\theta_L}: P\circ\bar{B}(P)\circ P\,\rightarrow\,P\circ\bar{B}(P)\circ P$
and consists of components of the differential $\beta: B(P,P,P)\,\rightarrow\,B(P,P,P)$
defined in paragraph~\ref{CoefficientBarDifferential}{\bf .b)}.

We have $B(L,P,R) = L\circ_P B(P,P,P)\circ_P R$
and the differential of $B(L,P,R)$ is induced by the differential of $B(P,P,P)$.
Let us mention
that the derivation $d_{\theta_R}: P\circ\bar{B}(P)\circ P\,\rightarrow\,P\circ\bar{B}(P)\circ P$ vanishes in $B(P,P,I)$
and $d_{\theta_L}: P\circ\bar{B}(P)\circ P\,\rightarrow\,P\circ\bar{B}(P)\circ P$ vanishes in $B(I,P,P)$.

\subsection{Twisting cochains}\label{TwistingCochains}
We would like to point out that the maps $\theta_R$ and $\theta_L$
are determined by a certain homogeneous morphism
$\phi: \bar{B}(P)\,\rightarrow\,P$.
Explicitly,
we consider the composite of the projection morphism
\[F^c(\Sigma\tilde{P}) = \bigoplus_{r=0}^\infty F^c_{(r)}(\Sigma\tilde{P})\,\rightarrow\,F^c_{(1)}(\Sigma\tilde{P})\]
with the canonical isomorphism
\[F^c_{(1)}(\Sigma\tilde{P})\,\simeq\,\Sigma\tilde{P},\]
so that we obtain a homogeneous morphism
$\phi: \bar{B}(P) = F^c(\Sigma\tilde{P})\,\rightarrow\,P$
of degree $-1$.
According to E. Getzler and J. Jones,
this definition gives a \emph{universal twisting cochain} $\phi: D\,\rightarrow\,P$
on the bar cooperad $D = \bar{B}(P)$
(\emph{cf}. \cite{GetzlerJones}).
We just explain how to recover morphisms
$\theta_R: D\,\rightarrow\,D\circ P$
and
$\theta_L: D\,\rightarrow\,P\circ D$
from $\phi: D\,\rightarrow\,P$.
We give the general construction.
It should be clear that we recover the components of the bar differential
defined in paragraph~\ref{CoefficientBarDifferential}
in the case $D = \bar{B}(P)$.

As in the definition of the reduced cobar construction,
the partial coproducts of a cooperad $D$ determine a morphism
$\nu: D\,\rightarrow\,F^c_{(2)}(D)$
which is homogeneous of degree $0$.
We observe that the map $\phi: D\,\rightarrow\,P$
induces morphisms
$1*\phi: F^c_{(2)}(D)\,\rightarrow\,D\circ P$
and
$\phi*1: F^c_{(2)}(D)\,\rightarrow\,P\circ D$.
The morphisms $\theta_R$ and $\theta_L$
are given by the composites
\[D\,\xrightarrow{\nu}\,F^c_{(2)}(D)\,\xrightarrow{1*\phi}\,D\circ P
\qquad\text{and}
\qquad D\,\xrightarrow{\nu}\,F^c_{(2)}(D)\,\xrightarrow{\phi*1}\,P\circ D.\]
To be explicit,
we recall that an element of $F^c_{(2)}(D)$
is represented by a labeled tree $\tau$
with two vertices $u$ and $v$
(\emph{cf}. figure~\ref{TwoVertexTree}).
We assume that the lower and upper vertices of $\tau$ are labeled by $x\in D(I_u)$ and $y\in D(I_v)$
respectively.
We set precisely
$(1*\phi)(x\otimes y) = x(1,\ldots,\phi(y),\ldots,1)\in D\circ P$
and
$(\phi*1)(x\otimes y) = \phi(x)(1,\ldots,y,\ldots,1)\in P\circ D$.
Graphically,
these elements are represented by the labeled trees with two levels
of figure~\ref{TwistingDifferentials}.
\begin{figure}
\[\xymatrix@M=0pt@R=6mm@C=6mm{ & *+<3mm>{i_{1}}\ar[d] & *+<3mm>{i_{2}}\ar[d] & *{\cdots} & *+<3mm>{i_{k-1}}\ar[d] &
*+<3mm>{j_{1}}\ar[dr] & *{\cdots} & *+<3mm>{j_{l}}\ar[dl] \\
*+<3mm>{1}\ar@{.}[r] &
*+<3mm>{1}\ar[drr]\ar@{.}[r] & *+<3mm>{1}\ar[dr]\ar@{.}[r] &
*+<3mm>{\cdots}\ar@{.}[r] & *+<3mm>{1}\ar[dl]\ar@{.}[rr] &
& *+<3mm>{\phi(y)}\ar[dlll]\ar@{.}[r] & \\
*+<3mm>{0}\ar@{.}[rrr] & & & *+<3mm>{x}\ar[d]\ar@{.}[rrrr] & & & & \\
& & & *+<3mm>{0} & & & & \\
& & & & & & & \\
& *+<3mm>{i_{1}}\ar[d] & *+<3mm>{i_{2}}\ar[d] & *{\cdots} & *+<3mm>{i_{k-1}}\ar[d] &
*+<3mm>{j_{1}}\ar[dr] & *{\cdots} & *+<3mm>{j_{l}}\ar[dl] \\
*+<3mm>{1}\ar@{.}[r] &
*+<3mm>{1}\ar[drr]\ar@{.}[r] & *+<3mm>{1}\ar[dr]\ar@{.}[r] &
*+<3mm>{\cdots}\ar@{.}[r] & *+<3mm>{1}\ar[dl]\ar@{.}[rr] &
& *+<3mm>{y}\ar[dlll]\ar@{.}[r] & \\
*+<3mm>{0}\ar@{.}[rrr] & & & *+<3mm>{\phi(x)}\ar[d]\ar@{.}[rrrr] & & & & \\
& & & *+<3mm>{0} & & & & \\ }\]
\caption{}\label{TwistingDifferentials}\end{figure}

In the case of the universal twisting cochain $\phi: \bar{B}(P)\,\rightarrow\,P$,
the composite morphisms
$\theta_R = 1*\phi\cdot\nu$ and $\theta_L = \phi*1\cdot\nu$
drop simply extremal vertices from labeled trees
(which represent elements of $\bar{B}(P)$).
Therefore, we recover the definition of paragraph~\ref{CoefficientBarDifferential}.

\section{The levelization morphism}\label{LevelizationProcess}

We define the \emph{levelization morphism}
$\phi(L,P,R): B_*(L,P,R)\,\rightarrow\,N_*(L,P,R)$
in this section.
We prove in lemma~\ref{LevelizationDifferential}
that the levelization morphism
is a morphism of chain complexes.
We assume that $P$ is a connected operad, $L$ is a right $P$ module and $R$ is a connected left $P$ module.

\subsection{Levelization of composite trees}\label{TreeLevelization}
We define a relation between composite trees and certain trees with levels.
We introduce the word \emph{levelization} for this process.

We assume that $\tau^l$ is a tree with $d+2$ levels such that, for $i = 1,\ldots,d$,
the set $V_i(\tau^l)$ has only one non-unital vertex $v_i\in V_i(\tau^l)$.
We recall that a unit vertex is a vertex
which has exactly one entry
and which can be labeled only by a unit element.
We associate a composite tree $\tau$ to the level tree $\tau^l$
by removing all unit vertices in $\tau^l$.
We set precisely $V_l(\tau) = V_0(\tau^l)$, $V_p(\tau) = \{v_1,\ldots,v_d\}$,
and $V_r(\tau) = V_{d+1}(\tau^l)$.
A sequence of egdes in $\tau^l$
\[\xymatrix@M=0pt@R=6mm@C=6mm{ *+<3mm>[o][F]{u}\ar[r]^{e_1} &
*+<3mm>[o][F]{1}\ar[r]^{e_2} & *+<3mm>{\cdots}\ar[r]^{e_{m-1}} &
*+<3mm>[o][F]{1}\ar[r]^{e_m} & *+<3mm>[o][F]{v} \\ },\]
where intermediate vertices are all unit vertices,
defines an edge from $u$ to $v$ in the composite tree $\tau$.

We call the level tree $\tau^l$
a levelization of the associated composite tree $\tau$.
We observe that a levelization $\tau^l$ of a given composite tree $\tau$
is determined by a bijection
$l: V_p(\tau)\,\rightarrow\,\{1,\ldots,d\}$
such that $l(t(e))<l(s(e))$
for all edges $e\in E(\tau)$
which have $s(e)\in V_p(\tau)$ and $t(e)\in V_p(\tau)$.
Precisely,
we can arrange the vertices of $\tau$ on levels of $\tau^l$
according to the map
$l: V_p(\tau)\,\rightarrow\,\{1,\ldots,d\}$.
As an example,
the composite tree of figure~\ref{DisplayCompositeTree} has $2$ levelizations $\tau^{l_1}$ and $\tau^{l_2}$,
which are represented in figure~\ref{LevelizedTree}.
\begin{figure}
\[\xymatrix@M=0pt@R=6mm@C=6mm{ *+<3mm>{2}\ar[dr] & & *+<3mm>{4}\ar[dl] &
*+<3mm>{6}\ar[d] & *+<3mm>{1}\ar[d] & *+<3mm>{3}\ar[d] & *+<3mm>{5}\ar[d] \\
*+<3mm>{3}\ar@{.}[r] & *+<3mm>[o][F]{w_1}\ar[d]\ar@{.}[rr] & &
*+<3mm>[o][F]{w_2}\ar[d]\ar@{.}[r] & *+<3mm>[o][F]{w_3}\ar[dr]\ar@{.}[r] &
*+<3mm>[o][F]{w_4}\ar[d]\ar@{.}[r] & *+<3mm>[o][F]{w_5}\ar[dl] \\
*+<3mm>{2}\ar@{.}[r] & *+<3mm>{1}\ar[dr]\ar@{.}[rr] & &
*+<3mm>{1}\ar[dl]\ar@{.}[rr] & & *+<3mm>[o][F]{v_2}\ar[d]\ar@{.}[r] & \\
*+<3mm>{1}\ar@{.}[rr] & & *+<3mm>[o][F]{v_1}\ar[dr]\ar@{.}[rrr] & & & *+<3mm>{1}\ar[dll]\ar@{.}[r] & \\
*+<3mm>{0}\ar@{.}[rrr] & & & *+<3mm>[o][F]{u_1}\ar[d]\ar@{.}[rrr] & & & \\
& & & *+<3mm>{0} & & & \\
& & & & & & \\
*+<3mm>{2}\ar[dr] & & *+<3mm>{4}\ar[dl] &
*+<3mm>{6}\ar[d] & *+<3mm>{1}\ar[d] & *+<3mm>{3}\ar[d] & *+<3mm>{5}\ar[d] \\
*+<3mm>{3}\ar@{.}[r] & *+<3mm>[o][F]{w_1}\ar[dr]\ar@{.}[rr] & &
*+<3mm>[o][F]{w_2}\ar[dl]\ar@{.}[r] & *+<3mm>[o][F]{w_3}\ar[d]\ar@{.}[r] &
*+<3mm>[o][F]{w_4}\ar[d]\ar@{.}[r] & *+<3mm>[o][F]{w_5}\ar[d] \\
*+<3mm>{2}\ar@{.}[rr] & & *+<3mm>[o][F]{v_1}\ar[d]\ar@{.}[rr] & &
*+<3mm>{1}\ar[dr]\ar@{.}[r] & *+<3mm>{1}\ar[d]\ar@{.}[r] & *+<3mm>{1}\ar[dl] \\
*+<3mm>{1}\ar@{.}[rr] & & *+<3mm>{1}\ar[dr]\ar@{.}[rrr] & & & *+<3mm>[o][F]{v_2}\ar[dll]\ar@{.}[r] & \\
*+<3mm>{0}\ar@{.}[rrr] & & & *+<3mm>[o][F]{u_1}\ar[d]\ar@{.}[rrr] & & & \\
& & & *+<3mm>{0} & & & \\
}\]
\caption{}\label{LevelizedTree}\end{figure}

\begin{fact}\label{LevelizedTensors}
We assume that $\tau^l$ is a levelization of a composite tree $\tau$.
We have a canonical isomorphism $\tau^l(L,P,R)\,\xrightarrow{\simeq}\,\tau(L,P,R)$.
\end{fact}

We let $v_i\in V(\tau)$ denote the vertex of $\tau$
which lies on level $i$ in $\tau^l$.
We have
\[\tau^l(L,P,R) = \Bigl\{\bigotimes_{v\in V_l(\tau)} L(I_v)\Bigr\}
\otimes\Bigl\{\bigotimes_{i=1}^d P(I_{v_i})\Bigr\}
\otimes\Bigl\{\bigotimes_{v\in V_r(\tau)} R(I_v)\Bigr\},\]
because operad units and unital vertices can be removed from the definition of the module $\tau^l(L,P,R)$.
We recall that this tensor product is ordered according to the level of vertices.
We consider the isomorphism
\[\Sigma^d\tau^l(L,P,R)\,\simeq\,\Bigl\{\bigotimes_{v\in V_l(\tau)} L(I_v)\Bigr\}
\otimes\Bigl\{\bigotimes_{i=1}^d \Sigma\tilde{P}(I_{v_i})\Bigr\}
\otimes\Bigl\{\bigotimes_{v\in V_r(\tau)} R(I_v)\Bigr\}\]
which moves suspensions from left to right.
For that purpose,
we generalize the definition of paragraph~\ref{SuspensionTensorProduct}.
In particular,
we assume that our isomorphism involves a sign
yielded by the permutation of suspension symbols
with homogeneous tensors.

As an example,
we consider the composite tree $\tau$ of figure~\ref{DisplayCompositeTree}
and the associated trees with levels $\tau^{l_1}$ and $\tau^{l_2}$,
represented in figure~\ref{LevelizedTree}.
Suppose given an element
\begin{multline*}
x_1\otimes\Sigma p_1\otimes\Sigma p_2\otimes y_1\otimes\cdots\otimes y_5 \\
\in\underbrace{L(I_{u_1})\otimes\Sigma\tilde{P}(I_{v_1})\otimes\Sigma\tilde{P}(I_{v_2})
\otimes R(I_{w_1})\otimes\cdots\otimes R(I_{w_5})}_{\displaystyle\tau(L,P,R)}.
\end{multline*}
The levelization process gives
\begin{align*}
& x_1\otimes\Sigma p_1\otimes\Sigma p_2\otimes y_1\otimes\cdots\otimes y_5 \\
& \qquad\begin{aligned}[t]
\mapsto\,\pm x_1\otimes\bigl\{p_1\otimes 1\bigr\}
\otimes\bigl\{1\otimes 1\otimes p_2\bigr\}
\otimes\bigl\{y_1\otimes\cdots\otimes y_5\bigr\} \\
\in\tau^{l_1}(L,P,R)
\end{aligned} \\
& x_1\otimes\Sigma p_1\otimes\Sigma p_2\otimes y_1\otimes\cdots\otimes y_5 \\
& \qquad\simeq\,\pm x_1\otimes\Sigma p_2\otimes\Sigma p_1\otimes y_1\otimes\cdots\otimes y_5 \\
& \qquad\begin{aligned}[t]
\mapsto\,\pm x_1\otimes\bigl\{1\otimes p_2\bigr\}
\otimes\bigl\{p_1\otimes 1\otimes 1\otimes 1\bigr\}
\otimes\bigl\{y_1\otimes\cdots\otimes y_5\bigr\} \\
\in\tau^{l_2}(L,P,R)
\end{aligned}
\end{align*}

On the other hand,
we have a canonical isomorphism $\tau^{l_1}(L,P,R)\,\simeq\,\tau^{l_2}(L,P,R)$
defined by tensor permutations
(after removing operad units from tensor products).
But,
as in paragraph~\ref{SuspensionTensorProduct},
this isomorphism makes the isomorphisms
\[\tau(L,P,R)\,\xrightarrow{\simeq}\,\tau^{l_1}(L,P,R)\qquad\text{and}\qquad\tau(L,P,R)\,\xrightarrow{\simeq}\,\tau^{l_2}(L,P,R)\]
differ by the sign $-1$,
since the levelization process involves a permutation of suspension symbols.

\subsection{The levelization morphism}
The levelization morphism
\[\phi(L,P,R): B(L,P,R)\,\rightarrow\,N(L,P,R)\]
is defined componentwise.
Namely,
for each composite tree $\tau$,
we add the isomorphisms
\[\tau(L,P,R)\,\xrightarrow{\simeq}\,\tau^l(L,P,R)\subset N(L,P,R)\]
associated to levelizations of $\tau$.

\begin{lem}\label{LevelizationDifferential}
The levelization process defines a natural morphism of dg-modules
\[\phi(L,P,R): B(L,P,R)\,\rightarrow\,N(L,P,R).\]
This morphism is an embedding provided that $P$ is a connected operad and $R$ is a connected $\Sigma_*$-module.
\end{lem}

\begin{proof}
We check that the levelization morphism preserves differentials.
We are given a composite tree $\tau$ of degree $d$.
We recall that the bar differential has a component
$\beta^{\tau/\sigma}: \tau(L,P,R)\,\rightarrow\,\tau/\sigma(L,P,R)$
for each subtree $\sigma\subset\tau$ such that $V(\sigma) = \{u\}\amalg I_u$,
where $u\in V_p(\sigma)$ and $I_u\subset V_r(\tau)$,
and for each subtree $\sigma = \xymatrix@M=0pt@R=6mm@C=6mm{ *+<3mm>[o][F]{v}\ar[r]^{e} & *+<3mm>[o][F]{u} \\ }$
where $u\in V_p(\sigma)$ or $u\in V_l(\sigma)$ and $v\in V_p(\sigma)$.
In the first case,
all levelizations of the differential $\beta^{\tau,\tau/\sigma}$
occur as the $d_d$-face of a uniquely determined levelization of $\tau$,
which has the vertex $u$ on level $d$.
In the second case,
all levelizations of the differential $\beta^{\tau,\tau/\sigma}$
occur as the $d_i$-face of a uniquely determined levelization of $\tau$,
which has the vertex $u$ on level $i$ and the vertex $v$ on level $i+1$.
We prove that the faces $d_i(\tau^l)$
which do not occur in this correspondence
cancel two by two.

We consider a level contraction $d_i(\tau^l)$,
where $\tau^l$ is a levelization of $\tau$.
We assume that $\tau^l$ is determined by a bijection $l: V_p(\tau)\,\rightarrow\,\{1,\ldots,d\}$.
We consider the case $i = d$ which gives a contraction of the top levels of $\tau^l$.
We let $u\in V(\tau)$ denote the vertex of $\tau$ which lies on level $d$ in $\tau^l$.
We have by definition $I_u\subset V_{d+1}(\tau^l) = V_r(\tau)$.
We consider the subtree $\sigma\subset\tau$ such that $V(\sigma) = \{u\}\amalg I_u$,
We observe that the contraction of levels $d$ and $d+1$ in $\tau^l$
is equivalent to the contraction of the tree $\sigma$
in $\tau$.
More precisely,
the tree $d_d(\tau^l)$ coincides with the levelization of $\tau/\sigma$
determined by the map $l': V_p(\tau/\sigma)\,\rightarrow\,\{1,\ldots,d-1\}$
such that $l'(w) = l(w)$ for all $w\in V_p(\tau)\setminus\{u\}$.
Moreover,
we observe that the isomorphisms involved in the definition of the levelization morphism
match in a commutative diagram
\[\xymatrix{\tau(L,P,R)\ar[r]^{\simeq}\ar[d]_{\beta^{\tau,\tau/\sigma}} & \tau^{l}(L,P,R)\ar[d]^{(-1)^d d_d} \\
\tau/\sigma(L,P,R)\ar[r]^{\simeq} & d_d(\tau^{l})(L,P,R) \\ }\]

We assume now $0\leq i\leq d-1$.
We let $u$ (respectively, $v$)
denotes the non-unital vertex of $\tau$
which lies on level $i$ (respectively, $i+1$) in $\tau^l$.
There are $2$ possibilities: the vertices $u$ and $v$ are either connected by an edge $e\in E(\tau)$
or not.
In the first case,
the contraction of levels $i$ and $i+1$ in $\tau^l$
is equivalent to the contraction of the subtree
$\sigma_e = \xymatrix@M=0pt@R=6mm@C=6mm{ *+<3mm>[o][F]{v}\ar[r]^{e} & *+<3mm>[o][F]{u} \\ }$
in $\tau$.
To be more precise,
the tree $d_i(\tau^l)$ coincides with the levelization of $\tau/\sigma_e$
determined by the map $l': V_p(\tau/\sigma_e)\,\rightarrow\,\{1,\ldots,d-1\}$
such that:
$l'(x) = l(x)$, if $x\in V_p(\tau)$ satisfies $l(x)<i$;
$l'(\sigma_e) = i$;
and $l'(x) = l(x)-i$, if $x\in V_p(\tau)$ satifies $l(x)>i+1$.
Moreover,
we have a commutative diagram
\[\xymatrix{\tau(L,P,R)\ar[r]^{\simeq}\ar[d]_{\beta^{\tau,\tau/\sigma_e}} & \tau^{l}(L,P,R)\ar[d]^{(-1)^i d_i} \\
\tau/\sigma_e(L,P,R)\ar[r]^{\simeq} & d_i(\tau^{l})(L,P,R) \\ }\]

In the second case,
we have a unique levelization $\tau^{l'}\not=\tau^{l}$
such that $d_i(\tau^{l}) = d_i(\tau^{l'})$.
As an example,
the levelizations of figure~\ref{LevelizedTree}
have the same $d_1$-face.
In general,
the tree $\tau^{l'}$ is determined by the map $l': V_p(\tau)\,\rightarrow\,\{1,\ldots,d\}$
which has $l'(u) = i+1$, $l'(v) = i$
and which agrees with $l(x)$ for vertices $x\not = u,v$.
We have a canonical isomorphism $\tau^{l'}(L,P,R)\,\simeq\,\tau^{l}(L,P,R)$,
since, after removing unit vertices,
the tensor products $\tau^{l'}(L,P,R)$ and $\tau^{l}(L,P,R)$
differ by a transposition $P(I_v)\otimes P(I_u)\,\simeq\,P(I_u)\otimes P(I_v)$.
On one hand,
this isomorphism $\tau^{l'}(L,P,R)\,\simeq\,\tau^{l}(L,P,R)$ fits in the commutative diagram
\[\xymatrix{ \tau^{l'}(L,P,R)\ar[r]^{d_i}\ar[d]^{\simeq} & d_i(\tau^{l'})(L,P,R)\ar[d]^{=} \\
\tau^{l}(L,P,R)\ar[r]^{d_i} & d_i(\tau^{l})(L,P,R) \\ }\]
On the other hand,
as in the example of paragraph~\ref{LevelizedTensors},
the isomorphism $\tau^{l'}(L,P,R)\,\simeq\,\tau^{l}(L,P,R)$
makes the levelizations
\[\tau(L,P,R)\,\xrightarrow{\simeq}\,\tau^{l'}(L,P,R)\qquad\text{and}\qquad\tau(L,P,R)\,\xrightarrow{\simeq}\,\tau^{l}(L,P,R)\]
differ by the sign $-1$,
because the levelization process involves permutations of suspension symbols.
We conclude that the $d_i$-face of $\tau^{l}(L,P,R)$ is cancelled by the $d_i$-face of $\tau^{l'}(L,P,R)$.
\end{proof}

\subsection{The example of the associative operad}\label{AssociativeLevelization}
We give a few indications about the example of the associative operad $\A$.
We know from \cite[E. Getzler and J. Jones]{GetzlerJones}
and \cite[V. Ginzburg and M. Kapranov]{GinzburgKapranov}
that the differential graded bar construction $\bar{B}(\A)$
is related to the cellular complex of Stasheff's \emph{associahedron}.
We recall that the associahedron $K_n$,
introduced by J. Stasheff in \cite{Stasheff},
is an $n-2$-dimensional polyhedron
whose $d$-dimensional faces $K_\tau\subset K_n$ are indexed by planar trees $\tau$
with $n-1-d$ vertices and $n$ entries.
We have a face inclusion $K_{\tau'}\subset K_{\tau}$
if and only if $\tau'$ is a quotient of $\tau$.
We let $C_*(K_n)$ denote the cell complex of $K_n$.

We have precisely
\[\bar{B}_d(\A)(n) = C_{n-1-d}(K_n)^{\vee}\otimes\K[\Sigma_n].\]
We recall that $\A(n)$ is the regular representation of $\Sigma_n$.
In the definition of $\bar{B}_d(\A)$,
the labeling of an abstract $n$-tree by elements of the associative operad $p_v\in \A(I_v)$
is equivalent to an embedding of $\tau$ in the plane.
In the definition of the associahedron,
the entries of a planar tree are indexed by $1,\ldots,n$
according to the orientation of the plane.
On the other hand,
in the expansion of $\bar{B}_d(\A)$,
we consider planar trees indexed by any permutation of $1,\ldots,n$.
Therefore,
in the relation above,
we insert a tensor product by the regular representation of $\Sigma_n$.

The simplicial bar construction $\bar{N}(\A)$
is related to the cellular complex of Milgram's \emph{permutohedron}.
We have precisely
\[\bar{N}_d(\A)(n) = C_{n-1-d}(P_{n-1})^{\vee}\otimes\K[\Sigma_n].\]
Classically,
the permutohedron $P_{n-1}$,
introduced by J. Milgram in \cite{Milgram},
is an $n-2$-dimensional polyhedron
whose faces are indexed by partitions of $\{1,\ldots,n-1\}$.
On observe that such partitions
are equivalent to planar trees $\tau$ with $n$ entries and level structures
(\emph{cf}. S. Saneblidze and R. Umble \cite{SaneblidzeUmble}).
Therefore,
as for associahedra,
we obtain a relation $\bar{N}_d(\A)(n) = C_{n-1-d}(P_{n-1})^{\vee}\otimes\K[\Sigma_n]$.
In fact,
it is classical that Milgram's permutohedron $P_{n-1}$
is Poincar\'e dual to a simplicial set,
namely to the coxeter complex of $\Sigma_{n-1}$
(\emph{cf}. M. Kapranov \cite{Kapranov}).
The relationship between the simplicial bar construction $\bar{N}_*(\A)$
and the normalized chain complex of this simplicial set
is deduced from the relationship between planar trees with levels
and partitions.

One has a cellular map $P_{n-1}\,\rightarrow\,K_n$ defined by using the relationship
between planar trees with levels and planar trees
(\emph{cf}. J.-L. Loday and M. Ronco \cite{LodayRonco},
S. Saneblidze and R. Umble \cite{SaneblidzeUmble},
A. Tonks \cite{Tonks}).
The levelization morphism $\bar{B}_*(\A)\,\rightarrow\,\bar{N}_*(A)$
is clearly the chain morphism induced by this cellular map.

\section{Proofs}\label{ProofsBarCoefficientsPpties}

We prove the results stated in section~\ref{BarCoefficientsPpties}.
We are given a connected dg-operad $P$.
We assume that $L$ is a right $P$-module and $R$ is a connected left $P$-module.
If the ground ring $\K$ is not a field,
then we assume in addition
that $P$ (respectively, $L$, $R$)
is projective as a $\K$-module.
We recall that the levelization morphism defines an embedding of chain complexes
$\phi(L,P,R): B_*(L,P,R)\,\hookrightarrow\,N_*(L,P,R)$
(\emph{cf}. lemma~\ref{LevelizationDifferential}).

\begin{lem}
The extra degeneracy of the simplicial bar construction
$s_{-1}: N_{d}(P,P,R)\,\rightarrow\,N_{d+1}(P,P,R)$
verifies
$s_{-1}(B_{d}(P,P,R))\subset B_{d+1}(P,P,R)$.
\end{lem}

\begin{proof}
By definition,
the extra degeneracy $s_{-1}: N_{d}(P,P,R)\,\rightarrow\,N_{d+1}(P,P,R)$
is obtained by inserting one unit vertex in the zero level of trees.
From this picture,
it should be clear that the extra degeneracy
preserves levelizations.
Formally,
if $\tau$ is a composite tree,
then we let $s_{-1}(\tau)$ denotes the composite tree
obtained by putting a unit vertex in the lower level of vertices $V_l(s_{-1}(\tau))$
and which has $V_p(s_{-1}(\tau)) = V_l(\tau)\amalg V_p(\tau)$
and $V_r(s_{-1}(\tau)) = V_r(\tau)$.
We have a canonical morphism $\tau(P,P,R)\,\rightarrow\,s_{-1}(\tau)(P,P,R)$.
Explicitly, we assume that the bottom vertex $u\in V_l(\tau)$ in the tree $\tau$
is labeled by the operation $p_u\in P(I_u)$.
Then, in the tree $s_{-1}(\tau)$,
we mark this vertex by the suspension $\Sigma p_u\in\Sigma\tilde{P}(I_u)$, if $p_u\in\tilde{P}(I_u)$,
and by $0$, if $p_u\in\K 1$.
The vertices $v\in V_p(\tau)$ have unchanged labels $\Sigma p_v\in\Sigma\tilde{P}(I_v)$
in $s_{-1}(\tau)$ and so do the vertices $v\in V_r(\tau)$.

Any levelization of $s_{-1}(\tau)$ keeps necessarily the vertex of $V_l(\tau)$ in level $1$
and, consequently, is equivalent to the $s_{-1}$-degeneracy of a levelization of $\tau$.
We conclude that the insertion process above corresponds to the $s_{-1}$-degeneracy in $N_*(L,P,R)$.
\end{proof}

The next statement is a direct consequence of the previous lemma:

\begin{lem}\label{LeftModuleBarResolution}
The augmentation $\epsilon(P,P,R): B(P,P,R)\,\rightarrow\,R$
is a quasi-isomorphism of left $P$-modules.
\end{lem}

In the case $R = P$,
we obtain the following result:

\begin{lem}\label{BiModuleBarResolution}
The augmentation $\epsilon(P,P,P): B(P,P,P)\,\rightarrow\,P$
is a quasi-isomorphism of $P$-bimodules.
\end{lem}

We deduce the next lemma from the result above and from theorem~\ref{QuasiFreeCompositionProduct}
about derived composition products.

\begin{lem}
The augmentation $\epsilon(L,P,P): B(L,P,P)\,\rightarrow\,L$
is a quasi-isomorphism of right $P$-modules.
\end{lem}

\begin{proof}
By lemma~\ref{BiModuleBarResolution},
we have quasi-isomorphism $\epsilon(P,P,P): B(P,P,P)\,\rightarrow\,P$.
The dg-module $B(P,P,P)$ is a quasi-free left $P$-module,
since we have by construction $B(P,P,P) = P\circ B(I,P,P)$.
The dg-operad $P$ is obviously a free left $P$-module.
We have explicitly $P = P\circ I$.
We deduce from theorem~\ref{QuasiFreeCompositionProduct} that the induced morphism
\[L\circ_P\epsilon(P,P,P): L\circ_P B(P,P,P)\,\rightarrow\,L\circ_P P\]
is a quasi-isomorphism.
This proves the lemma,
because this morphism is identified with the augmentation
$\epsilon(L,P,P): B(L,P,P)\,\rightarrow\,L$.
\end{proof}

We are now in position to achieve the proof of theorem~\ref{Levelization}:

\begin{lem}
The levelization morphism
\[\phi(L,P,R): B(L,P,R)\,\rightarrow\,N(L,P,R)\]
is a quasi-isomorphism.
\end{lem}

\begin{proof}
We deduce this result from theorem~\ref{QuasiFreeCompositionProduct}.
To be explicit,
we have a quasi-iso\-mor\-phism of right $P$-modules
$\phi(L,P,P): B(L,P,P)\,\rightarrow\,N(L,P,P)$,
because the chain complexes $B(L,P,P)$ and $N(L,P,P)$ are both quasi-isomorphic to $L$.
The right $P$-module $B(L,P,P)$ is quasi-free by definition.
Moreover,
we observe in paragraph ~\ref{NormalizedRightModule}
that the right $P$-module $N(L,P,P)$ is quasi-free.
We conclude from theorem~\ref{QuasiFreeCompositionProduct} that the induced morphism
\[\phi(L,P,P)\circ_P R: B(L,P,P)\circ_P R\,\rightarrow\,N(L,P,P)\circ_P R\]
is a quasi-isomorphism.
Namely,
the lemma follows from this assertion,
since we have the identity $\phi(L,P,R) = \phi(L,P,P)\circ_P R$.
\end{proof}

\section[Quasi-free resolutions of operads]{Quasi-free resolutions of operads and bar constructions}\label{ProofCobarBarResolution}

We give a proof of proposition~\ref{CobarBarResolution}.
More precisely,
one verifies easily
that the morphism of free operads
\[\epsilon: F(\Sigma^{-1}\tilde{F}^c(\Sigma\tilde{P}))\,\rightarrow\,P\]
introduced in proposition~\ref{CobarBarResolution}
defines a morphism of dg-operads $\epsilon: \bar{B}^c(\bar{B}(P))\,\rightarrow\,P$.
We prove that this morphism is a quasi-isomorphism.

\subsection{The cobar construction with coefficients}
We consider the cobar construction of a dg-cooperad $D$ with coefficients
in a right $D$-comodule $L$ and a left $D$-comodule $R$,
whose definition is dual to that of the bar construction of a dg-operad.
We obtain a complex $B^c(L,D,R)$ such that $B^c(L,D,R) = L\circ \bar{B}^c(D)\circ R$
and that satisfies the dual of the properties proved
in section~\ref{ProofsBarCoefficientsPpties}.
In particular, the chain complex $B^c(I,D,D)$ is acyclic.
The differential of $B^c(D,D,D)$ is also determined by a twisting cochain $\psi: D\,\rightarrow\,B^c(D)$
provided by the canonical inclusion
$\Sigma^{-1}\tilde{D}\,\hookrightarrow\,F(\Sigma^{-1}\tilde{D})$.

We consider the cobar construction $B^c(L,D,R)$
in the case $L = I$, $D = \bar{B}(P)$ and $R = \bar{B}(P)$.
We obtain then an acyclic complex $B^c(L,D,R) = B^c(I,\bar{B}(P),\bar{B}(P))$
such that $B^c(I,\bar{B}(P),\bar{B}(P)) = \bar{B}^c(\bar{B}(P))\circ\bar{B}(P)$.

\begin{lem}\label{TwistedCochainComparison}
The composition product
$\epsilon\circ\bar{B}(P): \bar{B}^c(\bar{B}(P))\circ\bar{B}(P)\,\rightarrow\,P\circ\bar{B}(P)$,
where $\epsilon: \bar{B}^c(\bar{B}(P))\,\rightarrow\,P$ is the morphism of dg-operads
provided by proposition~\ref{CobarBarResolution},
defines a morphism of chain complexes
from the cobar complex
$B^c(I,\bar{B}(P),\bar{B}(P)) = \bar{B}^c(\bar{B}(P))\circ\bar{B}(P)$
to the bar complex
$B(P,P,I) = P\circ\bar{B}(P)$.
\end{lem}

\begin{proof}
We deduce this lemma from the description of differentials by twisting cochains
(\emph{cf}. paragraph~\ref{TwistingCochains}).
We recall that the differential of $B(P,P,I) = P\circ\bar{B}(P)$
is determined by the internal differential of the reduced bar construction $\bar{B}(P)$
and by a homogeneous morphism
$\theta_L: \bar{B}(P)\,\rightarrow\,P\circ\bar{B}(P)$
induced by a certain twisting cochain $\phi: \bar{B}(P)\,\rightarrow\,P$.
Similarly, the differential of $B^c(I,\bar{B}(P),\bar{B}(P)) = \bar{B}^c(\bar{B}(P))\circ\bar{B}(P)$
is determined by the internal differentials of $\bar{B}^c(\bar{B}(P))$ and $\bar{B}(P)$
and by a homogeneous morphism
$\theta_R: \bar{B}^c(\bar{B}(P))\circ\bar{B}(P)\,\rightarrow\,\bar{B}^c(\bar{B}(P))$
induced by a twisting cochain $\psi: \bar{B}(P)\,\rightarrow\,\bar{B}^c(\bar{B}(P))$.

The morphism $\epsilon: \bar{B}^c(\bar{B}(P))\,\rightarrow\,P$ is compatible with operad differentials.
Therefore,
the associated map $\epsilon\circ\bar{B}(P): \bar{B}^c(\bar{B}(P))\circ\bar{B}(P)\,\rightarrow\,P\circ\bar{B}(P)$
preserves clearly internal differentials.
In addition,
the morphism $\epsilon: \bar{B}^c(\bar{B}(P))\,\rightarrow\,P$
is compatible with twisting cochains.
Explicitly,
the definitions imply immediately that the diagram
\[\xymatrix{ \bar{B}(P)\ar[r]^{\psi}\ar[dr]_{\phi} & \bar{B}^c(\bar{B}(P))\ar[d]^{\epsilon} \\ & P \\ }\]
is commutative.
We deduce readily from this property
that the map
$\epsilon\circ\bar{B}(P): \bar{B}^c(\bar{B}(P))\circ\bar{B}(P)\,\rightarrow\,P\circ\bar{B}(P)$
preserves the components of the differential determined by twisting cochains.
The conclusion follows.
\end{proof}

\begin{proof}[Proof of proposition~\ref{CobarBarResolution}]
We are now in position to deduce proposition~\ref{CobarBarResolution}
from our comparison theorems.
More properly,
the chain complexes $B^c(I,\bar{B}(P),\bar{B}(P)) = \bar{B}^c(\bar{B}(P))\circ\bar{B}(P)$
and $B(P,P,I) = P\circ\bar{B}(P)$
are quasi-cofree right $\bar{B}(P)$-como\-dules.
Therefore,
we consider the dual statement of theorem~\ref{RightModuleComparison}.
Since $\epsilon\circ\bar{B}(P): \bar{B}^c(\bar{B}(P))\circ\bar{B}(P)\,\rightarrow\,P\circ\bar{B}(P)$
is a morphism of acyclic quasi-cofree right $\bar{B}(P)$-comodules,
this theorem implies immediately
that $\epsilon: \bar{B}^c(\bar{B}(P))\,\rightarrow\,P$
is a quasi-isomorphism.
This achieves the proof of proposition~\ref{CobarBarResolution}.
\end{proof}

%\input{KoszulDuality}
% 30/1/2003

\chapter{Koszul duality for operads}\label{KoszulDuality}

We define the \emph{Koszul dual} of an operad in this chapter.
We follow closely the classical theory of Koszul algebras
set by S. Priddy in \cite{PriddyKoszul}
(see also A. Beilinson, V. Ginzburg and W. Soergel \cite{BeilinsonGinzburgSoergel}).
As mentioned in the prolog,
we generalize constructions of \cite[V. Ginzburg and M. Kapranov]{GinzburgKapranov}.
To be precise,
our results hold in positive characteristic
and are valid for a larger class of \emph{quadratic operads}
(see paragraphs~\ref{QuadraticOperads} and~\ref{BinaryOperadDuality}).

\section{Weight graded operads}

\subsection{Weight graded objects}\label{WeightGradings}
In this section,
we consider $\K$-modules $V$ equipped with a \emph{weight grading},
which is equivalent to a splitting
$V = \bigoplus_{s=0}^{\infty} V_{(s)}$.
In the case of a dg-module $V$,
we assume that the homogeneous components $V_{(s)}$ are sub-dg-modules of $V$.
A tensor product of weight graded modules is equipped with a canonical weight grading.
We have explicitly
$V\otimes W = \bigoplus_{n=0}^{\infty} (V\otimes W)_{(n)}$,
where
$(V\otimes W)_{(n)} = \bigoplus_{s+t=n} V_{(s)}\otimes W_{(t)}$.
Consequently,
the weight graded modules form a symmetric monoidal category.

In principle,
the sign $\pm$ which occurs in a symmetry isomorphism
$c(V,W)(v\otimes w) = \pm w\otimes v$, where $v\in (V_{(s)})_i$ and $w\in (W_{(t)})_j$,
is determined by differential degrees $|v| = i$ and $|w| = j$
(see paragraphs~\ref{DGTensor} and~\ref{GradedMod}).
Accordingly,
we assume that weights do not contribute to commutation signs.
In particular,
no sign occur in the symmetry isomorphism of weight graded $\K$-modules
which are not equipped with a differential graded structure.

\subsection{Weight graded operads}\label{WeightOperads}
The definitions of chapter~\ref{Operads}\linebreak make sense in the context of weight graded modules,
since weight graded modules form a symmetric monoidal category.
In particular,
a composite of weight graded $\Sigma_*$-modules $M\circ N$
is equipped with a canonical weight grading.
An operad $P$ equipped with a splitting $P = \bigoplus_{s=0}^{\infty} P_{(s)}$
is a weight graded operad
if the operad composition product
$P\circ P\,\rightarrow\,P$
preserves weight gradings.
Dually, a cooperad $D$ together with a splitting $D = \bigoplus_{s=0}^{\infty} D_{(s)}$
is a weight graded cooperad
if the coproduct
$D\,\rightarrow\,D\circ D$
preserves weight gradings.

More explicitly,
an operad $P$ is weight graded
provided that a partial composite of homogeneous elements
$p\in P_{(s)}(m)$ and $q\in P_{(t)}(n)$
satisfies
$p\circ_i q\in P_{(s+t)}(m+n-1)$.
Accordingly,
we observe that any operad $P$
is equipped with a canonical weight grading.
We set precisely
\[P_{(s)}(r) = \begin{cases} P(r), & \text{if $s = r-1$}, \\
0, & \text{otherwise}. \end{cases}\]
We note that,
according to this definition,
we have
$P_{(-1)}(r) = P(0)$, if $r = 0$, and $P_{(-1)}(r) = 0$, otherwise.
Consequently,
we have no component in negative weight
if and only if the operad satisfies $P(0) = 0$.

\subsection{On connected operads}\label{ConnectedWeightGradedOperad}
We observe that the operad unit $1\in P(1)$ of a weight graded operad
is necessarily homogeneous of weight $0$.
An operad $P$ is connected in regard to a weight grading
if we have $P_{(0)}(r) = \K\,1$, for $r = 1$, and $P_{(0)}(r) = 0$, for $r\not=1$.

Consider the canonical weight grading defined in paragraph~\ref{WeightOperads}.
We observe that the operad $P$ has no component in negative weight if and only if $P(0) = 0$.
We assume that this property holds for a connected operad.
We have in addition $P_{(0)}(r) = P(1)$, if $r = 1$, and $P_{(0)}(r) = 0$, otherwise.
Therefore,
an operad $P$ is connected in regard to the canonical weight grading
if and only if $P(0) = 0$ and $P(1) = \K\,1$.

The results of this chapter hold for all operads $P$
which are connected in regard to a weight grading
and
which have $P(0) = 0$,
because the theorems of sections~\ref{ModuleComplexes},~\ref{ReducedBar} and~\ref{BarCoefficients}
(from which we deduce our results)
can be generalized in this context.
For that purpose,
one has just to replace operadic gradings by weight gradings
in the spectral sequences of sections~\ref{ModuleSpectralSequences}
and~\ref{ProofsComparisonQuasiFreeOperads}.
In this article,
we assume that all connected operads verify $P(1) = \K\,1$ for simplicity
(so that the results of sections~\ref{ModuleComplexes},~\ref{ReducedBar} and~\ref{BarCoefficients} hold without changes),
but only the relation $P(0) = 0$ is crucial.

\subsection{Homogeneous ideals}
We consider homogeneous operad ideals $R\subset P$
which have a weight decomposition $R = \bigoplus_{s=0}^{\infty} R_{(s)}$
such that $R_{(s)} = R\cap P_{(s)}$.
As an example,
the augmentation ideal of a connected operad $\tilde{P}$ is homogeneous.
We have the following classical property:

\begin{fact}
The quotient of a weight graded operad by a homogeneous ideal
is equipped with a natural weight grading.
\end{fact}

We observe in section~\ref{FreeOperadStructure}
that a free operad is equipped with a natural weight decomposition
\[F(M) = \bigoplus_{r=0}^{\infty} F_{(r)}(M).\]
The next statement is a straightforward generalization of this property:

\begin{prp}\label{WeightedFreeOperad}
We assume that $M$ is a $\Sigma_*$-module equipped with a weight grading
$M = \bigoplus_{s=0}^{\infty} M_{(s)}$.
The homogeneous components of the free operad $F_{(r)}(M)$
are equipped with an internal splitting into weight components
\[F_{(r)}(M) = \bigoplus_{s=0}^{\infty} F_{(r)}(M)_{(s)}.\]
This internal weight grading is preserved by composition products
and is characterized by the relation $F_{(1)}(M)_{(s)} = M_{(s)}$.

If the $\Sigma_*$-module $M$ is equipped with the canonical weight grading
(explicitly, if we have $M_{(n-1)}(n) = M(n)$ and $M_{(s)}(n) = 0$ for $s\not=n-1$),
then so is the free operad.
Explicitly,
we have $F_{(r)}(M)_{(n-1)}(n) = F_{(r)}(M)(n)$
and $F_{(r)}(M)_{(s)}(n) = 0$ for $s\not=n-1$.

If the $\Sigma_*$-module $M$ is equipped with a trivial weight grading
(explicitly, if $M_{(1)}(n) = M(n)$ and $M_{(s)}(n) = 0$ for $s\not=1$),
then we have $F_{(r)}(M)_{(r)}(n) = F_{(r)}(M)(n)$
and $F_{(r)}(M)_{(s)}(n) = 0$ for $s\not=r$.
\end{prp}

\begin{proof}
We recall that the $\K$-module $F_{(r)}(M)(n)$ is generated by composite elements
\[(\cdots((x_1\circ_{i_2} x_2)\circ_{i_3}\cdots)\circ_{i_r} x_r\]
where $x_1\in M(n_1),\ldots,x_r\in M(n_r)$
(see paragraph~\ref{FreeOperad}).
If the elements $x_1,\ldots,x_r$
are all homogeneous of respective weights $s_1,\ldots,s_r$,
then the composite above is homogeneous of weight $s_1+\cdots+s_r$,
because an operad composition product is supposed to preserve weight gradings.
Therefore,
the free operad is equipped with an internal weight grading
which is characterized by the properties stated in the proposition.

Equivalently,
we consider the treewise expansion of the free operad
$F_{(r)}(M)(n) = \colim_{\tau\in T_r(n)} \tau(M)$
(see paragraph~\ref{WeightTreewiseTensors}).
We have by definition $\tau(M) = \bigotimes_{v\in V(\tau)} M(I_v)$.
Accordingly,
the $\K$-module $\tau(M)$ is equipped with the weight grading of a tensor product.
In addition,
an isomorphism of trees $f: \tau\,\rightarrow\,\tau'$
gives a homogeneous morphism of weight graded modules
$f_*: \tau(M)\,\rightarrow\,\tau'(M)$.
We conclude that the free operad $F_{(r)}(M)$ is equipped with a weight grading.
\end{proof}

The remainder assertions of the proposition are straighforward consequences of this construction.
Similarly, by analyzing the structure of the tensor product $\tau(M)$,
we obtain:

\begin{lem}\label{ConnectednessFreeOperad}
We assume that $M$ is a weight-graded $\Sigma_*$-module
which has $M_{(0)} = 0$.
In this situation,
the internal grading of the free operad $F(M)$
satisfies $F_{(r)}(M)_{(r)} = F_{(r)}(M_{(1)})$
and $F_{(r)}(M)_{(s)} = 0$ for $r>s$.
\end{lem}

\begin{proof}
We have a splitting
\[\tau(M) = \bigotimes_{v\in V(\tau)}\Bigl\{M_{(1)}(I_v)\oplus M_{(2)}(I_v)\oplus\cdots\Bigr\}
= \Bigl\{\bigotimes_{v\in V(\tau)} M_{(1)}(I_v)\Bigr\}\oplus\cdots.\]
If we assume $\tau\in T_{(r)}(n)$
(explicitly, if we assume that $\tau$ has $r$ vertices),
then the leading term of this expansion has weight $r$.
The conclusion follows.
\end{proof}

\section{Koszul operads}

We deduce the following proposition
from an analogue of theorem~\ref{WeightedFreeOperad}
in the context of a cofree cooperad $F^c(M)$:

\begin{lem}\label{WeightedReducedBar}
We assume that $P$ is a weight graded operad.
The reduced bar construction $\bar{B}(P)$ is equipped with a canonical weight grading
and forms a weight graded cooperad.
To be more precise,
each homogeneous component of the bar construction
is equipped with an internal splitting
$\bar{B}_d(P) = \bigoplus_{s=0}^{\infty}\bar{B}_d(P)_{(s)}$
such that
$\beta(\bar{B}_{d}(P))_{(s)}\subset\bar{B}_{d-1}(P)_{(s)}$.
\end{lem}

\begin{proof}
We have by definition $\bar{B}_d(P) = F^c_{(d)}(\Sigma\tilde{P})$.
Therefore, we set $\bar{B}_d(P)_{(s)} = F^c_{(d)}(\Sigma\tilde{P})_{(s)}$.
We observe that the homogeneous morphism
\[\theta: F_{(2)}(\Sigma\tilde{P})\,\rightarrow\,\Sigma\tilde{P},\]
which is determined by the partial composition products of $P$
(see lemma~\ref{BarCochain}),
is homogeneous with respect to weight gradings,
because so are the partial composition products.
We deduce readily that the associated coderivation $\beta = d_\theta$
(see paragraphs~\ref{TreeCoderivation} and~\ref{ReducedBarDifferential})
preserves homogeneous weight components.
\end{proof}

The next property is a direct corollary of lemma~\ref{ConnectednessFreeOperad}.

\begin{lem}\label{ConnectednessReducedBar}
If we assume that $P$ a connected weight graded operad,
then we have $\bar{B}_{s}(P)_{(s)} = F^c_{(s)}(\Sigma\tilde{P}_{(1)})$ and $\bar{B}_{d}(P)_{(s)} = 0$ for $d>s$.
\end{lem}

\subsection{The Koszul construction}\label{KoszulConstruction}
We are given a connected operad $P$ graded by weights (as in lemma~\ref{ConnectednessReducedBar}).
We consider the weight graded $\Sigma_*$-module $\bar{K}(P)$
such that
\[\bar{K}(P)_{(s)} = H_s(\bar{B}_{*}(P)_{(s)},\beta).\]
Since $\bar{B}_{d}(P)_{(s)} = 0$ for $d>s$,
we obtain equivalently
\[\bar{K}(P)_{(s)} = \ker\bigl(\beta: \bar{B}_{s}(P)_{(s)}\,\rightarrow\,\bar{B}_{s-1}(P)_{(s)}\bigr),\]
so that the $\Sigma_*$-module $\bar{K}(P)$ forms a subcomplex of the bar construction $\bar{B}(P)$.

We would like to mention that the $\K$-module $\bar{K}(P)_{(s)}$
has an internal differential graded structure
as long as $P$ is a differential graded operad.
In particular, the $\K$-module $\bar{K}(P)_{(s)}$
is equipped with an internal differential
$\delta: (\bar{K}(P)_{(s)})_*\,\rightarrow\,(\bar{K}(P)_{(s)})_{*-1}$
induced by the differential of $P$
and
deduced from the internal differential of the bar construction
$\delta: (\bar{B}_{s}(P)_{(s)})_*\,\rightarrow\,(\bar{B}_{s}(P)_{(s)})_{*-1}$.
In general, for simplicity,
we assume that the differential of $P$ is zero,
so that the module $\bar{K}(P)_{(s)}$ is equipped with an internal grading
but has no internal differential.
If we assume in addition that $P$ is a non graded operad
(if the modules $P(r)$ are concentrated in degree $0$),
then we observe that the modules $(\bar{B}_{d}(P)_{(s)})_*$ are concentrated in total degree $*=d$.
Consequently,
in this situation,
the modules $\bar{K}(P)_{(s)}$ are concentrated in degree $*=s$.

We have a dual construction for connected weight graded cooperads $D$.
In this context,
the cobar complex satisfies
$\bar{B}^c_{s}(D)_{(s)} = F_{(s)}(\Sigma^{-1}\tilde{D}_{(1)})$ and $\bar{B}^c_{d}(D)_{(s)} = 0$ for $d>s$.
We set
\[\bar{K}^c(D)_{(s)} = H^s(\bar{B}^c_{*}(D)_{(s)},\beta).\]
We have equivalently
\[\bar{K}^c(D)_{(s)} = \coker\bigl(\beta: \bar{B}^c_{s-1}(D)_{(s)}\,\rightarrow\,\bar{B}^c_{s}(D)_{(s)}\bigr),\]
so that the $\Sigma_*$-module $\bar{K}^c(D)$ forms a quotient complex of the cobar construction $\bar{B}^c(D)$.

\begin{lem}
Let $P$ be a connected weight graded operad $P$.
If the ground ring $\K$ is not a field, then we assume that $P$ is projective as a $\K$-module.
The $\Sigma_*$-module $\bar{K}(P)$, associated $P$,
is a subcooperad of the reduced bar construction $\bar{B}(P)$.

Dually, the $\Sigma_*$-module $\bar{K}^c(D)$,
associated to a connected weight graded cooperad $D$,
is a quotient operad of the reduced cobar construction $\bar{B}^c(D)$.
\end{lem}

\begin{proof}
We mention in paragraph~\ref{DGOperads} that the homology of a dg-operad
is equipped with the structure of a graded operad.
In the situation of the lemma,
we verify that the cokernel
$\bar{K}^c(D)_{(*)} = \coker\bigl(\beta: \bar{B}^c_{*-1}(D)_{(*)}\,\rightarrow\,\bar{B}^c_{*}(D)_{(*)}\bigr)$
has induced partial composition products
$\circ_i: \bar{K}^c(D)_{(s)}(m)\otimes\bar{K}^c(D)_{(t)}(n)\,\rightarrow\,\bar{K}^c(D)_{(s+t)}(m+n-1)$.
Explicitly,
we are given $p\in\bar{B}^c_{s}(D)_{(s)}(m)$ and $q\in\bar{B}^c_{t}(D)_{(t)}(n)$
such that $p = p' + \beta(p'')$ and $q = q' + \beta(q'')$.
Since $\beta(p) = \beta(q) = 0$,
we obtain $p\circ_i q = p\circ_i q' + \pm\beta(p\circ_i q'')$
and $p\circ_i q' = p'\circ_i q' + \beta(p'\circ_i q')$.
We conclude that the composites $p\circ_i q$ and $p'\circ_i q'$
represent the same element in $\bar{K}^c(D)_{(s+t)}$.

In the dual context,
we prove that the coproduct of the bar operad $\nu: \bar{B}(P)\,\rightarrow\,\bar{B}(P)\circ\bar{B}(P)$
preserves $\bar{K}(P)\subset\bar{B}(P)$.
According to the expansion of the composite $\Sigma_*$-module $\bar{B}(P)\circ\bar{B}(P)$,
we consider tensor products of the kernel
$\bar{K}(P)_{(s)} = \ker\bigl(\beta: \bar{B}_{s}(P)_{(s)}\,\rightarrow\,\bar{B}_{s-1}(P)_{(s)}\bigr)$
with $\K$-modules $\bar{B}_{d}(P)_{(s)}$.
Therefore,
the assertion about the $\Sigma_*$-module $\bar{K}(P)$
holds provided that the bar operad $\bar{B}(P)$ is projective as a $\K$-module.
But, if the operad $P$ is connected,
then the bar operad $\bar{B}(P)$ has a direct sum expansion
which consists of tensor products of $\K$-modules $P(n)$.
Consequently, if $P$ is projective as a $\K$-module,
then so is $\bar{B}(P)$.
\end{proof}

One observes that cooperads $\bar{K}(P)$ and operads $\bar{K}^c(D)$ can be described by generators and relations.
For this purpose,
we introduce the following notion
which generalizes a definition of V. Ginzburg and M. Kapranov
(\emph{cf}. \cite{GinzburgKapranov}).

\subsection{Quadratic operads}\label{QuadraticOperads}
We define a \emph{quadratic operad} $P$ to be a quotient $P = F(M)/(R)$,
where $(R)$ is an ideal of the free operad $F(M)$ generated by a $\Sigma_*$-module $R$
such that $R\subset F_{(2)}(M)$.
We have a dual notion of a quadratic cooperad.
According to this definition,
the ideal $(R)$ is homogeneous in the free operad $F(M)$,
provided $F(M)$ is equipped with its natural weight grading
$F(M) = \bigoplus_{r=0}^{\infty} F_{(r)}(M)$.
Consequently,
a quadratic operad $P$ has a natural weight grading
defined by the relation
$P_{(r)} = F_{(r)}(M)/F_{(r)}(M)\cap(R)$.
We note that $P_{(0)} = I$, $P_{(1)} = M$ and $P_{(2)} = F_{(2)}(M)/R$.
The relation $P_{(0)} = I$
implies that a quadratic operad is connected in regard to its weight grading.
In general, we assume that the $\Sigma_*$-module $M$ satisfies $M(0) = 0$,
so that the associated quadratic operad verifies also $P(0) = 0$ and $P(1) = \K\,1$.

We observe that the structure of a quadratic operad $P$
is determined by the weight grading
$P = \bigoplus_{r=0}^\infty P_{(r)}$
and
by the homogeneous components $P_{(1)}\subset P$ and $P_{(2)}\subset P$.
More precisely,
the generating module $M$ is determined by the relation $P_{(1)} = M$.
Moreover,
the quotient map
$F(M)\,\rightarrow\,P$
is necessarily induced by the morphism of $\Sigma_*$-modules
$M\,\xrightarrow{=}\,P_{(1)}\,\hookrightarrow\,P$.
Then,
the module of relations $R\subset F_{(2)}(M)$ of a quadratic operad $P$
is determined by the relation $P_{(2)} = F_{(2)}(M)/R$,
because this relation gives $R = \ker(F_{(2)}(M)\,\rightarrow\,P_{(2)})$.

\begin{lem}
We let $P$ be a connected weight graded operad.
The cooperad $\bar{K}(P)$ is quadratic and is determined by the relations
\begin{align*}
\bar{K}(P)_{(1)} & = \Sigma\tilde{P}_{(1)} \\
\text{and}\qquad\bar{K}(P)_{(2)}
& = \ker\bigl(\beta: F^c_{(2)}(\Sigma\tilde{P}_{(1)})\,\rightarrow\,\Sigma\tilde{P}_{(2)}\bigr).
\end{align*}
We have in particular $\bar{K}(P)_{(0)} = I$, so that $\bar{K}(P)$ is clearly connected.

Dually, we let $D$ be a connected weight graded cooperad.
The operad $\bar{K}^c(D)$ is quad\-ra\-tic and is determined by the relations
\begin{align*}
\bar{K}^c(D)_{(1)} & = \Sigma^{-1}\tilde{D}_{(1)} \\
\text{and}\qquad\bar{K}^c(D)_{(2)}
& = \coker\bigl(\beta: \Sigma^{-1}\tilde{D}_{(2)}\,\rightarrow\,F^c_{(2)}(\Sigma^{-1}\tilde{D}_{(1)})\bigr).
\end{align*}
We have in particular $\bar{K}^c(D)_{(0)} = I$, so that $\bar{K}^c(D)$ is connected.
\end{lem}

\begin{proof}
We generalize some arguments of Ginzburg-Kapranov (\emph{cf. loc. cit}.).
For instance, in the case of a cooperad,
we have clearly $\bar{K}^c(D)_{(0)} = I$ and $\bar{K}^c(D)_{(1)} = \Sigma^{-1}\tilde{D}_{(1)}$.
Moreover, by lemma~\ref{ConnectednessReducedBar},
the $\Sigma_*$-module $\bar{K}^c(D)_{(r)}$ is a quotient
of $\bar{K}^c_r(D)_{(r)} = F_{(r)}(\Sigma^{-1}\tilde{D}_{(1)})$.
Then,
it is straightforward to observe that the image of the cobar differential
\[\beta(F_{(r-1)}(\Sigma^{-1}\tilde{D})_{(r)})\subset F_{(r)}(\Sigma^{-1}\tilde{D})_{(r)}\]
defines the operad ideal generated
by
$R = \beta(\Sigma^{-1}\tilde{D}_{(2)})\subset F_{(2)}(\Sigma^{-1}\tilde{D}_{(1)})$,
because the $\Sigma_*$-module $F_{(r-1)}(\Sigma^{-1}\tilde{D})_{(r)}$ consists of treewise tensors
\[\bigotimes_{v\in V(\tau)}\Sigma^{-1} q_v\in\bigotimes_{v\in V(\tau)}\Sigma^{-1}\tilde{D}_{(*)}(I_v)\]
in which one factor
$\Sigma^{-1} q_{v_0}\in\Sigma^{-1}\tilde{D}_{(s)}(I_{v_0})$
has weight $s = 2$
and all other factors
$\Sigma^{-1} q_{v}\in\Sigma^{-1}\tilde{D}_{(s)}(I_{v})$, $v\not=v_0$, have weight $s = 1$.
By construction,
the cobar differential maps these tensors to elements
\[\beta(\Sigma^{-1} q_{v_0})\otimes\Bigl\{\bigotimes_{v\not=v_0}\Sigma^{-1} q_v\Bigr\}
\in F_{(r)}(\Sigma^{-1}\tilde{D}_{(1)})\]
and gives all composites of $F(\Sigma^{-1}\tilde{D}_{(1)})$
which have a factor in $R = \beta(\Sigma^{-1}\tilde{D}_{(2)})$.
The conclusion follows.
\end{proof}

\subsection{The quadratic duality of operads according to V. Ginzburg and M. Kapranov}\label{BinaryOperadDuality}
In this paragraph,
we give the relationship between the Koszul construction $P\mapsto\bar{K}(P)$
and the quadratic duality operation $P\mapsto P^!$
defined by V. Ginzburg and M. Kapranov in \cite{GinzburgKapranov}
for certain quadratic operads.
As alluded to in the introduction,
the Koszul construction $P\mapsto\bar{K}(P)$
is more general than the quadratic duality operation $P\mapsto P^!$.
To be precise,
one considers connected operads $P$
together with the canonical weight grading of operads,
introduced in paragraph~\ref{WeightOperads}:
\[P_{(s)}(r) = \begin{cases} P(r), & \text{if $s = r-1$}, \\
0, & \text{otherwise}. \end{cases}\]
In this context,
the construction of V. Ginzburg and  M. Kapranov gives an operad $P^!$
which can be deduced from the cooperad $\bar{K}(P)$
by the relation:
\[\bar{K}(P)(r) = \Sigma^{r-1} P^!(r)^{\vee}\otimes\sgn_r.\]
(To be more precise, this relation holds provided that the ground ring $\K$ is a field
and the operad $P$ consists of finitely generated $\K$-modules.)
Equivalently,
the Koszul construction $\bar{K}(P)$ coincides with the cooperad denoted by $P^{\bot}$
in \cite[E. Getzler and J. Jones]{GetzlerJones}.

Let us make explicit the presentation by generators and relations of the operad $\bar{K}(P)^{\vee}$
in order to explain this relationship.
First,
if the operad $P$ is equipped with the canonical weight grading of operads,
then so is the associated cooperad $\bar{K}(P)$.
In particular,
we obtain:
\begin{align*}
\bar{K}_{(1)}(P)(r) & = \left\{\begin{aligned} & \Sigma\tilde{P}(2),\ \text{if $r=2$}, \\
& 0,\ \text{otherwise}, \end{aligned}\right. \\
\text{and}\qquad\bar{K}_{(2)}(P)(r)
& = \left\{\begin{aligned} & \ker\bigl(\beta: F^c_{(2)}(\Sigma\tilde{P})(3)\,\rightarrow\,\Sigma\tilde{P}(3)\bigr),\ \text{if}\ r=3, \\
& 0,\ \text{otherwise}. \end{aligned}\right. \end{align*}
Consequently,
we have $\bar{K}(P)^{\vee} = F(M)/(R^\bot)$
where $M(2) = \bar{K}(P)(2)^{\vee}\linebreak = \Sigma^{-1}\tilde{P}(2)^{\vee}$ and $M(r) = 0$ for $r\not=2$.
The $\K$-module $R^\bot\subset F_{(2)}(M)$, which satisfies also $R^\bot(r) = 0$ for $r\not=3$,
is de termined by the short exact sequence
\[0\,\rightarrow\,R^\bot(3)\,\rightarrow\,F_{(2)}(M)(3)\rightarrow\,R(3)^{\vee}\,\rightarrow\,0,\]
where $R(3) = \ker\bigl(\beta: F^c_{(2)}(\Sigma\tilde{P})(3)\,\rightarrow\,\Sigma\tilde{P}(3)\bigr)$.

We refer to the article of V. Ginzburg and M. Kapranov
for calculations in the classical cases $P = \C,\A,\L$.
We recall these results in paragraph~\ref{ClassicalKoszulDuals}.

\subsection{Koszul operads}\label{DefinitionKoszulOperad}
We say that a connected and weight graded operad $P$
is \emph{Koszul}
if the inclusion morphism $\bar{K}(P)\,\hookrightarrow\,\bar{B}(P)$
is a quasi-isomorphism.
In this situation,
the cooperad $\bar{K}(P)$ is called the \emph{Koszul dual} of $P$.
If the ground ring $\K$ is not a field,
then a Koszul operad $P$ is also assumed to be projective as a $\K$-module
as well as the dual cooperad $\bar{K}(P)$.

Dually,
a connected weight graded cooperad $D$
is Koszul
if the quotient morphism $\bar{B}^c(D)\,\rightarrow\,\bar{K}^c(D)$
is a quasi-isomorphism.
If the ground ring $\K$ is not a field,
then a Koszul cooperad $D$ is also assumed to be projective as a $\K$-module
as well as the dual operad $\bar{K}^c(D)$.

\begin{lem}\label{LinearDualKoszulOperads}
We assume that $D$ is a connected weight graded cooperad
equipped with a trivial differential $\delta = 0$
and such that the sequence $D(r)$
consists of finitely generated projective $\K$-modules.
If $D$ is Koszul, then so is the dual operad $D^{\vee}$.
Moreover, we have $\bar{K}(D^{\vee}) = \bar{K}^c(D)^{\vee}$.

We assume that $P$ is a connected weight graded operad
equipped with a trivial differential $\delta = 0$
and such that the sequence $P(r)$
consists of finitely generated projective $\K$-modules.
If $P$ is Koszul, then so is the dual cooperad $P^{\vee}$.
Moreover, we have $\bar{K}^c(P^{\vee}) = \bar{K}(P)^{\vee}$.
\end{lem}

\begin{proof}
We observe that $\bar{B}^c(D)^{\vee} = \bar{B}(D^{\vee})$ (\emph{cf}. proposition~\ref{LinearDualBarCooperad}).
Consequently, by dualization, the short exact sequence
\[\bar{B}^c_{-s+1}(D)_{(s)}\,\xrightarrow{\beta}\,\bar{B}^c_{-s}(D)_{(s)}\,\xrightarrow{}\,\bar{K}^c(D)_{(s)}\,\xrightarrow{}\,0\]
gives
\[0\,\xrightarrow{}\,\bar{K}^c(D)_{(s)}^{\vee}\,\xrightarrow{}\,\bar{B}_{s}(D^{\vee})_{(s)}
\,\xrightarrow{\beta}\,\bar{B}_{s-1}(D^{\vee})_{(s)}.\]
Hence,
we obtain immediately $\bar{K}(D^{\vee}) = \bar{K}^c(D)^{\vee}$.
If $D$ is a Koszul cooperad,
then, for a fixed weight $s$, the quotient morphism
\[\bar{B}^c(D)_{(s)}(n)\,\rightarrow\,\bar{K}^c(D)_{(s)}(n)\]
is a quasi-isomorphism of finitely generated complexes of projective $\K$-modules.
We deduce from standard results of homological algebra
that the $\K$-dual morphism has also this property.
Therefore, the operad $D^{\vee}$ is Koszul.

Similarly, in the dual case of a Koszul operad $P$,
we obtain a quasi-isomorphism
\[\bar{B}(P)_{(s)}^{\vee}\,\xrightarrow{\sim}\,\bar{K}(P)_{(s)}^{\vee}.\]
We have also $\bar{B}(P)^{\vee} = \bar{B}^c(P^{\vee})$ (\emph{cf}. proposition~\ref{LinearDualBarCooperad}).
Consequently,
the cobar complex $\bar{B}^c_d(P^{\vee})_{(s)}$
satisfies
$H_s(\bar{B}^c_*(P^{\vee})_{(s)},\beta) = \bar{K}(P)_{(s)}^{\vee}$
and has
$H_d(\bar{B}^c_*(P^{\vee})_{(s)},\beta) = 0$
in degree $d\not=s$.
We conclude that $P^{\vee}$ is a Koszul cooperad
and the dual operad verifies $\bar{K}^c(P^{\vee}) = \bar{K}(P)^{\vee}$.
\end{proof}

\begin{lem}\label{KoszulDualityEquivalence}
We assume that $P$ is a connected weight graded operad
equipped with a trivial differential $\delta = 0$.
If $P$ is Koszul, then so is the cooperad $\bar{K}(P)$ and, moreover, we have $P = \bar{K}^c(\bar{K}(P))$.
(As a by-product, we deduce from this relation that a Koszul operad $P$ is necessarily quadratic.)

Dually, we assume that $D$ is a connected weight graded cooperad
equipped with a trivial differential $\delta = 0$.
If $D$ is Koszul, then so is the operad $\bar{K}^c(D)$ and, moreover, we have $D = \bar{K}(\bar{K}^c(D))$.
(As a by-product, we deduce from this relation that a Koszul cooperad $D$ is necessarily quadratic.)
\end{lem}

\begin{proof}
We deduce these assertions from lemma~\ref{KoszulResolution} below.
Explicitly, if $P$ is Koszul,
then we have a quasi-isomorphism
$\bar{B}^c(\bar{K}(P))\,\xrightarrow{\sim}\,P$.
Consequently, if $P$ is equipped with a trivial internal differential $\delta=0$,
then we obtain
\[H_{-d}(\bar{B}^c_{*}(\bar{K}(P))_{(s)},\beta)
= \begin{cases} P_{(s)}, & \text{if $d=s$}, \\
0, & \text{otherwise}. \end{cases}\]
The conclusion follows immediately.
The arguments are similar in the dual case of a Koszul cooperad $D$.
\end{proof}

\subsection{Koszul resolutions}\label{KoszulResolutionDefinition}
We have morphisms of dg-operads
\[\bar{B}^c(\bar{K}(P))\,\xrightarrow{}\,\bar{B}^c(\bar{B}(P))\,\xrightarrow{\sim}\,P.\]
The first morphism is induced by the inclusion of cooperads $\bar{K}(P)\,\hookrightarrow\,\bar{B}(P)$.
By the comparison theorem of quasi-free operads (\emph{cf}. theorem~\ref{ComparisonQuasiFreeOperads}),
this morphism is a quasi-isomorphism if and only if $P$ is Koszul.
The next morphism is provided by proposition~\ref{CobarBarResolution}
and is a quasi-isomorphism for all operads (which are projective as $\K$-modules).
As a conclusion,
we can record the following statement:

\begin{prp}\label{KoszulResolution}
We assume $P$ is a connected weight graded operad.
If the ground ring $\K$ is not a field, then $P$ is also supposed to be projective as a $\K$-module.
The canonical morphism
\[\bar{B}^c(\bar{K}(P))\,\xrightarrow{}\,P\]
is a quasi-isomorphism
if and only if $P$ is Koszul.
\end{prp}

In addition,
we have the following result,
which is a corollary of proposition~\ref{OperadCofibrantResolutions}:

\begin{prp}\label{KoszulCofibrantResolutions}
We assume $P$ is a Koszul operad.
We let $\tilde{K}(P)$ denote the coaugmentation coideal of the Koszul construction of $P$.
If $Q = F(M)$ is a cofibrant quasi-free resolution of $P$
such that $M(0) = M(1) = 0$,
then we have a quasi-isomorphism of dg-$\Sigma_*$-modules
$M\,\xrightarrow{\sim}\,\Sigma^{-1}\tilde{K}(P)$.
Conversely,
if we are given a resolution of $\Sigma^{-1}\tilde{K}(P)$ by projective $\Sigma_*$-modules,
let $M\,\xrightarrow{\sim}\,\Sigma^{-1}\tilde{K}(P)$,
then we have a cofibrant quasi-free resolution of $P$ such that $Q = F(M)$.
\end{prp}

\section[Koszul complexes]{Koszul complexes and characterization of Koszul operads}

\subsection{The Koszul construction with coefficients}\label{CoefficientKoszulConstruction}
We give a characterization of Koszul operads in theorem~\ref{KoszulCharacterization}.
For this purpose,
we introduce Koszul complexes with coefficients $K(L,P,R)$.
We set $K(L,P,R) = L\circ\bar{K}(P)\circ R$.
We have also $K(L,P,R) = L\circ_P K(P,P,P)\circ_P R$.
As in the definition of the differential of $B(L,P,R)$
by twisting cochains (\emph{cf}. paragraph~\ref{TwistingCochains}),
the differential of $K(P,P,P)$ is induced by morphisms
\[\theta_R: \bar{K}(P)\,\rightarrow\,\bar{K}(P)\circ P
\qquad\text{and}
\qquad\theta_L: \bar{K}(P)\,\rightarrow\,P\circ\bar{K}(P)\]
associated to a certain twisting cochain
$\phi': \bar{K}(P)\,\rightarrow\,P$.
To be precise,
we consider the composite of the inclusion morphism $\bar{K}(P)\,\hookrightarrow\,\bar{B}(P)$
with the twisting cochain of the bar construction $\phi: \bar{B}(P)\,\rightarrow\,P$.

Hence,
the maps $\theta_R$ and $\theta_L$
are determined by the process of paragraph~\ref{TwistingCochains}.
We consider the derivations of quasi-free modules
\begin{align*}
& d_{\theta_R}: K(I,P,P)\,\rightarrow\,K(I,P,P) \\
\text{and}\qquad & d_{\theta_L}: K(P,P,I)\,\rightarrow\,K(I,P,P),
\end{align*}
whose restrictions to $\bar{K}(P)$
agree with $\theta_R$ and $\theta_L$ respectively
(\emph{cf}. paragraphs~\ref{QuasiFreeRightModule} and~\ref{QuasiFreeLeftModule}).
The differential of $K(I,P,P)$ (respectively, $K(P,P,I)$)
reduce to the derivation $d_{\theta_R}$ (respectively, $d_{\theta_L}$).
Since $K(P,P,P) = P\circ K(I,P,P)$ and $K(P,P,P) = K(P,P,I)\circ P$,
we have also induced derivations of $P$-bimodules
\begin{align*}
& d_{\theta_R}: K(P,P,P)\,\rightarrow\,K(P,P,P) \\
\text{and}\qquad & d_{\theta_L}: K(P,P,P)\,\rightarrow\,K(P,P,P).
\end{align*}
The differential of $K(P,P,P)$ is the sum of these derivations.

Finally,
according to definitions,
the inclusion morphism $\bar{K}(P)\,\hookrightarrow\,\bar{B}(P)$
is a morphism of dg-cooperads and preserves the twisting co\-chains
which determine the differentials of $K(L,P,R)$ and $B(L,P,R)$.
Accordingly:

\begin{lem}\label{KoszulComplex}
We are given a connected weight graded operad $P$,
a right $P$-module $L$ and a connected left $P$-module $R$.
If the ground ring $\K$ is not a field,
then we assume in addition that $P$, $L$ and $R$ are projective as $\K$-modules.
The inclusion morphism $L\circ\bar{K}(P)\circ R\,\hookrightarrow\,L\circ\bar{B}(P)\circ R$
makes the Koszul construction with coefficients
$K(L,P,R) = L\circ\bar{K}(P)\circ R$,
a subcomplex of the bar construction with coefficients
$B(L,P,R) = L\circ\bar{B}(P)\circ R$.
\end{lem}

The next theorem gives the main result of the theory of Koszul operads.
In our framework,
this theorem is a straighforward corollary of the comparison theorems
of chapter~\ref{ModuleComplexes}.

\begin{thm}\label{KoszulCharacterization}
We let $P$ be a connected weight graded operad.
If the ground ring $\K$ is not a field,
then we assume that the operad $P$ is projective as a $\K$-module,
so that the reduced bar construction $\bar{B}(P)$ is projective as a $\K$-module.
We assume also that the Koszul construction $\bar{K}(P)$ is projective as a $\K$-module.

The following assertions are equivalent:
\textbf{(a)} the inclusion morphism $\bar{K}(P)\,\hookrightarrow\,\bar{B}(P)$
is a quasi-isomorphism;
\textbf{(b)} the chain complex $K(I,P,P)$ is acyclic;
\textbf{(c)} the chain complex $K(P,P,I)$ is acyclic.
If the operad $P$ has a trivial internal differential $\delta = 0$,
then property \textbf{(a)} holds if and only if $H_d(\bar{B}_{*}(P)_{(s)},\beta) = 0$
for $d\not=s$.

Moreover, if these properties are satisfied,
then the inclusion morphism $K(L,P,R)\,\hookrightarrow\,B(L,P,R)$
is a quasi-isomorphism for all coefficients $L$ and $R$,
such that
$L$ is a right $P$-module which is projective as a $\K$-module,
and
$R$ is a connected left $P$-module $R$ which is projective as a $\K$-module.
\end{thm}

\begin{proof}
The equivalence between properties \textbf{(a)}, \textbf{(b)} and \textbf{(c)}
is a corollary of the comparison theorems~\ref{RightModuleComparison} and~\ref{LeftModuleComparison}.
To be precise,
according to these results,
the inclusion morphism
$\bar{K}(P)\,\hookrightarrow\,\bar{B}(P)$
is a quasi-isomorphism
if and only if so are the induced inclusion morphisms
$K(I,P,P)\,\hookrightarrow\,B(I,P,P)$
and
$K(P,P,I)\,\hookrightarrow\,B(P,P,I)$.
But,
we know from lemma~\ref{BarResolution}
that the chain complexes $B(I,P,P)$ and $B(P,P,I)$ are acyclic.
Therefore,
properties \textbf{(a)}, \textbf{(b)} and \textbf{(c)}
are equivalent.

Similarly,
the last assertion of theorem~\ref{KoszulCharacterization} is a direct consequence
of theorem~\ref{QuasiFreeCompositionProduct},
because we have $K(L,P,R) = L\circ_P K(P,P,P)\circ_P R$.
\end{proof}

We have a dual result for the cobar construction of a connected weight graded cooperad.

\subsection{Koszul homology of algebras over an operad}\label{KoszulHomology}
We recall that an algebra over an operad $P$
is equivalent to a left $P$-module $R$
which has $R(r) = 0$ for $r\not=0$ (see paragraph~\ref{LeftModules}).
We consider the chain complexes
$K_*(I,P,R)$
and
$N_*(I,P,R)$
associated to such left $P$-modules $R$.
We have clearly $K_*(I,P,R)(r) = N_*(I,P,R)(r) = 0$ for $r\not=0$.
Finally,
we have chain complexes of $\K$-modules $K_*(I,P,A)$ and $N_*(I,P,A)$
associated to all $P$-algebras $R(0) = A$.
Furthermore,
our constructions (see theorem~\ref{Levelization} and lemma~\ref{KoszulComplex})
supply a natural morphism of chain complexes
\[K_*(I,P,A)\,\rightarrow\,N_*(I,P,A).\]
In positive characteristic,
the Koszul property does not imply that this comparison morphism is a quasi-isomorphism,
because the left $P$-module equivalent to a $P$-algebra
is not a connected $\Sigma_*$-module.
We can observe that $N_*(I,P,A)$ is nothing but the \emph{cotriple complex}
considered by J. Beck (\emph{cf}. \cite{Beck})
and by P. May (\emph{cf}. \cite{MayLoop}).
We define the \emph{Koszul homology} of a $P$-algebra
to be the homology of the chain complex $K_*(I,P,A)$.

We can assume that $A$ is a free $P$-algebra $A = S(P,V)$.
In this case,
we obtain readily
\begin{align*}
& K_*(I,P,S(P,V)) = S(K_*(I,P,P),V) \\
\text{and}\qquad & N_*(I,P,S(P,V)) = S(N_*(I,P,P),V).
\end{align*}
The cotriple homology of a free $P$-algebra vanishes
for general reasons.
But, on the other hand,
the Koszul property does not imply that the Koszul complex $K_*(I,P,S(P,V)) = S(K_*(I,P,P),V)$ is acyclic,
although the coefficients $K_*(I,P,P)$
form an acyclic complex of $\Sigma_*$-modules.
(This result holds if $K_*(I,P,P)$ is a chain complex of projective $\Sigma_*$-modules, but not in general.)
We have only the following implication:

\begin{prp}\label{VanishingKoszulHomology}
If the Koszul homology of a free $P$-algebra $A = S(P,V)$ vanishes,
then the Koszul complex $K_*(I,P,P)$ is acyclic.
\end{prp}

\begin{proof}
We observe that the functor $C: \Func\,\rightarrow\,\LMod{\Sigma_*}$,
left adjoint to $S: \LMod{\Sigma_*}\,\rightarrow\,\Func$ (see paragraph~\ref{CrossEffectModule}),
is exact (dislike $S: \LMod{\Sigma_*}\,\rightarrow\,\Func$).
This property holds because the $r$th homogeneous cross effect of a functor $F\in\Func$
is identified with a direct summand of the module $F(\K\oplus\cdots\oplus\K)$
(\emph{cf}. A.K. Bousfield \cite{BousfieldFunc}, E. Friedlander and A. Suslin \cite{FriedlanderSuslin}).
Since $K_*(I,P,S(P,V)) = S(K_*(I,P,P),V)$,
the chain complex of $\Sigma_*$-modules $K_*(I,P,P)$
is associated to the functor $V\mapsto K_*(I,P,S(P,V))$
(\emph{cf}. proposition~\ref{MonoidalCoef}).
The conclusion follows.
\end{proof}

\setcounter{subsection}{0}
\renewcommand{\thesubsection}{\thechapter.\arabic{subsection}}
%\input{Epilog}
% 30/1/2003

\chapter{Epilog: partition posets}\label{CommutativeDuality}

In this chapter,
we go back to partition posets and we prove the theorem announced in the prolog.
For that purpose,
we consider the example of the commutative operad $P = \C$.
On one hand,
we identify the simplicial bar construction of the commutative operad
with the normalized chain complex of the partition poset:

\begin{obv}\label{ClassicalKoszulConstructions}
The simplicial bar construction of the commutative operad $N_*(\C)(r)$
is isomorphic to the normalized chain complex $N_*(\bar{K}(r))$
of the simplicial set $\bar{K}(r)$.
\end{obv}

\begin{proof}
We recall that the $\K$-module $\C(r)$
is $1$ dimensional.
Consequently,
the tensor product $\tau(\C)$
associated to a tree is a $1$ dimensional $\K$-module.
Equivalently,
we can omit labels of vertices
in the definition of the bar constructions $\bar{B}_*(\C)$ and $\bar{N}_*(\C)$.
To be explicit,
the $\K$-module $\bar{B}_d(\C)$
is spanned by abstract trees $\tau$ with $d$ vertices
and
the $\K$-module $\bar{N}_d(\C)$
is spanned by abstract trees $\tau$ with $d$ levels.

We associate a sequence of partitions $\lambda_0\leq\cdots\leq\lambda_n$
to each abstract tree with $n$ levels $\tau$.
By assumption,
the last partition $\lambda_n$
is given by the collection
$\{1\},\ldots,\{r\}$.
If $i<n$,
then we assume that the components of the partition $\lambda_{i}$, denoted by $\lambda_{i}^{v}$,
are indexed by the vertices $v\in V_{i+1}(\tau)$
(which lie in level $i+1$).
We obtain the components of $\lambda_{i-1}$
by grouping the components of $\lambda_{i}$
according to the structure of the tree $\tau$.
Explicitly,
for $u\in V_i(\tau)$,
we set $\lambda_{i-1}^{u} = \coprod_{v\in I_u} \lambda_{i}^{v}$.
Since the first level of $\tau$ is reduced to a single vertex,
the partition $\lambda_0$ has a single component
and is the initial object of the partition poset.
Finally,
our construction gives a map $N_*(\C)(r)\,\rightarrow\,N_*(\bar{K}(r))$
which defines clearly an isomorphism of chain complexes.
\end{proof}

On the other hand,
we recall that the Koszul dual of the commutative operad
can be related to the Lie operad.
More precisely,
we have the following result, proved in \cite[V. Ginzburg and M. Kapranov]{GinzburgKapranov}:

\begin{fact}\label{ClassicalKoszulDuals}
The Koszul construction of the commutative operad $\bar{K}(\C)$
has the module
$\bar{K}(\C)(r)_* = \L(r)^{\vee}\otimes\sgn_r$
in degree $* = r-1$
and vanishes in degree $* \not= r-1$.
Dually, the Koszul construction of the Lie operad $\bar{K}(\L)$
has the module
$\bar{K}(\L)(r)_* = \C(r)^{\vee}\otimes\sgn_r$
in degree $* = r-1$
and vanishes in degree $* \not= r-1$.
\end{fact}

Let us mention that the Lie operad $\L$
consists of free $\Z$-modules (\emph{cf}. C. Reutenauer \cite{Reutenauer}).

\smallskip
Let us make explicit the presentation of the operads $\bar{K}(\C)^{\vee}$ and $\bar{K}(\L)^{\vee}$
which are dual to $\bar{K}(\C)$ and $\bar{K}(\L)$.
For that purpose,
we go back to the formulas of paragraph~\ref{BinaryOperadDuality}.
Since $\C(2)$ is the trivial representation of $\Sigma_2$,
the $\K$-module $\bar{K}(\C)(2)^{\vee}$
is generated by an operation $l' = l'(x_1,x_2)$ of degree $-1$
which verifies $l'(x_1,x_2) = l'(x_2,x_1)$.
Then,
one proves that the $\K$-module of relations $R^{\bot}(3)$
is generated by the composite element
\[r'(x_1,x_2,x_3) = l'(x_1,l'(x_2,x_3)) + l'(x_3,l'(x_1,x_2)) + l'(x_2,l'(x_3,x_1)).\]
(We refer to V. Ginzburg and M. Kapranov \cite{GinzburgKapranov} for details.)
In fact,
an algebra over this operad $\bar{K}(\C)^{\vee}$
is equivalent to the suspension of a Lie algebra $\Sigma\G$.
Explicitly,
the operation $l': \Sigma\G\otimes\Sigma\G\,\rightarrow\,\Sigma\G$ associated to $l'\in\bar{K}(\C)^{\vee}(2)$
is determined by the composite map
\[\Sigma\G\otimes\Sigma\G\,\xrightarrow{\simeq}\,\Sigma^2(\G\otimes\G)\,\xrightarrow{[-,-]}\,\Sigma^2\G.\]
The relation above is equivalent to the Jacobi identity.

Similarly,
the $\K$-module $\bar{K}(\L)^{\vee}(2)$
is generated by an operation $m' = m'(x_1,x_2)$ of degree $-1$
which verifies $m'(x_1,x_2) = - m'(x_2,x_1)$.
One proves that the $\K$-module of relations $R^{\bot}(3)$
is generated by the composite element
\[r'(x_1,x_2,x_3) = m'(m'(x_1,x_2),x_3) + m'(x_1,m'(x_2,x_3)).\]
In fact,
an algebra over this operad $\bar{K}(\L)^{\vee}$
is equivalent to the suspension of a commutative algebra $\Sigma A$.
Explicitly,
the operation $m': \Sigma A\otimes\Sigma A\,\rightarrow\,\Sigma A$ associated to $m'\in\bar{K}(\L)^{\vee}(2)$
is determined by the composite map
\[\Sigma A\otimes\Sigma A\,\xrightarrow{\simeq}\,\Sigma^2(A\otimes A)\,\xrightarrow{\cdot}\,\Sigma^2 A.\]
The relation above is equivalent to the associativity identity.

\subsection{Koszul homology of Lie algebras}
We observe that the Chevalley-Eilenberg homology of a Lie algebra
is a variant of the Koszul homology theory
introduced in paragraph~\ref{KoszulHomology}.
To be more precise,
the construction of paragraph~\ref{KoszulHomology} gives a chain complex $K_*(I,\L,\G)$
such that
\[K_*(I,\L,\G) = S(\bar{K}_*(\L),\G) = \bigoplus_{r=1}^{\infty} (\sgn_r\otimes\G^{\otimes r})_{\Sigma_r}.\]
Thus,
if the ground ring $\K$ is a field of characteristic $p\not=2$,
then we obtain exactly the (reduced) Chevalley-Eilenberg complex $C^{CE}_*(\G)$.
In the general case,
the Chevalley-Eilenberg complex is given by the functor $\Lambda(\bar{K}(\L))$
introduced in paragraph~\ref{DividedPowerFunctors} (see also proposition~\ref{DividedCommutativeAlgebras}).
We have explicitly:
\[C^{CE}_*(\G) = \Lambda(\bar{K}_*(\L),\G).\]

Similarly,
we observe in proposition~\ref{DividedLieAlgebras}
that the free Lie algebra corresponds to the functor $V\mapsto\Lambda(\L,V)$.
Consequently,
the Chevalley-Eilenberg complex of free Lie algebras
defines a functor $V\mapsto C^{CE}_*(\Lambda(\L,V))$
such that
\[C^{CE}_*(\Lambda(\L,V)) = \Lambda(\bar{K}_*(\L)\circ\L,V) = \Lambda_*(K_*(I,\L,\L),V).\]
The next assertion is classical:

\begin{lem}\label{VanishingChevalleyEilenbergHomology}
The Chevalley-Eilenberg complex of a free Lie algebra is acyclic.
\end{lem}

We deduce from this result that the Koszul complex $K_*(I,\L,\L)$ is acyclic,
because, as in the proof of proposition~\ref{VanishingKoszulHomology},
the dg-$\Sigma_*$-module $K_*(I,\L,\L)$
is associated to the functor $V\mapsto C^{CE}_*(\Lambda(\L,V))$.
Consequently,
we obtain the following theorem:

\begin{thm}
The Lie operad $\L$ is Koszul.
In particular, the comparison morphisms
\[\bar{K}(\L)\,\rightarrow\,\bar{B}(\L)\,\rightarrow\,\bar{N}(\L)\]
are all quasi-isomorphisms.
\end{thm}

\subsection{About the proof of lemma~\ref{VanishingChevalleyEilenbergHomology}}
Let us recall the arguments involved in the proof of lemma~\ref{VanishingChevalleyEilenbergHomology}.
We refer to the classical reference \cite[J.-L. Koszul]{Koszul}
for more details about the homology of Lie algebras.
In this paragraph,
we adopt the classical notation $\L(V)$ for the free Lie algebra,
$\Lambda(V)$ for the exterior algebra,
$S(V)$ for the symmetric algebra
and $T(V)$ for the tensor algebra.
First of all,
we recall that the free Lie algebra $\L(V)$
forms a free $\Z$-module
(\emph{cf}. C. Reutenauer \cite{Reutenauer}).
By adjunction,
the enveloping algebra of the free Lie algebra $U(\L(V))$
is isomorphic to the tensor algebra $T(V)$
(\emph{cf. loc. cit}.).

One observes that the Chevalley-Eilenberg complex with coefficients
\[C^{CE}_*(\L(V),U(\L(V)) = \Lambda(\L(V))\otimes U(\L(V))\]
forms an acyclic complex.
For that purpose,
one considers the filtration of the enveloping algebra $U(\L(V))$
which, according to the Poincar\'e-Birkhoff-Witt theorem,
satisfies the relation $S(\L(V)) = \gr_* U(\L(V))$.
Then,
we obtain a spectral sequence
\[E^0\,\Rightarrow\,H^{CE}_*(\L(V),U(\L(V)).\]
such that
$(E^0,d^0) = (\Lambda(\L(V))\otimes S(\L(V)),d^0)$
is the original Koszul complex
(\emph{cf}. J.-L. Koszul \cite{Koszul}, see also \cite[S. Priddy]{PriddyKoszul}
and \cite[P. May]{MayLie}).
This chain complex is acyclic
and the assertion above follows.

One can deduce from the vanishing result above
that the Chevalley-Eilenberg homology verifies
\[H^{CE}_*(\L(V)) = \Tor^{U(\L(V))}_*(\K,\K) = \Tor^{T(V)}_*(\K,\K).\]
In turn,
it is not difficult to determine the modules $\Tor^{T(V)}_*(\K,\K)$
(see for instance \cite[J.-L. Loday]{Loday}).
In fact,
the tensor algebra $T(V)$ is an obvious example of a Koszul algebra
(\emph{cf}. S. Priddy \cite{PriddyKoszul}).
Equivalently,
we have a natural morphism of right $U(\L(V))$-modules
\[\epsilon: \Lambda(\L(V))\otimes U(\L(V))\,\rightarrow\,V\otimes U(\L(V)).\]
According to the classical theory of Koszul algebras,
the source and the target of this morphism
are both acyclic chain complexes (see above).
Consequently,
we have an induced quasi-isomorphism
\[\bar{\epsilon}: \Lambda(\L(V))\,\rightarrow\,V.\]
This proves lemma~\ref{VanishingChevalleyEilenbergHomology}.

\subsection{Koszul homology of Commutative algebras}\label{KoszulCommutativeHomology}
We observe that the Koszul homology of a commutative algebra
coincides with the classical Harrison complex.
In fact,
the classical Harrison complex with trivial coefficients $C^{Harr}_*(A)$
is the indecomposable part of the Hochschild complex
(\emph{cf}. M. Barr \cite{BarrHarrison}, D. Harrison \cite{Harrison}, J.-L. Loday \cite{Loday}).
Consequently,
this chain complex is given by the dual $\K$-module of the free restricted Lie algebra
(\emph{cf}. C. Reutenauer \cite{Reutenauer}).
According to proposition~\ref{DividedLieAlgebras},
we obtain exactly
\[C^{Harr}_*(A) = S(\bar{K}_*(\C),A) = K_*(I,\C,A).\]
As mentioned in the prolog,
one observes that the Harrison homology of a symmetric algebra
does not vanish in positive characteristic
(\emph{cf}. M. Barr \cite{BarrHarrison}, D. Harrison \cite{Harrison}).
The same problem holds
for a variant of the Harrison complex
(\emph{cf}. S. Whitehouse \cite{Whitehouse}).
In fact,
the calculations of P. Goerss
(\emph{cf}. \cite{Goerss})
prove that the Andr\'e-Quillen homology (which is equivalent to the cotriple homology)
can not agree with any variant of Harrison homology.

Nevertheless,
by Koszul duality,
we obtain:

\begin{thm}\label{CommutativeKoszulness}
The commutative operad $\C$ is Koszul.
In particular, the comparison morphisms
\[\bar{K}(\C)\,\rightarrow\,\bar{B}(\C)\,\rightarrow\,\bar{N}(\C)\]
are all quasi-isomorphisms.
\end{thm}

\subsection{About the proof of theorem~\ref{CommutativeKoszulness}}
To be precise,
we deduce from lemma~\ref{KoszulDualityEquivalence} that the operad
\[\C' = \bar{K}(\L)^{\vee} = \bar{K}^c(\L^{\vee})\]
is Koszul,
because the Lie operad $\L$ is so.
We obtain the same result for the commutative operad $\C$
because $\C' = \bar{K}(\L)^{\vee}$ is essentially a suspension of $\C$.

Theorem~\ref{CommutativeKoszulness}
can also be deduced from the work of A. Bj\"orner
(\emph{cf}. \cite{Bjorner}).
Namely,
this author proves directly that the reduced homology of the partition posets $\tilde{H}_*(\tilde{K}(r))$
vanishes in degree $*\not=r-1$.
Consequently,
in regard to the bar construction,
we obtain
\[\tilde{H}_*(\bar{B}_*(\C)(r)) = \tilde{H}_*(\bar{N}_*(\C)(r)) = 0\]
if $*\not=r-1$.
This property implies immediately
that the commutative operad is Koszul
(\emph{cf}. theorem~\ref{KoszulCharacterization}).

Let us point out that the result of A. Bj\"orner reproves that the Lie operad consists of free $\Z$-modules.
For that purpose,
we recall that the module $\bar{K}(\C)(r)$
is defined by the kernel of $\beta: \tilde{B}_{r-1}(\C)(r)\,\rightarrow\,\tilde{B}_{r-2}(\C)(r)$
and that the bar construction $\tilde{B}_*(\C)(r)$ forms a projective $\Z$-module,
because the commutative operad $\C$ is so.
Since $\C$ is Koszul (according to Bj\"orner's theorem),
the $\Z$-module $\bar{K}^c(\C^{\vee})$,
given by the cokernel of $\beta: \tilde{B}^c_{-r+2}(\C^{\vee})(r)\,\rightarrow\,\tilde{B}^c_{-r+1}(\C^{\vee})(r)$,
is dual to $\bar{K}(\C)(r)$
(\emph{cf}. lemma~\ref{LinearDualKoszulOperads}).
The assertion about the Lie operad follows,
because one proves that this cokernel $\bar{K}^c(\C^{\vee})$
is essentially a suspension of the Lie operad $\L$
(see fact~\ref{ClassicalKoszulDuals}).

\smallskip
As mentioned in the prolog,
we would like to point out that the generalized Koszul construction $P\mapsto\bar{K}(P)$
allows to obtain the right homology modules for all partition posets
considered by P. Hanlon and M. Wachs in \cite{HanlonWachs}
and related to generalized Lie algebra structures
(see also A.V. Gnedbaye \cite{Gnedbaye}).

\subsection{On cofibrant resolutions of the commutative operad and $E_\infty$-algebras}
We deduce some results about cofibrant resolutions
of the commutative operad $\C$ from theorem~\ref{CommutativeKoszulness}.
Explicitly,
if $Q = F(M)$ is a cofibrant quasi-free resolution of the commutative operad $\C$
such that $M(0) = M(1) = 0$,
then we have a quasi-isomorphism of dg-$\Sigma_*$-modules
$M\,\xrightarrow{\sim}\,\Sigma^{-1}\tilde{K}(\C)$.
Conversely,
if we are given a resolution of $\Sigma^{-1}\tilde{K}(\C)$ by projective $\Sigma_*$-modules,
let $M\,\xrightarrow{\sim}\,\Sigma^{-1}\tilde{K}(\C)$,
then we have a cofibrant quasi-free resolution of $\C$ such that $Q = F(M)$
(\emph{cf}. proposition~\ref{KoszulCofibrantResolutions}).
For instance,
we can assume that $M(r)$ is the tensor product $M(r) = \Sigma^{-1}\tilde{K}(\C)(r)\otimes\E(r)$,
where $\E(r)$ denotes the classical bar construction
of the symmetric group $\Sigma_r$.

An algebra over a fixed cofibrant resolution $Q$ of the commutative operad $\C$ is called an \emph{$E_\infty$-algebra}.
The augmentation morphism $\epsilon: Q\,\rightarrow\,\C$
gives rise to a restriction functor $\epsilon^!: \LAlg{\C}\,\rightarrow\,\LAlg{Q}$
(see paragraph \ref{ExtensionRestrictionFunctors}).
Therefore,
commutative algebras are $E_\infty$-algebras.
In fact,
the definition of an $E_\infty$-algebra is motivated by the following observations.
Firstly,
in positive characteristic,
the dg-algebras over the commutative operad
do not form a model category
(see the argument given in paragraph \ref{LiftingProperty})
unlike algebras over a cofibrant operad.
Secondly,
cofibrant resolutions of the commutative operad
have equivalent homotopy categories of algebras
(see C. Berger and I. Moerdijk \cite{BergerMoerdijk}, V. Hinich \cite{Hinich}).

\subsection{On the Barratt-Eccles operad and $E_\infty$-algebras}
The classical bar constructions of the symmetric groups $\E(r)$ form an operad,
called the \emph{Barratt-Eccles operad},
and the augmentation morphisms $\E(r)\,\rightarrow\,\K$
define an operad quasi-isomorphism $\E\,\xrightarrow{\sim}\,\C$
(\emph{cf}. M. Barratt and P. Eccles \cite{BarrattEccles}, C. Berger and B. Fresse \cite{BergerFresse}).
We would like to recall that algebras over the Barratt-Eccles operad
are equivalent to $E_\infty$-algebras,
although the Barratt-Eccles operad is not cofibrant.
To be precise,
one observes that algebras over the Barratt-Eccles operad
form a model category (\emph{cf}. \cite{BergerFresse}).
Furthermore,
if $Q$ is a cofibrant resolution of the commutative operad $\C$,
then we have an operad quasi-isomorphism $\phi: Q\,\xrightarrow{\sim}\,\C$
which arises from the lifting diagram
\[\xymatrix{ & \E\ar[d]^{\sim} \\ Q\ar@{-->}[ur]\ar[r]^{\sim} & \C \\ }\]
One can prove that the derived extension functors
$L\phi_!: \Ho(\LAlg{Q})\,\rightarrow\,\Ho(\LAlg{\E})$
and the restriction functors
$R\phi^!: \Ho(\LAlg{\E})\,\rightarrow\,\Ho(\LAlg{Q})$
yield adjoint equivalences of homotopy categories
(\emph{cf}. \cite{BergerMoerdijk}).
In fact,
this equivalence of homotopy categories can be generalized to any $\Sigma_*$-projective
resolution of the commutative operad
(\emph{cf}. M. Mandell \cite{Mandell}, M. Spitzweck \cite{Spitzweck}).
But, in this context,
one has to introduce semi-model structures,
because algebras over a  $\Sigma_*$-projective operad
do not form a model category in general.

\subsection{On the homology of $E_\infty$-algebras}
As $E_\infty$-algebras form a model category,
there is a homology theory for $E_\infty$-algebras,
called \emph{$\Gamma$-homology}
(\emph{cf}. A. Robinson and S. Whitehouse \cite{RobinsonWhitehouseGammaH}, A. Robinson \cite{Robinson}),
which can be defined as the classical Quillen homology
in the category of simplicial commutative algebras
(\emph{cf}. D. Quillen \cite{QuillenProc}).
Let us mention that $\Gamma$-homology of commutative algebras differs from Quillen homology
(\emph{cf. loc. cit.} and B. Richter and A. Robinson \cite{RichterRobinson}).

According to A. Robinson (\emph{cf}. \cite{Robinson}),
the $\Gamma$-homology of a commutative algebra $A$
can be determined by a chain complex $C^{\Gamma}_*(A)$
such that
\[C^{\Gamma}_*(A) = \bigoplus_{r=0}^{\infty} (\Sigma^{-1}\tilde{K}(\C)(r)\otimes\E(r)\otimes A^{\otimes r})_{\Sigma_r}.\]
On the other hand,
according to M. Kontsevich and Y. Soibelman (\emph{cf}. \cite{KontsevichSoibelman}),
there is a similar complex
\[C^{\Gamma}_*(A) = \bigoplus_{r=0}^{\infty} (M(r)\otimes A^{\otimes r})_{\Sigma_r}.\]
for any cofibrant quasi-free resolution of the commutative operad $Q$
such that $Q = F(M)$
(at least in characteristic $0$).
In some sense,
the result of A. Robinson makes this construction
explicit for a cofibrant quasi-free resolution $Q = F(M)$
such that $M(r) = \Sigma^{-1}\tilde{K}(\C)(r)\otimes\E(r)$.

%\input{References}
% 18/1/2003

\backmatter
%\input{Index}
% 13/11/2002

\chapter*{Glossary}

\begin{itemize}

\renewcommand{\labelitemi}{}
\renewcommand{\labelitemii}{}

\item algebra
\begin{itemize}
\item algebra over a monad, \ref{MonadAlgebras}
\item algebra over an operad, \ref{OperadAlgebras},  \ref{LeftModules}
\end{itemize}

\item augmentation ideal of an operad,  \ref{ConnectedOperads},  \ref{ConnectedWeightGradedOperad}

\item bar construction
\begin{itemize}
\item differential graded bar construction with coefficients,  \ref{BarCoefficientsPpties},  \ref{DifferentialGradedBarConstruction}
\item simplicial bar construction with coefficients,  \ref{SimplicialBarConstructionPpties},  \ref{SimplicialBarConstruction}
\item reduced bar construction,  \ref{ReducedBarPpties},  \ref{ConstructionBarDifferential}
\end{itemize}

\item cobar construction
\begin{itemize}
\item reduced cobar construction,  \ref{ReducedCobarPpties}
\end{itemize}

\item cofibrant resolution
\begin{itemize}
\item cofibrant resolution of an operad,  \ref{CofibrantOperadResolutions}
\item cofibrant resolution of a right-module over an operad,  \ref{ModuleResolutions}
\end{itemize}

\item composition product
\begin{itemize}
\item composition product of $\Sigma_*$-modules,  \ref{MonoidalFunc},  \ref{CompositionProductDefinition},  \ref{ExpansionComposite},  \ref{TreewiseComposite}
\item relative composition product of modules over an operad,  \ref{RelativeCompositionProducts}
\item derived relative composition product of modules over an operad,  \ref{ModuleResolutions}
\item composition products of an operad,  \ref{ClassicalOperadDefinition}
\item partial composition products of an operad,  \ref{PartialComposites},  \ref{OperadIndexing}
\end{itemize}

\item connected
\begin{itemize}
\item connected operad,  \ref{ConnectedOperads},  \ref{ConnectedWeightGradedOperad}
\item connected $\Sigma_*$-module,  \ref{SymmetricFunctors}
\end{itemize}

\item cooperad,  \ref{Cooperads}
\begin{itemize}
\item bar cooperad,  \ref{ReducedBarPpties}
\item coalgebra over a cooperad,  \ref{Cooperads}
\item coderivation,  \ref{OperadDerivationDefinition},  \ref{FreeOperadDerivation},  \ref{TreeCoderivation}
\item cofree cooperad,  \ref{CofreeCooperad}
\item quasi-cofree cooperad,  \ref{QuasiFreeOperads}
\end{itemize}

\item cross effect,  \ref{PolynomialCrossEffect}

\item differential,  \ref{DGMod}
\begin{itemize}
\item bar differential,  \ref{ReducedBarDifferential},  \ref{CoefficientBarDifferential},  \ref{TwistingBarDifferential}
\item differential graded module (also dg-module),  \ref{DGMod}
\item differential graded $\Sigma_*$-module (also dg-$\Sigma_*$-module),  \ref{DGSigmaModule}
\item differential graded operad (also dg-operad),  \ref{DGOperads}
\item differential of a composite $\Sigma_*$-module,  \ref{DGSigmaModule},  \ref{CompositionBicomplex}
\item differential of a quasi-free left module over an operad,  \ref{QuasiFreeLeftModule}
\item differential of a quasi-free right module over an operad,  \ref{QuasiFreeRightModule},  \ref{RightModuleDifferential}
\item differential of a quasi-free operad,  \ref{QuasiFreeOperads}
\item differential of a quasi-cofree cooperad,  \ref{QuasiFreeOperads}
\end{itemize}

\item free
\begin{itemize}
\item free algebra,  \ref{MonadAlgebras}
\item free left module over an operad,  \ref{FreeModules}
\item free operad,  \ref{FreeOperad},  \ref{dgFreeOperad},  \ref{FreeOperadStructure}
\item free right module over an operad,  \ref{FreeModules}
\end{itemize}

\item functor associated to a $\Sigma_*$-module,  \ref{ClassicalOperadDefinition},  \ref{SymmetricFunctors},  \ref{DividedPowerFunctors}

\item indecomposable quotient
\begin{itemize}
\item indecomposable quotient of a left-module over an operad,  \ref{ModuleIndecomposableQuotient}
\item indecomposable quotient of an operad,  \ref{FreeOperad}
\item indecomposable quotient of a quasi-free left module over an operad,  \ref{QuasiFreeLeftModule}
\item indecomposable quotient of a quasi-free operad,  \ref{QuasiFreeOperads}
\item indecomposable quotient of a quasi-free right module over an operad,  \ref{QuasiFreeRightModule}
\item indecomposable quotient of a right-module over an operad,  \ref{ModuleIndecomposableQuotient}
\end{itemize}

\item Koszul
\begin{itemize}
\item Koszul complex with coefficients,  \ref{CoefficientKoszulConstruction}
\item Koszul construction,  \ref{KoszulConstruction}
\item Koszul homology,  \ref{KoszulHomology}
\item Koszul operad,  \ref{DefinitionKoszulOperad}
\item Koszul resolution of an operad,  \ref{KoszulResolutionDefinition}
\end{itemize}

\item levelization
\begin{itemize}
\item levelization morphism,  \ref{Levelization},  \ref{LevelizationProcess}
\item levelization of a composite tree,  \ref{TreeLevelization}
\end{itemize}

\item modules over an operad
\begin{itemize}
\item bimodule over an operad,  \ref{LeftModules}
\item free left module over an operad,  \ref{FreeModules}
\item free right module over an operad,  \ref{FreeModules}
\item left module over an operad,  \ref{LeftModules}
\item quasi-free left module over an operad,  \ref{QuasiFreeRightModule}
\item quasi-free right module over an operad,  \ref{QuasiFreeLeftModule}
\item right module over an operad,  \ref{RightModules}
\end{itemize}

\item monad,  \ref{MonoidMonads}
\begin{itemize}
\item algebra over a monad,  \ref{MonadAlgebras}
\item monad associated to an operad,  \ref{ClassicalOperadDefinition},  \ref{SymmetricMonads},  \ref{DividedPowerMonads}
\end{itemize}

\item operad,  \ref{ClassicalOperadDefinition},  \ref{MonoidOperads}
\begin{itemize}
\item algebra over an operad,  \ref{OperadAlgebras},  \ref{LeftModules}
\item augmented operad,  \ref{ConnectedOperads}
\item associative operad,  \ref{ClassicalOperadExamples}
\item commutative operad,  \ref{ClassicalOperadExamples},  \ref{DividedCommutativeAlgebras}
\item composition products of an operad,  \ref{ClassicalOperadDefinition}
\item cobar operad,  \ref{ReducedCobarPpties}
\item cofibrant resolution of an operad,  \ref{CofibrantOperadResolutions}
\item connected operad,  \ref{ConnectedOperads},  \ref{ConnectedWeightGradedOperad}
\item operad derivation,  \ref{OperadDerivationDefinition}
\item free operad,  \ref{FreeOperad},  \ref{dgFreeOperad},  \ref{FreeOperadStructure}
\item left module over an operad,  \ref{LeftModules}
\item Koszul operad,  \ref{DefinitionKoszulOperad}
\item Koszul resolution of an operad,  \ref{KoszulResolutionDefinition}
\item Lie operad,  \ref{ClassicalOperadExamples},  \ref{DividedLieAlgebras}
\item operad ideal,  \ref{OperadIdeals}
\item operad of trees,  \ref{TreeOperad}
\item partial composition products of an operad,  \ref{PartialComposites},  \ref{OperadIndexing}
\item quadratic operad,  \ref{QuadraticOperads},  \ref{BinaryOperadDuality}
\item quasi-free operad,  \ref{QuasiFreeOperads}
\item quasi-free resolution of an operad,  \ref{CofibrantOperadResolutions}
\item right module over an operad,  \ref{RightModules}
\item weight graded operad,  \ref{WeightOperads}
\end{itemize}

\item projective
\begin{itemize}
\item $\K$-projective $\Sigma_*$-module,  \ref{DGSigmaModule}
\item $\Sigma_*$-projective $\Sigma_*$-module,  \ref{DGSigmaModule}
\end{itemize}

\item quadratic operad,  \ref{QuadraticOperads},  \ref{BinaryOperadDuality}

\item quasi-free
\begin{itemize}
\item quasi-free left module over an operad,  \ref{QuasiFreeRightModule}
\item quasi-free operad,  \ref{QuasiFreeOperads}
\item quasi-free resolution of an operad,  \ref{CofibrantOperadResolutions}
\item quasi-free resolution of a right-module over an operad,  \ref{ModuleResolutions}
\item quasi-free right module over an operad,  \ref{QuasiFreeLeftModule}
\end{itemize}

\item $\Sigma_*$-module,  \ref{SymmetricFunctors}
\begin{itemize}
\item composition product of $\Sigma_*$-modules,  \ref{MonoidalFunc},  \ref{CompositionProductDefinition},  \ref{ExpansionComposite},  \ref{TreewiseComposite}
\item connected $\Sigma_*$-module,  \ref{SymmetricFunctors}
\item functors associated to a $\Sigma_*$-module,  \ref{SymmetricFunctors},  \ref{DividedPowerFunctors}
\item $\K$-projective $\Sigma_*$-module,  \ref{DGSigmaModule}
\item $\Sigma_*$-module associated to a functor,  \ref{CrossEffectModule},  \ref{MonoidalCoef}
\item $\Sigma_*$-projective $\Sigma_*$-module,  \ref{DGSigmaModule}
\item tensor product of $\Sigma_*$-modules,  \ref{MonoidalFunc},  \ref{SymmetricTensorProduct}
\end{itemize}

\item suspension of a dg-module,  \ref{Suspension}

\item tree,  \ref{TreeStructure}
\begin{itemize}
\item composite tree,  \ref{CompositeTrees}
\item isomorphism of trees,  \ref{TreeIsomorphisms}
\item levelization of a composite tree,  \ref{TreeLevelization}
\item operad of trees,  \ref{TreeOperad}
\item quotient tree,  \ref{QuotientTrees}
\item reduced tree,  \ref{ReducedFreeOperad}
\item subtree,  \ref{SubTrees}
\item tree with levels,  \ref{LevelTrees}
\end{itemize}

\item treewise tensor product,  \ref{TreewiseTensors},  \ref{OperadBarSimplices},  \ref{ExpansionCoefficientBar}

\item twisting cochain,  \ref{TwistingCochains},  \ref{CoefficientKoszulConstruction}

\item weight
\begin{itemize}
\item weight grading,  \ref{WeightGradings}
\item weight graded operad,  \ref{WeightOperads}
\end{itemize}

\end{itemize}

\chapter*{Notation}

\begin{itemize}

\newcommand{\notationentry}[1]{\makebox[6em][l]{#1}}

\item[\notationentry{$\K$}] the ground ring, \ref{Conventions}
\item[\notationentry{$\Sigma_r$}] the symmetric group on $r$ elements, \ref{Conventions}
\item[\notationentry{$S(P)$}] the monad associated to an operad,  \ref{ClassicalOperadDefinition}
\item[\notationentry{$\tilde{P}$}] the augmentation ideal of an operad,  \ref{ConnectedOperads}
\item[\notationentry{$F(M)$}] the free operad,  \ref{FreeOperad}
\item[\notationentry{$\C$}] the commutative operad,  \ref{ClassicalOperadExamples}
\item[\notationentry{$\A$}] the associative operad,  \ref{ClassicalOperadExamples}
\item[\notationentry{$\L$}] the Lie operad,  \ref{ClassicalOperadExamples}
\item[\notationentry{$S(M)$}] the coinvariant functor associated to a $\Sigma_*$-module,  \ref{SymmetricFunctors}
\item[\notationentry{$\Lambda(M)$}] the norm functor associated to a $\Sigma_*$-module,  \ref{DividedPowerFunctors}
\item[\notationentry{$\Gamma(M)$}] the invariant functor associated to a $\Sigma_*$-module,  \ref{DividedPowerFunctors}
\item[\notationentry{$M\otimes N$}] the tensor product of $\Sigma_*$-modules,  \ref{SymmetricTensorProduct}
\item[\notationentry{$M\circ N$}] the composition product of $\Sigma_*$-modules,  \ref{CompositionProductDefinition}
\item[\notationentry{$V^{\vee}$}] the dual of a $\K$-module, \ref{Conventions}
\item[\notationentry{$M^{\vee}$}] the dual of a $\Sigma_*$-module, \ref{Conventions}
\item[\notationentry{$L\circ_P R$}] the relative composition product over an operad,  \ref{RelativeCompositionProducts}
\item[\notationentry{$\bar{L}$}] the indecomposable quotient of a right module,  \ref{ModuleIndecomposableQuotient}
\item[\notationentry{$\bar{R}$}] the indecomposable quotient of a left module,  \ref{ModuleIndecomposableQuotient}
\item[\notationentry{$d_\theta$}] a derivation of a free right module over an operad,  \ref{QuasiFreeRightModule},\newline
or a derivation of a free left module over an operad,  \ref{QuasiFreeLeftModule},\newline
or a derivation of a free operad,  \ref{FreeOperadDerivation},\newline
or a coderivation of a cofree cooperad,  \ref{FreeOperadDerivation}
\item[\notationentry{$F^c(M)$}] the cofree cooperad,  \ref{CofreeCooperad}
\item[\notationentry{$\bar{B}_*(P)$}] the reduced bar construction of an operad,  \ref{ReducedBarPpties}
\item[\notationentry{$\bar{B}^c_*(D)$}] the reduced bar construction of a cooperad,  \ref{ReducedCobarPpties}
\item[\notationentry{$V(\tau)$}] the set of vertices of a tree,  \ref{TreeStructure}
\item[\notationentry{$E(\tau)$}] the set of edges of a tree,  \ref{TreeStructure}
\item[\notationentry{$I_v$}] the entries of a vertex $v$ in a tree,  \ref{TreeStructure}
\item[\notationentry{$\tau(M)$}] the treewise tensor product,  \ref{TreewiseTensors}
\item[\notationentry{$B_*(L,P,R)$}] the differential graded bar construction with coefficients,  \ref{BarCoefficientsPpties}
\item[\notationentry{$C_*(L,P,R)$}] the simplicial bar construction,  \ref{SimplicialBarConstructionPpties}
\item[\notationentry{$N_*(L,P,R)$}] the normalized chain complex of the simplicial bar construction,  \ref{Levelization},  \ref{NormalizedRightModule}
\item[\notationentry{$\bar{N}_*(P)$}] the normalized chain complex of the simplicial bar construction with trivial coefficients,  \ref{Levelization}
\item[\notationentry{$V_i(\tau)$}] the set of vertices of level $i$ in a tree with levels,  \ref{LevelTrees},\newline
or the set of vertices of level $i = l,p,r$ in a composite tree,  \ref{CompositeTrees}
\item[\notationentry{$\tau(L,P,R)$}] the treewise tensor product with coefficients,  \ref{OperadBarSimplices},  \ref{ExpansionCoefficientBar}
\item[\notationentry{$\bar{K}_*(P)$}] the Koszul construction of an operad,  \ref{KoszulConstruction}
\item[\notationentry{$K_*(L,P,R)$}] the Koszul complex with coefficients,  \ref{CoefficientKoszulConstruction}

\end{itemize}

\end{document}